\documentclass[11pt,leqno,a4paper]{book}
\usepackage{amsmath}
\usepackage{amsthm}
\usepackage{amssymb}
\usepackage{amscd}
\usepackage{mathabx}
\usepackage{accents}
\usepackage{makeidx}
\usepackage{array}

\usepackage{tikz, tikz-cd}
\usetikzlibrary{matrix,arrows,decorations.pathmorphing}

\usepackage{enumerate}
\usepackage{verbatim}

\usepackage{hyperref}

\usepackage[english]{babel}
\usepackage[babel]{csquotes}
\usepackage{guit}

\usepackage[totalheight=22 true cm, totalwidth=14 true cm]{geometry}
\usepackage{mathrsfs}
\usepackage[all]{xy}

\newtheorem{thm}{Theorem}[section]
\newtheorem{cor}[thm]{Corollary}
\newtheorem{lemma}[thm]{Lemma}

\newtheorem{prop}[thm]{Proposition}

\newtheorem{defn}[thm]{Definition}
\newtheorem{conj}[thm]{Conjecture}

\theoremstyle{remark}

\theoremstyle{definition}

\newtheorem{rmk}[thm]{Remark}
\newtheorem{exa}[thm]{Example}

\newtheorem{notation}[thm]{Notation}
\numberwithin{equation}{thm}

\def\beq{\begin{equation}}
\def\eeq{\end{equation}}

\def\beqn{\begin{equation*}}
\def\eeqn{\end{equation*}}

\def\ben{\begin{enumerate}}
\def\een{\end{enumerate}}

\def\crash#1{}
\def\A{{\mathbb A}}

\def\C{{\mathbb C}}
\def\D{{\mathbb D}}

\def\F{{\mathbb F}}

\def\N{{\mathbb N}}

\def\Q{{\mathbb Q}}
\def\R{{\mathbb R}}
\def\T{{\mathbb T}}

\def\Z{{\mathbb Z}}

\def\l{\left}
\def\r{\right}

\def\an{{\rm an}}

\def\End{{\rm End}}
\def\Im{{\rm Im \,}}
\def\Ker{{\rm Ker \,}}
\def\Max{{\rm Max \,}}
\def\Coim{{\rm Coim \,}}
\def\Coker{{\rm Coker \,}}
\def\cf{\emph{cf.}~}
\def\ie{\emph{i.e.}~}

\def\cB{{\mathcal B}}
\def\cC{{\mathcal C}}
\def\cD{{\mathcal D}}
\def\cE{{\mathcal E}}
\def\cF{{\mathcal F}}

\def\cI{{\mathcal I}}
\def\cM{{\mathcal M}}
\def\cN{{\mathcal N}}
\def\cO{{\mathcal O}}
\def\cH{{\mathcal H}}

\def\cP{{\mathcal P}}
\def\cR{{\mathcal R}}
\def\cS{{\mathcal S}}
\def\cT{{\mathcal T}}
\def\cV{{\mathcal V}}
\def\cU{{\mathcal U}}
\def\cX{{\mathcal X}}
\def\cY{{\mathcal Y}}

\def\sC{{\mathscr C}}

\def\sF{{\mathscr F}}
\def\sG{{\mathscr G}}

\def\sI{{\mathscr I}}

\def\sP{{\mathscr P}}

\def\sU{{\mathscr U}}
\def\sV{{\mathscr V}}

\def\bC{{\mathbf C}}

\def\bF{{\mathbf F}}

\def\bR{{\mathbf R}}
\def\bS{{\mathbf S}}

\def\bone{{\mathbf 1}}
\def\btwo{{\mathbf 2}}

\def\bee{{\mathbf e}}
\def\bk{{\mathbf k}}
\def\bn{{\mathbf n}}
\def\bm{{\mathbf m}}

\def\fm{{\mathfrak m}}

\def\wtilde{\widetilde}

\def\what{\widehat}

\def\a{\alpha}
\def\be{\beta}

\def\la{\lambda}

\def\an{{\rm an}}

\def\Spec{{\rm Spec \,}}
\def\Hom{{\rm Hom \,}}

\def\Gal{{\rm Gal\,}}

\def\id{{\rm id\,}}
\def\ker{{\rm Ker\,}}

\def\ul{\underline}

\def\ol{\overline}

\def\limpro{\mathop{\lim\limits_{\displaystyle\leftarrow}}}

\def\limind{\mathop{\lim\limits_{\displaystyle\rightarrow}}}

\def\acb{{A^\circ_{\cB}}}
\def\then{\Rightarrow}
\def\Ct{\C^\times}

\def\ob{{\rm ob\,}}

\def\lt{\langle}
\def\gt{\rangle}

\def\void{{\rm \varnothing}}
\def\Frac{{\rm Frac\,}}

\def\rhook{{\hookrightarrow}}

\def\ov{{\rm \diamond\,}}       
\def\cc{{\rm \circ \circ\,}}    
\def\Id{{\rm Id\,}}

\def\Int{{\rm Int\,}}    
\def\Supp{{\rm Supp\,}}   
\def\an{{\rm an\,}}

\def\dim{{\rm dim\,}}  

\def\bMon{{\mathbf{Mon} \,}}
\def\bSet{\mathbf{Sets}}
\def\bBorn{\mathbf{Born}}
\def\bSBorn{\mathbf{SBorn}}
\def\bCBorn{\mathbf{CBorn}}
\def\bCat{\mathbf{Cat}}
\def\bHalos{\mathbf{Halos}}
\def\bAb{\mathbf{Ab}}
\def\bRings{\mathbf{Rings}}
\def\bDom{\mathbf{Dom}}
\def\bFields{\mathbf{Fields}}
\def\bMonoids{\mathbf{Monoids}}
\def\bAlg{\mathbf{Alg}}
\def\bSn{\mathbf{Sn}}
\def\bNrm{\mathbf{Nrm}}
\def\bBan{\mathbf{Ban}}
\def\bAff{\mathbf{Aff}}
\def\bStein{\mathbf{Stein}}

\def\bGerms{\mathbf{Germs}}

\def\bMod{\mathbf{Mod}}
\def\bCoh{\mathbf{Coh}}
\def\bAn{\mathbf{An}}
\def\bPic{\mathbf{Pic}}
\def\bPro{\mathbf{Pro}}
\def\bInd{\mathbf{Ind}}

\def\bOuv{\mathbf{Ouv}}
\def\bRig{\mathbf{Rig}}

\def\bTop{\mathbf{Top}}

\author{Federico Bambozzi \thanks{Universit\`{a} di Padova,
Dipartimento di matematica pura e applicata, Via Trieste, 63, 35121 Padova, Italy.}}

\title{}

\makeindex

\begin{document}

\thispagestyle{empty}

\linespread{1.2} 

\begin{center}

\begin{tabular}{m{2cm}l}
\includegraphics[width=20mm]{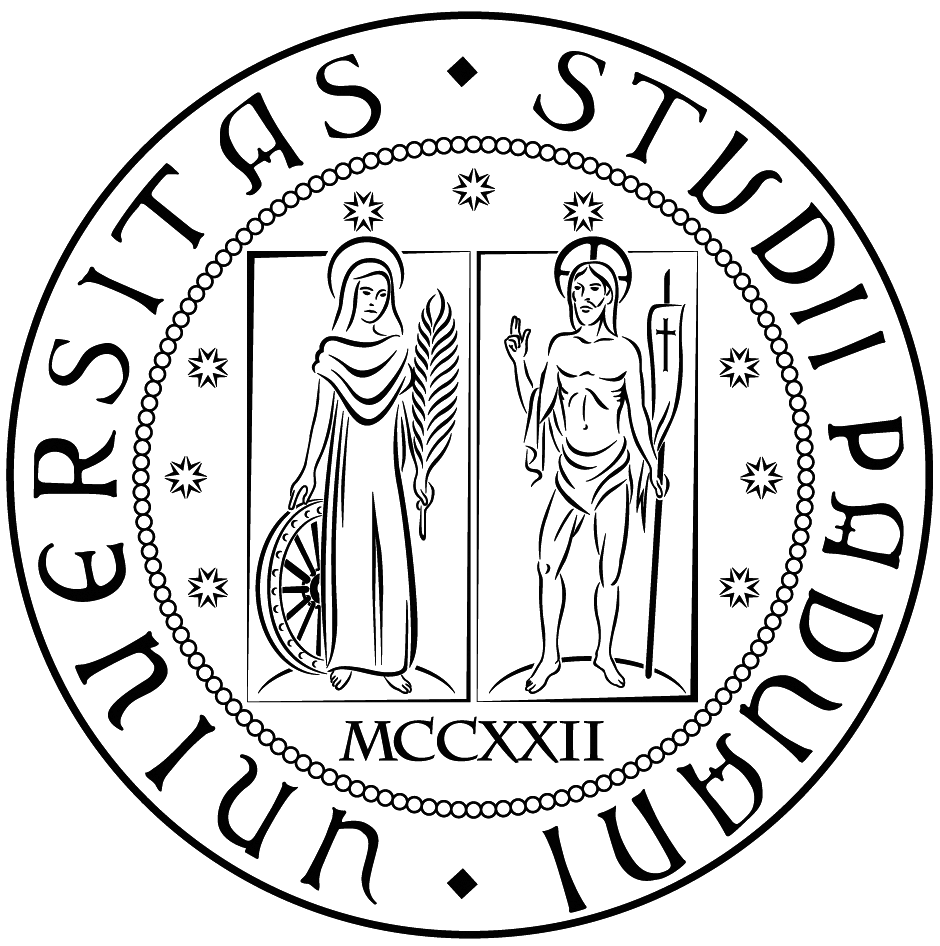} & \Large{\scshape{Universit\`{a} degli Studi di Padova}}
\end{tabular}

\end{center}

\vspace{\stretch{.3}}

\noindent Sede Amministrativa: Universit\`a degli Studi di Padova\\
\noindent Dipartimento di Matematica \\

\vspace{\stretch{.1}}

\noindent DOTTORATO DI RICERCA IN MATEMATICA \\
\noindent INDIRIZZO: MATEMATICA \\
\noindent CICLO XXVI

\vspace{\stretch{.6}}

\begin{center}

\textbf{\huge On a generalization of affinoid varieties }

\end{center}

\vspace{\stretch{1.2}}

\noindent Direttore della Scuola: Ch.mo Prof. Paolo Dai Pra \\
\noindent Coordinatore d'indirizzo: Ch.mo Prof. Franco Cardin \\
\noindent Supervisore: Ch.mo Prof. Francesco Baldassarri

\vspace{\stretch{.4}}

\begin{flushright}
Dottorando: Federico Bambozzi
\end{flushright}

\vspace{\stretch{1}}

\begin{center}
 \large{Dicembre 2013}
\end{center}

\linespread{1} 

\tableofcontents

\newpage

\section*{Abstract}
\thispagestyle{empty}
\markboth{}{}

In this thesis we develop the foundations for a theory of analytic geometry over a non-trivially valued field, uniformly encompassing the case 
when the base field is equipped with a non-archimedean valuation and that in which it has an archimedean one. We will use the theory of
bornological algebras to reach this goal. Since, at our knowledge, there is not a vast literature on bornological algebras on which we can base our results,
we will start from scratch with the theory of bornological algebras and we will develop it as far as we need for our scopes. In this way
we can construct our theory of analytic spaces taking as building blocks dagger affinoid algebras, \ie bornological algebras isomorphic to quotients of the algebras of germs of analytic
functions on polycylinders. It turns out that for the category of dagger affinoid algebras we can show the analogous of all the main results of the category of classical affinoid algebras and hence
we obtain a good theory of dagger affinoid spaces. Indeed, these spaces behave very similarly to affinoid spaces used to construct rigid analytic spaces.
We emphasize that the archimedean case has special features, that we study in the first section of the fifth chapter, and that these special properties allow
us to prove the generalization, in archimedean context, of the main result of the affinoid algebras theory: the Gerritzen-Grauert theorem 
(we remark that our proof is inspired by the new proof of Gerritzen-Grauert theorem in Berkovich geometry given by Temkin in \cite{TEM}).

After proving that we have a good affine theory, we end this thesis by constructing the global theory of dagger analytic spaces. We use the theory of Berkovich
nets and hence a Berkovich-like approach to the globalization, as done in \cite{BER4}. In this way we obtain the category of dagger analytic spaces and we
study the relations of the spaces we have found with the ones already present in literature. In particular, we see that our spaces are very strongly related to dagger spaces of Grosse-Kl\"onne, \cite{GK}, and that, in the archimedean case, the category of classical complex analytic spaces embeds fully faithfully in the category of complex dagger analytic spaces.

In conclusion, in this work we obtain a complete affinoid theory of dagger spaces over any valued field and we show that on this theory a good
theory of dagger analytic spaces can be developed, which deserves to be strengthened in the future.

\newpage

\section*{Riassunto}
\thispagestyle{empty}
\markboth{}{}

In questa tesi vengono sviluppate le basi per una teoria degli spazi analitici su campi valutati con valutazione non-banale, che comprende in modo uniforme il caso in cui il campo
base sia archimedeo o non-archimedeo. Per raggiungere questo obiettivo useremo la teoria delle algebre bornologiche. Siccome non sembra essere
disponibile una ampia letteratura riguardo le algebre bornologiche, la prima parte della tesi si occupa di stabilire alcuni risultati fondazionali a riguardo.
In questo modo \`{e} possibile costruire una teoria degli spazi analitici prendendo come mattoni fondamentali le algebre dagger affinoidi, \ie
algebre bornologiche isomorfe a quozienti delle algebre dei germi di funzioni analitiche su policilindri. Si ottiene quindi che la categoria delle algebre  dagger affinoidi
cosi definite soddisfa l'analogo di tutte le propriet\`a pi\`{u} importanti della categoria delle algebre affinoidi classiche e che di conseguenza 
si ottiene una buona teoria degli spazi dagger affinoidi la cui categoria \`{e} analoga alla categoria degli spazi affinoidi usati in geometria rigida. C'\`e da sottolineare che il caso in cui il campo base \`e archimedeo 
gli spazi dagger affinoidi hanno propriet\`a particolari, studiate all'inizio del quinto capitolo, e che queste propriet\`a permettono di 
ottenere anche per spazi dagger affinodi archimedei il principale risultato
strutturale della teoria delle algebre affinoidi: il teorema di Gerritzen-Grauert (la dimostrazione qui data di questo teorema \`{e} ispirata dalla nuova dimostrazione data da
Temkin in \cite{TEM} valida in geometria di Berkovich).

Dopo aver ottenuto una buona teoria affine, la tesi si conclude discutendo la teoria globale degli spazi dagger analitici. Verr\`{a} utilizzato un approccio analogo
a quello usato da Berkovich in \cite{BER4} per la globalizzazione. Dopo aver definito la categoria degli spazi dagger analitici verranno studiate le relazioni tra questi spazi e gli spazi analitici presenti in letteratura. In particolare si vedr\`a che questi spazi sono fortemente correlati agli spazi dagger di Grosse-Kl\"onne e che nel caso archimedeo la categoria degli spazi
analitici complessi classica si immerge in modo pienamente fedele nella categoria degli spazi dagger analitici complessi.

In conclusione, in questo lavoro si ottiene una teoria affinoide completa per spazi dagger su ogni campo valutato sulla base della quale viene proposta una teoria degli spazi dagger globali la quale dovr\`a essere approfondita in futuro.

\newpage

\section*{Acknowledgements}
\thispagestyle{empty}
\markboth{}{}

My utmost gratitude goes to my advisor, prof. Francesco Baldassari who gave the fascinating topic of my thesis. It has been an honour and pleasure to work with him, and it would not have been possible for me to do this thesis without the useful discussions that we had during the three years of my Ph.D.

A special thanks goes to Velibor Bojkovic for his help in checking arguments of proofs on earlier versions of this work and the help he gave me in order 
to understand Berkovich geometry.
I also wish to thank the referees of my Ph.D. dissertation, Oren Ben-Bassat and Elmar Grosse-Kl\"onne, for their comments (especially Elmar's critics which spurred me to give a better form and presentation to the text) and their encouragements without which I probably would not have had the energies to finish the last chapter in time for the bureaucratic deadlines. I am also very greatful to J\'er\^ome Poineau for pointing out several small mistakes and for useful comments. I also would like to thank Florent Martin for useful and insightful discussions about non-archimedean analytic geometry. I also would like to acknowledge Mauro Porta, Fr\'ederic Pauguam and Kobi Kremnizer for several interesting discussions.

Last, but not the least, a special thanks goes to Walter Gubler and Klaus K\"unnemann for their support and for giving me the time to finish improving the material I wrote during my Ph.D during my post-doc in Regensburg.

Carrying out this project would not have been possible without the financial support of the University of Padova by MIUR PRIN2010-11 ``Arithmetic Algebraic Geometry and Number Theory", and the University of Regensburg with the DFG funded CRC 1085 ``Higher Invariants. Interactions between Arithmetic Geometry and Global Analysis". I am very grateful for the support I received by both institutions.

\newpage

\section*{Notation}

\begin{itemize}

\item We shall follow standard notations: $\N$ is the set of natural numbers (with zero), $\Z$ the set of integers,  $\Q$ the set of rational numbers, $\R$ the set of real numbers, $\C$ the set of complex numbers, $\Q_p$ the set of $p$-adic numbers and we consider on them the usual operations. We shall use the notation $\R_{\ge 0} = \{ x \in \R | x \ge 0 \}$ and $\R_+ = \{ x \in \R | x > 0 \}$.

\item $k$ will always denote a fixed base field which is always supposed to be complete with respect to a non-trivial absolute value, archimedean or non-archimedean.

\item By an extension of valued fields $K/k$ we mean a valued field $K$ equipped with an embedding $k \rhook K$ which is an isometry onto its image.

\item With the term \emph{polyradius} we refer to a $n$-tuple $\rho = (\rho_1, \ldots, \rho_n) \in \R_+^n$. We define a partial order on the set of polyradii by setting $\rho < \rho'$ if and only if $\rho$ and $\rho'$ have the same number of components and $\rho_i < \rho_i'$ for every $1 \le i \le n$.

\item Given a polyradius $\rho = (\rho_1, \ldots, \rho_n)$, then:
\begin{itemize}
	\item $\D_k(\rho^+) = \{ x \in k^n | |x|_i \le \rho_i, \forall i \}$ is the \emph{closed polycylinder} of polyradius $\rho$;
	\item $\D_k(\rho^-) = \{ x \in k^n | |x|_i < \rho_i, \forall i \}$ is the \emph{open polycylinder} of polyradius $\rho$;
	\item if $\rho_i = r$, for some $r \in \R_+$ and all $i$, we say that $\D_k(\rho^+)$ and $\D_k(\rho^-)$ are \emph{polydisks}.
\end{itemize} 

\item For any ring $A$ we denote with $A^\times \subset A$ the subset of units of $A$.

\item For a commutative ring $A$ we denote with $\bMod_A$ the category of modules over $A$ and with $\bMod_A^f$ the category of finite modules over $A$.

\item Given a category $\bC$, we denote with $\ob(\bC)$ the class of objects of $\bC$.

\item $\bSet$ denotes the category of sets.

\item $\bSn_k$ denotes the category of semi-normed spaces over $k$, $\bNrm_k$ denotes the category of normed spaces over $k$, $\bBan_k$ denotes the category of Banach spaces over $k$ where the morphisms are bounded linear maps.

\item Given a category $\bC$, we denote by $\bInd(\bC)$ the category of ind-objects of $\bC$ and by $\bPro(\bC)$ the category of pro-objects of $\bC$, see appendix \ref{pro_appendix} for precise definitions.

\item If $\bC$ is a site we denote by $\wtilde{\bC}$ its associated topos.

\item Given a diagram $\sF:I \to \bC$, where $I$ is an index category, we say that the diagram is \emph{monomorphic} if for any morphism $\phi: i \to j$ in $I$ the morphism $\sF(\phi)$ is a monomorphism in $\bC$. In the same way, we say that the diagram is epimorphic if for all $\phi: i \to j$ the morphism $\sF(\phi)$ is an epimorphism in $\bC$.

\item We fix notation for inductive systems once for all, because we will use them many times in the text. For any ordered diagram $\{E_i\}_{i \in I}$ of objects of $\bC$ we denote by $\varphi_{i, j}$ the system morphisms $\varphi_{i,j}: F(i) \to F(j)$ for any $i, j \in I$ with $i \le j$. The notation is unambiguous because we will only consider diagrams for which there is at most one morphism between any couple of objects $i,j \in I$. We denote by $\a_i: A_i \to \underset{i \in I} \limind A_i$ the canonical morphisms to the colimit.

We say that the colimit $\underset{i \in I} \limind A_i$ is \emph{filtered} if the diagram $I \to \bC$ is directed, and we say that is \emph{filtrant} if $I \to \bC$ is linearly ordered and each map $\varphi_{i,j}$ is a monomorphism.

\end{itemize}

\chapter{Introduction}

The aim of this work is to build the foundations for a theory of analytic geometry over a non-trivially valued field, uniformly encompassing the case 
in which the base field is equipped with a non-archimedean valuation and that in which it has an archimedean one. The main tools we are going to use are the dagger affinoid algebras (or algebras of germs of analytic functions), analogous to the dagger affinoid
 algebras defined by Grosse-Kl\"onne \cite{GK}, and the theory of bornological vector spaces.

Historically, one of the first attempts to develop a theory of analytic spaces over non-archimedean valued fields in a systematic way, in particular analytic manifolds over such fields, was made by Serre in his notes on Lie algebras and Lie groups, \cite{SERRE}. 
The idea was to mimic the theory of analytic spaces over the field of complex numbers, but it was soon clear that this point of view did not give the expected 
results because of the pathological
topologies non-archimedean valuations induces on fields. Some years later, a new approach by John Tate, first published in \cite{TATE}, overcame these problems with 
new and highly original ideas. His theory has the power to give a non-archimedean version of a big part of the main results of complex analytic geometry and it
has also proved to be so much useful for solving problems in arithmetic geometry that now it is a fundamental tool for any arithmetic geometer. 
But beside these good features, the theory of rigid analytic spaces introduces an awkward asymmetry
between two theories: archimedean analytic geometry and non-archimedean analytic geometry. These theories start from (apparently) very different foundation, use different techniques and they arrive to work out two parallel theories with analogous main results as for example the GAGA principle and direct image theorem. In this work we will see that these differences can disappear to an accurate inspection and with the use germs of analytic functions on closed polydisks as starting point.

Let's see in detail the contents of the chapters. In the first chapter we introduce the concept of a bornology on a set: a bornology
on a set $X$ is a collection of subsets $\cB \subset \sP(X)$ which is an ideal of the boolean algebra $\sP(X)$ and that covers $X$, cf. definition \ref{defn_bornology}. 
We will give a structure of category to the class of bornological sets by defining bounded maps as the set-theoretical functions which
map all the bounded subsets of the domain to bounded subsets of the codomain.
Then, we study the basic properties of algebraic bornological structures: bornological groups, bornological rings and bornological modules which
we define as algebraic structures over bornological sets for which the structural maps are required to be bounded maps.
The main results of this chapter are the proof that the category of abelian bornological groups 
is a quasi-abelian category (proposition \ref{prop_quasi_ab}), generalizing the result of Prosman and Schneiders given in 
\cite{PS} for bornological vector spaces over $\C$, and the introduction of the notion of convexity over bornological rings and modules (definition \ref{convbal}). 
In order to define convexity we attach to the submonoid of power-bounded 
elements $A^\circ \subset A$ of a bornological ring a set $\cR^n(A^\circ) \subset A^n$, for any $n \in \N$, which encodes the information of how the monoid $A^\circ$ is embedded in $A$. It turns out that the collection $\{  \cR^n(A^\circ) \}_{n \in \N}$ defines in a natural way a
``generalized ring'' in the sense of Durov's thesis, cf. \cite{DUR}. In this context absolutely convex subsets of an $A$-module are interpreted as
$\Sigma_{A^\circ}$-submodules (where $\Sigma_{A^\circ}$ is the monad, \ie the generalized ring, that we canonically associate to the immersion $A^\circ \rhook A$) in
total analogy with non-archimedean analysis where a subset of a $k$-vector space is called absolutely convex if it is a $k^\circ$-module.
As it appears, the ``generalized ring'' $\Sigma_{A^\circ}$ is a more interesting object to study than simply $A^\circ$ if one wants to develop, for example, a theory of formal 
models for archimedean dagger analytic spaces. We end this chapter by recalling the basic definitions of the Halos theory, as developed by Paugam in \cite{PAU}. Our use of Halos theory is bounded to a simplification of notation and formulas, not main results of \cite{PAU} are needed.

In the second chapter we develop the tools needed to study our dagger affinoid algebras. We start by recalling the definition of the adjoint pair 
of functors ${}^b$ and ${}^t$, given in \cite{H1} and their basic properties. The functor ${}^b$ goes from the category of locally convex
topological vector spaces over a non-trivially valued field to the category of bornological vector spaces of convex type, while ${}^t$ goes in the other direction.
These functors correspond to the classical well-known constructions called the canonical or Von Neumann bornology, for the functor ${}^b$, and
the topology of bornivorous subsets, for the functor ${}^t$.
Then, we focus on complete bornological algebras $A$ over $k$ (cf. definition \ref{complete_born_alg}) which is the
more interesting class of bornological algebras for us, since dagger affinoid algebras belong to this class which is a broad generalization of the class of Banach algebras.
We define the bornological spectrum $\cM(A)$ of a bornological algebra $A$ using a
Berkovich-like approach: points of the spectrum are equivalence classes of bounded characters from $A$ to complete valued fields or equivalently
bounded multiplicative seminorms, see definition \ref{born_spectrum}. 
The association $A \mapsto \cM(A)$ has some problem if it is considered over the full category of bornological algebras over $k$, for example there are bornological algebras with empty spectrum. We show that if $A$ is a bornological m-algebra (cf. definition \ref{defn:bornological_m}) $\cM(A)$ is functorially associated to $A$ and it is a non-empty compact Hausdorff topological space, cf. theorem \ref{thm_born_spectrum}.
After that, we define the ring of bornological strictly convergent power-series $A \lt X \gt$ for any bornological
m-algebra and we show results analogous to the classical theory for normed algebras (cf. theorems \ref{thm_born_ps_wds}, \ref{thm_dense_poly} and \ref{thm_univ_prop_born_ps}): 
\begin{itemize}
 \item $A \lt X \gt$ is a bornological m-algebra;
 \item if $A$ is complete then also $A \lt X \gt$ is complete;
 \item $A[X] \subset A \lt X \gt \subset A \ldbrack X \rdbrack$ and moreover $A[X]$ is bornologically dense in $A \lt X \gt$,
       \ie the bornological closure of $A[X]$ in $A \lt X \gt$ is equal to all $A \lt X \gt$;
 \item $A \lt X \gt$ satisfies the following universal property in the category of bornological m-algebras:
       for every bounded morphism of algebras $\phi: A \to B$, with $B$ complete, and any power-bounded element $b \in B^\circ$ there exists a unique bounded morphism $\Phi: A \lt X \gt \to B$ such that $\Phi|_A = \phi$ and $\Phi(X) = b$. 
\end{itemize}

Then, we introduce the ring of bornological overconvergent power-series $A \lt X \gt^\dagger$, and analogous properties to the ones just listed are proved for $A \lt X \gt^\dagger$ (see theorems \ref{thm:over_conv_ring}, \ref{thm:over_conv_ring_complete} and \ref{thm:over_conv_ring_dense}). We would like to point out that to 
in order to express correctly the universal property characterising $A \lt X \gt^\dagger$ we introduce the notion of weakly power-bounded elements of a bornological algebra. In this way, $A \lt X \gt^\dagger$ satisfies the same universal property of $A \lt X \gt$ providing that we substitute the concept of power-bounded element
with the concept of weakly power-bounded element. This allows us to state the universal property which characterizes dagger affinoid algebras, and
in particular the algebras of overconvergent analytic functions on polydisks, in a more satisfying way than what is stated in \cite{GK}. Moreover, if $A$ is a regular bornological algebra (cf. definition \ref{defn_born_regular})
the weakly power-bounded elements coincide with the power-bounded elements of $A$ equipped with its spectral seminorm (cf. definition \ref{spec_seminorm}), 
with respect to which $A$ is not a Banach $k$-algebra in general. 
Hence, the concept of weakly power-bounded elements permits us to study some kinds of non-complete normed algebras as if they were complete (as they are, but with respect to a finer notion of completeness). We end this chapter by showing the functoriality property of the map $A \mapsto \Sigma_{A^\ov}$ which associates to a regular bornological algebra
the algebraic monad associated to the submonoid of weakly power-bounded elements and then by studying
the relations of the classical notion of spectra on topological $k$-algebras with the notion of bornological spectrum we have introduced.

In the next chapter we start studying dagger affinoid algebras, dagger affinoid spaces and dagger affinoid subdomains. We introduce the ring of overconvergent
analytic functions over the polydisk of polyradius $1$ as the bornological algebra 
\[ W_k^n \doteq k \lt X_1, \dots, X_n \gt^\dagger \doteq \limind_{\rho > 1} T_k^n(\rho) \]
where $\rho$ is a polyradius bigger than $(1, \dots, 1)$, \ie every component is strictly bigger than $1$, and 
$T_k^n(\rho)$ is the algebra of summable power-series in the polydisk of polyradius $\rho$ (see the first paragraph of the third chapter for the 
definition of $T_k^n(\rho)$ and the last section of the first chapter for our conventions on summations), and we equip $W_k^n$ with the direct 
limit bornology.
We then recall and discuss some property of  $k \lt X_1, \dots, X_n \gt^\dagger$: the fact that is a Noetherian, factorial domain and a regular LF-space, cf.
theorem \ref{thm_W_regular}. The main result of the first section is the 
fact that any ideal $I \subset W_k^n$ is bornologically closed, cf. theorem \ref{thm_W}. We underline that it is usual to meet this property in the theory of affinoid algebras,
but this does not have an immediate archimedean analog: for example it is known that every $\R$-Banach algebra for which all principal ideals are closed must be a division algebra. 
Hence, every Banach algebra $k \lt \rho_1^{-1} X_1, \dots, \rho_n^{-1} X_n \gt$ has non-closed ideals and also non-finitely generated ideals because it is known that every Noetherian $\C$-Banach algebra must be finite dimensional as $\C$-vector space. 
Beside these facts, $\C \lt X_1, \dots, X_n \gt^\dagger$ is Noetherian and all its ideals are bornologically closed. Although this might looks strange for the reasons explained so far, it becomes reasonable once one notices the bornological isomorphism $\C \lt X_1, \dots, X_n \gt^\dagger \cong \underset{\rho > 1}\limind \cO(\D(\rho^-))$ (where $\cO(\D(\rho^-))$ is the Stein algebra of holomorphic functions in the open polycylinder $\D(\rho^-)$), when $\cO(\D(\rho^-))$ is equipped with its canonical structure of Fr\'echet algebra. It is a classical result of complex analytic geometry that all finitely generated ideals of $\cO(\D(\rho^-))$ are closed (we are going to study in details this isomorphism in the third section of the third chapter).
We also notice that $W_k^n$ is characterized by a universal property which is the same for any $k$, archimedean or not, fixing the asymmetry that one finds by the lacking of an archimedean analogue of affinoid algebras theory (see remark \ref{rmk_universal} for a more detailed discussion of this topic). 
In the second section of the third chapter we define the category of dagger affinoid algebras as
the subcategory of the category of bornological algebras which are isomorphic to quotients of the form $\frac{W_k^n}{I}$. We will show that dagger affinoid algebras are regular bornological algebras (in the sense of definition \ref{defn_born_regular}) and that all algebra morphisms between dagger affinoid algebras are bounded (generalizing the results of the theory of affinoid algebras). We also show that all morphisms are strongly bounded, with respect to the ``sub-generalized ring'' $\Sigma_{A^\ov} \rhook \Sigma_A$.  We describe basic constructions like tensor products
and direct products of dagger affinoid algebras and this machinery permits us to define the ``generalized ring of fractions'' as done 
in rigid geometry, with the difference that now the universal properties must be stated in term of weakly power-bounded elements.
Next, we define the category of dagger affinoid spaces as the dual of the category of dagger affinoid algebras and we represent the objects of this category by means of the bornological spectrum that we defined in the previous chapter. We see that
in the non-archimedean case the category of strict dagger affinoid spaces is equivalent to the category of wide strict $k$-affinoid space
defined by Berkovich (proposition \ref{prop_strict_germs}), which is an already-known result, and that this equivalence gives an homeomorphism of the 
underlying topological spaces between the underlying space of the germ and the bornological spectrum of the associated dagger affinoid algebra. We also generalize this result to the case in which the dagger affinoid algebras and the germs of analytic spaces are not strict in
proposition \ref{prop_germs}. Then, we see that in the archimedean case the bornological spectrum of a dagger affinoid algebra coincides
with the set of its maximal ideals, giving an identification of archimedean dagger affinoid spaces with certain compact Stein subsets of $\C^n$. We conclude this chapter by defining the notion of dagger affinoid subdomain by requiring the usual universal property characterizing open embeddings. We check that most of the basic properties of the subdomains of classical rigid geometry still hold, and the proofs usually rely on the same ideas. We underline, again, that the machinery developed 
so far permits us to handle both the archimedean and the non-archimedean case as base field.

The next chapter is devoted to prove the Tate's acyclicity theorem for coverings made of strict dagger rational subdomains. The proof of the theorem is classical and uses the same reduction argument to the case of Laurent coverings used by Tate himself. It is quite obvious that Tate's proof can be used to prove the
acyclicity of the structural sheaf for dagger affinoid spaces when $k$ is non-archimedean. It is not so obvious that the same argument works with
$k$ archimedean. So, this chapter is devoted to check that this is true, although the arguments used in this chapter are mainly the same that can be found in chapter 8 of \cite{BGR}.

The fifth chapter contains the main structural result in the theory of dagger affinoid spaces: the Gerritzen-Grauert theorem. 
We start this chapter by studying in detail the notion of dagger affinoid subdomain when $k$ is archimedean. We show that in this case a subdomain $U \subset X = \cM(A)$ is a compact Stein subset that may have empty topological interior, when considered as a compact Stein subset of $\C^n$. In theorem \ref{thm_steinness} we show that the dagger affinoid algebra $A_U$ we associate to $U$ by requiring the universal property characterizing subdomains is in fact, for any rational subdomain of $X$, isomorphic to the algebra of germs of analytic functions $\wtilde{\cO}_X(U)$
that the composition of immersions $U \rhook X \rhook \C^n$ induces on $U$. This permits us to give another proof of Tate's acyclicity for rational coverings based on the Theorem B of Cartan, when $k$ is archimedean (theorem \ref{thm_arch_tate}), by exploiting the fact that the rational dagger site of $X$ is contained in what we call the 
compact Stein site of $X$: the site given by the $G$-topology for which the coverings of $X$ are given by finite coverings by compact Stein subsets.
We will then see another peculiar feature of archimedean subdomains: the failure of the transitivity property of Weierstrass subdomains, \ie 
if $U \subset V \subset X$ are two Weierstrass subdomain embeddings then, in rigid geometry, we can conclude that also the composition $U \subset X$
is a Weierstrass subdomain embedding. But the proof of this result relies on the non-archimedean character of the base field and we will see that 
in the archimedean case one can only assert that the inclusion $U \subset X$ is a rational subdomain embedding. We end the chapter with our main result in the theory
of dagger affinoid spaces: our generalization of the Gerritzen-Grauert theorem. We remark that in the non-archimedean case this result easily follows from the classical one, instead in the archimedean case one has to work out a new proof. Our proof is inspired 
by the new Temkin's proof of the Gerritzen-Grauert theorem \cite{TEM} in Berkovich geometry, and it will exploit the peculiar properties of archimedean subdomains we proved.

In the last chapter we construct the category of $k$-dagger analytic spaces. With the aid of the dagger affinoid theory developed
so far, we can use the same methods of Berkovich, from \cite{BER4}, to define our analytic spaces. 
The first section is devoted to some auxiliary results on the category of finite modules
over a dagger affinoid algebra and on the category of coherent sheaves over a dagger affinoid space. In particular we endow any finite module with a canonical filtration induced by the 
filtration of the dagger affinoid algebra over which it is defined and we shows results analog to the ones of \cite{BER2}: the category of finite modules over a dagger affinoid algebra is equivalent to the
category of finite dagger modules (proposition \ref{prop_dagger_modules}), \ie every finite module admits a unique canonical structure of complete bornological module over $A$
and every $A$-linear morphism between finite dagger modules is bounded. Then, we deduce the general version of Tate's acyclicity theorem for the weak G-topology and for any module 
over any (also non-strict) $k$-dagger algebra. We end this section with a proof of Kiehl's theorem which is an adaptation of the classical proof in rigid geometry, cf. theorem \ref{thm:kiehl}.
In the next section we introduce the notion of a $k$-dagger analytic space using the methods of Berkovich nets developed in
\cite{BER4} and we essentially follow the discussion of the first chapter of \cite{BER4} to define our spaces. Therefore, according to our definition, a $k$-dagger analytic space is a triple $(X, A, \tau)$ where $X$ is a locally Hausdorff topological space, $\tau$ a Berkovich net on $X$ and $A$ an atlas of dagger affinoid algebras for the net $\tau$. We see that we can perform
all the basic constructions of \cite{BER4}, underlining the differences of the two approaches and the relations of dagger analytic spaces with the pro-analytic structure on the 
pro-site (as defined in appendix \ref{pro_appendix}) of $X$. 
The second-last section of the sixth chapter is devoted to study the relations between our dagger analytic spaces and the other definitions of analytic spaces already present in literature.
First we deal with the non-archimedean case, in which we are mainly interested to compare our dagger analytic spaces with Grosse-Kl\"onne ones and with Berkovich's analytic spaces. Then, we show that in the complex
case, the category of classical complex analytic spaces embeds in a fully faithful way in the category of $\C$-dagger analytic spaces. In the last section we show how the notion of flatness behaves better, and more naturally, for a morphism of $k$-dagger analytic spaces with respect to a morphism of Berkovich $k$-analytic spaces.

\chapter{Bornological structures}

In this chapter we introduce bornologies and bornological algebraic structures in a general fashion. The idea of studying abstract bornological structures is not new but there is no systematic study in literature (at our knowledge) of notions like bornological groups, bornological rings, bornological modules, etc... This chapter has the aim to perform a detailed study of these structures and to fix some notations and terminology that will be used in subsequent chapters. We underline that even if this chapter is mainly of introductory nature, some of its results are new.

We start by recalling the notion of bornological set and by studying some properties of the category they form. This category has several analogies with the category of topological spaces but it seems to have no direct geometrical meaning. That is why, apparently, bornologies are not interesting to be studied outside functional analysis. As a further step, we introduce bornological groups and we prove that the category of bornological abelian groups is quasi-abelian, generalizing the analogous result of Prosmans and Schneiders \cite{PS} for bornological vector spaces. Then, we pass to bornological rings and to study a notion of convexity which can be defined over any bornological ring. We end the chapter by recalling some definitions and notations of Paugam's Halos theory that will be used later on.

\section{Bornological sets}

\index{bornological set}
\index{bornology}
\begin{defn} \label{defn_bornology}
Let $X$ be a set. A \emph{bornology} on $X$ is a collection $\cB \subset \sP(X)$ such that
\begin{enumerate}
\item $\cB$ is a covering of $X$, \ie $\forall x \in X, \exists B \in \cB$ such that $x \in B$;
\item $\cB$ is stable under inclusions, \ie $A \subset B \in \cB \then A \in \cB$;
\item $\cB$ is stable under finite unions, \ie for each $n \in \N$ and $B_1, ..., B_n \in \cB$, $\underset{i = 1}{\overset{n}\bigcup} B_i \in \cB$.
\end{enumerate}

A pair $(X, \cB)$ is called a \emph{bornological set}, and the elements of $\cB$ are called \emph{bounded subsets} of $X$ 
(with respect to $\cB$, if a specification is needed).  
\end{defn}

We recall that $\sP(X)$ has a canonical structure of boolean algebra with respect to the operations $\cup, \cap$ and the complementation $Y \mapsto \ol{Y} = X - Y$,
for $Y \in \sP(X)$. Thus the definition of a bornology $\cB$ on $X$ is equivalent to the request that the collection $\cB \subset \sP(X)$ is an ideal of $\sP(X)$ 
that contains the ideal $\sF(X)$ defined by finite subsets of $X$. So, the family of bornologies over a set $X$ is in a one-to-one correspondence 
with the ideals of the boolean algebra $\sP(X)/\sF(X)$ and forms a complete lattice with respect to unions (more precisely, the ideal generated by the unions) and intersections. 

\index{bornological basis}
\index{bornological pre-basis}
For any collection $\cS \subset \sP(X)$ there is a minimal bornology $\cB$ on $X$ such that $\cS \subset \cB$, that is the smallest ideal of $\sP(X)$ 
that contains the ideal generated by $\cS$ and $\sF(X)$. This bornology is called the bornology \emph{generated} by $\cS$ and is denoted $\langle \cS \rangle$, 
while $\cS$ is said to be a \emph{bornological pre-basis} of $\langle \cS \rangle$. A family $\cS \subset \sP(X)$ is called a \emph{bornological basis} of 
$\langle \cS \rangle$ if it is a covering of $X$ and for any $B_1, B_2 \in \lt \cS \gt$ there exists $B \in \cS$ such that $B_1 \cup B_2 \subset B$.

\begin{exa}
\begin{enumerate}
\index{trivial bornology}
\index{chaotic bornology}
\item Each set $X$ has two distinguished bornologies, which coincide if and only if $X$ is finite. They are the maximum and the minimum elements 
      of the lattice of the family of bornologies over $X$. The former is called the \emph{chaotic bornology}, for which $\cB = \sP(X)$, and the 
      latter \emph{trivial bornology}, for which $\cB = \sF(X)$.
\item The fact that $\cB$ is an ideal of the boolean algebra $\sP(X)$ implies that to check that a bornology is chaotic it is enough to know 
      that for one subset $Y \subset X$ both $Y, X-Y \in \cB$. In fact, the only ideal for which this can happen is the improper ideal.
\item The bounded subsets in a metric space form a bornology.
\end{enumerate}
\end{exa}

\index{morphism of bornological sets}

A \emph{morphism} of bornological sets $\varphi: (X, \cB_X) \to (Y, \cB_Y)$ is a bounded map $\varphi: X \to Y$, \ie a map of sets such that 
$\varphi(B) \in \cB_Y$ for all $B \in \cB_X$. To check that a map is bounded is sufficient to check that $\varphi(B) \in \cB_Y$ for $B$ varying on a basis of $\cB_X$. 
It is clear that the set theoretic composition of two bounded maps is a bounded map, so we can define the category of bornological 
sets, denoted $\bBorn$, whose objects are bornological sets and arrows are  bounded maps.

\begin{rmk} \label{weak_rmk}
If $\cB_1, \cB_2$ are two bornologies on a set $X$ then $\Id_X : (X, \cB_1) \to (X, \cB_2)$ is a bounded map if and only if
$\cB_1 \subset \cB_2$, which is the opposite of the condition that one finds for topologies over $X$. Indeed, the strongest bornology on $X$
is the trivial bornology and the weakest one is the chaotic bornology.
\end{rmk}

Let $(X, \cB)$ be a bornological set and $f:Y \to X$ a map of sets. The family $f^{-1}(\cB) \doteq \{f^{-1}(B)\}_{B \in \cB}$ is a basis for a 
bornology on $Y$ denoted $f^* \cB$; it is the maximal bornology on $Y$ such that $f$ is a bounded map. If $f$ is an injective map, 
we write $\cB|_Y = f^* \cB$ and call it the \emph{induced bornology} on $Y$. 

\index{induced bornology}

For a set theoretic inclusion $i: Y \to (X, \cB)$ the induced bornology on $Y$ has the usual universal property: In the diagram

\beqn
\begin{tikzpicture}
\matrix(m)[matrix of math nodes,
row sep=2.6em, column sep=2.8em,
text height=1.5ex, text depth=0.25ex]
{Y&X\\
Z \\};
\path[right hook->,font=\scriptsize]
(m-1-1) edge node[auto] {$i$} (m-1-2);
\path[->,font=\scriptsize]
(m-2-1) edge node[auto] {$f$} (m-1-1);
\path[->,font=\scriptsize]
(m-2-1) edge node[below] {$i \circ f$} (m-1-2);
\end{tikzpicture}
\eeqn

$f$ is bounded if and only if $i \circ f$ is bounded. In fact, if $f$ is bounded then $i \circ f$ is bounded because it is a composition of two bounded morphisms. Now suppose that $f$ is not bounded, then $f$ must map a bounded subset of $Z$ on a subset of $Y$ which is not bounded with respect to the bornology induced on $Y$ by $X$, hence also $i \circ f$ will map the same subset of $Z$ to a subset of $X$ which is not bounded. Therefore $i \circ f$ is not bounded.
This universal property characterizes the induced bornology uniquely up to isomorphism.

\index{regular monomorphism between bornological sets}
\begin{defn} \label{defn:regular_mono}
We say that a monomorphism $(Y, \cB_Y) \to (X, \cB_X)$ is \emph{regular} if $Y$ has the bornology induced by $X$.
A subobject of $(X, \cB_X)$ is an equivalence class of regular monomorphisms.
\end{defn}

In a dual way, if $g: X \to Z$ is another map of sets, and $\cB$ a bornology on $X$, the family $g(\cB) \doteq \{ g(B) \}_{B \in \cB}$ is a pre-basis for a bornology
$g_* \cB \doteq \lt g(\cB) \gt$ on $Z$, minimal among the bornologies on $Z$ for which $g$ is a bounded map. 
If $g$ is surjective, we call $g_*(\cB)$ the \emph{quotient bornology} of $\cB$ under $g$ and moreover in this case $g_*(\cB) = g(\cB)$.

\index{quotient bornology}

For a surjection of sets $\pi: (X, \cB) \to Z$ the quotient bornology on $Z$ has the following universal property: In the diagram

\beqn
\begin{tikzpicture}
\matrix(m)[matrix of math nodes,
row sep=2.6em, column sep=2.8em,
text height=1.5ex, text depth=0.25ex]
{X&\\
Z& Y\\};
\path[->>,font=\scriptsize]
(m-1-1) edge node[auto] {$\pi$} (m-2-1);
\path[->,font=\scriptsize]
(m-2-1) edge node[auto] {$f$} (m-2-2);
\path[->,font=\scriptsize]
(m-1-1) edge node[auto] {$f \circ \pi$} (m-2-2);
\end{tikzpicture}
\eeqn

$f$ is bounded if and only if $f \circ \pi$ is bounded. Indeed, if $f$ is bounded then $f \circ \pi$ is bounded. 
Instead, if $f$ is not bounded then there exists a bounded subset for the quotient bornology $B \subset Z$ such that $f(B)$ is not bounded. 
But by definition there exists a bounded subset $B' \subset X$ such that $\pi(B') = B$, so $f \circ \pi$ is not bounded.

\index{regular epimorphism between bornological sets}
\begin{defn} \label{defn:regular_epi}
We say that an epimorphism $f: (X, \cB_X) \to (Y, \cB_Y)$ is \emph{regular} if $Y$ has the quotient bornology given by $f$.
A quotient object of $(X, \cB_X)$ is an equivalence class of regular epimorphisms.
\end{defn}

\index{product bornology}
The category of bornological sets admits small limits, that we are going to describe now. Let $(X_i, \cB_i)_{i \in I}$ be a small family of 
bornological sets, then the direct product $\underset{i \in I}\prod (X_i, \cB_i)$ is defined to be the set $\underset{i \in I}\prod X_i$ equipped with 
the bornology generated by the family
\beqn
\l \{ \prod_{i \in I} B_i | B_i \in \cB_i, \forall i \in I \r\}.
\eeqn
This bornology is the biggest bornology on $\underset{i \in I}\prod X_i$ such that all the projections 
\beqn
\pi_j: \prod_{i \in I} X_i \to X_j
\eeqn
are bounded maps.

Moreover, given two bounded maps of bornological sets $\varphi: (X, \cB_X) \to (Z, \cB_Z), \psi: (Y, \cB_Y) \to (Z, \cB_Z)$ the fiber product 
$(X, \cB_X) \times_{(Z, \cB_Z)} (Y, \cB_Y)$ is defined as the fiber product of sets $X \times_Z Y = \{ (x,y) \in X \times Y | \varphi(x) = \psi(y) \}$, 
equipped with the bornology generated by the family
\beqn
\l \{ B \times_Z C | B \in \cB_X, C \in \cB_Y \r \}.
\eeqn
Hence $(X, \cB_X) \times_{(Z, \cB_Z)} (Y, \cB_Y)$ is equipped with the bornology induced by the inclusion in $(X, \cB_X) \times (Y, \cB_Y)$. 
Notice that, as in the case of topological spaces, the bornology on $(X, \cB_X) \times_{(Z, \cB_Z)} (Y, \cB_Y)$ does not depend on 
the bornology on $Z$.

\index{coproduct bornology}
For coproducts, let $(X_i, \cB_i)_{i \in I}$ be a small family of bornological sets, then the 
coproduct $\underset{i \in I}\coprod (X_i, \cB_i)$ is defined to be the set $\underset{i \in I}\coprod X_i$ equipped with the bornology generated by the family
\beqn
\l \{ \coprod_{i \in I} B_i | B_i \in \cB_i, \text{ and } B_i = \varnothing, \forall \text{ but a finite number of } i \r\}.
\eeqn
This bornology is the smallest bornology on $\underset{i \in I}\coprod X_i$ such that all the injections 
\beqn
\a_j: X_j \to \coprod_{i \in I} X_i
\eeqn
are bounded maps.

We can describe pushouts considering two bounded maps $\varphi: (X, \cB_X) \to (Y, \cB_Y), \psi: (X, \cB_X) \to (Z, \cB_Z)$ and giving to the 
set theoretic pushout $Y \coprod_X Z$ the bornology generated by the family
\beqn
\l \{ C \coprod_{X} D | C \in \cB_Y, D \in \cB_Z \r\}.
\eeqn
Hence $(Y, \cB_Y) \coprod_{(X, \cB_X)} (Z, \cB_Z)$ is equipped with the bornology induced by the surjection from $(Y, \cB_Y) \coprod (Z, \cB_Z)$. 
Notice that, as in the case of topological spaces, the bornology on $(Y, \cB_Y) \coprod_{(X, \cB_X)} (Z, \cB_Z)$ does not depend on 
the bornology on $X$.

So, we can summarize what we said with the next proposition.

\begin{prop}
The category of bornological sets is complete and cocomplete and the forgetful functor $U: \bBorn \to \bSet$, to the category of sets, 
commutes with limits and colimits.
\end{prop}

The next proposition explains our choice of names in definitions \ref{defn:regular_mono} and \ref{defn:regular_epi}.

\begin{prop}
Let $f: A \to B$ be a morphism in $\bBorn$, then 
\begin{itemize}
\item $f$ is a regular monomorphism in the sense of definition \ref{defn:regular_mono} if and only if it is regular in the sense of category theory;
\item $f$ is a regular epimorphism in the sense of definition \ref{defn:regular_epi} if and only if it is regular in the sense of category theory.
\end{itemize}
Moreover, in $\bBorn$ the categorical notions of regular and strict monomorphism (and respectively regular and strict epimorphism)
coincide.
\end{prop}
\begin{proof}

In proposition 5.1.5 of \cite{KAS} one can find the categorical definitions of strict and regular epimorphism. More precisely, condition (1)(d) of proposition 5.1.5 of ibid. defines regular epimorphisms and condition (1)(e) of proposition 5.1.5 of ibid. defines strict epimorphisms. $\bBorn_k$ satisfies the condition of that proposition, therefore an epimorphism between bornological sets is regular (in categorical sense) if and only if it is strict. Then, proposition 5.1.5 of \cite{KAS} has a clear dual statement which settles the case of monomorphisms. 

Then, consider a regular monomorphism (in the sense of definition \ref{defn:regular_mono}) of bornological sets $i: (X, \cB_X) \rhook (Y, \cB_Y)$ and consider the push-out of the diagram
\beqn
\begin{tikzpicture}
\matrix(m)[matrix of math nodes,
row sep=2.6em, column sep=2.8em,
text height=1.5ex, text depth=0.25ex]
{X & Y\\
 Y & \\};
\path[->,font=\scriptsize]
(m-1-1) edge node[auto] {$i$} (m-2-1);
\path[->,font=\scriptsize]
(m-1-1) edge node[auto] {$i$} (m-1-2);
\end{tikzpicture}.
\eeqn
Then, $(X, \cB_X) \cong \ker(Y \rightrightarrows Y \coprod_X Y)$, hence $i$ is a regular monomorphism. The dual argument works for regular epimorphisms.

\end{proof}

\begin{prop}
$\bBorn$ is a concrete topological category.
\end{prop}
\begin{proof}
A concrete category $U: \bBorn \to \bSet$ is said topological if $U$ is a topological functor\footnote[1]{Here we use the notion of topological functor of \cite{ACC}, definition 21.1, but in literature the term topological functor is sometimes used with different meanings.}. This means that for any diagram of the form
\[ (Y \stackrel{f_i}{\to} U X_i)_{i \in I} \]
there exists a unique initial lift 
\[ X \stackrel{\ol{f}_i}{\to} X_i. \]
In the category of topological spaces this universal problem is solved by the weak topology defined by the maps $f_i$, and in the category
of bornological sets we can give the same description of the ``weak bornology'' which solves the same problem. Consider the family of subsets of $Y$ given by the intersection
\[ \bigcap_{i \in I} f_i^*(\cB_i) \subset \sP(Y) \]
of all the bornologies induced by the $f_i: Y \to X_i$. This family is a bornology on $Y$ and it is the biggest bornology (hence the weakest, 
cf. remark \ref{weak_rmk}) for which all the maps $f_i$ are bounded. 
To give a lift of $Y$ is the same to give a bornology on $Y$. If $X'$ is a lift of $Y$ such that
\[ X' \stackrel{\ol{f}_i'}{\to} X_i \]
is well defined, \ie all maps are bounded, then the bornology of $X'$ is smaller than the bornology of $X$ and so the identity map
\[ X' \to X \]
is bounded, showing that $X$ is the unique initial lift.
\end{proof}

We conclude this section introducing an internal hom functor on $\bBorn$.

\begin{defn}
Let $(X, \cB_X)$ and $(Y, \cB_Y)$ be two bornological sets. We define $\ul{\Hom}_{\bBorn}(X, Y)$ to be the set $\Hom_{\bBorn}(X, Y)$ equipped with the 
\emph{equiboundedness bornology}, \ie $H \subset \Hom_{\bBorn}(X, Y)$ is bounded if for any bounded subset $B \in \cB_X$ the set
\[  H(B) = \bigcup_{h \in H} h(B) \]
is bounded in $(Y, \cB_Y)$.
\end{defn}

\begin{prop} \label{prop:f_1_born}
$\bBorn$ equipped with the internal hom functor given by the equiboundedness bornology and the direct product of bornological sets is a closed symmetric monoidal category.
\end{prop}
\begin{proof}
The only non-trivial condition to check is the bijection
\[ \Hom_\bBorn(X \times Y, Z) \cong \Hom_\bBorn(X, \ul{\Hom}_\bBorn(Y, Z)) \]
for any triple of bornological sets $(X, \cB_X), (Y, \cB_Y), (Z, \cB_Z)$. Let $f: X \times Y \to Z$ be a bounded map. This map induces a map of sets $\widehat{f}: X \to \Hom_\bBorn(Y, Z)$
by the currying operation and we have to check that $\widehat{f}$ is a bounded map when we equip the target with the equiboundednees bornology. Let $B_X \in \cB_X$ then
\[ \widehat{f}(B_X) = \{ f(b, \cdot) | b \in B_X \} \]
and for any $B_Y \in \cB_Y$,
\[ (\widehat{f}(B_X))(B_Y) = \bigcup_{b \in B_X} f(b, B_Y) = f(B_X, B_Y) \in \cB_Z \]
hence $\widehat{f}$ is a bounded map. It is clear that the previous reasoning can be reversed, therefore starting from any bounded map $\widehat{f}$ we can recover uniquely the bounded map $f$ from which $\widehat{f}$ comes from.
\end{proof}

\begin{rmk}
Last proposition marks a fundamental difference between $\bBorn$ and the category of topological space $\bTop$. It is well known 
that $\bTop$ cannot be equipped with a structure of a closed monoidal category.
\end{rmk}

\begin{exa}
\begin{enumerate}
	\item Applying to $\bBorn$ the machinery of relative algebraic geometry of T\"oen-Vezzosi we can define the category of bornological $\F_1$-schemes. In the same way affine
	      $\F_1$-schemes can be seen as the dual category of the category of commutative monoids, which is the category of commutative algebras (in the sense of tensor categories) of $\bSet$, the category
				of affine bornological $\F_1$-schemes can be defined as the dual category of bornological monoids.
	\item We will see that $\bBorn(\bAb)$, the category of bornological abelian groups, is quasi-abelian and proposition \ref{prop:f_1_born} easily implies that the bornology of 
	      equiboundedness induces a structure of closed symmetric monoidal category also on $\bBorn(\bAb)$. This implies that the category of bornological commutative rings can be seen as the category of commutative monoids of $\bBorn(\bAb)$, allowing the application of the machinery of \cite{BEKR} to its study. We remark that also in this case the category $\bTop(\bAb)$ does not share these properties: although $\bTop(\bAb)$ is quasi-abelian closed category it has no structure which renders it a closed symmetric monoidal category.
\end{enumerate}
\end{exa}

\section{Bornological groups}

\index{bornological group}
\begin{defn} \label{defn_bornological_group}
A group $(G, \cdot)$ is said to be a \emph{bornological group} if the set $G$ is equipped with a bornology such that the map 
$(g,h) \mapsto g \cdot h$ (considering on $G \times G$ the product bornology) and the map $g \mapsto g^{-1}$ are bounded. If $(G, \cB)$ is a bornological group, we also say that $\cB$ 
is a \emph{group bornology}. A morphism of bornological groups is a group homomorphism which is also a bounded map.
\end{defn}

\begin{exa}

\begin{enumerate}
\item On $\Q$ each valuation, archimedean $|\cdot|_{\infty}$ or non-archimedean $|\cdot|_p$, give rise to a bornology on the additive group $(\Q, +)$. 
A basis for the bornology is the set of closed disks 
\[ D_{|\cdot|_*}(0, q^+) = \{x \in \Q | |x|_* \le q \} \]
with $* \in \{\infty, 2, 3, 5, ...\}$.

\item The additive group ($\Z$, +) can be equipped with the bornology generated by the half-lines $\{n \in \Z | n \le a\}$, for $a \in \Z$. 
This bornology is not a group bornology because the map $x \mapsto -x$ is not bounded.

\item Consider the multiplicative group $(\Ct, \cdot)$, endowed with the bornology induced by the usual bornology on $\C$. 
This is not a group bornology on $\Ct$, because the inversion is not bounded. The standard way to make $\Ct$ into a bornological group is to equip it with the 
bornology induced by the embedding $(id_\C, inv_\C): \Ct \to \C \times \C$ (where $inv_\C(x) = x^{-1}$), where the target space is equipped with the product bornology.

\end{enumerate}
\end{exa}

Let $(G, \cB)$ be a bornological group and let $f: H \to G$ be a homomorphism of groups, then $(H, f^* \cB)$ is a bornological group. 
In fact, if $B_1, B_2 \in f^* \cB$, then $B_1 \subset f^{-1}(B_1^G)$ and $B_2 \subset f^{-1}(B_2^G)$ for some $B_1^G, B_2^G \in \cB$.
So 
\[ B_1 + B_2 \subset f^{-1}(B_1^G) + f^{-1}(B_2^G) = f^{-1}(B_1^G + B_2^G) \]
is a bounded subset.

On the other hand, given a morphism $f: G \to H$ then $f_* \cB$ is a group bornology on $H$ if $f$ is surjective. Indeed, in this case, 
given $B_1, B_2 \in f_* \cB$ then there exist $B_1^G, B_2^G \in G$ such that $f(B_1^G) = B_1$ and $f(B_2^G) = B_2$ so
\[ f(B_1^G + B_2^G) = f(B_1^G) + f(B_2^G) = B_1 + B_2 \]
which is bounded. But this argument works only if $f$ is surjective, and it is easy to construct examples where $f_* \cB$ fails to be a group
bornology (for example take the immersion $\R \to \R^2$ with the canonical metric bornology on $\R$).

In any case there exists a smallest group bornology which contains $f_* \cB$, thanks to the next lemma.

\begin{lemma} \label{borngrop_int}
Let $G$ be a group and $(\cB_i)_{i \in I}$ be any family of group bornologies on $G$, then $\underset{i \in I} \bigcap \cB_i$ is a group bornology on $G$.
\end{lemma}
\begin{proof}
We know that $\cB = \underset{i \in I}\bigcap \cB_i$ is a bornology, we must show that it is a group bornology. Let $B, C \in \cB$, then 
$B, C \in \cB_i$ for any $i \in I$ and so
\beqn
B + C \in \cB_i
\eeqn
because $\cB_i$ are group bornologies. Therefore $B + C \in \cB$.

To check that the map $x \mapsto - x$ is bounded, we can proceed in a similar fashion.
\end{proof}

\begin{lemma} \label{borngrop_char}
Let $G$ be a group and $\cB$ a bornology on $G$. Then $\cB$ is a group bornology if and only if there exists a basis $\{B_i\}_{i \in I}$ 
for $\cB$ such that
\[ B_i + B_j \subset B_{i_1} \cup \dots \cup B_{i_n} \]
for any $i,j \in I$ and 
\[ -B_i \subset B_{j_1} \cup \dots \cup B_{j_m} \]
for any $i \in I$ and some $i_1, \dots, i_n, j_1, \dots, j_m \in I$.
\end{lemma}
\begin{proof}
If $\cB$ is a group bornology then $B_i + B_j$ is bounded and hence must be contained in a finite union of the elements of the basis and 
the same is true for $-B_i$.

On the other hand, if $\cB$ has such a basis then for any couple of bounded subsets $R_1, R_2 \in \cB$ we have that $R_1 \subset B_{i_1}$
and $R_2 \subset B_{i_2}$ for some $i_1, i_2 \in I$. Then 
\[ -R_1 \subset -B_{i_1} \]
and
\[ R_1 + R_2 \subset B_{i_1} + B_{i_2} \]
which shows that multiplication and inversion are bounded maps.
\end{proof}

We denote $\bBorn(\bAb)$ the category of bornological abelian groups. 
It is easy to describe products and coproducts for any small family $(G_i, \cB_i)_{i \in I} \subset \ob(\bBorn(\bAb))$; indeed the group $\underset{i \in I}\prod G_i$ endowed with the bornology defined by the family

\beqn
\prod_{i \in I} \cB_i = \l \{ B  | B \subset \prod_{i \in I} B_i, B_i \in \cB_i \r \}
\eeqn
is the product of $(G_i, \cB_i)_{i \in I}$ in $\bBorn(\bAb)$. Whereas endowing $\underset{i \in I}\bigoplus G_i$ with the bornology
\beqn
\bigoplus_{i \in I} \cB_i = \l \{ B | B \subset \bigoplus_{i \in I} B_i, B_i \in \cB_i, B_i = 0 \mbox{ for almost all } i \in I \r \}
\eeqn
we get the coproducts in $\bBorn(\bAb)$.

In these settings we have that the canonical projections 
\beqn
\pi_j: \prod_{i \in I} G_i \to G_j
\eeqn
and the canonical injections
\beqn
i_j: G_j \to \bigoplus_{i \in I} G_i
\eeqn
are morphisms of bornological groups. So, the next proposition  permits us to calculate any limit and colimit in $\bBorn(\bAb)$. In next pages we follow very closely the discussion of Prosmans-Schneiders \cite{PS} on bornological vector spaces, adapting it to bornological abelian groups.

\begin{prop} \label{prop:kernel}
The category $\bBorn(\bAb)$ is additive and moreover, if $f: G \to H$ is a morphism of bornological groups, then
\ben
\item $\Ker f$ is the subgroup $f^{-1}(0)$ of $G$ endowed with the induced bornology;
\item $\Coker f$ is the quotient $H/f(G)$ endowed with the quotient bornology;
\item $\Im f$ is the subgroup $f(G)$ of $H$ endowed with the induced bornology;
\item $\Coim f$ is the quotient $G/f^{-1}(0)$ endowed with the quotient bornology.
\een
\end{prop}
\begin{proof}
$\bBorn(\bAb)$ is additive because for any $g_1,g_2 \in \Hom(G, H)$, $h_1 \in \Hom(H_2, G)$, $ h_2 \in \Hom(H, H_1) $ and $B$ a bounded subset of $G$
\beqn
(g_1+g_2)(B) = \l \{g_1(x) + g_2(x) | x \in B \r \} \subset g_1(B) + g_2(B) 
\eeqn
is bounded in $H$ and the equalities $(g_1+g_2) \circ h_1 = g_1 \circ h_1 + g_2 \circ h_1$, $h_2 \circ (g_1+g_2) = h_2 \circ g_1 + h_2 \circ g_2$ 
hold because they hold for abelian groups and all maps are bounded.

The other statements have similar proofs, so we only give explicit proof of the first one as an example. Let $f:G \to H$ be a bounded homomorphism of bornological groups
and $\Ker f$ be the kernel of $f$ in the algebraic sense. Then, for any
group homomorphism $\phi: G' \to G$ such that $f \circ \phi = 0$ there exists a unique map $\psi: G' \to \Ker(f)$ such that the diagram
\[
\begin{tikzpicture}
\matrix(m)[matrix of math nodes,
row sep=2.6em, column sep=2.8em,
text height=1.5ex, text depth=0.25ex]
{G' & G & H \\
 \Ker(f) & \\};
\path[->,font=\scriptsize]
(m-1-1) edge node[auto] {$\phi$} (m-1-2);
\path[->,font=\scriptsize]
(m-1-1) edge node[auto] {$\psi$} (m-2-1);
\path[->,font=\scriptsize]
(m-2-1) edge node[auto] {$k$} (m-1-2);
\path[->,font=\scriptsize]
(m-1-2) edge node[auto] {$f$} (m-1-3);
\end{tikzpicture}
\]
commutes, where $k$ is the canonical morphism. Now, if we ask that all maps are bounded we see that the request that 
$k \circ \psi$ is bounded for any bounded homomorphism $\phi$ is satisfied precisely when $\Ker f$ has the bornology induced by the inclusion 
$\Ker f \rhook G$ because of the characterization of regular monomorphism of bornological sets we gave so far, before definition \ref{defn:regular_mono}.
\end{proof}

Thus, $\bBorn(\bAb)$ is complete and cocomplete and the forgetful functor $\bBorn(\bAb) \to \bAb$, from the category of bornological abelian groups to the category of abelian groups, commutes with all limits and colimits.

\index{strict morphism}
We recall that, in general, a morphism $f: G \to H$ in an additive category with kernels and cokernels is called 
\emph{strict} if the induced morphism $\overline{f}: \Coim u \to \Im u$  is an isomorphism.

\begin{prop}
A morphism of bornological abelian groups $f: G \to H$ is strict if and only if for any bounded subset $B$ of $H$ there exists a bounded subset $B'$ of $G$ such that
\beqn
B \cap f(G) = f(B').
\eeqn
\begin{proof}
Thanks to proposition \ref{prop:kernel} we know that $f$ is strict if and only if the morphism
\beqn
 \overline{f}: G/f^{-1}(0) \to f(G) 
\eeqn
is an isomorphism in $\bBorn(\bAb)$. $\ol{f}$ is an algebraic isomorphism, so it is an isomorphism of bornological groups if and only if it identifies the 
quotient bornology on $G/f^{-1}(0)$ with the bornology induced on $f(G)$ by its inclusion in $H$. This happens precisely when the conditions of the proposition
are verified because, by definition, each bounded subset in $G/f^{-1}(0)$ is the image of a bounded subset $B' \subset G$ and any bounded subset in $f(G)$ is of the 
form $B \cap f(G)$ for a bounded subset $B \subset H$.
\end{proof}
\end{prop}

\begin{cor} \label{stric_epi_mono}
Let $f: G \to H$ be a morphism of bornological abelian groups. $f$ is a strict monomorphism (resp. a strict epimorphism) if and only if it is a regular 
monomorphism (resp. regular epimorphism) of the underlying bornological sets.
\end{cor}

\index{quasi-abelian category}
\begin{defn}
An additive category is called \emph{quasi-abelian} if it admits kernels and cokernels, and if the class of strict epimorphisms (resp. strict monomorphisms)
 is stable under pullbacks (resp. pushouts). This means that given a cartesian square
\beqn
\begin{CD}
P @>v>> B\\
@VVfV @VVgV\\
A @>u>> C
\end{CD}
\eeqn
where $u$ is a strict epimorphism then $v$ is a strict epimorphism. In a dual way for monomorphisms and cocartesian squares.
\end{defn}

We record the following fact for the clearness sake.

\begin{prop} 
Let $G, H$ be two bornological groups and $G \oplus H$ be their direct sum, then the canonical maps 
\[ \pi_G: G \oplus H \to G, \ \ \pi_H: G \oplus H \to H \]
and
\[ i_G: G \to G \oplus H, \ \ i_H: H \to G \oplus H \]
are strict morphisms.
\end{prop}
\begin{proof}
By corollary \ref{stric_epi_mono}, for a map to be a strict mono(resp. epi)morphism is equivalent to be a regular mono(resp. epi)morphism as map of the underlying bornological sets. Hence, the proposition is readily deduced by the description of the bornology of $G \oplus H$ given so far.
\end{proof}

\begin{prop} \label{prop_quasi_ab}
The category $\bBorn(\bAb)$ is quasi-abelian.
\end{prop}
\begin{proof}
We have only to show that pullbacks and pushouts preserve, respectively, strict epimorphisms and strict monomorphisms.

Consider first the cartesian square
\beqn
\begin{CD}
P @>v>> B\\
@VVfV @VVgV\\
A @>u>> C
\end{CD}
\eeqn
where $u$ is a strict epimorphism. By the universal properties of the direct product and of the pullback if we define the map
\beqn
\a = (u, -g): A \oplus B \to C
\eeqn
then we have that
\beqn
P \cong \Ker \a = \{(a, b) \in A \oplus B | u(a) = g(b) \}.
\eeqn
We also have the canonical epimorphisms $\pi_A: A \oplus B \to A, \pi_B: A \oplus B \to B$, the canonical monomorphism $i_\a: \Ker \a \to A \oplus B$ and we can write
\beqn
f = \pi_A \circ i_\a, v = \pi_B \circ i_\a.
\eeqn
The map $v$ is surjective, as can be deduced, for example, applying proposition 11.18 of \cite{ACC} and the fact that the forgetful functor $\bBorn(\bAb) \to \bSet$ commutes with limits; we have to show that $v$ is strict. Let $B_1 \subset B$ be a bounded subset, so $g(B_1)$ is bounded and since $u$ is 
strict we can find a bounded subset $B_2 \subset A$ such that $u(B_2) = g(B_1)$. Hence, $\pi_A^{-1}(B_2) \cap \pi_B^{-1}(B_1) = B_2 \times B_1$ 
is a bounded subset of $A \oplus B$, therefore
\beqn
P \cap\pi_A^{-1}(B_2) \cap \pi_B^{-1}(B_1) = P \cap (B_2 \times B_1)
\eeqn
is a bounded subset of $P$. Moreover, since $u(B_2) = g(B_1)$ we have that
\[ B_2 \times B_1 \subset P, \]
thus
\beqn
v ( B_2 \times B_1 ) = \pi_B \circ i_\a (B_2 \times B_1 ) = \pi_B(B_2 \times B_1) = B_1,
\eeqn
proving that $v$ is strict.

Now consider the cocartesian square
\beq
\begin{CD}
C @>u>> B\\
@VVgV @VVfV\\
A @>v>> S
\end{CD}
\eeq
where $u$ is a strict monomorphism. If we set
\beqn
\a = (g, -u): C \to A \oplus B
\eeqn
then we have
\beqn
S \cong \Coker \a = \frac{A \oplus B}{\a(C)}
\eeqn
with the canonical morphisms $i_A: A \to A \oplus B$, $i_B: B \to A \oplus B$, $q_\a: A \oplus B \to S$ and
\beqn
v = q_\a \circ i_A,  f = q_\a \circ i_B.
\eeqn
$v$ is injective because 
\[ v(x) = 0 \then (q_\a \circ i_A)(x) = 0 \then q_\a((x, 0)) = 0 \then  \] 
\[ \then (x, 0) \in \a(C) \then \exists c \in C | g(c) = x, -u(c) = 0 \then c = 0 \then x = 0. \]
We have to show that $v$ is strict. Let $R \subset S$ be a bounded subset and $R_1 \subset A \oplus B$ another one such that $q_\a(R_1) = R$. 
We must show that $v^{-1}(q_\a(R_1))$ is bounded in $A$. By definition we can find 
$R_A \subset A$ and $R_B \subset B$ such that $R_1 \subset R_A \times R_B$, hence
\[ v^{-1}(q_\a(R_1)) \subset v^{-1}(q_\a(R_A \times R_B)) =  v^{-1}( q_\a(R_A \times \{0\}) +  q_\a(\{0\} \times R_B)) = \]
\[  = v^{-1}( q_\a \circ i_A(R_A)) + v^{-1} (q_\a \circ i_B (R_B)) \]
because $R_A \times R_B = (R_A \times \{0\}) + (\{0\} \times R_B)$ and $v^{-1}(R + R') = v^{-1}(R) + v^{-1}(R')$, and
\beqn
 v^{-1} ( q_\a \circ i_A(R_A)) + v^{-1} ( q_\a \circ i_B (R_B)) = v^{-1} ( v (R_A)) + v^{-1} ( f (R_B)).
\eeqn
Since $u$ is a monomorphism it has a set-theoretical left inverse $u^{-1}$ and 
\[ u^{-1}(R_B) = R_B \cap C \]
which is bounded in $C$ because $u$ is a strict monomorphism and so $C$ can be (bornologically) identified with $\Im(C) \subset B$.
It is enough to show that $v^{-1}(f(R_B)) \subset g(u^{-1}(R_B))$, because $g$ is a bounded map. So, since $v$ is injective
\[ x \in v^{-1}(f(R_B)) \iff x \in f(R_B) \cap \Im(A) \iff x \in f(R_B) \cap \Im(v \circ g) \iff x \in \Im(f \circ u). \]
Indeed, the inclusion $f(R_B) \cap \Im(v \circ g) \subset f(R_B) \cap \Im(A)$ is obvious. So, if $x \in f(R_B) \cap \Im(A)$ it means that there exist $a \in A$ and $b \in B$ such that $f(b) = x$ and $v(a) = x$, hence
\[ f(b) = (q_\a \circ i_B)(b) = q_\a((0, b)) = x = v(a) = (q_\a \circ i_A)(a) = q_\a((a, 0)) \]
hence
\[ q_\a((a, 0)) - q_\a((0, b)) = q_\a((a,-b)) = 0 \then (a, -b) \in \a(C)  \]
and this implies that there exists $c \in C$ such that $g(c) = a$ and $u(c) = b$. Therefore, we proved that $\exists y \in R_B \cap C$ such that $f(u(y)) = x$ \ie $v(g(y)) \in f(R_B)$ which implies that $x = g(y) \in g(u^{-1}(R_B))$.

So, we have that
\beqn
v^{-1}(q_\a(R_1)) \subset R_A + v^{-1} ( f (R_B)) 
\eeqn
is bounded because $(x,y) \mapsto x+y$ is a bounded map. And this proves the proposition.
\end{proof}

\begin{rmk}
The proof of proposition \ref{prop_quasi_ab} is in total analogy with the proof proposition 1.8 of \cite{PS}, which works for bornological vector spaces of convex type over a non-trivially valued field.
We remark that, in order to show that cocartesian squares preserve strict monomorphisms, our generalization of the proof uses the boundedness of
the addition instead of using the absolute convex envelopes as done in \cite{PS}. Thus, it can be applied to bornological vector space without any convexity assumption.
\end{rmk}

\begin{prop}
In the category $\bBorn(\bAb)$ direct sums and direct products are kernel and cokernel preserving.
\end{prop}
\begin{proof}
We show that direct sums are cokernel preserving. If $(f_i: E_i \to F_i)_{i \in I}$ is a family of morphisms of $\bBorn(\bAb)$  
it is well-known that algebraically $\Coker (\underset{i \in I}\oplus f_i) \cong \underset{i \in I}\bigoplus \Coker f_i$. We denote by $\a: \underset{i \in I}\bigoplus F_i \to \Coker(\underset{i \in I}\oplus f_i)$
and $\a_i: F_i \to \Coker(f_i)$ the canonical maps, so that we have the following diagram
\[
\begin{tikzpicture}
\matrix(m)[matrix of math nodes,
row sep=2.6em, column sep=2.8em,
text height=1.5ex, text depth=0.25ex]
{ \underset{i \in I}\bigoplus E_i & \underset{i \in I}\bigoplus F_i & \underset{i \in I}\bigoplus \Coker(f_i)  \\
  & \Coker(\underset{i \in I}\oplus f_i) \\};
\path[->,font=\scriptsize]
(m-1-1) edge node[auto] {$\underset{i \in I}\bigoplus f_i$} (m-1-2);
\path[->,font=\scriptsize]
(m-1-2) edge node[auto] {$\underset{i \in I}\bigoplus \a_i$} (m-1-3);
\path[->,font=\scriptsize]
(m-1-2) edge node[auto] {$\a$} (m-2-2);
\path[->,font=\scriptsize]
(m-2-2) edge node[auto] {$\gamma$} (m-1-3);
\end{tikzpicture}.
\]
$\Coker(\underset{i \in I}\oplus f_i)$ carries the quotient bornology given by the map $\a$ and this bornology
can be described by the image of bounded subsets of $\underset{i \in I}\bigoplus F_i$. Thus, the bounded subsets of $\Coker(\underset{i \in I} \oplus f_i)$ are of the form
$\a( B_{i_1} \times ... \times B_{i_n})$ for bounded subsets $B_{i_j} \subset F_{i_j}$, hence the bornologies on $\Coker (\underset{i \in I}\oplus f_i)$ and on
$\underset{i \in I}\bigoplus \Coker f_i$ coincide.

Let us show that direct sums are kernel preserving. Consider a family $(u_i: E_i \to F_i)_{i \in I}$ of morphisms of $\bBorn(\bAb)$  
and for any $u_i$ denote $K_i$ the kernel of $u_i$ and $k_i: K_i \to E_i$ the canonical morphisms. We have the diagram 
\[
\begin{tikzpicture}
\matrix(m)[matrix of math nodes,
row sep=2.6em, column sep=2.8em,
text height=1.5ex, text depth=0.25ex]
{ 0 & \underset{i \in I}\bigoplus K_i & \underset{i \in I}\bigoplus E_i & \underset{i \in I}\bigoplus F_i  \\
  & \Ker(\underset{i \in I}\oplus u_i) \\};
\path[->,font=\scriptsize]
(m-1-1) edge node[auto] {} (m-1-2);
\path[->,font=\scriptsize]
(m-1-2) edge node[auto] {$\underset{i \in I}\bigoplus k_i$} (m-1-3);
\path[->,font=\scriptsize]
(m-1-3) edge node[auto] {$\underset{i \in I}\bigoplus u_i$} (m-1-4);
\path[->,font=\scriptsize]
(m-1-2) edge node[auto] {$\gamma$} (m-2-2);
\path[->,font=\scriptsize]
(m-2-2) edge node[auto] {$\a$} (m-1-3);
\end{tikzpicture}.
\]
where $\gamma$ is an isomorphism of the underlying abelian groups. We show that it is also an isomorphism of bornological abelian groups.

Let $B$ be a bounded subset of $\underset{i \in I}\bigoplus E_i$. By corollary \ref{stric_epi_mono}, it is sufficient to show that
\beqn
 (\bigoplus_{i \in I} k_i  )^{-1}  (B)
\eeqn
is a bounded subset of $\underset{i \in I}\bigoplus K_i$. Bounded subsets of $\underset{i \in I}\bigoplus E_i$ are of the form $B \subset \underset{i \in I}\bigoplus B_i$ for some bounded subset $B_i \subset E_i$ where 
$B_i = \{0 \}$ for almost all $i$. Moreover, $k_i^{-1}(B_i)$ is a bounded subset of $K_i$ for any $i \in I$ and $k_i^{-1}(\{0\}) = \{ 0 \} \subset K_i$
because $k_i$ is injective. It follows that 
\beqn
(\bigoplus_{i \in I} k_i  )^{-1}  (B) = \bigoplus_{i \in I} ( k_i^{-1})  (B_i)
\eeqn
is a bounded subset of $\underset{i \in I}\bigoplus K_i$.

Analogous (dual) arguments work for direct products. 
\end{proof}

\begin{rmk}
The property of being quasi-abelian of the category of bornological abelian groups is also important when one wants to study sheaves of bornological groups
(especially bornological vector spaces or bornological modules) over a topological space because it allows to use the theory of 
quasi-abelian categories and quasi-abelian sheaves developed by Schneiders \cite{SN}.
\end{rmk}

\section{Bornological rings}

\index{bornological ring}
In this work the word \emph{ring} will always mean a commutative ring with $1$. Given a ring $(A, +, \cdot)$ and a bornology $\cB$ on $A$, the quadruple $(A, +, \cdot, \cB)$ is said to be a \emph{bornological ring} if

\ben
\item $(A, +, \cB)$ is a bornological group;
\item the map $(x,y) \mapsto x \cdot y$ is bounded with respect to the product bornology on $A \times A$.
\een

Bornological rings (resp. domains, resp. fields) form a category, with bounded homomorphisms as morphisms, which we denote $\bBorn(\bRings)$ 
(resp. $\bBorn(\bDom)$, resp. $\bBorn(\bFields)$). Let $(A, \cB)$ be a bornological ring and let $f: C \to A$ (resp. $f:A \to C$) 
be a homomorphism of rings, then $(C, f^* \cB)$ (resp. $(C, f_* \cB)$ with $f$ surjective) is a bornological ring and $f$ 
is a morphism in $\bBorn(\bRings)$.

We say that a morphism $f: A \to B$ in $\bBorn(\bRings)$ is a \emph{strict} morphism if it is strict as morphism of the underlying abelian groups.

\begin{lemma} \label{ring_born_init}
Let $A$ be a ring and $(\cB_i)_{i \in I}$ be a family of ring bornologies on $A$, then $\underset{i \in I}\bigcap \cB_i$ is a ring bornology on $A$.
\begin{proof}
Similar to the proof of lemma \ref{borngrop_int}.
\end{proof}
\end{lemma}

\begin{lemma}
Let $A$ be a ring and $\cB$ a bornology on $A$, then $\cB$ is a ring bornology if and only if $\cB$ has a basis $\{B_i\}_{i \in I}$ such that $\{B_i\}_{i \in I}$
is a basis for a group bornology (for the additive group of $A$) and for every $i,j \in I$
\[ B_i B_j \subset B_{i_1} \cup \ldots \cup B_{i_n} \]
for some $i_1, \ldots, i_n \in I$.
\end{lemma}
\begin{proof}
Similar to the proof of lemma \ref{borngrop_char}.
\end{proof}

$\bBorn(\bRings)$ admits small cofiltered limits and all small colimits which we now describe. The products in $\bBorn(\bRings)$
are just the set theoretic products equipped with the product bornology.  
In the same way one can see that the forgetful functors $\bBorn(\bRings) \to \bBorn(\bAb)$ and $\bBorn(\bRings) \to \bBorn$ commute with all small cofiltered limits.

On the other hand, in the category of commutative rings coproducts are given by tensor products and we can make them into bornological rings in the 
following way. Given two morphisms $R \to A$ and $R \to B$, their pushout is $A \otimes_R B$ and in the case $R = \Z$ this is the coproduct
(this definition extends naturally also for an arbitrary family of morphisms, \ie the coproduct of a small family of rings is represented by the filtered colimit of all the finite tensor products over $\Z$ of all the finite sub-families). Hence, it is enough to describe the bornology to put on the tensor product of two rings that ensure the universal property of push-out is satisfied, because the forgetful functor $\bRings \to \bSet$ commutes with filtered colimits which implies that the same is true for $\bBorn(\bRings) \to \bBorn$. By lemma \ref{ring_born_init} the smallest ring bornology which contains the images
of the bornologies $A$ and $B$ by the canonical maps to $A \otimes_R B$ is the bornology we are looking for. More explicitly, denoting $i_A: A \to A \otimes_R B$ defined $a \mapsto a \otimes 1$ and $i_B: B \to A \otimes_R B$ defined $b \mapsto 1 \otimes b$. Then, we take as the bornology
of $A \otimes_R B$ the bornology generated by the family of all subsets of $A \otimes_R B$ obtained by applying a finite number of algebraic operations
to the image of the bounded subsets of $A$ and $B$ with respect to $i_A$ and $i_B$.

So, $\bBorn(\bRings)$ is cocomplete but none of the forgetful functors $\bBorn(\bRings) \to \bBorn(\bAb)$ nor $\bBorn(\bRings) \to \bBorn$
commute with colimits.

\index{non-archimedean bornological ring}
\index{archimedean bornological ring}
\begin{defn}
A bornological ring $(A, \cB)$ is said \emph{non-archimedean} if its subring $\Z \cdot 1_A$ is $\cB$-bounded and it is said \emph{archimedean} otherwise.
\end{defn}

An archimedean bornological ring therefore contains $\Z$ as sub-object equipped with a bornology for which it is an unbounded subset, and so it is necessarily of characteristic $0$.

\index{power-bounded element}
Let $(A, \cB)$ be a bornological ring. We define
\[ \acb = \l \{x \in A | \{x^n\}_{n \in \N} \in \cB \r \} \]
and we call it the subset of $\cB$\emph{-power bounded} elements of $A$. Notice that $\acb$ is a submonoid of $(A, \cdot)$ containing
 $\pm 1$ and $0$, since for $x, y \in \acb$ we have

\beqn
\{(xy)^n\}_{n \in \N} \subset \{x^n\}_{n \in \N} \cdot \{y^n\}_{n \in \N} \in \cB .
\eeqn

Let $U_{\cB}(A)$ denotes the subgroup of $\acb$ consisting of its invertible elements, \ie $x \in U_{\cB}(A)$ if and only if $x^{-1} \in \acb$. It is clear that $U_\cB(A)$ is the 
unique minimal face\footnote[1]{A \emph{face} in a monoid $(M, \cdot)$, commutative with $1$, is a 
submonoid $F \subset M$ such that, $\forall x,y \in M, xy \in F \iff x \in F$ and $y \in F$. 
An \emph{ ideal} of $M$ is a subset $I \subset M$ such that $M \cdot I \subset I$.} of the monoid $\acb$, so the complement 
$A_\cB^{\circ \circ} = \acb \backslash U_\cB(A)$ is the unique maximal ideal of the monoid. 

We note also that if $\varphi: (A_1, \cB_1) \to (A_2, \cB_2)$ is a morphism of bornological rings then 
$\varphi(A_1^\circ) \subset A_2^\circ$, because if $\{x^n\}_{n \in \N} \subset A_1$ is a bounded subset
then $\varphi(\{x^n\}_{n \in \N}) = \{\varphi(x)^n\}_{n \in \N}$ must be bounded in $A_2$.
Hence, we get a functor ${}^\circ: \bBorn(\bRings) \to \bMonoids$ from the category of bornological rings to the category of abelian monoids.

\begin{prop}
Let $(A, \cB)$ be a bornological ring then $\acb$ is a \emph{saturated} submonoid of $(A, \cdot)$, in the sense that if $x \in A$
 and there exists $n > 0$ such that $x^n \in \acb$ then $x \in \acb$.
\end{prop}
\begin{proof}
Suppose that $x \in A$ and that there exists $n > 0$ such that $x^n \in \acb$. Then $P = \{(x^n)^m\}_{m \in \N}$ is a bounded 
subset of $(A, \cB)$. Hence, $x P$, $x^2 P$, $\ldots$, $x^{n-1} P$ are $\cB$-bounded subsets and the subset
\beqn
P \cup xP \cup \ldots \cup x^{n-1} P = \{ x^n\}_{n \in \N} \in \cB
\eeqn 
\ie $x \in \acb$. 
\end{proof}

\subsection{Convexity on bornological rings}

We introduce the notion of convexity over a bornological ring.

\begin{notation}
Let $(A, \cB)$ be a bornological ring. For all $n = 1, 2, 3, \ldots$ we define $\cR^n(A_\cB^\circ) \subset A^n$ as the set of $n$-tuples 
$(\la_1, \ldots, \la_n) \in A^n$ such that for all $x_1, \dots, x_n \in \acb$
\beqn
\sum_{i = 1}^n \la_i x_i \in \acb.
\eeqn
We define
\beqn
\cR_\cB(A) = \bigcup_{n=1}^\infty \cR_\cB^n (A)
\eeqn
where, if $m \ge n$, we identify $\cR_\cB^n(A)$ with the subset of $\cR_\cB^m(A)$ of $m$-tuples with the last $m - n$ zero coordinates.
\end{notation}

\index{absolutely convex subset}
\index{absolutely convex hull}
\begin{defn} \label{convbal}
Let $(A, \cB)$ be a bornological ring and $V$ an $A$-module. A subset $S \subset V$ is said \emph{absolutely convex} if for all $n \in \N$,
for all $x_1, \ldots, x_n \in S$ and all $(\la_1, \ldots, \la_n) \in \cR^n(\acb)$ one has that 
\[ \sum_{i = 1}^n \la_i x_i \in S. \]
For $S \subset V$, the \emph{absolutely convex hull} $\Gamma_\cB(S)$ of $S$ in $V$ is defined by the intersection of all absolutely
convex subsets of $V$ containing $S$.
\end{defn}

\begin{prop}
Let $(A, \cB)$ be a bornological ring and $V$ an $A$-module. The intersection of any family of absolutely convex subsets of $V$
is absolutely convex, hence the definition of absolutely convex hull is well posed.
\end{prop}
\begin{proof}
Let $\{ S_i \subset V \}_{i \in I}$ be a family of absolutely convex subsets of $V$. Let $x_1, \ldots , x_n \in \underset{i \in I}\bigcap S_i$ and $(\la_1, \ldots , \la_n) \in \cR^n(\acb)$
then by hypothesis 
\[ \sum_{j = 1}^n \la_j x_j \in S_i \]
for all $S_i$'s, hence 
\[ \sum_{j = 1}^n \la_j x_j \in \bigcap_{i \in I} S_i \]
so $\underset{i \in I}\bigcap S_i$ is absolutely convex.
\end{proof}

The following, easy-to-show properties, are useful.

\begin{lemma} \label{lemma:R_ideal}
Let $(A, \cB)$ be a bornological ring then
\ben
\item $\cR^n(\acb) \subset (\acb)^n$;
\item let $(\la_1, \ldots, \la_n) \in \cR^n(\acb)$, then for any $(a_1, \ldots, a_n) \in \acb$ we have that $(\la_1 a_1, \ldots, \la_n a_n) \in \cR^n(\acb)$,
      \ie $\cR^n(\acb)$ is an ideal in $(\acb)^n$ with respect to the direct product monoid structure.
\een
\end{lemma}
\begin{proof}
\ben
\item Since $0,1 \in \acb$ by taking $x_1 = 1, x_2 = 0, \ldots , x_n = 0$
      we see that 
\beqn
\sum_{i = 1}^n \la_i x_i \in \acb
\eeqn
 implies $\cR^n(\acb) \subset (\acb)^n$.
\item If $(\la_1, \ldots, \la_n) \in \cR^n(\acb)$, $(a_1, \ldots, a_n) \in \acb$ and it is given any other $n$-tuple 
      $(b_1, \ldots, b_n) \in \acb$ then
      \[ \sum_{i = 1}^n \la_i a_i b_i = \sum_{i = 1}^n \la_i (a_i b_i) \]
      and since $\acb$ is a monoid then $a_i b_i \in \acb$. Thus, $\underset{i = 1}{\overset{n}\sum} \la_i (a_i b_i) = \underset{i = 1}{\overset{n}\sum} (\la_i a_i) b_i \in \acb$, which means that 
      $(\la_1 a_1, \ldots, \la_n a_n) \in \cR^n(\acb)$.
\een
\end{proof}

\begin{prop}
Let $(A, \cB)$ be a bornological ring and $V$ an $A$-module. For all $S \subset V$ one has
\beqn
 \Gamma_\cB(S) = \l \{ \sum_{i = 1}^n \la_i x_i | n = 1, 2, \ldots, (\la_1, \ldots, \la_n) \in \cR^n(\acb), x_1, \ldots, x_n \in S \r \}.
\eeqn
\end{prop}
\begin{proof}
Let $T$ be an absolutely convex subset of $V$ containing $S$, then
\[ \Gamma_\cB(S) \subset T \]
by the definition of absolutely convex subset and the fact that $S \subset T$. Thus, we only need to check that $\Gamma_\cB(S)$ is itself
absolutely convex. Let $x_1, \ldots, x_n \in \Gamma_\cB(S)$ and $(\la_1, \ldots, \la_n) \in \cR^n(\acb)$, so
\[ \sum_{i = 1}^n \la_i x_i = \sum_{i = 1}^n \la_i \sum_{j = 1}^{m_i} \mu_{i, j} s_{i, j} \]
with $(\mu_{i, 1}, \ldots, \mu_{i, m_i}) \in \cR^{m_i}(\acb)$ and $s_{i,j} \in S$. So, the only thing left to show is that 
$(\la_1 \mu_{1, 1}, \la_1 \mu_{1, 2}, \ldots, \la_n \mu_{n, m_n}) \in \cR^m(\acb)$ where $m = \underset{i = 1}{\overset{n}\sum} m_i$. Let $(a_{i ,j}) \in (\acb)^m$, then for all $i$
\[ \tilde{a}_i = \sum_{j = 1}^{m_i} \mu_{i, j} a_{i, j} \in \acb \]
by the choice of the $\mu_{i, j}$. Since by the lemma \ref{lemma:R_ideal} $\cR^n(\acb)$ is an ideal of the monoid $(\acb)^n$ we have that
$(\la_1 \tilde{a}_1, \ldots, \la_n \tilde{a}_n) \in \cR^n(\acb)$, which implies that $\underset{i = 1}{\overset{n}\sum} \la_i x_i \in \Gamma_\cB(S)$. So, the proposition is proved.
\end{proof}

\begin{prop}
The definition of absolutely convex subset of an $A$-module coincides with the usual definition for $A = \R, \C$ or for $A$ a 
non-archimedean field, equipped with the bornologies induced by the valuations.
\end{prop}
\begin{proof}
For $A= \R, \C$ and $V$ an $A$-vector space, a subset $S$ of $V$ is absolutely convex in the classical sense if for any 
$\la_1, \ldots, \la_n \in A$ with $\underset{i = 1}{\overset{n}\sum} |\la_i| \le 1$, and any $x_1, \ldots, x_n \in S$ it happens that $\underset{i = 1}{\overset{n}\sum} \la_i x_i \in S$. 
So, we have to show that, if $A=\R, \C$ then
\beqn
 (\la_1, \ldots, \la_n) \in \cR^n(\acb) \iff \sum_{i = 1}^n |\la_i| \le 1.
\eeqn
If $ (\la_1, \ldots, \la_n) \in \cR_\cB^n(A)$, we consider the elements
$x_i = \frac{|\la_i|}{\la_i} \in \acb = D_\C(0, 1^+) \doteq \{x \in \C | |x| \le 1\}$, for $i = 1, 2, \ldots, n$. So, by definition
\beqn
 \l | \sum_{i = 1}^n \la_i x_i \r | = \l |\sum_{i = 1}^n \la_i \frac{|\la_i|}{\la_i} \r | = \sum_{i = 1}^n |\la_i| \le 1.
\eeqn
For the converse, if $\la_1, \ldots, \la_n \in \C$ are such that $\underset{i = 1}{\overset{n}\sum} |\la_i| \le 1$, then for any $x_1, \ldots, x_n \in D_\C(0, 1^+)$, we have 
\[ \l | \sum_{i=1}^n \la_i x_i \r | \le \sum_{i=1}^n |\la_i||x_i| \le  \sum_{i=1}^n |\la_i| \le 1. \]
If $A$ is a non-archimedean field, we call $R$ its ring of integers, and let $V$ be an $A$-vector space. A subset $S$ of $V$ is usually
said to be absolutely convex if it is an $R$-submodule of $V$. On the other hand
\beqn
(\la_1, \ldots, \la_n) \in \cR^n(\acb) \iff \la_i \in R, \forall i = 1, \ldots, n,
\eeqn
is an immediate consequence of the non-archimedean nature of $A$ which implies that $R = \acb$ is a subring of $A$.
\end{proof}

The data $\{ \cR^n(\acb) \}_{n \in \N}$ is equivalent to the generalized ring $\Sigma_{\acb}$ associated to the 
inclusion of monoids $\acb \rhook A$, defined definition \ref{defn:power_bounded_monad}. Moreover, using the language of generalized rings, a subset $X \subset A$ is absolutely convex 
if and only if it is a $\Sigma_\acb$-submodule of $A$.

So, up to now, we have showed that given any bounded morphism of bornological rings $f: (A, \cB_A) \to (B, \cB_B)$ we get the following commutative diagram of maps of monoids
\[
\begin{tikzpicture}
\matrix(m)[matrix of math nodes,
row sep=2.6em, column sep=2.8em,
text height=1.5ex, text depth=0.25ex]
{ A       & B   \\
  A_{\cB_A}^\circ & B_{\cB_B}^\circ \\};
\path[->,font=\scriptsize]
(m-1-1) edge node[auto] {$f$} (m-1-2);
\path[->,font=\scriptsize]
(m-2-1) edge node[auto] {$f^\circ$} (m-2-2);
\path[right hook->,font=\scriptsize]
(m-2-1) edge node[auto] {} (m-1-1);
\path[right hook->,font=\scriptsize]
(m-2-2) edge node[auto] {} (m-1-2);
\end{tikzpicture}.
\]
Hence, it is natural to ask whether this map gives rise to a map of generalized rings making the following diagram commutative
\[
\begin{tikzpicture}
\matrix(m)[matrix of math nodes,
row sep=2.6em, column sep=2.8em,
text height=1.5ex, text depth=0.25ex]
{ \Sigma_A  & \Sigma_B   \\
  \Sigma_{A_{\cB_A}^\circ} & \Sigma_{B_{\cB_B}^\circ} \\};
\path[->,font=\scriptsize]
(m-1-1) edge node[auto] {$f$} (m-1-2);
\path[->,font=\scriptsize]
(m-2-1) edge node[auto] {$f^\circ$} (m-2-2);
\path[right hook->,font=\scriptsize]
(m-2-1) edge node[auto] {} (m-1-1);
\path[right hook->,font=\scriptsize]
(m-2-2) edge node[auto] {} (m-1-2);
\end{tikzpicture}.
\]
Thanks to proposition \ref{prop:monad_functorial} we know that this is always the case.

\section{Bornological modules}

In this section $A$ denotes a bornological ring (always commutative with identity) and $E$ denotes an $A$-module.

\index{bornological module}
\begin{defn}
We say that $E$ is a bornological module over $A$ if $E$ is equipped with a group bornology for which the map 
$(\la, x): (A \times E) \to E$ defined $x \mapsto \la x$ is bounded.
\end{defn}

In total analogy with previous sections we can define the category of bornological modules over $A$ whose morphism
are bounded homomorphism of modules and it will be denote with $\bBorn(\bMod_A)$.
The results that we discussed for the category of bornological groups easily extend to  $\bBorn(\bMod_A)$. We only list these properties, without giving proves.

\ben
\item $\bBorn(\bMod_A)$ is complete and cocomplete and the forgetful functor to the category of $A$-modules commutes with both limits and colimits;
\item $\bBorn(\bMod_A)$ is additive and admits kernel and cokernel;
\item Given a morphism of bornological modules $f: E \to F$, we have that 
      \ben
      \item $\Ker f = f^{-1}(0)$ endowed with the subspace bornology;
      \item $\Coker f = F/f(E)$ endowed with the quotient bornology;
      \item $\Im f = f(E)$ endowed with the bornology induced by the inclusion $f(E) \rhook F$;
      \item $\Coim f = E/f^{-1}(0)$ endowed with the quotient bornology.
      \een
\item $\bBorn(\bMod_A)$ is quasi-abelian.
\een

In the next chapter we recall the rich theory of bornological vector spaces over a valued field.

\section{The language of halos and the tropical addition} \label{sec:halos}

We recall here some definitions used by Paugam in \cite{PAU}, which will be useful in the next sections.

\index{halo}
\begin{defn}
A \emph{halo} is a semiring $A$, with an identity, whose underlying set is equipped with a partial order $\le$ which is compatible with its operations:
$x \le z$ and $y \le t$ implies $x y \le z t$ and $x + y \le z + t$. A \emph{morphism} between two halos is an increasing map $f : A \to B$
which is submultiplicative, \ie $f(1) = 1$ and
\[ f(ab) \le f(a)f(b), \forall a, b \in A \]
and subadditive, \ie $f(0) = 0$ and
\[ f(a + b) \le f(a) + f(b), \forall a, b \in A. \]

A halo morphism is called 
\begin{itemize}
\item \emph{square-multiplicative} if $f(a^2) = f(a)^2$ for all $a \in A$;
\item \emph{power-multiplicative} if $f(a^n) = f(a)^n$ for all $a \in A$ and $n \in \N$;
\item \emph{multiplicative} if $f(a b) = f(a) f(b)$ for all $a, b \in A$.
\end{itemize}

The class of halos with halo morphisms form a category that is denoted $\bHalos$.
\end{defn}

\begin{lemma}
On a ring $A$ there exists only one halo structure, up to isomorphisms. This structure is given by the trivial order on $A$.
\end{lemma}
\begin{proof}
\cite{PAU}, Lemma 1.
\end{proof}

\begin{rmk}
This lemma not only implies that there exists a canonical functor that embeds the category of rings in the category of halos,
but also that this functor is fully faithful. Indeed, the inequalities in the subadditivity and submultiplicativity conditions for halo morphisms become equalities when the halos are equipped with the trivial order.
\end{rmk}

In the next definition we change the names used in \cite{PAU} a little, for our convenience.

\index{aura}
\begin{defn}
Let $R$ be a halo. We say that $R$ is
\begin{itemize}
\item \emph{positive}, if $0 < 1$;
\item a \emph{pre-aura}, if $R$ is a semifield, \ie if $R - \{0\}$ is a group for the multiplication operation;
\item an \emph{aura}, if $R$ is a positive pre-aura whose order is total;
\item \emph{tropical}, if it is totally ordered and $\forall a,b \in R$ we have that $a + b = \max \{ a, b \}$.
\end{itemize}
\end{defn}

\index{generalized seminorm}
\begin{defn} \label{defn:gen_seminomr}
Let $A$ be a ring and $R$ an aura. A (generalized) \emph{seminorm} on $A$ with values on $R$ is a halo morphism $|\cdot|: A \to R$.
\end{defn}

In the following chapters we use the language of halos as a shortcut to deal at the same time with archimedean and non-archimedean seminorms over algebras over valued fields.
Our use of the concepts we have introduced in this section aims at having a better notation for dealing with normed algebras over archimedean and non-archimedean base fields at the same time. In particular, when we deal with an algebra $A$ over a non-archimedean valued field $k$, we (tacitly) equip $\R_{\ge 0}$ with its structure of tropical aura in order to obtain from the triangle inequality of definition \ref{defn:gen_seminomr}, formally written $|a + b| \le |a| + |b|$, the ultrametric triangle inequality $|a + b| \le \max \{ |a|, |b| \}$.
Moreover, this will be handy in dealing with power-series over $k$ and the summation norm on $\underset{i = 0}{\overset{\infty}\sum} a_i X^i \in k \ldbrack X \rdbrack$
which will be formally written
\[ \| \sum_{i = 0}^\infty a_i X^i \| = \sum_{i = 0}^\infty |a_i| \]
but its meaning is different with respect to the halos structure considered on $\R_{\ge 0}$. If we are working on an archimedean valued field
the summation symbol stands for the usual addition between real numbers whereas if we are dealing with a non-archimedean one it stands for tropical addition. From now on we will adopt this convention (we will also have the care of remind the reader about this fact in its most important uses).

Finally, we remark that when we deal with algebras over a non-archimedean valued field, this convention implies that we work in a full subcategory of the category of all seminormed algebras, \ie precisely the one endowed with an ultrametric seminorm. Indeed, we are avoiding to consider algebras which are endowed with a seminorm which does not satisfy the ultrametric triangle inequality when they are defined over a non-archimedean valued field (considered as pathological objects by some authors). This is a point we must keep in mind especially when we state the universal properties of the algebras that we will study in the following chapters. In particular, from this choice it depends the correct statement of the universal properties characterizing the Tate's algebras over non-archimedean valued fields.

\chapter{Bornological spectra and power-series} \label{sec:spectra}

This chapter is devoted to the study of bornological algebras over a non-trivially valued complete field, archimedean or non-archimedean. Some results of this chapter can be found in some form in literature, in different contexts or less generality. This chapter is not meant to be an encyclopaedic collection of results on bornological algebras. Its aim is to study bornological algebra as far as needed in next chapters, mainly to have a well-established theory over which we can base the theory of dagger affinoid algebras.

In the first section we define the notion of spectrum of a bornological algebra. This notion naturally extends Berkovich's notion of the spectrum of a seminormed algebra and it is functorially, contravariantly, attached to bornological algebras. Thus, for a bornological algebra $A$, its spectrum $\cM(A)$ is defined as a suitable set of multiplicative and bounded seminorms over $A$ and equipped with the weak topology obtained pulling back the topology of $\R_{\ge 0}$. The most interesting fact about this spectra is that we can show that, for a large class of bornological algebras, $\cM(A)$ is a non-empty, Hausdorff and compact topological space. After this, we show some results about power-series rings of bornological algebras. In particular we characterize the rings $A \lt X \gt$ and $A \lt X \gt^\dagger$ of strictly convergent and overconvergent power-series over a bornological algebra by a universal property and study some of their basic properties. In particular, we show that the associations $A \mapsto A \lt X \gt$ and $A \mapsto A \lt X \gt^\dagger$ preserve several properties of the algebra $A$ we care about.

We end this chapter by recalling some results on spectra of locally convex topological algebras which are often found in literature explained only over $\C$, but which are valid over any valued field, and we relate these with our results on bornological algebras.

\section{Spectra of bornological algebras}

Let $(k, |\cdot|)$ be a fixed field complete with respect to a non-trivial absolute value $|\cdot|: k \to \R_{\ge 0}$.
From now on, all vector spaces and algebras are defined over $k$. It is clear that the valuation on $k$ induces a bornology given by subsets $B \subset k$ such that $|B| < C$ for some $C \in \R_{\ge 0}$. We call this bornology \emph{canonical} and we denote it by $\cB_{|\cdot|}$.

\index{bornological vector space}
\begin{defn}
A vector space $E$ over $k$ is said to be \emph{a bornological vector space} if on $E$ is defined a bornology compatible with the structure of vector space,
\ie the vector sum and the multiplication by scalars are bounded maps, where on $k$ we consider the canonical bornology given by $|\cdot|$. 
In other words, $E$ is a bornological module over the bornological ring $(k, \cB_{|\cdot|})$.
\end{defn}

As we saw in the previous chapter the category of bornological vector spaces is complete, cocomplete and quasi-abelian. We are interested in studying bornologies over $k$-vector spaces which can be defined by mean of seminorms. Thus, we recall now some results and definitions about the theory of bornological vector spaces of convex type.

\index{Minkowsi functional}
\begin{defn}
Let $E$ be a $k$-vector space and $A \subset E$ an \emph{absolutely convex subset}, also called \emph{disks} (see \cite{H1}, page 80 for the definition of disks). We associate to $A$ the \emph{gauge} or \emph{Minkowski functional} $\mu_A: E \to \R_{\ge 0}$ defined by
\[ \mu_A(x) \doteq \inf_{\la \in k} \l \{ |\la| | x \in \la A \r \}, \] 
which is a seminorm. We denote $E_A$ the vector subspace of $E$ generated by $A$ equipped with the seminorm given by the gauge $\mu_A$.
\end{defn}

The unit ball of a seminormed $k$-vector space $E$ is an absolutely convex subset of $E$ whose gauge coincides with the given seminorm. It is therefore possible to show that there is a bijective correspondence between seminorms on $E$ and absorbing absolutely convex subsets of $E$ modulo a suitable equivalence relation, that we will not describe here since is not interesting for our scopes. It is also clear that the family of all bounded subsets of a seminormed space induces a structure of bornological $(k, \cB_{|\cdot|})$-module on $E$. The family of balls of all possible radii is a base for this bornology which hence has a base made of absolutely convex subsets. So, in next definition we want to generalize these properties of seminormed spaces.

\index{bornological vector space of convex type}
\begin{defn}
A bornological vector space $E$ over $k$ is said of \emph{convex type} if there exists a base for the bornology consisting of absolutely convex subsets of $E$.
\end{defn}

This is the same as saying that $E$ is a filtered inductive limit, calculated in the category of bornological vector spaces, of 
seminormed vector spaces (equipped with their canonical bornology induced by the seminorms) as shown by Houzel in \cite{H1} proposition 1, page 92. Explicitly, let $E$ be a bornological vector space of convex type over $k$ and let $\cB_E$ denote the family of bounded disks of $E$, then
\begin{equation} \label{eqn:born_space}
 E \cong \limind_{B \in \cB_E} E_B 
\end{equation}
in $\bBorn(\bMod_k)$. On the other hand, it is easy to verify that given any inductive system $\{E_i\}_{i \in I}$ of seminormed spaces, then the canonical bornologies of the seminorms of $E_i$  induce on $E = \underset{i \in I}\limind E_i$ a bornology of convex type. More precisely, one can show that the category of bornological vector spaces of convex type is equivalent to the full subcategory defined by the essentially monomorphic objects in the category of $\bInd$-objects of the category of seminormed spaces, cf. \cite{PS} corollary 3.5. So, from now on, we will make no difference between a bornological vector space and a monomorphic inductive system that represents it in the category of $\bInd(\bSn_k)$. We will denote with $\bBorn_k$ the full subcategory of $\bBorn(\bMod_k)$ identified by the bornological vector spaces of convex type.

\index{complete bornological vector space}
\index{separated bornological vector space}
\begin{defn} \label{defn:sep_born_space}
Let $E$ be a bornological vector space of convex type and consider the system 
$\underset{B \in \cB_E} \limind E_B$ of equation \ref{eqn:born_space} which represents $E$. Then,
\begin{itemize}
\item $E$ is said to be \emph{separated} if each $E_B$ is normed;
\item $E$ is said to be \emph{complete} if $\underset{B \in \cB_E} \limind E_B$ admits a final subsystem for which each $E_B$ is a $k$-Banach space.
\end{itemize} 
\end{defn}

We will denote by $\bSBorn_k$ (resp. $\bCBorn_k$) the full subcategories of $\bBorn_k$ identified by the separated (resp. complete) bornological vector spaces.

\begin{prop}
A bornological vector space of convex type is complete if and only if 
\[ E \cong \limind_{i \in I} E_i \]
for a monomorphic filtered system where all $E_i$ are $k$-Banach spaces.
\end{prop}
\begin{proof}
\cite{H1} page 96, proposition 7.
\end{proof}

It is not difficult to check that this last proposition implies that $\bCBorn_k$ is equivalent to the full subcategory of $\bInd(\bBan_k)$ identified by essentially monomorphic inductive systems, see proposition 5.15 of \cite{PS}. The analogous result also hold for $\bSBorn_k$ which is equivalent to the full full subcategory of $\bInd(\bNrm_k)$ made of essentially monomorphic objects, cf. corollary 4.20 of \cite{PS}.

\begin{notation}
We denote by $\bTop_k$ the category of locally convex topological vector spaces over $k$ whose morphisms are continuous linear maps of vector spaces.
\end{notation}

\index{canonical bornology}
We recall that Houzel in \cite{H1} has defined two functors ${}^b: \bTop_k \to \bBorn_k$, 
${}^t: \bBorn_k \to \bTop_k$ and that these functors are adjoints. More precisely, ${}^t$ is the left adjoint functor of ${}^b$. 
We recall briefly how these functors are defined: Let $(E, \cT)$ be a  topological vector space, the \emph{canonical bornology} on $E$
is defined by the family of subsets $B \subset E$ such that for any open neighborhood of $0$ there exists
a $\la \in k^\times$ with $B \subset \la U$; \ie the bounded subsets are the subsets absorbed by any open neighborhood of $0 \in E$. It is easy to check that
this family induces on $E$ the structure of a bornological vector space which is of convex type if $E$ is a locally convex vector space. This association is functorial and is precisely the definition of the functor ${}^b$, \ie $(E, \cT)^b$ is the bornological vector space obtained equipping $E$ with the canonical (also called \emph{von Neumann}) bornology.

\index{bornivorous topology}
On the other hand, if we have a bornological vector space $(E, \cB)$, then the family of subsets which absorb all bounded subsets of
$E$ forms a filter of neighborhoods of $0$ which induces a topology on $E$ which endows $E$ with a structure of topological vector space which turns out to be a locally convex space if $(E, \cB)$ is of convex type.
The association of this topology is functorial, is denoted by $(E, \cB)^t$ and is called the \emph{topology of bornivorous subsets}.

\index{normal bornological vector space}
\index{normal topological vector space}
\begin{defn} \label{defn:normal_space}
A bornological vector space $E$ over $k$ is said to be \emph{normal} if $(E^t)^b \cong E$. The same definition is given for a topological vector space 
$E$, which is \emph{normal} if $(E^b)^t \cong E$.
\end{defn}

The full subcategories of normal bornological and topological vector spaces are equivalent by mean of the functors ${}^b$ and ${}^t$. 
These subcategories are very useful since they are the categories where the concepts of boundedness and of continuity of linear maps are equivalent.

\begin{exa}
\ben
\item Metrizable, bornological or topological, vector spaces are normal, cf. \cite{H1} page 109.
\item Normal topological vector spaces are what Bourbaki calls ``bornological locally convex spaces'', cf. definition 1 in page III.12 of \cite{Bour}.
\item In next chapters we will work with LF-spaces\footnote[1]{What we call LF-spaces are sometimes called generalized LF-spaces by some authors. This happens mainly when the term LF-space is reserved for what other authors call strict LF-spaces; we follow the nomenclature of calling strict LF-space the ones defined by inductive systems with strict monomorphism and use LF-space otherwise.}, 
      defined as monomorphic countable filtrant direct limits of Fr\'echet spaces, and very often LB-space \ie monomorphic countable filtrant direct limits
      of Banach spaces. On these spaces one can consider the locally convex direct limit topology which has the main drawback of being not easy to deal with. For example, it is known that no LB-space (our main object of study in the following) is 
      metrizable as topological vector space. This implies that the family of neighborhoods of zero is not countable. 
      We can simplify the handling of such spaces by equipping them
      with the direct limit bornology, which has a countable base, it is easy to describe and often contains the same amount of information of the direct limit locally convex topology. LF-spaces are not always normal and there are many works in literature studying normality of some classes of LF-spaces.
\item 
\index{Washnitzer algebra}
Our main example of LB-space is the Washnitzer algebra defined by Grosse-Kl\"onne \cite{GK} as the colimit
      \[ W_k^n = \limind_{\stackrel{\rho > 1}{\rho \in \sqrt{k^\times}}} T_k^n(\rho) \]
      where $T_k^n(\rho)$ are the Tate algebras for the polydisk of radius $\rho$ for a non-archimedean field $k$. We will see that $W_k^n$ is a normal and LF-regular LB-space.
\een
\end{exa}

So, on a normal bornological vector space we do not lose any information by considering the canonical topology associated to it and vice-versa. It is known that if $E$ is a complete locally convex vector space then $E^b$ is a complete bornological vector space. But it is possible
that $E^b$ is bornologically complete even if $E$ is not complete.
A common issue when dealing with the functors ${}^b$ and ${}^t$ is whether they commute with some particular kind of limits or colimits.

\begin{prop} \label{prop_commute_limits}
Let $\{E_i\}_{i \in I}$ be any directed family of bornological vector spaces of convex type then
\ben
\item $(\underset{i \in I}\limind E_i)^t \cong \underset{i \in I}\limind E_i^t$;
\item if $I = \N$, with the canonical order, and the morphisms $E_{n +1} \to E_n$ are strict epimorphisms then $(\underset{i \in I}\limpro E_i)^t \cong \underset{i \in I}\limpro E_i^t$.
\een
Let $\{E_i\}_{i \in I}$ be any directed family of locally convex topological vector spaces then
\ben
\item $(\underset{i \in I} \limpro E_i)^b \cong \underset{i \in I}\limpro E_i^b$;
\item if $I = \N$, with the canonical order, and the morphisms $E_n \to E_{n + 1}$ are closed strict monomorphisms 
      then $(\underset{i \in I} \limind E_i)^b \cong \underset{i \in I}\limind E_i^b$.
\een
\end{prop}
\begin{proof}
The commutation of ${}^t$ with colimits and ${}^b$ with limits follows from the fact that they form an adjoint pair of functors. For the other claims we refer to \cite{H1} pages 104-105.
\end{proof}

Although last proposition will be very useful, it is not enough for all our purposes and we will need stronger results that we will show in next chapters. The following property gives an idea of why bornological vector spaces are suitable for studying direct limits of seminormed spaces.

\begin{prop} \label{prop:complete_born}
The forgetful functor $\bBorn_k \to \bBorn(\bMod_k)$ from the category of bornological vector spaces of convex type to the category of bornological vector spaces commutes with all colimits.
\end{prop}
\begin{proof}
\cite{H1} page 95.
\end{proof}

The previous property is in contrast with what one meets in the theory of locally convex topological vector space, for which the forgetful functor
to the category of topological vector spaces commutes with limits but not with colimits, making it hard to understand locally convex direct limit topologies.

So, after this recall of the theory of bornological vector spaces we start the study of the main objects of this section: Bornological algebras over $k$.

\index{bornological algebra}
\begin{defn}
A bornological vector space $A$ over $k$ equipped with a bilinear associative function $A \times A \to A$, called multiplication map, 
is said a \emph{bornological $k$-algebra} if the multiplication map is bounded. 
We always suppose that $A$ has an identity and that the multiplication is commutative.
A morphism of bornological algebras is a bounded linear map that preserves multiplication and maps $1$ to $1$.
We denote by $\bBorn(\bAlg_k)$ the category of bornological algebras over $k$.
\end{defn}

Thus, there is an inclusion functor $\bBorn(\bAlg_k) \rhook \bBorn(\bRings)$. It is easy to check that this functor commutes with any kind of limit (for the ones that exist) and colimit.
We can compose the functors $\bBorn(\bAlg_k) \rhook \bBorn(\bRings) \stackrel{{}^\circ}{\to} \bMonoids$, where the functor ${}^\circ$ is the one defined so far. So, given an element 
$A \in \ob(\bBorn(\bAlg_k))$ we will write
\[ A^\circ \doteq \l \{ f \in A | \{f, f^2, ... \} \text{ is bounded} \r \} \]
for the subset of \emph{power-bounded} elements of $A$.

\index{multiplicatively convex bornological algebra}
\begin{defn} \label{defn:bornological_m}
A bornological algebra $A \in \ob(\bBorn(\bAlg_k))$ is said \emph{multiplicatively convex} (or simply \emph{m-algebra}) if 
\[ A \cong \limind_{i \in I} A_i \]
in the category of bornological algebras, with $A_i$ seminormed algebras.
\end{defn}

\begin{rmk}
The important request in the previous definition is that the maps of the system $\underset{i \in I}\limind A_i$ are required to be algebras morphisms and not only linear maps!
\end{rmk}

\index{complete bornological algebra}
\begin{defn} \label{complete_born_alg}
Let $A$ be a multiplicatively convex bornological algebra we say that $A$ is \emph{complete} if there is an isomorphism
\[ A \cong \limind_{i \in I} A_i \]
where the system is filtered, the system morphisms $\varphi_{i,j}$ are monomorphisms of the underlying bornological spaces and $A_i$ are $k$-Banach algebras.
\end{defn}

We denote by $\widehat{\bBorn}(\bAlg_k)$ the full subcategory of $\bBorn(\bAlg_k)$ identified by complete bornological algebras.

\begin{prop}
The underlying bornological vector space of a complete bornological algebra is a complete bornological vector space.
\end{prop}
\begin{proof}
Immediate consequence of proposition \ref{prop:complete_born}. 
\end{proof}

The inclusion $\widehat{\bBorn}(\bAlg_k) \rhook \bBorn(\bAlg_k)$ commutes with monomorphic filtered colimits. 
Notice that given a multiplicatively convex bornological algebra $A \cong \underset{i \in I}\limind A_i$ such that the system is filtered and monomorphic it is not clear that the completion $\what{A} \cong \underset{i \in I}\limind \what{A}_i$ is well defined, where $\what{A}_i$ is the separated completion of $A_i$. This is due to the fact that the completion functor, also considered on the category of normed $k$-algebras, does not preserve monomorphisms in general. Since we only work with complete algebras we do not discuss the issue of defining the completion of bornological m-algebras.

The notion of multiplicatively convex bornological algebra that we defined, following Houzel's ideas from \cite{H2}, is equivalent to the notion of locally multiplicative bornological algebra introduced by Meyer in \cite{MEYER}, cf. theorem 3.10 and definition 3.11. Meyer discusses only the archimedean case of the theory but his arguments easily extend to any non-trivially valued complete base field. Following Meyer we have the next proposition.

\begin{prop} \label{prop:m_algebras}
	The category of multiplicatively convex bornological algebras is closed under the following operations:
	\begin{enumerate}
		\item subalgebras, equipped with the induced bornology, of bornological m-algebras are m-algebras;
		\item all kind of small colimits;
		\item finite limits (when they exist);
	\end{enumerate}
	where all limits and colimits are calculated in $\bBorn(\bAlg_k)$.
	\begin{proof} 
		\cite{MEYER} lemma 3.18, page 113.
	\end{proof}
\end{prop}

As immediate corollary we have the following result.

\begin{cor} \label{cor:fitered_m_algebra}
	The following are equivalent:
	\begin{enumerate}
		\item $A$ is a multiplicatively convex bornological algebra;
		\item $A$ can be written as a monomorphic filtered direct limit of seminormed algebras.
	\end{enumerate}
\end{cor}

\index{bounded seminorm}
\begin{defn} \label{defn:bounded_seminorm}
Let $A$ be a bornological vector space. We say that a seminorm $\|\cdot \|: A \to \R_{\ge 0}$ is \emph{bounded} if it is bounded as map in the 
category of bornological sets, with $\R_{\ge 0}$ equipped with the canonical metric bornology. 
\end{defn}

We will be interested to study a particular kind of bounded seminorms defined over bornological algebras.

\begin{defn}
Let $A$ be a bornological algebra over $k$. A seminorm $\|\cdot \|: A \to \R_{\ge 0}$ is said to be \emph{compatible} with the vatuation on $k$ if the restriction of $\|\cdot\|$ to $k$ is equal to the valuation on $k$, \ie that the following diagram 
\[
\begin{tikzpicture}
\matrix(m)[matrix of math nodes,
row sep=2.6em, column sep=2.8em,
text height=1.5ex, text depth=0.25ex]
{ k  & A \\
   & \R_{\ge 0} \\};
\path[->,font=\scriptsize]
(m-1-1) edge node[auto] {$i$} (m-1-2);
\path[->,font=\scriptsize]
(m-1-2) edge node[auto] {$\|\cdot \|$} (m-2-2);
\path[->,font=\scriptsize]
(m-1-1) edge node[auto] {$|\cdot|$} (m-2-2);
\end{tikzpicture}
\]
is commmutative.
\end{defn}

\index{bornological spectrum}
\begin{defn} \label{born_spectrum}
Let $A$ be a bornological algebra. The \emph{bornological spectrum} (or simply the spectrum) of $A$ is the set $\cM(A)$ 
of all bounded multiplicative seminorms compatible with the valuation on $k$. $\cM(A)$ is equipped with the weakest topology such that all maps 
$\cM(A) \to \R_+$ defined $\|\cdot\| \mapsto \|f\|$, for $f \in A$, are continuous.
\end{defn}

The classical definition of boundedness for seminorms, for example as defined in the first chapter of \cite{BER2}, says that $|\cdot|_2: A \to \R_+$ is bounded with respect to $|\cdot|_1: A \to \R_+$ if there exist a constant $C > 0$ such that
\[ |f|_2 \le C |f|_1 \]
for all $f \in A$, where $A$ can be a general Banach (or normed) ring, not necessarily a $k$-algebra. We now check that our definitions naturally extend the classical ones for $k$-algebras, while this is no longer true for a general Banach ring.

\begin{prop} \label{prop:compatible_seminorm}
Let $(A, |\cdot|)$ be a seminormed $k$-algebra, in classical sense. If $|\cdot|$ is non-trivial and multiplicative then it is also compatible with the valuation on $k$.
\end{prop}
\begin{proof}
Quotienting by the kernel of $|\cdot|$ we can reduce the proof to the case when $A$ is normed, since by hypothesis $\ker(|\cdot|) \ne A$. Then, since $k$ is a complete valued field and $|\cdot|$ is a norm this implies that $|\cdot|$ is compatible with the valuation on $k$ because the restriction of $|\cdot|$ is a valuation on $k$ and therefore uniquely determined. 
\end{proof}

\begin{prop} \label{prop:bounded_coincides}
Let $(A, |\cdot|_1)$ be a seminormed algebra over $k$, in the classical sense, and $|\cdot|_2: A \to \R_+$ a multiplicative seminorm. Then, $|\cdot|_2$ is bounded with respect to $|\cdot|_1$ in the sense of \cite{BER2} if and only if it is in the sense of definition \ref{defn:bounded_seminorm}.
\end{prop}
\begin{proof} 
Suppose that $|\cdot|_2$ is bounded with respect to $|\cdot|_1$ in the sense of \cite{BER2}. This means that if $B \subset A$ is bounded for $|\cdot|_1$, \ie there exists a $\rho \in \R_+$ such that 
\[ |f|_1 \le \rho \]
for each $f \in B$, then 
\[ |f|_2 \le \rho C \]
for some constant $C$. This means that each bounded subset, with respect to $|\cdot|_1$, is mapped to a bounded subset in $\R_{\ge 0}$ by $|\cdot|_2$,
\ie $|\cdot|_2$ is a bounded map of bornological sets, when $A$ is equipped with the bornology induced by the subsets which are bounded for the seminorm $|\cdot|_1$. This shows one implication.

On the other hand, suppose that the bornology induced by $|\cdot|_2$ is coarser than the bornology induced by $|\cdot|_1$. 
This means that each bounded subset for $|\cdot|_1$ is bounded for $|\cdot|_2$. Then, we have
\[ |D_{|\cdot|_1}(0, 1^+)|_2 \subset [0, C] \]
for some $C \in \R_+$ (where $D_{|\cdot|_1}(0, 1^+)$ denotes the closed ball of radius $1$, centered in $0$ for the seminorm $|\cdot|_1$), which is equivalent to say that
\[ |f|_2 \le C \]
for each $f \in D_{|\cdot|_1}(0, 1^+)$. Given an arbitrary $f \in A$, then we can always find a $\la \in k$ such that\footnote{For non-archimedean base fields, finding such a $\la \in k$ is always possible if the norm $|\cdot|_1$ is solid, in sense of definition 2.1.10 of \cite{PGS}. By Theorem 2.1.11 of ibid. we can always find an equivalent seminorm to $|\cdot|_1$ which is solid and hence we can always suppose that $|\cdot|_1$ is solid, because the boundedness of $|\cdot|_2$ with respect to $|\cdot|_1$ only depends on the equivalence class of $|\cdot|_1$.} for any $\epsilon > 0$
\[ 1 - \epsilon \le |\la^{-1} f|_1 \le 1 \then |\la^{-1} f|_2 \le C. \]
So
\[ |f|_2 \le C |\la|_2 = C |\la| \le \frac{C}{ (1 - \epsilon)} |f|_1 \]
where we use the fact that $|\cdot|_2$ is multiplicative and proposition \ref{prop:compatible_seminorm}. Last inequality implies that $|\cdot|_2$ is bounded with respect to $|\cdot|_1$, in classical sense.
\end{proof}

\begin{cor}
Let $(A, |\cdot|_1)$ be a seminormed algebra over $k$, in the classical sense, then the bornological spectrum of $A$ coincides with the spectrum of $A$ in the sense of Berkovich. 
\end{cor}
\begin{proof}
Let $|\cdot|_2$ be an element of the Berkovich's spectrum of $(A, |\cdot|_1)$. Then $(A, |\cdot|_2)$ is a multiplicative seminormed algebra, thus by proposition \ref{prop:compatible_seminorm} it is compatible with the valuation of $k$ and by proposition \ref{prop:bounded_coincides} it is bounded in bornological sense. Therefore it is an element of the bornological spectrum of $(A, |\cdot|_1)$. On the other hand, by proposition \ref{prop:bounded_coincides} every element of the bornological spectrum of $(A, |\cdot|_1)$ defines an element of the Berkovich's spectrum, which hence coincides with the bornological spectrum.
\end{proof}

\begin{rmk}
From now on, when we talk of seminorms on a bornological algebra we always assume them to be compatible with the valuation of $k$.
\end{rmk}

\begin{rmk} \label{rmk_bounded_seminorm}
The statement of proposition \ref{prop:bounded_coincides} is false over a general normed ring.
Consider for example the ring $\Z$ equipped with the trivial valuation $|\cdot|_0$ and a $p$-adic valuation $|\cdot|_p$ for a prime $p$.
If we look at them as maps of bornological sets with respect to the bornology they pullback from $\R_{\ge 0}$, we see that both $|\cdot|_0$
and $|\cdot|_p$ induce the chaotic bornology on $\Z$. Thus, they induces the same bornology on $\Z$ and hence they are equivalent from a bornological point of view. 
Instead, if now we apply the classical definition of boundedness, in one case we see that $|n|_p \le 1$ for every $n \ne 0$ which shows that 
\[ |n|_p \le |n|_0 \]
for every $n$, hence $|\cdot|_p$ is bounded with respect to $|\cdot|_0$ in the usual sense. But on the other hand 
\[ |p^n|_p = \frac{1}{p^n} \]
can be arbitrarily near to $0$ so there does not exists a constant $C > 0$ such that
\[ |n|_0 \le C |n|_p \] 
for each $n \in \Z$.
\end{rmk}

We denote by $\cU: \bBorn(\bAlg_k) \to \bBorn$ the forgetful functor from the category of bornological algebras to the category of bornological sets.

\begin{lemma} \label{lemma:bounded_colimit}
Let $A = \underset{i \in I}\limind (A_i, |\cdot|_i)$ be a bornological multiplicatively convex algebra which satisfies the condition 
\[ \limind_{i \in I} \cU(A_i) \cong \cU(\limind_{i \in I} A_i). \] 
Then, a seminorm $\|\cdot\|: A \to \R_{\ge 0}$ is bounded if and only if $\|\cdot\| \circ \a_i$ is bounded with respect to $|\cdot|_i$ for each $i \in I$ (where $\a_i: A_i \to A$ are the canonical morphisms). 
\end{lemma}
\begin{proof} 
Let $\|\cdot\|: A \to \R_{\ge 0}$ be a bounded seminorm and let $\a_i: A_i \to A$ be the canonical morphisms. Each $\| \cdot \| \circ \a_i$ is a seminorm on $A_i$ and since it is a composition of bounded maps, it is bounded.

Now let $\|\cdot\|: A \to \R_{\ge 0}$ be a seminorm such that the map $\| \cdot \| \circ \a_i$ is bounded for every $i$. 
We have the following canonical maps of bornological sets
\[ \cU(A_i) \stackrel{\cU(\a_i)}{\to} \limind_{i \in I} \cU(A_i) \stackrel{J}{\to} \cU(\limind_{i \in I} A_i) \stackrel{\|\cdot\|}{\to} \R_{\ge 0}. \]
By the universal property of the direct limit $J \circ \|\cdot\|$ is bounded if and only if $\cU(\a_i) \circ J \circ \|\cdot\|$ are bounded for each $i$. 
Thus, the same is true for $\|\cdot\|$ if $J$ is an isomorphism of plain bornological sets.
\end{proof}

It is known that the property stated in last lemma may fail if $\underset{i \in I}\limind \cU(A_i) \not\cong \cU(\underset{i \in I}\limind A_i)$, for example one can check that $\cM(T_k^2) \not\cong \cM(T_k) \times \cM(T_k)$ where $T_k$ is the Tate algebra over the non-archimedean field $k$. Nevertheless, this lemma applies in the following interesting case.

\begin{prop} \label{prop:commutation}
Let $A = \underset{i \in I}\limind A_i$ be bornological multiplicatively convex algebra, then 
\[ \limind_{i \in I} \cU(A_i) \cong \cU(\limind_{i \in I} A_i), \] 
if the limit $\underset{i \in I}\limind A_i$ is filtered.
\end{prop}
\begin{proof} 
The forgetful functor $\bBorn(\bAlg_k) \to \bAlg_k$ and the forgetful functor $\bBorn \to \bSet$, from bornological sets to sets, commute with all limits and colimits. Thanks to the commutative square of forgetful functors
\[
\begin{tikzpicture}
\matrix(m)[matrix of math nodes,
row sep=2.6em, column sep=2.8em,
text height=1.5ex, text depth=0.25ex]
{ \bBorn(\bAlg_k) & \bAlg_k \\
  \bBorn & \bSet \\};
\path[->,font=\scriptsize]
(m-1-1) edge node[auto] {$$} (m-2-1);
\path[->,font=\scriptsize]
(m-2-1) edge node[below] {$$} (m-2-2);
\path[->,font=\scriptsize]
(m-1-1) edge node[auto] {$$} (m-1-2);
\path[->,font=\scriptsize]
(m-1-2) edge node[below] {$$} (m-2-2);
\end{tikzpicture}
\]
 it is enough to check that the functor $\bAlg_k \to \bSet$ commutes with filtered colimits. This is proved in \cite{ACC}, Example 13.2 (2).
\end{proof}

\begin{cor} \label{cor:filtered_bounded_seminorms}
Let $A = \underset{i \in I}\limind (A_i, |\cdot|_i) $
be a bornological multiplicatively convex algebra, where the inductive limit is filtered. Then, a seminorm $\|\cdot\|: A \to \R_{\ge 0}$ 
is bounded if and only if $\|\cdot\| \circ \a_i$ is bounded with respect to $|\cdot|_i$ for each $i \in I$ (where $\a_i: A_i \to A$ are the canonical morphisms). 
\end{cor}
\begin{proof} 
Thanks to proposition \ref{prop:commutation} we can apply lemma \ref{lemma:bounded_colimit} to the system $\underset{i \in I}\limind (A_i, |\cdot|_i)$ to deduce the corollary.
\end{proof}

\begin{lemma} \label{lemma:filtered_multiplicative_seminorm}
Let $A = \underset{i \in I}\limind (A_i, |\cdot|_i)$
be a multiplicatively convex bornological algebra such that $\underset{i \in I}\limind \cU(A_i) \cong \cU(\underset{i \in I}\limind A_i)$. 
A seminorm $\|\cdot\|: A \to \R_{\ge 0}$ is multiplicative if and only if $\|\cdot\| \circ \a_i$ is multiplicative for any $i$. 
\end{lemma}
\begin{proof} 
Let $\|\cdot\|: A \to \R_+$ be a multiplicative seminorm, then $\| \cdot \| \circ \a_i$ is multiplicative because $\a_i$ is an algebra morphism.

Suppose now that each $\|\cdot\| \circ \a_i$ is multiplicative. This means that
\[ (\|\cdot\| \circ \a_i)(x y) = (\|\cdot\| \circ \a_i)(x) (\|\cdot\|\circ \a_i)(y) \]
for all $x,y \in A_i$, \ie
\[ \|\a_i(x y)\| = \|\a_i(x)\| \|\a_i(y)\| \then \|\a_i(x) \a_i(y)\| = \|\a_i(x)\| \|\a_i(y)\|. \]
The condition $\underset{i \in I}\limind \cU(A_i) \cong \cU(\underset{i \in I}\limind A_i)$ ensure that any couple of elements $x,y \in A$ is in the set-theoretic image of some $A_i$, and hence the assertion of the lemma is proved.
\end{proof}

Next one is our main result about spectra of bornological algebras.

\begin{thm} \label{thm_born_spectrum}
Let $A = \underset{i \in I}\limind A_i$ be a multiplicatively convex bornological algebra for which 
\[ \limind_{i \in I} \cU(A_i) \cong \cU(\limind_{i \in I} A_i). \]
Then, $\cM(A)$ is a compact Hausdorff topological space, which is always non-empty.
\end{thm}
\begin{proof} 

Reasoning as in lemma \ref{lemma:filtered_multiplicative_seminorm} it is easy to see that a seminorm on $A$ is compatible with the valuation of $k$ if and only if its pullback on $A_i$ is a compatible seminorm over $k$ for every $i \in I$. 
Combining lemma \ref{lemma:filtered_multiplicative_seminorm} and corollary \ref{cor:filtered_bounded_seminorms} we see that to give an element $\|\cdot\| \in \cM(A)$ is equivalent to give a system of bounded multiplicative seminorms $| \cdot |_i: A_i \to \R_{\ge 0}$ such that
\[ |\cdot|_j \circ \varphi_{i,j} = |\cdot|_i \]
where $\varphi_{i,j}$ are the morphisms of the system that defines $A$. The association $A_i \mapsto \cM(A_i)$ is contravariantly functorial (see the first chapter of \cite{BER2}) hence to the system $\underset{i \in I}\limind A_i$ we can associate the projective system $\underset{i \in I}\limpro \cM(A_i)$, where the maps of the system are the continuous maps induced 
by the pullback by $\varphi_{i,j}$, as above.
So, we showed that to give a bounded multiplicative seminorm on $A$ is equivalent to give a an element 
\[ \|\cdot\| \in \limpro_{i \in I} \cM(A_i) \]
therefore there is a bijection of sets $\underset{i \in I}\limpro \cM(A_i) \cong \cM(A)$.

We endow $\underset{i \in I}\limpro \cM(A_i)$ with the projective limit topology and we show that this topology coincides with the topology of $\cM(A)$. The maps that associate  $\|\cdot\| \mapsto \|f\|$ for each $\|\cdot\| \in \underset{i \in I}\limpro \cM(A_i)$ are continuous for all $f \in A$. More precisely, fixed an $f \in A$ the map $\|\cdot\| \mapsto \|f\|$ factors
\[
\begin{tikzpicture}
\matrix(m)[matrix of math nodes,
row sep=2.6em, column sep=2.8em,
text height=1.5ex, text depth=0.25ex]
{\underset{i \in I}\limpro \cM(A_i)  \\
  \cM(A_i) & \R \\};
\path[->,font=\scriptsize]
(m-1-1) edge node[auto] {$\pi_i$} (m-2-1);
\path[->,font=\scriptsize]
(m-2-1) edge node[below] {$|\cdot|_i \mapsto |f|_i$} (m-2-2);
\path[->,font=\scriptsize]
(m-1-1) edge node[auto] {$\|\cdot\|  \mapsto \|f \|$} (m-2-2);
\end{tikzpicture}
\]
in a cofinal part of the projective limit, where $i$ is such that $f \in A_i$. The map $(|\cdot|_i \mapsto |f|_i)\circ \pi_i$ is continuous because $\cM(A_i)$ agrees with the Berkovich spectrum of $A_i$. Therefore, the identity $\underset{i \in I}\limpro \cM(A_i) \to \cM(A)$ is continuous because $\cM(A)$ is endowed with weakest topology for which all the maps $\|\cdot\| \mapsto \|f\|$ are continuous. On the other hand, the association $\| \cdot \| \mapsto |\cdot|_i \circ \a_i$, where $\a_i: A_i \to A$ are the canonical morphisms, is a continuous map $\cM(A) \to \cM(A_i)$ for each $i \in I$ and so the identity $\cM(A) \to \underset{i \in I}\limpro \cM(A_i)$ is a continuous map, by the universal property characterizing the projective limit.

To prove last claim, notice that $\cM(A_i)$ are compact Hausdorff spaces (because they are Berkovich spectra and we can apply theorem 1.2.1 of \cite{BER2}) and non-empty because $A_i$ are $k$-algebras. It follows that
$\cM(A)$ is compact Hausdorff too, because it can be described as a closed subset of a compact Haudorff space (namely the direct product $\underset{i \in I}\prod \cM(A_i)$), and it is non-empty because of a general result of A.H. Stone, cf. \cite{STO} Theorem 5 (to apply the theorem of Stone, notice that every continuous map between compact Hausdorff spaces is closed).
\end{proof}

\begin{cor} \label{cor_homeo_spectra}
Let $A$ be a multiplicatively convex bornological algebra and $B$ another algebra for which there exists
a bornological isomorphism $A \cong B$. Then, $\cM(A)$ is homeomorphic to $\cM(B)$; \ie $\cM(A)$ depends only on the isomorphism class of $A$ and not on the representation of $A$ as a direct limit of seminormed algebras.
\end{cor}
\begin{proof} 
It is an immediate consequence of the functoriality property of spectrum showed in last theorem, that it depends only on the isomorphism class of $A$.
\end{proof}

In general the association $A \mapsto \cM(A)$ is functor from the category of all bornological algebras to the category of topological spaces. Indeed, given a bounded homomorphism of bornological algebras $\phi: A \to B$ and a bounded multiplicative seminorm $|\cdot|: B \to \R_{\ge 0}$ then $|\cdot| \circ \phi$ is a bounded seminorm on $A$ because composition of bounded maps is a bounded map. Notice that $\cM$ is a contravariant functor. Next proposition describes a useful equivalent characterization of the bornological spectrum.

\begin{prop} \label{prop:born_spectrum_char}
Let $A$ be a bornological algebra. Then, $\cM(A)$ is in bijection with the equivalence classes of characters of $A$ to valued extensions of $k$.
\end{prop}
\begin{proof} 
Consider a point $x = |\cdot| \in \cM(A)$. Since $\Ker(|\cdot|)$ is a prime ideal of 
$A$, we can extend $|\cdot|$ to $A / \Ker(|\cdot|)$ and then to its field of fractions. In this way we obtain a field with a valuation whose completion is denoted $\cH(x)$.
Composing the morphisms involved, we get a bounded homomorphism $A \to \cH(x)$ which is a complete valued field which extends $k$ because $|\cdot|$ is supposed to be compatible with the valuation on $k$.

On the other hand, a character $A \to \cH$ to complete valued field, must factor through the quotient of $A$ by a prime ideal. Thus $A \to \cH$ gives an element of $\cM(A)$ by composing the valuation of $\cH$ with the character homomorphism. It is clear that the multiplicative seminorm so obtained on $A$ is compatible with the valuation on $k$ if and only if $\cH$ is a valued extension of $k$. 
Finally, two points $x = |\cdot|, x' = |\cdot|' \in \cM(A)$ are equal if and only if they define the same equivalence class of characters 
\ie if and only if they have the same kernel and there is a complete valued field $K/k$ and isometric embeddings $K \rhook \cH(x)$ and $K \rhook \cH(x')$ for which the diagram
\[
\begin{tikzpicture}
\matrix(m)[matrix of math nodes,
row sep=2.6em, column sep=2.8em,
text height=1.5ex, text depth=0.25ex]
{  & & \cH(x)  \\
 A & K \\
   & & \cH(x') \\};
\path[->,font=\scriptsize]
(m-2-1) edge node[auto] {$\chi_x$} (m-1-3);
\path[->,font=\scriptsize]
(m-2-1) edge node[below] {$\chi_{x'}$} (m-3-3);
\path[->,font=\scriptsize]
(m-2-1) edge node[auto] {} (m-2-2);
\path[->,font=\scriptsize]
(m-2-2) edge node[auto] {} (m-1-3);
\path[->,font=\scriptsize]
(m-2-2) edge node[auto] {} (m-3-3);
\end{tikzpicture}
\]
commutes.
\end{proof}

The last proposition shows that we can interpret the bornological spectrum as the set of equivalence classes of bounded character of a bornological algebra in an analogous way as one can do with the Berkovich spectrum of a ring.

\begin{exa} \label{exa_born_spectrum}
\ben
\item The Washnitzer algebra 
      \[ W_k^n = \lim_{\stackrel{\rho > 1}{\rho \in \sqrt{k^\times}}} T_k^n(\rho), \]
      fulfils all the hypothesis of theorem \ref{thm_born_spectrum} hence its bornological spectrum is non-empty, compact and Hausdorff. We will prove in the following pages that it is homeomorphic to the Berkovich spectrum of $T_k^n(1)$.
\item Using the theory of dagger affinoid algebras, developed by Grosse-Kl\"onne in \cite{GK}, we will see that to every dagger affinoid algebra $A$ we can associate its bornological spectrum.
      In this association the spectrum of $A$ is a Berkovich-like spectrum which is homeomorphic to the topological space of the germ of analytic spaces associated to $A$
      as done by Berkovich in \cite{BER}, \S 3.4. In next chapters we will study precisely the relations between dagger affinoid spaces and germs of analytic spaces in the sense of Berkovich.
\item Fr\'echet algebras are bornological algebras but they are not of multiplicatively convex type, in general.
\item The spectrum of a bornological algebra can be empty.
      We borrow an example of such an algebra from \cite{AKNA}. For any $n \in \N$, let $A_n$ be the set of analytic functions on $\C - \D(n^+)$, the complement in $\C$ of the closed disk of radius $n$. We equip $A_n$ with the usual Fr\'echet
      structure given by the family of norms 
      \[ |f|_{k, n} \doteq \sup_{n + \frac{1}{k} \le |z| \le n + k} |f(z)| \]
      for any $k \in \N$. In \cite{AKNA} is shown that the LF-algebra
      \[ A \doteq \limind_{n \to \infty} A_n. \]
       is such $\cM^t(A)$ is empty, where $\cM^t(A)$ is the topological spectrum of $A$ \ie the equivalences classes of continuous characters to $\C$. Proposition 15 of \cite{AKNA}, shows that $A$ is a normal bornological/topological vector space, so a character $A \to \C$ is bounded if and only if it is continuous which implies that the topological spectrum and the bornological spectrum of $A$ coincide and hence they are both empty.
\een
\end{exa}

In \cite{H2} Houzel introduced the following definition (cf. definition 6).

\begin{defn}
Let $A$ be a bornological algebra, the \emph{spectral radius} of $x \in A$ is defined to be the number
\[ q_A(x) \doteq \inf \l \{ |\la| | \la \in k, x \in \la A^\circ \r \}. \]
\end{defn}

Then, in \cite{H2} is shown the following proposition (cf. proposition 1).

\begin{prop}
Let $A$ be a (not necessarily commutative) multiplicatively convex bornological algebra and $x,y \in A$ such that $x y = y x$, then
\[ q_A(x + y) \le q_A(x) + q_A(y) \]
\[ q_A(x)q_A(y) \le q_A(x)q_A(y) \]
\[ q_A(x^n) = q_A(x)^n \]
and $q_A(x) < \infty$ for all $x \in A$, \ie $q_A$ is a seminorm.
\end{prop}
\begin{proof} 
\end{proof}

We would like to conclude that, at least under some suitable hypothesis, for all $f \in A$ one has $q_A(f) = \underset{x \in \cM(A)}\sup |f(x)|$.
But this is false in the following interesting case. Let 
\[ W_{\Q_p}^n = \limind_{\rho > 1} T_{\Q_p}^n(\rho) \]
be the $n$-dimensional Washnitzer over $\Q_p$ equipped with the direct limit bornology. We will see in next sections that $\cM(W_{\Q_p}^n) \cong \cM(T_{\Q_p}^n)$, for now let's take it for granted. For the sake of simplicity suppose $n = 1$ and
consider $X \in W_{\Q_p}$, the coordinate variable, then we have that
\[ |X|_{\sup} = \sup_{x \in \cM(W_{\Q_p}^n)} |X(x)| = 1. \]
But since the valuation group of $\Q_p$ is a discrete subgroup of $\R_+$ then 
\[ q_{W_{\Q_p}}(X) = \inf \{ |\la| | \la \in \Q_p, x \in \la W_{\Q_p}^\circ \} > 1 \]
because $X \notin W_{\Q_p}^\circ$ (since $X \notin T_{\Q_p}(\rho)^\circ$ for any $\rho > 1$, see lemma \ref{lemma_ind_pb}), hence $q_{W_{\Q_p}}(X) = p$. 

Thus, we introduce our version of the definition of spectral norm of a bornological m-algebra.

\index{spectral radius}
\index{spectral seminorm}
\begin{defn} \label{spec_seminorm}
Let $A$ be a bornological m-algebra then we define for all $f \in A$
\[ \rho_A(f) = \sup_{x \in \cM(A)} |f(x)| \]
and we call it the \emph{spectral radius} of $f$ and the function $\rho_A: A \to \R_{\ge 0}$ the \emph{spectral seminorm}.
\end{defn}

\begin{rmk}
Since we showed that $\cM(A)$ is compact and non-empty then $\rho_A(f) = \underset{x \in \cM(A)}\max |f(x)| < \infty$, is well defined.
\end{rmk}

\index{spectrally power-bounded element}
\begin{defn} \label{spec_power_bounded}
Let $A$ be a bornological m-algebra, then $f \in A$ is said \emph{spectrally power-bounded} if $f \in (A, \rho_A)^\circ$.
We denote the set of spectrally power-bounded elements by $A^s$.
\end{defn}

Notice that it is easy to check that $A^\circ \subset A^s$, but there is no reason to expect that an equality between the two sets holds in general. Later on we will study an important case when this inclusion is a proper inclusion.

\begin{thm} \label{thm_spec_functorial}
The association $A \mapsto (A, \rho_A)$ is a functor from the category of bornological m-algebras to the category of seminormed algebras.
\end{thm}
\begin{proof} 
Given a bounded map between complete bornological algebras which satisfy the above hypothesis
$\phi: A \to B$ then for any $f \in A$ we have that $\rho_B(\phi(f)) \le \rho_A(f)$. Because $\phi$ induces a map of spectra $\phi^*: \cM(B) \to \cM(A)$ and so 
\[ \rho_B(\phi(f)) = \sup_{x \in \cM(B)} |(f \circ \phi)(x)| = \sup_{x \in \phi^*(\cM(B))\subset \cM(A)} |f(x)| \le \sup_{x \in \cM(A)} |f(x)| = \rho_A(f), \]
which shows that $\phi: (A, \rho_A) \to (B, \rho_B)$ is bounded.
\end{proof}

\subsection{Bornological convergence and topologically nilpotent elements}

We need now to recall some definitions from \cite{H1} for introducing the concept of convergence on bornological vector spaces and related notions. Since there is no need to restrict the discussion to the case of bornological vector space of convex type, in this subsection we drop this hypothesis, although later the results will be applied only to bornological spaces of convex type.

\index{Mackey convergence}
\begin{defn}
Let $E$ be a bornological $k$-vector space and $\Phi$ a filter of subsets of $E$. We say that $\Phi$ \emph{converges to $0$ in the sense of Mackey} if
there exists a bounded subset $B \subset E$ such that for every $\la \in k^\times$ we have that $\la B \in \Phi$. We say that $\Phi$ converges to
$a \in E$ if $\Phi - a$ converges to $0$.
\end{defn}

The previous definition applies in particular to the filter associated to a sequence. 

\index{Cauchy-Mackey sequence}
\begin{defn}
Let $E$ be a bornological $k$-vector space, we say that a sequence $\{x_n\}_{n \in \N} \subset E$ is \emph{Cauchy-Mackey} if the double sequence
$\{x_n - x_m\}_{n,m \in \N}$ converges to zero in the sense of Mackey. We say that $E$ is \emph{semi-complete} if all the sequences of Cauchy-Mackey of $E$ have
a limit in $E$.
\end{defn}

A complete bornological vector space is semi-complete but the converse is not true in general, cf. \cite{H1} page 98. 

\index{topologically nilpotent element of a bornological algebra}
\begin{defn}
Let $A$ be a bornological algebra, we say that $a \in A$ is \emph{topologically nilpotent} if the sequence
\[ \{a, a^2, a^3, \ldots \} \]
converges to $0$ in the sense of Mackey.
\end{defn}

We denote by $A^{\circ \circ}$ the set of topologically nilpotent elements of $A$. We notice that $A^{\circ \circ}$ is a multiplicative sub-semigroup
of $A$ because given $a, b \in A^{\circ \circ}$ then 
\[ \lim_{n \to \infty} (a b)^n = \lim_{n \to \infty} a^n \lim_{n \to \infty} b^n = 0 \]
because the multiplication is a bounded map. We also notice that $A^{\circ \circ} \subset A^\circ$ because if $a \in A^{\circ \circ}$ it follows that there exists an $N \in \N$ such that 
\[ \{a^n\}_{n \ge N} \subset B \]
for a bounded subset $B \subset A$: Hence $\{a^n\}_{n \in \N}$ is bounded, because it is the union of a bounded subset and a finite subset. The image of a topologically nilpotent element by a bounded map is topologically nilpotent, because in general if $\phi: A \to B$
is a bounded map of bornological algebras and $\{a_n\}_{n \in \N} \subset A$ is a bornologically convergent sequence, then
\[ \lim_{n \to \infty} \phi(a_n) = \phi(\lim_{n \to \infty} a_n). \]
This shows that the association $A \mapsto A^{\circ \circ}$ is a functor from the category of bornological algebras to the category of
commutative semigroups.

\index{bornologically closed subset}
\begin{defn} \label{defn:born_closed}
Let $E$ be a bornological $k$-vector space and $X \subset E$ a subset. We say that $X$ is \emph{(bornologically) closed} in $E$ if for any 
sequence $\{f_1, f_2, ...\}$ of elements of $X$ which have a bornological limit in $E$, the limit is an element of $X$.
The closure of a subset $X \subset E$ is defined to be the intersection of all closed subsets that contain $X$.
\end{defn}

The closure of a subset is a closed subset and the intersection of any family of closed subsets is a closed subset, as usual. Indeed, one can show that the bornological closed subsets of $E$ define a topology on it, but $E$ endowed with this topology may not be a topological vector space and neither a topological group. The notion of separated bornological vector space (cf. definition \ref{defn:sep_born_space}) can be generalized to spaces of non-convex type in the following way.

\begin{defn}
Let $E$ be a bornological $k$-vector space, $E$ is said to be \emph{separated} if $\{0\}$ is the only bounded subspace of $E$.
\end{defn} 

It is easy to check that if $E$ is of convex type the definition of separatedness just given is equivalent to the one given so far.

\begin{prop} \label{prop:born_closed_subspace}
Let $E$ be a bornological $k$-vector space and $F \subset E$ a vector subspace, then $F$ is closed if and only if $E/F$ is separated.
\end{prop}
\begin{proof} 
\cite{H1}, page 50.
\end{proof}

We can restate the previous proposition by saying that closed subspaces of $E$ are in bijection with kernels of homomorphisms from $E$ to separated
bornological vector spaces. The following proposition easily follows.

\begin{prop} \label{prop_closed_subspaces}
Let $\phi: E \to F$ be a bounded morphism of bornological vector spaces and $H \subset F$ a closed subspace then $\phi^{-1}(H)$ is closed.
\end{prop}
\begin{proof}
It is enough to consider the composition $\psi: E \to F \to \frac{F}{H}$ and to notice that $\Ker(\psi) = \phi^{-1}(H)$.
\end{proof}

We end this section by giving some properties of nilpotent and power-bounded elements.

\begin{prop}
Let $A$ be a semi-complete bornological algebra and $x = 1 - y$ with $y \in A^\cc$, then $x$ is a unit in $A$ and
\[ x^{-1} = \sum_{n = 0}^\infty y^n. \]
\end{prop}
\begin{proof}
The series $\underset{n = 0}{\overset{\infty}\sum} y^n$ is convergent. Indeed, $\forall m > 0$
\[ \lim_{n \to \infty} y^n + \cdots +  y^{n + m} = \lim_{n \to \infty} y^n (1 + \cdots +  y^m) = 0 \]
because $y^n \to 0$ and the other factor is bounded. Since $A$ is semi-complete then the series defines an element of $A$. The usual argument of Neumann's series applies and it implies that the series $\underset{n = 0}{\overset{\infty}\sum} y^n$ defines an inverse of $x$.
\end{proof}

\begin{lemma} \label{lemma_ind_pb}
Let $A = \underset{i \in I}\limind A_i$ be a bornological m-algebra. Suppose that $I$ is filtered and the morphisms of the system are injective, then 
\[ A^\circ = \limind_{i \in I} A_i^\circ \]
as a set.
\end{lemma}
\begin{proof}
The image of a power-bounded element by a bounded map is a power-bounded element, as we remarked so far, therefore $\underset{i \in I} \limind A_i^\circ \subset A^\circ$. Let now $a \in A^\circ$, this means that $\{a, a^2, \cdots \}$ is a bounded subset of $A$, hence there are some $A_{i_1}, ..., A_{i_n}$ such that 
\[ \{a, a^2, \cdots \} \subset \Im(B_{i_1}) \cup \cdots \cup \Im(B_{i_n}) \]
for $B_{i_j}$ bounded in $A_{i_j}$. The limit is filtered, so we can find an index $k$ such that $i_1, ..., i_n \le k$, and this
means that there exists $a' \in A_k$ such that $\a_k(a') = a$ and $B_1', \cdots, B_n' \subset A_k$, bounded such that
$B_j' = \Im(B_{i_j})$. Since the system morphisms are injective then
\[ \{a, a^2, \cdots \} \subset B_1' \cup \cdots \cup B_n' \]
proving that $a \in A_k \subset \underset{i \in I} \limind A_i^\circ$.
\end{proof}

\begin{cor} 
Let $A$ be a bornological m-algebra, then there exists a representation of $A \cong \underset{i \in I}\limind A_i$, with $A_i$ seminormed, such that $A^\circ = \underset{i \in I}\limind A_i^\circ$.
\end{cor}
\begin{proof}
By corollary \ref{cor:fitered_m_algebra} we can always find a representation of $A$ such that \ref{lemma_ind_pb} applies. 
\end{proof}

\begin{lemma} \label{lemma:closed_subset}
Let $A = \underset{i \in I}\limind A_i$ be a bornological m-algebra which satisfies the same hypothesis of lemma \ref{lemma_ind_pb}, let $\{ X_i \subset A_i \}_{i \in I}$ be a system of closed subsets then $\underset{i \in I}\limind X_i$ is a closed subset of $A$.
\end{lemma}
\begin{proof}
Let $\{x_n\}_{n \in \N}$ be a bornologically convergent sequence in $A$, with limit $x \in A$ and with $x_n \in X_{i_n}$. By the definition of bornological convergence
this means that there exists a bounded subset $B \subset A$ such that for any $\la \in k^\times$
\[ x - x_n \in \la B \]
for all $n > N$, for some $N = N(\la) \in \N$. But to be bounded in $A$ means that there exists an $A_i$ such that $B \subset A_i$
and $B$ is bounded in $A_i$, hence $x - x_n \in A_i$ for $n > N$. Moreover, the condition 
\[ x - x_n \in \la B \]
turns out to be equivalent to Mackey convergence in $A_i$, which is equivalent to the convergence for the seminorm on $A_i$.

Now, since $x$ is an element of $A$ there must exist an index $k$ such that $x \in A_k$. Hence taking a $j \ge i, k$ we have
that in $A_j$
\[ x - x_n \in \la B \]
for all $n > N$ and since $x \in A_j$ and $x - x_n \in A_j$ then $x_n \in A_j$, and so $x_n \in X_j$. Since $X_j$ is closed in $A_j$ we have that $x \in X_j \subset \underset{i \in I}\limind X_i$.
\end{proof}

\begin{prop} \label{prop_closed_pb}
Let $A = \underset{i \in I} \limind A_i$ be a bornological m-algebra which satisfies the same hypothesis of lemma \ref{lemma_ind_pb}, then $A^\circ$ is a closed submonoid of $A$. So, if $A$ is semi-complete then also $A^\circ$ is semi-complete.
\end{prop}
\begin{proof}
It is enough to combine last two lemmas.
\end{proof}

\section{Bornological strictly convergent power-series} \label{sec:conv_pw_series}

\begin{notation}
From this section on, we will use the language halos and the conventions explained in section \ref{sec:halos}.
\end{notation}

Let $(A, |\cdot|)$ be a seminormed algebra over $k$ and consider the elements $f = \underset{j = 0}{\overset{\infty}\sum} a_j X^j \in A \ldbrack X \rdbrack$ such that
\begin{equation} \label{eqn:strict_conv}
 \lim_{i \to \infty} \sum_{n = i}^{i + m} |a_n| = 0, \forall m > 0 
\end{equation}
where $|\cdot|$ is the seminorm of $A$. Specializing this formula in the archimedean case, the request is that
\[ \sum_{i = 0}^\infty |a_n| < \infty. \]
In the non-archimedean case we have to interpret the summation symbol as a tropical sum, obtaining
\[ \lim_{n \to \infty} |a_n| = 0. \]
We denote by $A \lt X \gt \subset A \ldbrack X \rdbrack$ the subset of elements that satisfy equation (\ref{eqn:strict_conv}) and call this ring the \emph{ring of strictly convergent power series} over $A$. $A \lt X \gt$ is a seminormed ring whose seminorm is defined 
\begin{equation} \label{eqn:strict_norm} 
\|f\| \doteq \sum_{j = 0}^\infty | a_j | < \infty 
\end{equation}
for $f = \underset{j = 0}{\overset{\infty}\sum} a_j X^j \in A \lt X \gt$. Notice that following our convention on summations, formula (\ref{eqn:strict_norm}) means
\[ \|f\| = \sum_{j = 0}^\infty | a_j | < \infty  \] 
in the usual sense, if $k$ is archimedean, and 
\[ \|f\| = \max_{j \in \N} | a_j | \] 
if $k$ is non-archimedean.

This seminorm is well-defined and equips $A \lt X \gt$ with a structure of a seminormed subalgebra of $A \ldbrack X \rdbrack$. In fact,
for every $f = \underset{j = 0}{\overset{\infty}\sum} a_j X^j, \underset{k = 0}{\overset{\infty}\sum} b_k X^k \in A \lt X \gt$ we have that
\[ \|f \pm g \| \le \|f \| + \|g\| = \sum_{j = 0}^\infty |a_j| + \sum_{k = 0}^\infty |b_k| < \infty \]
showing that $f \pm g \in A \lt X \gt$. In a similar way
\[ \| f g \| = \| \sum_{ l = 0}^\infty \sum_{j + k =l} a_j b_k \| \le  \sum_{ l = 0}^\infty \sum_{j + k =l}  |a_j | |b_k | 
 = \l ( \sum_{j = 0}^\infty |a_j| \r) \l ( \sum_{k = 0}^\infty |b_k| \r ) < \infty \]
so $f g \in A \lt X \gt$.
Next is the most important (for our discussion) property of the ring of strictly convergent power series.

\begin{prop}
Let $(A, |\cdot|_A)$ and $(B, |\cdot|_B)$ be seminormed $k$-algebras with $B$ complete. For every bounded ring homomorphism
$\phi: A \to B$ and every $b \in B^\circ$ there exists a unique bounded homomorphism $\Phi: A \lt X \gt \to B$ such that $\Phi|_A = \phi$
and $\Phi(X) = b$, \ie the map 
\[ \Hom (A \lt X \gt, B) \to \Hom(A, B) \times B^\circ \]
defined by $\Phi \mapsto (\Phi|_A, \Phi(X))$ is bijective.
\end{prop}
\begin{proof} 
Let $\{a_i\}$ be a sequence in $A$ such that
\[ \sum_{i = 0}^\infty |a_i|_A < \infty \]
and $b \in B^\circ$, then the sequence 
\[ \lim_{n \to \infty} \sum_{i = 0}^n \phi(a_i) b^i \]
converges to an element in $B$ because $\phi$ is bounded, $b$ is power-bounded and $B$ is supposed to be complete. Therefore, the association
\[ \sum_{i = 0}^\infty a_i X^i  \mapsto \lim_{n \to \infty} \sum_{i = 0}^n \phi(a_i) b^i \]
is a well-defined map from $A \lt X \gt \to B$ which we denote by $\Phi$. By definition $\Phi|_A = \phi$ and $\Phi(X) = b$. Next, we show that $\Phi$ is a bounded ring homomorphism. 
We have that
\[ |\Phi(\sum_{i = 0}^\infty a_i X^i )|_B = | \sum_{i = 0}^\infty \phi(a_i) b^i |_B \le \sum_{i = 0}^\infty |\phi(a_i)|_B|b^i|_B. \]
Choose a $\sigma \in \R$ such that $|\phi(a)|_B \le \sigma |a|_A$ for all $a \in A$ and $|b^i|_B \le \sigma$ for all $i$. This choice is always possible and it follows that
\[ |\Phi(\sum_{i = 0}^\infty a_i X^i )|_B \le \sum_{i = 0}^\infty \sigma |a_i|_A \sigma = \sigma^2 \sum_{i = 0}^\infty |a_i|_A = \sigma^2 \|\sum_{i = 0}^\infty a_i X^i\|_A \]
hence $\Phi$ is bounded. The restriction of $\Phi$ to $A[X]$ is a ring homomorphism, so $\Phi$ is a ring homomorphism by extension 
by continuity.

If $\Psi: A \lt X \gt \to B$ is another bounded ring homomorphism, with $\Psi|_A = \phi$ and $\Psi(X) = b$, then by continuity
\[ \Psi(\sum_{i = 0}^\infty a_i X^i ) = \Psi(\lim_{n \to \infty}\sum_{i = 0}^n a_i X^i ) = \lim_{n \to \infty} \Psi (\sum_{i = 0}^n a_i X^i) = 
 \lim_{n \to \infty} \Phi (\sum_{i = 0}^n a_i X^i) = \Phi(\sum_{i = 0}^\infty a_i X^i). \]
\end{proof}

\begin{rmk}
In the last proposition we could have worked with the summation norm (\ie the norm $\|f\| = \underset{j = 0}{\overset{\infty}\sum} | a_j |$, with the summation symbol denoting the sum of real numbers) also for non-archimedean base fields. Indeed, this kind of algebras have been already discussed in \cite{GR}, especially in the first chapter.
The proof still works, but the algebra $A \lt X \gt$ will be equipped with a norm which does not satisfy the ultrametric triangle inequality. Even when $A = k$ then $k \lt X \gt$ would not be the Tate algebra used in non-archimedean geometry. It satisfies the analogous 
universal property of the Tate algebra in the category of all seminormed algebras over $k$, but in the category of algebras over $k$ equipped with non-archimedean seminorm the universal object described so far is the Tate algebra. This is important to point out mainly for easily compare our construction with the literature and to ensure that the dagger affinoid algebras we will introduce later will be characterized by the correct universal properties. 

A posteriori one can show that for overconvergent analytic functions one does not need to pay attention to this issue, as it is explained in section 6 of \cite{BABE}. But this is beyond the scope of the current work.
\end{rmk}

\begin{rmk}
Another fact to remark is that in the non-archimedean case the ``summation'' norm that we defined, specialized in the case when $A = k$, coincides with the 
Gauss norm, which agrees with the spectral norm of $k \lt X \gt$ and therefore it is multiplicative. If $k$ is archimedean, this is no longer true and on $k \lt X \gt$ 
the summation norm is neither multiplicative nor power-multiplicative. On $k \lt X \gt$ one can consider the spectral seminorm,
for which $k \lt X \gt$ is not complete, and then take its completion $\wtilde{k \lt X \gt}$. This new algebra satisfies the following universal property: given any
seminormed ring $A$ whose seminorm is power-multiplicative and $a \in A^\circ$, then there exists a unique map $\phi: \wtilde{k \lt X \gt} \to A$
such that $\phi(X) = a$.
We will see in section \ref{sec_ind_stein} how our approach fixes this asymmetrical behaviour obtaining a dagger affinoid algebra theory which
works naturally in archimedean and in non-archimedean settings and whose overconvergent function algebra on the unit polydisk satisfies 
the ``right'' universal property regardless on the choices that can be made to define it. 
\end{rmk}

\begin{prop} \label{prop_complete_series}
Let $A$ be a complete normed algebra, then $A \lt X \gt$ is complete.
\end{prop}
\begin{proof} 
Let 
\[ \{ f_j \} = \l \{ \sum_{i = 0}^\infty a_{i, j} X^i \r \} \]
be a Cauchy sequence in $A \lt X \gt$. By the estimate 
\[ |a_{i, j + m} - a_{i, j}| \le \| f_{j + m} - f_j \| \to 0 \]
we get that the sequences of coefficients $\{ a_{i, j}\}$ are Cauchy sequences for $j \to \infty$. Since $A$ is complete
these sequences have limits in $A$ and we denote by them $a_i$, for $i \in \N$. Consider the series
\[ f = \sum_{i = 0}^\infty a_i X^i \]
then
\[ \| f \| = \sum_{i = 0}^\infty |a_i| = \lim_{j \to \infty } \sum_{i = 0}^\infty |a_{i, j}| = \lim_{j \to \infty} \|f_j \| \]
and hence by the reverse triangle inequality and the fact that $\{ f_j \}$ is a Cauchy sequence we have that
\[ | \|f_{j + m} \| - \|f_j \|| \le \| f_{j + m} - f_j \| \to 0  \] 
so $\underset{j \to \infty}\lim \|f_j\|$ is a Cauchy sequence which converges to $\| f \|$ which is therefore finite. So, $f \in A \lt X \gt$ and $f = \underset{j \to \infty}\lim f_j$.
\end{proof}

Now, we want to generalize these facts about seminormed algebras to multiplicatively convex bornological algebras over $k$.

\index{ring of bornological strictly convergent power-series}
\begin{defn} \label{born_strict_convergent_ps}
Let $A \cong \underset{i \in I}\limind A_i$ be a multiplicatively convex bornological algebra. We define the \emph{ring of strictly convergent power-series} over $A$ as 
the ring
\[ A \lt X \gt \doteq \limind_{i \in I} A_i \lt X_i \gt \]
where the morphisms $A_i \lt X_i \gt \to A_j \lt X_j \gt$ are the natural extensions of $A_i \to A_j$ mapping $X_i \to X_j$.
\end{defn}

More clearly, let $\varphi_{i,j}: A_i \to A_j$ be a map of the system defining $A = \underset{i \in I}\limind A_i$. This map can be composed with the 
canonical injection $A_j \rhook A_j \lt X_j \gt$ obtaining a bounded map $A_i \to A_j \lt X_j \gt$. By the universal property of
$A_i \lt X_i \gt$ shown above, there exists a unique bounded map $A_i \lt X_i \gt \to A_j \lt X_j \gt$ extending $\varphi_{i,j}$ 
mapping $X_i$ to $X_j$, and we use it to define the system of definition \ref{born_strict_convergent_ps}.

\begin{rmk} 
Since the association $A \mapsto A \lt X \gt$ is manifestly functorial, definition \ref{born_strict_convergent_ps} does not depend (up to isomorphism) on the choice of the isomorphism $A \cong \underset{i \in I}\limind A_i$.
\end{rmk}

What follows is a simple lemma on filtered colimits, that we put here for the sake of clarity because we will use this argument several times in the next proofs.

\begin{lemma} \label{system_lemma}
Let $A = \underset{i \in I}\limind A_i$ be a bornological m-algebra indexed by the directed set $I$ and let $i_0 \in I$. Suppose also that the system morphism of the direct limit are injective, then
\[ \limind_{i \ge i_0} A_i \cong \limind_{i \in I} A_i. \]
\end{lemma}
\begin{proof}
Clearly there is a bounded injection 
\[ \limind_{i \ge i_0} A_i \rhook \limind_{i \in I} A_i. \]
We show that this map is also surjective. Let $a \in \underset{i \in I}\limind A_i$, then there exists $a_j \in A_j$ for some $j \in I$, such that $\a_j(a_j) = a$,
where $\a_j: A_j \to \underset{i \in I}\limind A_i$ is the canonical map. Since the colimit is filtered we have that there exists an element $k \in I$ such that $k \ge i_0$ and $k \ge j$.
This means that $A_k$ is in the subsystem obtained by taking all $i \ge i_0$ and moreover
\[ \a_j(a_j) = \a_k(\varphi_{j, k}(a_j)) = a \]
proving that $a \in \underset{i \ge i_0}\limind A_i$ as claimed.
\end{proof}

\begin{thm} \label{thm_born_ps_wds}
Let $A$ be a bornological m-algebra, then
$A \lt X \gt$ is a subalgebra of $A \ldbrack X \rdbrack$ and
\[ A \subset A[X] \subset A \lt X \gt. \]
Moreover, $A \lt X \gt$ is a bornological m-algebra, and if $A$ is complete then also $A \lt X \gt$ is complete.
\end{thm}
\begin{proof} 
We fix an isomorphism $A \cong \underset{i \in I}\limind A_i$ for a filtered monomorphic inductive system of seminormed algebras. Since all the constructions that we will use in the proof are functorial, nothing will depend on the choice of this isomorphism.
Then, $A \lt X \gt$ is clearly multiplicatively convex and it is complete if the $A_i$ are complete because then $\underset{i \in I}\limind A_i \lt X_i \gt$ are complete by proposition \ref{prop_complete_series} and the functor $A \mapsto A \lt X \gt$ preserve monomorphisms.

By definition for each $i$ we have
\[ A_i \lt X_i \gt \subset A_i \ldbrack X_i \rdbrack \] 
so we have that
\[  A \lt X \gt = \limind_{i \in I}  A_i \lt X_i \gt \subset \limind_{i \in I}  A_i \ldbrack X_i \rdbrack \]
where the system morphisms of $\underset{i \in I}\limind  A_i \ldbrack X_i \rdbrack$ are canonically defined as did above for $\underset{i \in I}\limind A_i \lt X_i \gt$ 
(and also in this case are uniquely determined). So, we are reduced to show the isomorphism
\[ \limind_{i \in I}  (A_i \ldbrack X_i \rdbrack) \cong (\limind_{i \in I} A_i) \ldbrack X \rdbrack = A \ldbrack X \rdbrack. \]
To show this isomorphism we check that $\underset{i \in I}\limind  A_i \ldbrack X_i \rdbrack$ has the same universal property of $A \ldbrack X \rdbrack$ in the 
category of $A$-algebras.  Indeed, the universal property characterizing $A \ldbrack X \rdbrack$ is the following: given any commutative $A$-algebra $A \to S$ 
and any ideal $I \subset S$ such that $S$ is complete in the $I$-adic topology and an element $x \in I$ there exists a unique morphism 
$\Phi: A \ldbrack X \rdbrack \to S$ such that
\ben
\item $\Phi$ is an $A$-algebra morphism;
\item $\Phi$ is continuous for the $(X)$-adic and $I$-adic topologies, respectively;
\item $\Phi(X) = x$. 
\een
Now, to give an $A$-algebra structure on $S$ is equivalent to give a system of $A_i$-algebra structures such that if $\varphi_{i, j}: A_i \to A_j$ is 
a morphism of the inductive system which defines $\underset{i \in I}\limind A_i$, we have the commutative diagram
\[
\begin{tikzpicture}
\matrix(m)[matrix of math nodes,
row sep=2.6em, column sep=2.8em,
text height=1.5ex, text depth=0.25ex]
{ A_i &  \\
  A_j & S \\};
\path[->,font=\scriptsize]
(m-1-1) edge node[auto] {} (m-2-2);
\path[->,font=\scriptsize]
(m-1-1) edge node[auto] {$\varphi_{i, j}$} (m-2-1);
\path[->,font=\scriptsize]
(m-2-1) edge node[auto] {} (m-2-2);
\end{tikzpicture}.
\]
This follows immediately from the universal property that characterizes $A = \underset{i \in I}\limind A_i$. By the universal property of the $A_i$-algebras $A_i \ldbrack X_i \rdbrack$, we have a unique system
of morphisms $\Phi_i: A_i \ldbrack X_i \rdbrack \to S$, with the three properties listed so far. This system of morphisms gives a map 
\[ \Phi = \limind_{i \in I} \Phi_i: \limind_{i \in I}  A_i \ldbrack X_i \rdbrack \to S \]
which is well-defined (by compatibility of the $\Phi_i$ with the system morphisms) and that is an $A$-algebra morphism.
Then, take $X_i \in A_i \ldbrack X_i \rdbrack$ and denote its image in $\underset{i \in I} \limind A_i \ldbrack X_i \rdbrack$ by the canonical morphisms 
$\a_i: A_i \ldbrack X_i \rdbrack \to \underset{i \in I}\limind A_i \ldbrack X_i \rdbrack$ with $X$. This element is independent of $i$ because if we take another
$X_k \in A_k \ldbrack X_k \rdbrack$ then there exists a $X_l \in A_l \ldbrack X_l \rdbrack$ such that $i \le l$ and $k \le l$ and so
\[ \a_i = \a_l \circ \varphi_{i,l}, \a_k = \a_l \circ \varphi_{k,l}  \]
and so 
\[ X = \a_i(X_i) = \a_l \circ \varphi_{i, l}(X_i) = \a_l (X_l) \]
and 
\[ \a_k(X_k) = \a_l \circ \varphi_{i, k}(X_k) = \a_l (X_l) = \a_i(X_i) = X. \]
Thus, we see that if $x = \Phi_i(X_i)$ for $i \in I$ then from the universal property of $A_i \ldbrack X_i \rdbrack$, we have
\[ x = \Phi_i (X_i) = \limind_{i \in I} \Phi_i( \a_i(X_i)) = \Phi(X). \]
To show that $\Phi$ is continuous it is enough to show that for each $n > 0$ there exists an $m = m(n) > 0$ such that
\[ (X^m) \subset \Phi^{-1}(I^n). \]
This follows from the fact that for each $i$ we have 
\[ (X_i^{m_i}) \subset \Phi_i^{-1}(I^n) \]
combined with the fact that $\varphi_{i_0,i}(X_{i_0}) = X_i$, for $i \ge i_0$. We can then write
\[ (X_i^{m_{i_0}}) \subset \Phi_i^{-1}(I^n) \]
for any $i \ge i_0$ for a fixed $i_0$, because
\[ \Phi_i^{-1}(I^n) \supset \varphi_{i_0,i}(\Phi_{i_0}^{-1}(I^n)) \supset \varphi_{i_0,i}( (X_{i_0}^{m_{i_0}}) ) \]
and since $\Phi_i^{-1}(I^n)$ is an ideal then $\Phi_i^{-1}(I^n) \supset (X_i^{m_{i_0}})$. By lemma \ref{system_lemma} we have that 
\[ \limind_{i \ge i_0} A_i \ldbrack X_i \rdbrack \cong \limind_{i \in I} A_i \ldbrack X_i \rdbrack \]
hence $\Phi$ is continuous.
Finally, the uniqueness of the morphism is ensured by the universal property of the direct limit.

To show that $A [X] \subset A \lt X \gt$ we proceed by showing an isomorphism similar to the previous one, \ie the following
\[ \limind_{i \in I} A_i [ X_i ] \cong (\limind_{i \in I} A_i) [ X ] = A [ X ] \]
defining the inductive system for $\underset{i \in I}\limind  A_i [ X_i ]$ always mapping $X_i \to X_j$.
So, we consider the universal property that characterizes $A [X]$, which is the following: given any ring $B$ with a ring homomorphism $\phi: A \to B$ and any 
$b \in B$ there exists a unique morphism $\Phi: A [X] \to A$ such that $\Phi(X) = b$ and $\Phi|_A = \phi$.
To give a morphism $A \to B$ is equivalent to give a compatible system of morphism $A_i \to B$ which gives a system of morphism
$\Phi_i: A_i[X_i] \to B$. These morphisms are compatible, hence we get a morphism $\Phi = \underset{i \in I}\limind \Phi_i: \underset{i \in I}\limind A_i [ X_i ] \to B$
whose restriction to $A$ is equal to $\phi$ by definition and $\Phi(X) = b$ since $\Phi_i(X_i) = b$ for any $i$.
As before the universal property of the direct limit ensure that this morphism is unique, and this proves that $A [X]$ and
$\underset{i \in I}\limind A_i [ X_i ]$ have the same universal property.
\end{proof}

\begin{thm} \label{thm_dense_poly}
If $A$ is a complete bornological m-algebra, then $A [X]$ is (bornologically) dense in $A \lt X \gt$.
\end{thm}
\begin{proof} 
We have, using the same notation of previous theorem, the isomorphism 
\[ \limind_{i \in I} A_i [X_i] \cong A [X], \]
where the system is monomorphic. For each $i$ 
\[ \ol{A_i [X_i]} \cong A_i \lt X_i \gt \]
and this induces an isomorphism
\[ \limind_{i \in I} \ol{A_i [X_i]} \cong \limind_{i \in I} A_i \lt X_i \gt = A \lt X \gt. \]
So, it is enough to check the equality
\[ \limind_{i \in I} \ol{A_i [X_i]} = \ol{ \limind_{i \in I} A_i [X_i]}. \]
In proposition \ref{prop:born_closed_subspace} we saw that a subspace of a bornological vector space $A$ is closed if and only if it is the kernel of a morphism
\[ \phi: A \to B \]
where $B$ is a separated bornological vector space, \ie $\ol{\{0\}} = \{0\}$ in $B$. In our case, we have that $A \lt X \gt$ is
a complete bornological m-algebra, which is to say that it is a filtered limit of 
$k$-Banach algebras where the morphisms $\varphi_{i, j}$ are monomorphisms. To give a morphism from $A \lt X \gt$ to a separated bornological
vector space is the same to give a system of morphisms $\phi_i: A_i \lt X_i \gt \to B$, such that $\phi_j \circ \varphi_{i, j} = \phi_i$.
Hence, all kernels $\Ker(\phi_i)$ must be closed subspaces of $A_i \lt X_i \gt$ and also
\[ \Ker(\phi) \cong \limind_{i \in I} \Ker(\phi_i), \]
because with filtered colimits commute with kernels in $\bCBorn_k$\footnote{Filtered colimits are exact and strongly exact in any elementary quasi-abelian category, see \cite{SN} proposition 2.1.16. For a proof that the category of complete bornological vector spaces is elementary quasi-abelian see lemma 3.45 of \cite{BABE}. 
}.

Now suppose that 
\[ \ol{\limind_{i \in I} A_i [X_i]} \subset A \lt X \gt = \limind_{i \in I} \ol{A_i [X_i]} \]
is a strict inclusion. Then, there must exist a linear map of bornological vector spaces
\[ f: A \lt X \gt \to B \]
with $B$ separable and $\Ker f = \ol{\underset{i \in I}\limind A_i [X_i]}$. This gives a system of maps
\[ f_i = f \circ \a_i: A_i \lt X_i \gt \to B \]
whose kernel is a closed subset that contains $\ol{A_i [X_i]}$, hence $\Ker (f_i) = A_i \lt X_i \gt$, \ie $f_i$ is the null map for all $i \in I$. 
This implies that also $f$ is the null map and finally that the closure of $\underset{i \in I}\limind A_i [X_i]$ is $A \lt X \gt$.
\end{proof}

\begin{defn}
Let $A$ be a bornological m-algebra, by induction on $n > 0$ we define
\[ A \lt X_1, ..., X_n \gt \doteq A \lt X_1, ..., X_{n-1} \gt \lt X_n \gt. \]
\end{defn}

Notice that theorem \ref{thm_born_ps_wds} immediately implies that $A \lt X_1, ..., X_n \gt$ is a bornological m-algebra and that $A \lt X_1, ..., X_n \gt$ is complete if $A$ is complete. Furthermore, it does not matter the order of the variables, for example
\[ A \lt X, Y \gt \cong A \lt Y, X \gt \]
\ie
\[ A \lt X \gt \lt Y \gt \cong A \lt Y \gt \lt X \gt. \]
We can see it using the isomorphisms
\[ A_i \lt X_i \gt \lt Y_i \gt \cong A_i \lt Y_i \gt \lt X_i \gt \]
for each $i \in I$, which give a well-defined map
\[ \limind_{i \in I} A_i \lt X_i \gt \lt Y_i \gt \cong \limind_{i \in I} A_i \lt Y_i \gt \lt X_i \gt \]
that induces the isomorphism
\[ A \lt X, Y \gt \cong A \lt Y, X \gt \]
because it is an isomorphism of systems.

\begin{prop}
Let $A \cong \underset{i \in I} \limind A_i$ be a complete bornological algebra, then
\[ ( A \lt X \gt )^\circ = \limind_{i \in I} (A_i \lt X_i \gt)^\circ \]
as monoids.
\end{prop}
\begin{proof}
We are in the hypothesis of lemma \ref{lemma_ind_pb}.
\end{proof}

For the next theorem we need a couple of lemmas.

\begin{lemma} \label{lemma:pwseries}
Let $A$ be a $k$-seminormed algebra such that $A \cong \limind_{i \in I} A_i$, with $A_i$ seminormed subalgebras of $A_i$ (\ie the morphism of the system are strict inclusions), then
\[ \limind_{i \in I} A_i \lt X_i \gt \cong A \lt X \gt \]
\ie the two definitions of $A \lt X \gt$ as a seminormed or as a bornological m-algebra coincide.
\end{lemma}
\begin{proof} 
Knowing the universal property of $A \lt X \gt$, it is enough to check that $\underset{i \in I}\limind A_i \lt X_i \gt$ satisfies the same one.
Let $\phi: A \to B$ be a bounded morphism of $k$-seminormed algebras, with $B$ complete, we obtain $\phi_i = \phi \circ \a_i: A_i \to B$,
by composing with the canonical morphisms induced by the inclusions. We get canonical morphisms
$\Phi_i: A_i \lt X_i \gt \to B$ which extends all the $\phi_i$. From this we obtain a unique morphism $\Phi: \underset{i \in I}\limind A_i \lt X_i \gt \to B$ with the required universal property.
\end{proof}

\begin{lemma}
Let $A$ be a seminormed algebra and $B$ a complete bornological m-algebra.
For every bounded morphism of rings $\phi: A \to B$ and any
power-bounded element $b \in B^\circ$ there exists a unique bounded morphism of rings $\Phi: A \lt X \gt \to B$ such that $\Phi|_A = \phi$ and $\Phi(X) = b$.
In other words, the map $\Hom (A \lt X \gt, B) \to \Hom(A, B) \times B^\circ$ defined by $\Phi \mapsto (\Phi|_A, \Phi(X))$ is bijective.
\end{lemma}
\begin{proof} 
We write $B \cong \underset{i \in I}\limind B_i$ as a monomorphic inductive limit of $k$-Banach algebras. Thus, we have that
\[ \Im(\phi) = \underset{i \in I}\limind (B_i \cap \phi(A)) \]
with the bornology induced by the inclusion $\Im(\phi) \rhook B$.
Each $B_i \cap \phi(A)$ is a normed subalgebra of $B_i$ and we have the commutative diagram
\[
\begin{tikzpicture}
\matrix(m)[matrix of math nodes,
row sep=2.6em, column sep=2.8em,
text height=1.5ex, text depth=0.25ex]
{A& & \limind_{i \in I} B_i\\
 & \limind_{i \in I} (B_i \cap \phi(A)) \\};
\path[->,font=\scriptsize]
(m-1-1) edge node[auto] {$\phi$} (m-1-3);
\path[->,font=\scriptsize]
(m-1-1) edge node[auto] {} (m-2-2);
\path[->,font=\scriptsize]
(m-2-2) edge node[auto] {} (m-1-3);
\end{tikzpicture}.
\]
From this we deduce a system of bounded homomorphisms of seminormed algebras
\[ \pi_i: \phi^{-1}(B_i \cap \phi(A)) \to B_i \cap \phi(A). \]
Moreover, $B_i \cap \phi(A)$ is a subalgebra of $\phi(A)$, hence $\phi^{-1}(B_i \cap \phi(A))$ is a subalgebra of $A$.
Now suppose that the element $b \in B^\circ$, chosen in the hypothesis, is such that $b \in B_{i_0}^\circ$, for a fixed $i_0 \in I$. 
This implies that $b \in B_i^\circ$ for all $i > i_0$, but this part of the inductive system is final in the system $\{ B_i \}_{i \in I}$ by lemma \ref{system_lemma},
therefore we can suppose $b \in B_i^\circ$ for each $i$. It follows that we can always find a unique system of maps of $k$-seminormed algebras
\[ \Pi_i: \phi^{-1}(B_i \cap \phi(A)) \lt X_i \gt \to B_i \]
with $\Pi_i|_{B_i \cap \phi(A)} = \pi_i$ and $\Pi_i(X_i) = b$.
We can apply the functor $\limind$, because we have the systems with the same diagrams, and obtain a map
\[ \Phi: \limind_{i \in I} \phi^{-1}(B_i \cap \phi(A)) \lt X_i \gt \to \limind_{i \in I} B_i \]
and, by lemma \ref{lemma:pwseries}
\[ \limind_{i \in I} \phi^{-1}(B_i \cap \phi(A)) \lt X_i \gt \cong ( \limind_{i \in I} \phi^{-1}(B_i \cap \phi(A))) \lt X \gt = A \lt X \gt. \]
Therefore, we have obtained a well-defined map
\[ \Phi: A \lt X \gt \to \limind_{i \in I} B_i \]
which is unique by construction, because all the maps of the system are uniquely determined. Moreover, $\Phi|_A = \phi$ and $\Phi(X) = b$ are clear 
by the definition of $\Phi$.
\end{proof}

\begin{thm} \label{thm_univ_prop_born_ps}
Let $A$ and $B$ be bornological m-algebras with $B$ complete.
For every bounded morphism of algebras $\phi: A \to B$ and any power-bounded element $b \in B^\circ$ there exists a unique bounded morphism 
of algebras $\Phi: A \lt X \gt \to B$ such that $\Phi|_A = \phi$ and $\Phi(X) = b$. In other words, the map 
\[ \Hom (A \lt X \gt, B) \to \Hom(A, B) \times B^\circ \]
defined by $\Phi \mapsto (\Phi|_A, \Phi(X))$ is bijective.
\end{thm}
\begin{proof} 
We can write $A = \underset{i \in I}\limind A_i$ and $B = \underset{j \in J}\limind B_j$ where the systems are monomorphic and $B_j$ are $k$-Banach algebras.
To give a morphism $\phi: A \to B$ is tantamount to give a system of morphisms $\phi_i: A_i \to B$, so by the previous lemma we get a system
of maps $\Phi_i: A_i \lt X_i \gt \to B$ of which we can take the direct limit map 
\[ \Phi: \limind_{i \in I} \Phi_i: \limind_{i \in I} A_i \lt X_i \gt \cong A \lt X \gt \to B. \]
The property $\Phi|_A = \phi$ and $\Phi(X) = b$ are inherited from the maps $\Phi_i$
\end{proof}
In next sections we will be mainly interested to the following generalization:

\begin{cor} 
Let $A$ and $B$ be bornological m-algebras with $B$ complete.
For every bounded morphism of rings $\phi: A \to B$ and any finite set of
power-bounded elements $b_1,...,b_n \in B^\circ$ there exists a unique bounded morphism of rings $\Phi: A \lt X_1, ..., X_n \gt \to B$ 
such that $\Phi|_A = \phi$ and $\Phi(X_i) = b_i$.
\end{cor}
\begin{proof}
Direct generalization of the proof of last theorem.
\end{proof}

\section{The ring of bornological overconvergent power-series} \label{sec:over_pw_series}

In order to define the ring of overconvergent power-series we describe the following generalization of definition \ref{born_strict_convergent_ps}. 
Given any seminormed $k$-algebra $A$ and a real number $\rho > 0$, we define the ring
\begin{equation} \label{eqn:polydisk_alg}
A \lt \rho^{-1} X \gt \doteq \l \{ \sum_{i = 0}^\infty a_i X^i \in A \ldbrack X \rdbrack | 
   \lim_{i \to \infty} \sum_{n = i}^{i + m} |a_n| \rho^i = 0, \forall m > 0 \r \}.
\end{equation}
This algebra satisfies the following universal property:

\begin{prop} \label{prop:rho_univ_prop}
Given any complete normed algebra $B$, any $\rho > 0$, a morphism $\phi: A \to B$ and a $f \in B$ such that 
$\la^{-1} f \in (B \otimes_k K)^\circ$ (where $K/k$ is any valued extension of $k$ with $\rho = |\la| \in |K|$), for $\la \in K$, then there exists a unique morphism
\[ \Phi: A \lt \rho^{-1} X \gt \to B \]
such that $\Phi(X) = f$ and $\Phi|_A = \phi$.
\end{prop} 
\begin{proof} 

The only thing to check is that for any sequence $(a_i)_{i \in \N}$ such that
\[ \sum_{i \in \N} |a_i| \rho^i < \infty \]
the series
\[ \sum_{i \in \N} |a_i| f^i \]
converges in $B$. But this is equivalent to say that 
\[ \max_{|\cdot| \in \cM(B)} |f| \le \rho \]
which holds by hypothesis.

\end{proof}

We can generalize (\ref{eqn:polydisk_alg}) to any bornological m-algebra, simply defining, 
for $A \cong \underset{i \in I}\limind A_i$
\[ A \lt \rho^{-1} X \gt \doteq \limind_{i \in I} A_i \lt \rho^{-1} X \gt \]
with obvious system morphisms. In the same vein we can define
\[ A \lt \rho_1^{-1} X_1, ..., \rho_n^{-1} X_n \gt \]
for any $n$-tuple of real numbers such that $\rho_i > 0$.

\begin{thm} 
Let $A$ be a bornological m-algebra, then
\ben
\item $A \lt \rho^{-1} X \gt$ is a bornological m-algebra;
\item $A[X] \subset A \lt \rho^{-1} X \gt \subset A \ldbrack X \rdbrack$;
\item $A[X]$ is bornologically dense in $A \lt \rho^{-1} X \gt$;
\item if $A$ is complete then also $A \lt \rho^{-1} X \gt$ is complete.
\een
\end{thm}
\begin{proof}
The proofs are very similar to the ones given in the previous section for the case $\rho = 1$.
\end{proof}

The next definition introduces the main object of this section. 

\index{ring of bornological overconvergent power-series}
\begin{defn} \label{defn:over_conv}
Given a bornological m-algebra $A \cong \underset{i \in I}\limind A_i$ we define the \emph{ring of overconvergent
power-series} as the bornological algebra
\begin{equation} \label{eqn:over_conv}
A \lt X \gt^\dagger \doteq \limind_{\rho > 1} A \lt \rho^{-1} X \gt
\end{equation} 
where the morphisms $A \lt \rho^{-1} X \gt \to A \lt (\rho')^{-1} X \gt$ are given by the canonical embeddings of $A \lt \rho^{-1} X \gt$
in $A \lt (\rho')^{-1} X \gt$ if $\rho > \rho'$.
Since we are calculating the colimit (\ref{eqn:over_conv}) in the category of bornoloigcal algebras, we consider on it the direct limit bornology.
\end{defn}

\begin{rmk} 
Again, the functorialiy of the association $A \mapsto A \lt X \gt^\dagger$ ensure that definition \ref{defn:over_conv} does not depend (up to isomorphism) on the choice of the isomorphism $A \cong \underset{i \in I}\limind A_i$.
\end{rmk}

\begin{thm} \label{thm:over_conv_ring}
Let $A$ be a bornological m-algebra, then $A \lt X \gt^\dagger$ is a bornological subalgebra of $A \ldbrack X \rdbrack$ and
\[ A[X] \subset A \lt X \gt^\dagger \subset A \lt X \gt. \]
Moreover, also $A \lt X \gt^\dagger$ is a bornological m-algebra.
\end{thm}
\begin{proof} 
The inclusion $ A \lt X \gt^\dagger \subset A \lt X \gt$ is deduced by the inclusions $A \lt \rho^{-1} X \gt \subset A \lt X \gt$ for any
$\rho > 1$ and also $A[X] \subset A \lt X \gt^\dagger$ is obtained by the inclusions $A[X] \subset A \lt \rho^{-1} X \gt$, for any $\rho > 1$. Finally, $A \lt X \gt^\dagger$ is by a bornological m-algebra because it is by definition a colimit of m-algebras and hence, by proposition \ref{prop:m_algebras}, it is a bornological m-algebra.

\end{proof}

\begin{thm} \label{thm:over_conv_ring_complete}
If $A$ is a complete bornological algebra, then also $A \lt X \gt^\dagger$ is complete.
\end{thm}
\begin{proof} 
Let $A \cong \underset{i \in I}\limind A_i$ with $A_i$ $k$-Banach algebras and the system morphism monomorphic. Then, by definition 
\[ A \lt \rho^{-1} X \gt = \limind_{i \in I} A_i \lt \rho^{-1} X_i \gt \]
where also $A_i \lt \rho^{-1} X_i \gt$ are $k$-Banach algebras. So,
\[ A \lt X \gt^\dagger = \limind_{\rho > 1} \limind_{i \in I} A_i \lt \rho^{-1} X_i \gt \cong \limind_{\stackrel{\rho > 1}{i \in I}} A_i \lt \rho^{-1} X_i \gt \]
where the last colimit is the composition of the two colimits. The fact that these two systems have the same colimit is a general result of category theory,
and the system is clearly filtered. So, $A \lt X \gt^\dagger$ is a filtered colimit of $k$-Banach algebras with injective system maps,
 hence it is bornologically complete.
\end{proof}

\begin{thm} \label{thm:over_conv_ring_dense}
If $A$ is a complete bornological algebra, then $A [X]$ is dense in $A \lt X \gt^\dagger$.
\end{thm}
\begin{proof} 
For any $\rho$ we have that $A [X]$ is dense in $A \lt \rho^{-1} X \gt$, hence the same reasoning we have used in the case of $A \lt X \gt$ 
(see the proof of theorem \ref{thm_dense_poly}) applies to show that $A [X]$ is bornologically dense in $A \lt X \gt^\dagger$. We omit the details.
\end{proof}

\begin{defn} \label{defn:over_conv_2}
Let $A$ be a bornological m-algebra, by induction for each $n > 0$ we define
\[ A \lt X_1, ..., X_n \gt^\dagger \doteq A \lt X_1, ..., X_{n-1} \gt^\dagger \lt X_n \gt^\dagger. \]
\end{defn}

Notice that by previous propositions $A \lt X_1, ..., X_n \gt^\dagger$ is a complete bornological algebra when $A$ is complete.
The order of the variables does not matter, for example
\[ A \lt X, Y \gt^\dagger \cong A \lt Y, X \gt^\dagger \]
\ie that
\[ A \lt X \gt^\dagger \lt Y \gt^\dagger \cong A \lt Y \gt^\dagger \lt X \gt^\dagger. \]
We can see it using the isomorphisms
\[ A_i \lt \rho^{-1} X_i \gt \lt \rho^{-1} Y_i \gt \cong A_i \lt \rho^{-1} Y_i \gt \lt \rho^{-1} X_i \gt \]
which give a well-defined map
\[ \limind_{\stackrel{\rho > 1}{i \in I}} A \lt \rho^{-1} X_i \gt \lt \rho^{-1} Y_i \gt \cong \limind_{\stackrel{\rho > 1}{i \in I}} A \lt \rho^{-1} Y_i \gt \lt \rho^{-1} X_i \gt \]
which gives the isomorphism
\[ A \lt X, Y \gt^\dagger \cong A \lt Y, X \gt^\dagger \]
since the two inductive systems are isomorphic.

The following notion will be of crucial importance for us. It will be our replacement of the notion of power-boundedness when we are working with bornological
algebras of overconvergent analytic functions. So, in the next chapter all the universal properties which will characterize our dagger affinoid algebras will be stated in terms of the concept of 
weak-power-boundedness, which we are now going to explain.

\index{weak power-bounded element}
\begin{defn}
Let $A$ be a bornological algebra over $k$ and $\ol{k}$ a completion of an algebraic closure of $k$.
We say that $f \in A$ \emph{satisfies the weak power-boundedness condition} if for any 
$\la \in \ol{k}$ with $|\la| > 1$, $\la^{-1} f \in (A \otimes_k \ol{k})^\circ$. We also say that $f$ is \emph{weakly power-bounded}.
We denote the set of weakly power-bounded elements of $A$ with $A^\ov$.
\end{defn}

\begin{rmk}
We need to use $\ol{k}$ to deal with the case when $|k^\times|$ is a discrete subgroup of $\R_+$. In that case, we would like to avoid the problems explained before definition \ref{spec_seminorm}. The reason why we use this definition of weak power-bounded element instead of using the spectral seminorm to define them is explained in remark \ref{rmk:why_def_weak_pw_bounded}.
\end{rmk}

\begin{thm} \label{thm:over_univ_prop}
Let $A$ and $B$ be bornological m-algebras and let $B$ be complete.
For every bounded morphism of algebras $\phi: A \to B$ and any element $b \in B^\ov$,
there exists a unique bounded morphism of rings $\Phi: A \lt X \gt^\dagger \to B$ such that $\Phi|_A = \phi$ and $\Phi(X) = b$. 
In other words, the map 
\[ \Hom (A \lt X \gt^\dagger, B) \to \Hom(A, B) \times B^\ov \]
defined by $\Phi \mapsto (\Phi|_A, \Phi(X))$ is bijective.
\end{thm}
\begin{proof} 
Direct consequence of the definition of $A \lt X \gt^\dagger$ and the universal properties that characterize the algebras $A \lt \rho^{-1} X \gt$ (cf. proposition \ref{prop:rho_univ_prop}).
\end{proof}

\begin{cor} 
Let $A$ and $B$ be bornological m-algebras and let $B$ be complete.
For every bounded morphism of algebras $\phi: A \to B$ and any finite set of elements $b_1,...,b_n \in B^\ov$,
there exists a unique bounded morphism of rings $\Phi: A \lt X_1, ..., X_n \gt^\dagger \to B$ such that $\Phi|_A = \phi$ and $\Phi(X_i) = b_i$.
\end{cor}
\begin{proof}
\end{proof}

\begin{rmk}
We notice that $X \notin (A \lt X \gt^\dagger)^\circ$ but $\la^{-1} X \in (A \lt X \gt^\dagger \otimes_k \ol{k})^\circ$ for any $\la \in \ol{k}$ with $|\la| > 1$.
\end{rmk}

The last remark shows that for a general bornological algebra $A^\circ \subset A^\ov$ in a strict way, and we have the following characterization.

\begin{prop} \label{prop:diamond}
Let $A$ be a bornological algebra, then 
\[ A^\ov = \bigcap_{\la \in k, |\la| > 1} \la (A \otimes_k \ol{k})^\circ \cap A = (A \otimes_k \ol{k})^\ov \cap A. \]
\end{prop}
\begin{proof}
By definition $f \in A^\ov$ if and only if $\la^{-1} f \in (A \otimes_k \ol{k})^\circ$ for all $\la \in \ol{k}, |\la| > 1$,
hence $f \in \la (A \otimes_k \ol{k})^\circ \cap A$ for all $\la \in \ol{k}, |\la| > 1$.
\end{proof}

\begin{prop}
Let $A$ be a bornological algebra, then $A^\ov$ is a multiplicative submonoid of $A$ containing $0$.
\begin{proof}
Let $f, g \in A^\ov$, then 
\[ f g \in (A \otimes_k \ol{k})^\ov(A \otimes_k \ol{k})^\ov \cap A \]
but $(A \otimes_k \ol{k})^\ov (A \otimes_k \ol{k})^\ov = (A \otimes_k \ol{k})^\ov$ since $|\ol{k}^\times|$ is dense in $\R_+$.
\end{proof}
\end{prop}

\begin{prop}
Let $\phi: A \to B$ be a morphism of bornological algebras, then $\phi(A^\ov) \subset B^\ov$.
\begin{proof}
It is easy to check that $\phi(A^\circ) \subset B^\circ$, so $\la \phi(A^\circ) \subset \la B^\circ$, and the proposition follows. 
\end{proof}
\end{prop}

\begin{cor}
The association $A \mapsto A^\ov$ is a functor from the category of bornological algebras to the category of multiplicative monoids.
\end{cor}
\begin{proof}
\end{proof}

\begin{prop} \label{prop:weak_pw_seminormed}
Let $A$ be a seminormed algebra, then $A^\circ = A^\ov$.
\begin{proof}
Suppose $k = \ol{k}$. We have to show that in a $k$-seminormed algebra 
\[ \bigcap_{\la \in k, |\la| > 1} \la A^\circ = A^\circ. \]
The inclusion $\underset{\la \in k, |\la| > 1}\bigcap \la A^\circ \supset A^\circ$ is clear because $\la A^\circ \supset A^\circ$ for any $\la \in \ol{k}$ with $|\la| > 1$. For the other inclusion, we consider 
$f \in \underset{\la \in k, |\la| > 1}\bigcap \la A^\circ$. This means that $|f|_{\sup} \le |\la|$ for any $|\la| > 1$ hence $|f|_{\sup} \le 1 \then f \in A^\circ$.

If $k \ne \ol{k}$ then by proposition \ref{prop:diamond}
\[ A^\ov = (A \otimes_k \ol{k})^\ov \cap A = \{ f \in A | |f|_{\sup} \le 1 \} = A^\circ. \]
\end{proof}
\end{prop}

\begin{prop}
Let $A$ be a bornological m-algebra, then $A^\ov$ is closed in $A$.
\begin{proof}
By proposition \ref{prop_closed_pb} we know that $A^\circ$ is closed. Hence $\la A^\circ$ is closed for any $\la \in k^\times$, because $f \mapsto \la f$ is an isomorphism of underlying bornological sets. It follows that $(A \otimes_k \ol{k})^\ov$ is closed, by applying proposition \ref{prop:diamond}. Since $A \to A \otimes_k \ol{k}$ is a bounded map then $A^\ov$ is bounded because it is the preimage of a bornologically closed subset.
\end{proof}
\end{prop}

\begin{lemma} \label{lemma:dense_injective}
Let $\phi: A \to B$ be a morphism of bornological m-algebras such that the set of elements of the form $\frac{\phi(f)}{\phi(g)}$ with $f, g \in A$
and $\phi(g)$ invertible, is dense in $B$. Then, the induced map on spectra $\phi^*: \cM(B) \to \cM(A)$ is injective.
\end{lemma}

We give two proofs of this lemma. One using a purely bornological reasoning, the other using the spectral semi-norm of $A$ and $B$. We start with the bornological proof.

\begin{proof}
First of all, if both $A$ and $B$ are seminormed $k$-algebras the result is known, see \cite{BER2} remark 1.2.2 (iii), where the result is discussed for Banach rings but it is easy to see that it holds for any seminormed $k$-algebra. We will reduce the general case to the case of seminormed algebras. We define the set
\[ A^\phi \doteq \{ f \in A | \phi(f) \in B^\times \}. \]
$A^\phi$ is a multiplicative system of $A$, so we can form the localization of $A$ with respect to $A^\phi$, which will be denoted by $A_{loc_\phi}$. By hypothesis $\phi$ has a factorization of the form
\[ A \to A_{loc_\phi} \stackrel{\tilde{\phi}}{\to} B \]
as a bare map of $k$-algebras. We can always suppose to have two representations of $A$ and $B$ as monomorphic filtered colimits of seminormed $k$-algebras indexed by the same set and that $\phi$ can be written as a system map between these systems, see \cite{BABE} remark 2.3 and proposition 2.5. Explicitly, we will write $A \cong \underset{i \in I}\limind A_i$, $B \cong \underset{i \in I}\limind B_i$, $\phi_i: A_i \to B_i$ and $\phi = \underset{i \in I}\limind \phi_i$, meaning that there exists a commutative diagram
\[
\begin{tikzpicture}
\matrix(m)[matrix of math nodes,
row sep=2.6em, column sep=2.8em,
text height=1.5ex, text depth=0.25ex]
{ A & B   \\
  \underset{i \in I}\limind A_i & \underset{i \in I}\limind B_i \\};
\path[->,font=\scriptsize]
(m-1-1) edge node[auto] {$\phi$} (m-1-2);
\path[->,font=\scriptsize]
(m-1-1) edge node[auto] {$\cong$} (m-2-1);
\path[->,font=\scriptsize]
(m-1-2) edge node[auto] {$\cong$} (m-2-2);
\path[->,font=\scriptsize]
(m-2-1) edge node[auto] {$\limind \phi_i$} (m-2-2);
\end{tikzpicture}.
\]
For each $i$ we define 
\[ (A_i)^{\phi_i} \doteq \{ f \in A_i | \phi(f) \in B_i^\times \} \]
and $(A_i)_{loc_{\phi_i}}$ the localization of $A_i$ with respect to $(A_i)^{\phi_i}$. We notice that as a monoids
\[ A^\phi \cong \limind_{i \in I} (A_i)^{\phi_i}. \]
The inclusion $\underset{i \in I}\limind (A_i)^{\phi_i} \subset A^\phi$ is obvious. To show the reverse inclusion, we pick an element $g \in A^\phi$. This means that there exists a $b \in B$ such that $\phi(g) b = 1$. Since $I$ is filtered, there must exist $i, j \in I$ such that $g \in A_i$, $b \in B_j$ and a $k \in I$ such that $k \ge i,j$. It then follows that $g \in (A_k)^{\phi_k}$ and hence $g \in \underset{i \in I}\limind (A_i)^{\phi_i}$. From this, we can deduce that
\[ A_{loc_\phi} = \l \{ \frac{f}{g} | f \in A, g \in A^\phi \r \} = \l \{ \frac{f}{g} | f \in \limind_{i \in I} A_i , g \in  \limind_{i \in I} (A_i)^{\phi_i} \r \} = \limind_{i \in I} (A_i)_{loc_{\phi_i}} \]
so we see that, as a morphism of $k$-algebras, we have a factorization
\[
\begin{tikzpicture}
\matrix(m)[matrix of math nodes,
row sep=2.6em, column sep=2.8em,
text height=1.5ex, text depth=0.25ex]
{ A_i & (A_i)_{loc_{\phi_i}} & B_i   \\};
\path[->,font=\scriptsize]
(m-1-1) edge node[auto] {$$} (m-1-2);
\path[->,font=\scriptsize]
(m-1-2) edge node[auto] {$\tilde{\phi}_i$} (m-1-3);
\path[->,font=\scriptsize, bend left=45]
(m-1-1) edge node[auto] {$\phi_i$} (m-1-3);
\end{tikzpicture}.
\]
In order to end the proof we just need to show that the map $(A_i)_{loc_{\phi_i}} \to B_i$ has dense image, because then we can apply remark 1.2.2 (iii) of \cite{BER2} to deduce that the map on spectra $\cM(B_i) \to \cM(A_i)$ is injective and so $\cM(B) \to \cM(A)$ turns out to be injective because it is a projective limit of injective maps. Hence, 
\[ B_i = \ol{\tilde{\phi}(A_{loc_\phi})} \cap B_i = \ol{ \{ x \in B | x \in B_i, x \in \tilde{\phi}(A_{loc_\phi}) \} } = \ol{ \{ x \in B_i | x \in \tilde{\phi}(\limind_{i \in I} (A_i)_{loc_{\phi_i}}) \} } = \]
\[ = \ol{ \{ x \in B_i | x \in \tilde{\phi}_i((A_i)_{loc_{\phi_i}}) \} } \]
proving the claim and the lemma.

\end{proof}

Then the argument with the spectral semi-norms.
 
\begin{proof}
 	By Theorem \ref{thm_spec_functorial} the spectral seminorm is functorial, so $\phi: A \to B$ induces a bounded morphism $\phi: (A, \rho_A) \to (B, \rho_B)$. In the previous proof we saw that as map of $k$-algebras we have a factorization
 	\[ A \to A_{loc_\phi} \stackrel{\tilde{\phi}}{\to} B. \]
 	Consider the sub-algebra $C = \tilde{\phi}(A_{loc_\phi}) \subset B$ equipped with the bornology induced by the inclusion. By hypothesis $C$ is dense in $B$ and the canonical map $B \to (B, \rho_B)$ is bounded. In particular this implies that $B \to (B, \rho_B)$ is continuous for the topology induced by the bornological convergence, so for any subset $X \subset B$ we have $\ol{X} \subset \ol{X}^{(\rho_B)}$, where $\ol{X}^{(\rho_B)}$ denotes the closure of $X$ for the topology induced on $B$ by $\rho_B$. Applying this reasoning to $C \subset B$ we obtain that the closure of $C$ with respect to the topology induced by $\rho_B$ is the whole $B$, and the lemma is proved, again.
 
\end{proof}

\begin{prop} \label{prop_dagger_spectrum}
Let $A \lt X_1, \dots, X_n \gt^\dagger$ be the ring of overconvergent power series over a complete bornological m-algebra $A \cong \underset{i \in I}\limind A_i$, then 
\[ \cM(A \lt X_1, \dots, X_n \gt^\dagger) = \bigcap_{\rho_1, ..., \rho_n > 1} \cM(A \lt \rho_1^{-1} X_1, \dots, \rho_n^{-1} X_n \gt). \]
\end{prop}
\begin{proof} 
Given two polyradii $\rho > \rho'$, the maps of the system $A \lt \rho_1^{-1} X_1, \dots, \rho_n^{-1} X_n \gt \to A \lt (\rho_1')^{-1} X_1, \dots, (\rho_n')^{-1} X_n \gt$ are bounded monomorphisms and as a consequence of theorem \ref{thm_dense_poly} they have bornologically dense image, because they contain the polynomials. Applying lemma \ref{lemma:dense_injective} we get that the induced map of spectra 
\begin{equation} \label{eqn:inj_maps}
\cM(A \lt (\rho_1')^{-1} X_1, \dots, (\rho_n')^{-1} X_n) \gt \to \cM(A \lt \rho_1^{-1} X_1, \dots, \rho_n^{-1} X_n \gt)
\end{equation}
is injective. Then, the projective system is in fact an intersection of topological spaces. The topology on $\cM(A \lt X_1, \dots, X_n \gt^\dagger)$ coincides with the topology of the intersection because the maps of equation (\ref{eqn:inj_maps}) are continuous maps between compact Hausdorff spaces, hence they are closed (and injective) maps.
\end{proof}

The following theorem gives the geometrical interpretation of weak-powerbounded elements.

\begin{thm} \label{thm:dagger_spectrum}
Let $A \lt X \gt^\dagger$ be the ring of overconvergent power series over a seminormed $k$-algebra, then 
\ben
\item $f \in (A \lt X \gt^\dagger)^\circ \iff |f(x)| \le 1$ for all $x \in \cM(A \lt \rho^{-1} X \gt)$ for a suitable $\rho > 1$;
\item $f$ is topologically nilpotent if and only if $|f(x)| < 1$ for all $x \in \cM(A \lt X \gt^\dagger)$;
\item $f \in (A \lt X \gt^\dagger)^\ov \iff |f(x)| \le 1$ for all $x \in \cM(A \lt X \gt^\dagger)$.
\een
\end{thm}
\begin{proof} 
\ben
\item By lemma \ref{lemma_ind_pb}, we know that $(A \lt X \gt^\dagger)^\circ = \underset{\rho > 1}\limind A^\circ \lt \rho^{-1} X \gt$, because $A \lt \rho^{-1} X \gt$ are seminormed 
      algebras, the system is monomorphic and $A \lt X \gt^\dagger$ is equipped with the direct limit bornology. Hence, $f \in (A \lt X \gt^\dagger)^\circ$
      is equivalent to $f \in A^\circ \lt \rho^{-1} X \gt$ for some $\rho > 1$ which is equivalent to $|f(x)| \le 1$ for all $x \in \cM(A \lt \rho^{-1} X \gt)$.
\item $f$ is topologically nilpotent in $A \lt X \gt^\dagger$ if and only if is topologically nilpotent in $A \lt \rho^{-1} X \gt$ for some $\rho > 1$, because of the definition of bornological convergence.
\item Suppose $f \in (A \lt X \gt^\dagger)^\ov$, this means that $f \in (A \lt \rho^{-1} X \gt)^\circ$ for all $\rho > 1$,
      and so $|f(x)| \le 1$ for any $x \in \cM(A \lt \rho^{-1} X \gt)$.
      This is true for any $\rho > 1$ hence $|f(x)| \le 1$ for all $x \in \cM(A \lt X \gt^\dagger)$.
\een
\end{proof}

\begin{exa} 
The main example of an element which is weakly-power-bounded but not power-bounded is $X \in A \lt X \gt^\dagger$: 
$X \in (A \lt X \gt^\dagger)^\ov$ because for any $\la \in k$ with $|\rho| > 1$ we have that $\rho^{-1} X \in (A \lt X \gt^\dagger)^\circ$.
But $X \notin (A \lt X \gt^\dagger)^\circ$ because $X \notin A \lt \rho^{-1} X \gt$ for any $\rho > 1$.
\end{exa}

\begin{rmk} \label{rmk:why_def_weak_pw_bounded}
The condition $f \in (A \lt X \gt^\dagger)^\ov$ is equivalent to $|f|_{\sup} \le 1$, where the sup is calculated on the point of the spectrum.
The difference between thinking $f$ as a weakly power-bounded element of $A \lt X \gt^\dagger$, endowed with the direct limit bornology
or as power-bounded element of $A \lt X \gt^\dagger$ endowed with its spectral norm, is that in the latter case there is no clear universal property 
that characterized $A \lt X \gt^\dagger$ in the category of seminormed algebras. 
On the contrary, working in our bornological settings we have a pretty good universal property to work with and a bornological notion of completeness to exploit.
\end{rmk}

The next proposition shows the strict link between the direct limit bornology of $A \lt X \gt^\dagger$ and the spectral norm.

\begin{prop} \label{prop:seminorm_dagger_nondagger_spectrum}
Let $A$ be a seminormed $k$-algebra, then
\[ \cM(A \lt X \gt^\dagger) \cong \cM((A \lt X \gt^\dagger, |\cdot|_{\sup})) \cong \cM(A \lt X \gt). \]
\end{prop}
\begin{proof} 
The second isomorphism is classical. $\cM(A \lt X \gt^\dagger) \cong \cM(A \lt X \gt)$ is a consequence of proposition \ref{prop_dagger_spectrum}, noticing that the norm of $A \lt X \gt$ is precisely the infimum of the family of norms of $A \lt \rho^{-1} X \gt$ on $A \lt X \gt^\dagger$. More precisely, for any $f \in A \lt X \gt^\dagger$ we can consider the norm
\[ \| f \| = \inf_{\rho > 1} |f|_\rho \]
where $|\cdot|_\rho$ is the norm of $A \lt \rho^{-1} X \gt$. Since $A \lt X \gt^\dagger = \underset{\rho > 1}\limind A \lt \rho^{-1} X \gt$ the value $\|f\|$ is well-defined and finite for all $f \in A \lt X \gt^\dagger$. The norm $\| \cdot \|$ concides with the restriction of the norm of $A \lt X \gt$ to $A \lt X \gt^\dagger$ and since $\| \cdot \|$ is the minimum of the family of the seminorms that defines the bornology of $A \lt X \gt^\dagger$, all bounded seminorms (and multiplicative in the case of the elements of the spectrum) for the bornology of $A \lt X \gt^\dagger$ must be also bounded for $\|\cdot\|$.
\end{proof}

Analogous statements of the previous proposition and theorem can be given for $A \lt X_1, \ldots, X_n \gt^\dagger$. We omit a detailed study of these cases that can be obtained by an easy induction argument.
We conclude this section by noticing that we can also define the algebras $A \lt \rho_1^{-1} X_1, \ldots, \rho_n^{-1} X_n \gt^\dagger$ whose explicit description is easy to deduce from what done up to here and therefore it is omitted. We need a concluding remark.

\begin{rmk} \label{rmk:LB_spectral_power_bounded}
The proof of theorem \ref{thm:dagger_spectrum} easily adapts to every LB-algebra. Therefore, for any LB-algebra $A$ we have that $A^s = A^\ov$.
\end{rmk}

\section{Weak power-bounded elements and algebraic monads}

Let $\phi: A \to B$ be a morphism bornological m-algebras. We saw so far that $\phi$ induces two maps of monoids
\[ \phi^\ov: A^\ov \to B^\ov, \]
\[ \phi^\circ: A^\circ \to B^\circ. \]
We study the first one that is more interesting for the geometrical study we will do in next chapters. In the aim of finding a geometrical meaning to the functor
${}^\ov$, we associate to $A$ not simply the monoid $A^\ov$ but the algebraic monad defined
\[ \Sigma_{A^\ov}(X) \doteq \l \{ \sum_{x \in X} m_x x | (m_x) \in \cR(A^\ov) \r \} \]
for any set $X$\footnote{See the appendix \ref{app_durov} for notations about monads and in particular definition \ref{defn:power_bounded_monad} for the monad associated to any multiplicative submonoid of a ring.}. We notice that $|\Sigma_{A^\ov}| = \Sigma_{A^\ov}(\bone) \cong A^\ov$ and that 
$\Sigma_{A^\ov}$ is an object which carry more information than $A^\ov$, in fact its structure depends also on the embedding $A^\ov \rhook A$ and on the additive structure of $A$. Recall that in proposition \ref{prop:weak_pw_seminormed} we showed that for seminormed $k$-algebra $A^\ov = A^\circ$. We will need the following characterization of power-bounded elements.

\begin{lemma} \label{lemma_char_power_bounded}
Let $A$ be a seminormed $\C$-algebra then 
\[ (f_i) \in \Sigma_{A^\ov}(\bn) \iff \sup_{x \in \cM(A)} \sum_{i = 1}^n |f_i(x)| \le 1. \]
\end{lemma}
\begin{proof} 
First we check that $\underset{x \in \cM(A)}\sup \underset{i = 1}{\overset{n}\sum} |f_i (x)| \le 1 \then (f_i) \in \Sigma_{A^\ov}(\bn)$. Given any $g_1, ..., g_n \in A^\ov$,
then $|g_i|_{\sup} \le 1$, hence
\[ \sup_{x \in \cM(A)} |\sum_{i = 1}^n (f_i g_i)(x) | \le \sup_{x \in \cM(A)} \sum_{i = 1}^n |f_i(x)| |g_i(x)| 
   \le \sup_{x \in \cM(A)} \sum_{i = 1}^n |f_i (x)| \le 1  \]
so $\underset{i = 1}{\overset{n}\sum} f_i g_i \in A^\ov$.

Now we suppose that $(f_i) \in \Sigma_{A^\ov}(\bn)$. Consider the elements
\[ z_i = \frac{|f_i(x)|}{f_i(x)} \in k^\circ \rhook A^\ov \]
for a fixed $x \in \cM(A)$. Since $(f_i) \in \Sigma_{A^\ov}(\bn)$,
\[ |\sum_{i = 1}^n (f_i z_i)(x)| \le 1 \]
so 
\[ | \sum_{i = 1}^n f_i(x) \frac{|f_i(x)|}{f_i(x)} | = | \sum_{i = 1}^n |f_i(x)| | = \sum_{i = 1}^n |f_i(x)| \le 1. \]
We can vary $x \in \cM(A)$ freely, so we have that
\[ \sum_{i = 1}^n |f_i(x)| \le 1, \forall x \in \cM(A) \]
which ends the proof.
\end{proof}

\begin{rmk}
If $k$ is non-archimedean, the statement of the previous lemma holds trivially.
Moreover, if $k = \R$ we can tensor $A$ with $\C$, applying the lemma and then pull back the result on $A$ because 
$A^\ov = (A \otimes_\R \C)^\ov \cap A$.
\end{rmk}

\begin{prop} \label{prop_functor_pb_monad}
The association $A \mapsto \Sigma_{A^\ov}$ is a functor from the category of $k$-seminormed algebras to the category of algebraic monads.
\begin{proof}
We have to check that a bounded $k$-algebra homomorphism $\phi: A \to B$ induces a morphism of algebraic monads 
$\phi^\ov: \Sigma_{A^\ov} \to \Sigma_{B^\ov}$. Since $\phi$ is bounded then for any $f \in A$ we have that
$|\phi(f)(x)| = |f(\phi^*(x))|$, for $x \in \cM(B)$. For $f_1, ..., f_n \in A$ we have the following equivalence of statements
\[ \sum_{i = 1}^n |f_i(x)| \le 1, \forall x \in \cM(A) \iff (f_i) \in \Sigma_{A^\ov}(\bn) \]
thanks to lemma \ref{lemma_char_power_bounded}. Hence
\[ \sum_{i = 1}^n |\phi(f_i)(x)| = \sum_{i = 1}^n |f_i(\phi^*(x))| \le 1, \forall x \in \cM(B) \]
which, again applying lemma \ref{lemma_char_power_bounded}, implies $(\phi(f_i)) \in \Sigma_{B^\ov}(\bn)$.
Thus, we have a well-defined map $\phi^\ov$ which is a monad morphism because $\phi$ is a morphism of rings and $\Sigma_{A^\ov}$,
$\Sigma_{B^\ov}$ are submonads of $\Sigma_A$ and $\Sigma_B$ respectively.
\end{proof}
\end{prop}

\begin{rmk}
If $k$ is non-archimedean the previous proposition is almost trivial and well-known.
\end{rmk}

\index{regular bornological algebra}
\begin{defn} \label{defn_born_regular}
Let $A$ be a bornological m-algebra, we say that $A$ is \emph{regular} if 
\[ A^\ov = \{ f \in A | |f|_{\sup} \le 1 \} = A^s. \]
\end{defn}

\ie a bornological m-algebra $A$ is regular if for the elements of $A$ the concept of spectrally power-boundedness coincides with the concept of 
weakly power-boundedness.
\begin{exa}
\begin{enumerate}
\item All seminormed algebras are regular, because of proposition \ref{prop:weak_pw_seminormed};
\item The main example of "non-trivial" regular bornological algebra is given by $A \lt X \gt^\dagger$ as a consequence of theorem \ref{thm:dagger_spectrum}.
\item Remark \ref{rmk:LB_spectral_power_bounded} shows that all LB-algebras are regular.
\end{enumerate}
\end{exa}

\begin{thm}
The association $A \mapsto \Sigma_{A^\ov}$ is a functor from the category of regular bornological m-algebras to the category of algebraic monads.
\begin{proof}
We write $A = \underset{i \in I}\limind A_i$, $B = \underset{j \in J}\limind B_j$ with $A_i, B_j$ seminormed algebras and $\phi: A \to B$ a bounded map. Thanks to proposition \ref{prop_functor_pb_monad} to show
that $\phi$ induces a map of algebraic monads $\phi^\ov: \Sigma_{A^\ov} \to \Sigma_{B^\ov}$ it is enough to show that $\phi$ is 
bounded with respect to the spectral seminorm of $A$ and $B$. But this is shown in theorem \ref{thm_spec_functorial}
\end{proof}
\end{thm}

\begin{rmk}
We showed in theorem \ref{thm:dagger_spectrum} that the ring of overconvergent power-series $A \lt X \gt^\dagger$ is regular, if $A$ is a $k$-seminormed algebra, because we characterized
the elements of $(A \lt X \gt^\dagger)^\ov$ as the power-bounded elements with respect to the spectral norm of $A \lt X \gt^\dagger$. Hence, we can restate
this observation by saying that the diagram of functors
\[
\begin{tikzpicture}
\matrix(m)[matrix of math nodes,
row sep=2.6em, column sep=2.8em,
text height=1.5ex, text depth=0.25ex]
{ \bR \bBorn(\bAlg_k)   & \bMon(\bSet)   \\
  \bSn(\bAlg_k) &  \\};
\path[->,font=\scriptsize]
(m-1-1) edge node[auto] {$\diamond$} (m-1-2);
\path[->,font=\scriptsize]
(m-1-1) edge node[auto] {$$} (m-2-1);
\path[->,font=\scriptsize]
(m-2-1) edge node[auto] {$\circ$} (m-1-2);
\end{tikzpicture}
\]
commutes, where $\bR \bBorn(\bAlg_k)$ denotes the category of regular bornological m-algebras over $k$, $\bSn(\bAlg_k)$ the category of seminormed algebras
over $k$ and the functor $\bR \bBorn(\bAlg_k) \to \bSn(\bAlg_k)$ is the functor which associates to any regular bornological algebra $A$
the seminormed algebra $(A, |\cdot|_{\sup})$, where $|\cdot|_{\sup}$ is the spectral seminorm of $A$.
\end{rmk}

\section{Spectra of locally convex algebras}

The aim of this section is to discuss spectra of topological algebras over any $k$ similarly as did in \cite{GOLD} and \cite{MAL} for 
topological algebras over $\C$, and to compare these spectra with our bornological spectra.
We do not go in a deep and detailed study of general topological algebras over $k$, restricting our attention on the 
multiplicatively convex. 

\index{locally convex topological algebra}
\begin{defn}
Let $A$ be a locally convex topological $k$-vector space. We call $A$ a \emph{locally convex topological algebra} if $A$ is equipped with a bilinear map
$A \times A \to A$, which is continuous and associative. We will always suppose that the multiplication is commutative and that $A$ has an identity element. A morphism of locally convex algebras is a continuous $k$-algebra homomorphism which preserves identities.
\end{defn}

\index{multiplicatively convex topological algebra}
\begin{defn}
Let $A$ be a locally convex algebra, we say that $A$ is \emph{multiplicatively convex} if
\[ A \cong \limpro_{i \in I} A_i \]
with $A_i$ seminormed algebras, and where the system morphisms are algebra morphisms.
\end{defn}

\index{spectrum of a topological algebra}
\begin{defn}
Let $A$ be a topological algebra, we define the \emph{spectrum} $\cM(A)$ of $A$ as the set of continuous multiplicative seminorms on $A$ which are compatible
with the valuation on $k$, \ie whose restriction on $k$ is equal to its valuation. We endow it with the weak topology induced by the evaluation morphism $|\cdot| \mapsto |f|$, for $|\cdot| \in \cM(A)$ and $f \in A$.
\end{defn}

\begin{prop} \label{prop:char_top_spectrum}
Let $\cS(A)$ be the set equivalence classes of continuous characters\footnote{A character $\chi: A \to K$ is a continuous $k$-algebra morphism from $A$ to a valued field extension $K/k$.}. Then, we have a bijection
\[ \cS(A) \cong \cM(A). \]
\end{prop}
\begin{proof}
The proof follows from classical argumentations. We consider a continuous character $\chi: A \to K$ to a complete valued field $K/k$. The valuation $|\cdot|: K \to \R_{\ge 0}$ induces
on $A$ a seminorm $|\cdot| \circ \chi$ which is a continuous multiplicative seminorm on $A$. Given two characters $\chi_1: A \to K_1$ and
$\chi_2: A \to K_2$, they induces two multiplicative seminorms $|\cdot|_1, |\cdot|_2: A \to \R_{\ge 0}$ and $|\cdot|_1 = |\cdot|_2$ if and only if
there exist two isometric embeddings $K \rhook K_1$, $K \rhook K_2$ and a continuous character $\chi: A \to K$ such that the following diagram commutes
\[
\begin{tikzpicture}
\matrix(m)[matrix of math nodes,
row sep=2.6em, column sep=2.8em,
text height=1.5ex, text depth=0.25ex]
{   &   & K_1 \\
  A & K &  \\
    &   & K_2 \\};
\path[->,font=\scriptsize]
(m-2-1) edge node[auto] {$\chi$} (m-2-2);
\path[->,font=\scriptsize]
(m-2-1) edge node[auto] {$\chi_1$} (m-1-3);
\path[->,font=\scriptsize]
(m-2-1) edge node[auto] {$\chi_2$} (m-3-3);
\path[->,font=\scriptsize]
(m-2-2) edge node[auto] {} (m-1-3);
\path[->,font=\scriptsize]
(m-2-2) edge node[auto] {} (m-3-3);
\end{tikzpicture},
\]
\ie if and only if $\chi_1$ and $\chi_2$ are in the same equivalence class of continuous characters. Moreover, since $A$ is a $k$-algebra we have that
$k \to A \to K$ is a morphism of complete valued field, hence it must be an isometric embedding, since the spectrum of a complete valued field is a point. 
Therefore, any continuous character $\chi$ induces a seminorm on $A$, $|\cdot| \circ \chi$, that is compatible with the valuation on $k$.

On the other hand, given a continuous multiplicative seminorm $|\cdot|: A \to \R_{\ge 0}$, then $|\cdot|^{-1}(0)$ is a closed prime ideal of $A$ and
we can extend $|\cdot|$ to a multiplicative norm on $K = \Frac \l ( \frac{A}{|\cdot|^{-1}(0)} \r )$. Taking the completion of $K$ we get
a complete valued field $\widehat{K}$ and a character $\chi: A \to \widehat{K}$. In general this character is continuous
but the embedding $k \rhook \widehat{K}$ may not be an isometric embedding. To ensure this to be true, the seminorm $|\cdot|$ must be compatible with the valuation on $k$, as it is for the elements of $\cM(A)$.

It is easy to check that the two associations described so far are inverse of each other set-theoretically.
\end{proof}

\begin{rmk}
We would like to think of points of spectra as seminorms, as in Berkovich's approach, but in this section we find it more convenient 
to think of points as characters.
\end{rmk}

Our next objective is to generalize lemma 7.1 of page 175 of \cite{MAL} to the case of a general valued field, in order to be able to express the spectrum of a multiplicatively convex topological algebra by mean of spectra of $k$-Banach algebras in a dual way compared to what we did in the case of bornological m-algebras in theorem \ref{thm_born_spectrum}. But the theory of multiplicatively convex topological algebras is far more complicated than the theory of bornological m-algebras. So, we bound ourself in discussing the few results we need. 

We say that a projective limit $A = \underset{i \in I}\limpro A_i$ is \emph{cofiltered} if the set $I$ is directed, and we say that
is \emph{cofiltrant} if $I$ is linearly ordered and maps are injective. We also use the notation $\varphi_{i, j}: A_i \to A_j$ for the 
system morphisms of the projective limit, for $i,j \in I$ with $j \le i$. We also write $\pi_i: \underset{i \in I}\limpro A_i \to A_i$ for the canonical projections from the limit.

\begin{defn}
Let $A = \underset{i \in I}\limpro A_i$ be a cofiltered locally convex topological algebra, we say that the limit is a \emph{dense} if for any
$i, j \in I$ with $i \le j$ we have that 
\[ \ol{\varphi_{j, i} (A_j)} = A_i. \]
We say that the limit is \emph{strictly dense} if it is a dense and for any $i$ the canonical map $\pi_i: \limpro A_i \to A_i$ has dense image.
\end{defn}

\begin{thm} \label{spectrum_loc_conv}
Let $A = \underset{i \in I}\limpro A_i$ be a strictly dense (and hence cofiltered) projective limit of multiplicatively convex topological algebras. Then, there is a continuous bijection of topological spaces
\[ \limind_{i \in I} \cM(A_i) \to \cM(A). \]
\end{thm}
\begin{proof}
The following proof is an adaptation of lemma 7.1, chapter V of \cite{MAL}.
By the functoriality of the spectrum, the topological space $\underset{i \in I} \limind \cM(A_i)$ is well-defined.
Given any $\phi \in \cM(A_i)$, by composition we get $\phi \circ \pi_i \in \cM(A)$. Moreover, for any $j \le i$ we also get that
\[ \phi \circ \pi_i = \phi \circ \varphi_{j, i} \circ \pi_j \]
with $\phi \circ \varphi_{j, i} \circ \pi_j \in \cM(A_j)$. This shows that the diagram
\[
\begin{tikzpicture}
\matrix(m)[matrix of math nodes,
row sep=2.6em, column sep=2.8em,
text height=1.5ex, text depth=0.25ex]
{ \cM(A) &  \\
 \cM(A_j) & \cM(A_i) \\};
\path[->,font=\scriptsize]
(m-2-1) edge node[auto] {$\pi_j^*$} (m-1-1);
\path[->,font=\scriptsize]
(m-2-2) edge node[auto] {$\varphi_i^*$} (m-2-1);
\path[->,font=\scriptsize]
(m-2-2) edge node[above] {$\pi_i^*$} (m-1-1);
\end{tikzpicture},
\]
commutes, where the $*$ stands for the opposite maps on spectra. Hence, by the definition of colimit we have a unique continuous map $h: \underset{i \in I}\limind \cM(A_i) \to \cM(A)$. We show that this map is injective.
We consider $u, v \in \underset{i \in I}\limind \cM(A_i)$ and suppose that $h(u) = h(v)$. We can assume that $u = u_i$ and $v = v_i$ for $u_i, v_i \in \cM(A_i)$, since
the limit $\underset{i \in I}\limind \cM(A_i)$ is filtered. Given a $x = (x_i) \in \underset{i \in I}\limpro A_i$ one has that
\[ h(u)(x) = u_i(x_i), \ \ h(v)(x) = v_i(x_i). \]
Therefore, we get that $u_i = v_i$ on $\pi_i(A)$ and since $\ol{\pi_i(A)} = A_i$, by the continuity of $u_i$ and $v_i$ we deduce that $u_i = v_i$ everywhere, obtaining injectivity claim.

The surjectivity is the only thing left to show. Let $u \in \cM(A)$, \ie a continuous character of $A$ compatible with $k$. Then, $u: A \to K$ is a continuous linear map between $k$-locally convex spaces, and in particular $K$ is a $k$-Banach space. This is equivalent to say (cf. proposition 1.1 of page 75 of  \cite{SCH}) that, if $\cP$ is the family of seminorms that defines the 
topology on $A$, there exists a finite set $p_1, ..., p_n \in \cP$ and a constant $c > 0$ such that
\[ \|u(x)\| \le c \sup p_i(x) \]
where $\|\cdot\|$ is the valuation on $K$. Consider the topology generated by $p_1, ..., p_n$ on $A$, and denote this topological vector space
by $\tilde{A}$. Since this is a finite family of seminorms then there exists a seminorm which induces the same topology of $p_1, ..., p_n$
and the identity $A \to \tilde{A}$ is a continuous map. Moreover, the projective limit topology of $A$ is obtained by pulling back by the canonical morphisms $\pi_i: A \to A_i$ the topologies of $A_i$. Since this system is cofiltered for any finite sets of seminorm in $\cP$ there exists a $k \in I$ such that $A_k \to \tilde{A}$ is continuous. Hence, the map $u: A \to K$ factors through $A_k$, proving the surjectivity and the theorem.
\end{proof}

In general the continuous bijection of the last theorem is not a homeomorphism, although it is in a particular, important case that we now explain.

\begin{lemma} \label{lemma:Aren_michael_frechet}
Let $A = \underset{i \in \N}\limpro A_i$ be a strictly dense (and hence cofiltered) projective limit of Banach algebras. Let $x \in \cM(A)$ and $\pi_n: A \to A_n$ the canonical projections, then
\[ x \in \cM(A_n) \iff |f(x)| \le |\pi_n(f(x))| \]
for all $f \in A$.
\end{lemma}
\begin{proof}
Let $x \in \cM(A)$ and let $\cH_x$ be the residue field of $x$. By the same argument of the last part of the proof of theorem \ref{spectrum_loc_conv}, the character
\[ f \mapsto f(x) \]
is continuous as a linear map $A \to \cH_x$ if and only if there exist a $C > 0$ and an $n \in \N$ such that
\[ |f(x)| \le C |\pi_n(f(x))|, \]
again as an application of proposition 1.1 of page 75 of \cite{SCH}. Since the seminorm which corresponds to $x$ is multiplicative we get that
\[ |f(x)|^m \le C |\pi_n(f^m(x))| \le C |\pi_n(f(x))|^m \then |f(x)| \le \sqrt[m]{C} |\pi_n(f(x))|. \]
Taking the limit for $m \to \infty$ yields $|f(x)| \le |\pi_n(f(x))|$ for all $f \in A$.

Conversely, suppose that $x \in \cM(A_n) \subset \cM(A)$. It is then obvious that $|f(x)| \le |\pi_n(f(x))|$ for all $f \in A$, concluding the proof. 
\end{proof}

\begin{prop} \label{prop:Aren_michael_frechet}
Let $A = \underset{i \in \N}\limpro A_i$ be a strictly dense (and hence cofiltered) projective limit of Banach algebras. Then, the map
\[ \limind_{i \in \N} \cM(A_i) \to \cM(A) \]
is a homeomorphism.
\end{prop}
\begin{proof}
By theorem \ref{spectrum_loc_conv} we only need to show that the inverse map $\cM(A) \to \underset{i \in \N}\limind \cM(A_i)$ is continuous. By corollary V.5.1 of \cite{MAL} it is enough to show that
\[ \limind_{i \in \N} \cM(A_i) = \bigcup_{i \in \N} \cM(A_i) \]
is an exhaustion of $\cM(A)$, \ie that $\cM(A)$ is hemicompact. More concretely, we need to show that given any compact subset $K \subset \cM(A)$ there exists an $i \in \N$ such that $K \subset \cM(A_i)$.

Consider on $A$ the family of norms 
\[ p_i \doteq \pi_i^*(|\cdot|_i), \ \ \ i \in \N \]
where $\pi_i: A \to A_i$ is the canonical morphism and $|\cdot|_i$ is the norm of $A_i$. The family $\{p_i\}_{i \in \N}$ induces on $A$ its Fr\'echet structure (the one given by the projective limit that defines $A$) and for each $i \in \N$ we have that
\[ p_i \le p_{i + 1}. \]
Let $K \subset \cM(A)$ be a compact subset. For each $n \in \N$ consider the set
\[ n K^{\circ} \doteq \bigcap_{x \in K} \{ f \in A | |f(x)| \le n \}. \]
$n K^{\circ}$ is closed in $A$ because each $\{ f \in A | |f(x)| \le n \}$ is closed. Since the map
\[ f \mapsto f(x) \]
is continuous and $K$ is compact in $\cM(A)$, then $f(K)$ is compact in $\R$ and hence
\[ A = \bigcup_{n \in \N} n K^\circ. \]
$A$ is a Baire space, and by Baire's theorem we can deduce that there exists an $n$ such that $n K^{\circ}$ has non-empty interior. Since the topology of $A$ is given by the norms $\{p_i\}_{i \in \N}$, we can then find a $j \in \N$, an $\epsilon \in |k^\times|$ and a $g \in n K^{\circ}$ such that
\[ \{ f \in A | p_j(f - g) < \epsilon \} \subset n K^\circ \iff \{ f \in A | p_j(f - g) < \frac{\epsilon}{n} \} \subset K^\circ. \]
Let $h \in A$ be such that $p_j(h) \ne 0$. Then, for any $\la \in k$ such that $|\la| = \frac{\epsilon}{ 2 p_j(h)}$
\[ g + \la h \in K^\circ \iff |(g + \la h)(x)| \le 1, \ \ \forall x \in K. \] 
The relation $|\phi(g)| \le 1$ yields
\[ |(\la h)(x)| \le 2 \then |h(x)| \le \frac{2}{\la} \iff |h(x)| \le \frac{4}{\epsilon} |p_j(h)|. \]
Since the semi-norm associated to $x$ is multiplicative
\[ |(h^n)(x)| = |h(x)|^n \le \frac{4}{\epsilon} |p_j(h^n)| \]
so taking the $n$th square root we get that
\[ |h(x)| \le \sqrt[n]{\frac{4}{\epsilon}} \sqrt[n]{ |p_j(h^n)|} \le \sqrt[n]{\frac{4}{\epsilon}}  |p_j(h)| \]
which implies
\begin{equation} \label{eqn:h_le_1}
|h(x)| \le |p_j(h)|, \forall x \in K.
\end{equation} 
On the other hand, if $p_j(h) = 0$, then for any $c \in k$
\[ g + c h \in K^\circ \]
because $g + c h$ is contained in any neighborhood of $g$ in $A$. Thus
\begin{equation} \label{eqn:h_0}
|h(x)| \le \frac{2}{|c|}, \ \ \forall c \in k \then |h(x)| = 0, \forall x \in K.
\end{equation} 
Together equations (\ref{eqn:h_le_1}) and (\ref{eqn:h_0}) imply that we can apply Lemma \ref{lemma:Aren_michael_frechet} to deduce that $x \in \cM(A_j) \subset \cM(A)$, which shows that $K \subset \cM(A_j)$ proving that $\cM(A)$ is hemicompact.
\end{proof}

\begin{cor} \label{cor:frechet}
Let $A$ be a topological multiplicatively convex Fr\'echet algebra whose topology is not given by a single norm, then $\cM(A)$ is hemicompact but not compact.
So, as bornological algebra there is no presentation of a multiplicatively convex Fr\'echet algebra of the form
\[ A \cong \limind_{i \in I} A_i \]
with $A_i$ seminormed algebras: \ie multiplicatively convex Fr\'echet algebras (in the topological sense), in general, are not multiplicatively convex in bornological sense.
\end{cor}
\begin{proof} 
The only thing to check is that for Fr\'echet algebras the notion of bornological spectrum and topological spectrum coincide. This follows from the fact that Fr\'echet spaces are normal spaces, in the sense of definition \ref{defn:normal_space}, so the boundedness of a character is equivalent to its continuity. The corollary is proved using proposition \ref{prop:char_top_spectrum} and proposition \ref{prop:born_spectrum_char} which characterize the bornological spectrum as the equivalence classes of bounded characters and the topological spectrum as the equivalence classes of continuous characters.
\end{proof}

\begin{rmk}
The same conclusions of corollary \ref{cor:frechet} were drawn in warning 3.12 of \cite{MEYER} for the case $k = \C$.
\end{rmk}

\begin{exa}
\ben
\item In complex analytic geometry the Stein spaces can be thought as spectra of Stein algebras. Indeed, Stein algebras have a canonical structure of Fr\'echet algebras as, for example, the space of analytic functions over the open disk of radius $1$ in $\C$ do have. Its spectrum has a canonical complex analytic space structure which is isomorphic to the open disk of radius $1$ in $\C$ itself. One can even show that the association of Stein algebras to Stein spaces gives
     an anti-equivalence between the category of Stein spaces and that of Stein algebras (the so-called Forster's theorem).
\item In \cite{BER2} Berkovich introduced the affine space $\A_k^n$ over a non-archimedean field $k$ as the set of bounded multiplicative seminorms over
     $k [X_1, ..., X_n]$ compatible with the valuation of $k$. He then showed that this space is homeomorphic to the space 
     \[ \bigcup_{\rho > 0} \cM(k \lt \rho^{-1} X_1, ..., \rho^{-1} X_n \gt) = \limind_{\rho > 0} \cM(k \lt \rho^{-1} X_1, ..., \rho^{-1} X_n \gt). \] 
     It seems more satisfactory, from the prospective of analytic geometry, to be able to describe this space as the spectrum of some topological algebra, equipped with some nice topology.
     In fact, we can do this noticing that we can define
     \[ \cO(\A_k^n) \doteq \limpro_{\rho > 0} k \lt \rho^{-1} X_1, ..., \rho^{-1} X_n \gt \]
     and call this algebra the \emph{algebra of entire functions} on $\A_k^n$. This algebra is a Fr\'echet algebra over $k$ and satisfies the hypothesis
     of theorem \ref{spectrum_loc_conv} hence 
     \[ \cM(\cO(\A_k^n)) \cong \limind_{\rho > 0} \cM(k \lt \rho^{-1} X_1, ..., \rho^{-1} X_n \gt) \cong \A_k^n \]
     coincides with Berkovich definition. Note also the analogy of this algebra with the Fr\'echet algebra $\cO(\C^n)$ of entire 
     holomorphic functions on $\C^n$. We shall analyse these analogies more in depth in the following chapters.
\een
\end{exa}

\index{LF-algebra}
\begin{defn}
Let $\cdots \to A_i \to A_{i + 1} \to \cdots$ be an inductive system of Fr\'echet algebras. We say that
\[ A \cong \limind_{i \in \N} A_i \]
is an \emph{LF-algebra} (where calculated in the category of topological algebras) if the colimit is filtrant.
We say that $A$ is a \emph{strict LF-algebra} if the system morphisms are  strict. We use the term LB if the Fr\'echet spaces $A_i$ are Banach.
\end{defn}

\index{LF-regular Lf-algebra}
\begin{defn}
Let $A$ be an LF-algebra, we say that $A$ is \emph{LF-regular} if the underlying topological vector space is LF-regular (\ie normal).
\end{defn}

\begin{prop} \label{prop_spectra_LB}
Let $A$ be an LF-regular LF-algebra. If $A$ can be presented as a LB-algebra $A \cong \underset{i \in \N}\limind A_i$, then there is a homeomorphism 
\[ \cM(A^b) \cong \cM(A) \cong \bigcap_{i \in \N} \cM(A_i) \]
between the bornological and topological spectrum.
\end{prop}
\begin{proof}
If $A \cong \underset{i \in I}\limind A_i$, with $A_i$ Banach and $A$ is LF-regular, then canonically $A^b \cong \underset{i \in \N}\limind A_i^b$. Moreover, for normal bornological/topological vector spaces boundedness of linear maps is equivalent to continuity, hence a character from $A$ to a valued field is continuous if and only if is bounded hence $\cM(A^b) \cong \cM(A)$. We can then apply theorem \ref{thm_born_spectrum} to the system $\underset{i \in I}\limind A_i^b$ and deduce the proposition.
\end{proof}

\begin{rmk}
We saw in example \ref{exa_born_spectrum} (4) that there exist LF-regular LF-algebras with empty spectrum. So, we cannot
relax the hypothesis that $A$ possesses a presentation as an LB-algebra in the last proposition.
\end{rmk}

\chapter{Dagger affinoid spaces} \label{chap_dagger}

In this chapter we introduce the basic ingredients of dagger analytic geometry: Dagger affinoid algebras and dagger affinoid spaces. We start the first section by reviewing some properties of the algebras of germs of analytic functions over polydisks and recasting some classical results using the language of bornological vector spaces and bornological algebras. Our main result of this first section is Theorem \ref{thm_W} which asserts that all the ideals of the algebra $W_k^n(\rho)$ are bornologically closed for the "dagger bornology" we are considering on it. This result permits  to develop the theory of dagger affinoid algebras, in the subsequent section, and to prove in our context some of the main basic properties these algebras are expected to satisfy: for example we show that every algebra morphism between dagger affinoid algebras is bounded. Then, we will devote a section to compare the algebras we defined with the classical algebras of germs of analytic functions over compact Stein subset when the base field is $\C$. After that we will introduce dagger affinoid spaces and compare them with analogous objects already present in literature: for $k$ a non-archimedean base field we compare them with the notion of germs of analytic spaces given by Berkovich in \cite{BER4} and when the base field is $\C$ we show that our dagger affinoid spaces have a canonical structure of compact Stein subset of $\C^n$, for some $n$. We conclude this chapter introducing the notion of dagger affinoid subdomain and discussing the basic properties of dagger affinoid subdomains.

\section{Overconvergent analytic functions on polydisks}

Let $k$ be an arbitrary complete non-trivially valued field. Given any polyradius $\rho = (\rho_1, ..., \rho_n) \in \R_+^n$ we define
\[ T_k^n(\rho) \doteq k \lt \rho_1^{-1} X_1, ..., \rho_n^{-1} X_n \gt \]
with $j = (j_1, ..., j_n)$ a multi-index, where we are using notation introduced in section \ref{sec:conv_pw_series}.
\index{Tate algebra}
We call it the \emph{Tate algebra} on the polycylinder of polyradius $\rho = (\rho_1, ..., \rho_n)$.
In the case of $k= \R, \C$ we can identify any element of $T_k^n(\rho)$ with an analytic function on the open polycylinder of polyradius $\rho$, 
whose algebra is denoted by $\cO(\D_k(\rho^-))$, and by the identity theorem we have that the map
\[ T_k^n(\rho) \rhook \cO(\D_k(\rho^-)) \]
is injective.
Moreover, we can consider on $T_k^n(\rho)$ the norm
\[ |\sum_{i \in \N^n} a_i X^i |_\rho \doteq \sum_{i \in \N^n} |a_i| \rho_1^{i_1} ... \rho_n^{i_n} \]
and $T_k^n(\rho)$ is a $k$-Banach algebra with respect to this norm (remember our conventions on summation symbols from section \ref{sec:halos}). 
For $f = \underset{i \in \N^n}\sum a_i X^i \in T_k^n(\rho)$ we have the estimate
\[ \|f\|_{\sup} = \sup_{x \in \D_k(\rho^+)} |f(x)| \le | f |_\rho \]
hence given a Cauchy sequence $\{f_n\} \to f$ for $|\cdot|_\rho$ we see that
\[ \|f_n - f_m \|_{\sup} \le |f_n - f_m|_\rho \to 0, \]
therefore the convergence for the norm $|\cdot|_\rho$ implies the uniform convergence. Therefore, we can deduce that, for $k = \R, \C$, any element in $T_k^n(\rho)$ defines an 
analytic function in $\D_k(\rho^-)$ which is continuous on the boundary of the polycylinder. But in general, not all analytic functions in $\D_k(\rho^-)$
which have a continuous extension on the boundary can be represented by a series which belongs to $T_k^n(\rho)$.  

\index{ring of overconvergent analytic function}
\begin{defn} \label{defn:over_convergent_W}
For any polyradius $\rho = (\rho_1, ..., \rho_n) \in \R_+^n$ we define
\[ W_k^n (\rho) = k \lt \rho_1^{-1} X_1, ..., \rho_n^{-1} X_n \gt^\dagger \doteq \limind_{r > \rho} T_k^n(r) = \bigcup_{r > \rho} T_k^n(r) \]
and we call it the ring of \emph{overconvergent analytic functions} on the polycylinder of polyradius $\rho = (\rho_1, ..., \rho_n)$ centered in zero.
Moreover, if $\rho = (1, ..., 1)$ we simply write
\[  W_k^n = k \lt X_1, ..., X_n \gt^\dagger \]
and we call it the ring of \emph{overconvergent analytic functions} on the polydisk of radius $1$.

We put on $W_k^n (\rho)$ the direct limit bornology induced by the $k$-Banach algebra structure on $T_k^n(r)$, and we call it the 
\emph{canonical (or dagger) bornology} on $W_k^n(\rho)$.
\end{defn}

\begin{rmk}
	Definition \ref{defn:over_convergent_W} is just definition \ref{defn:over_conv_2} applied with $A = k$.
\end{rmk}

We use the notation $W_k^n$ when we do not need an explicit reference to the variables, otherwise we use the notation $k \lt X_1, ..., X_n\gt^\dagger$.

\begin{rmk}
Our main objects of study will be the algebras $W_k^n(\rho)$ and their quotients. If $k$ is trivially valued then $W_k^n = T_k^n(1) = k [X_1, ..., X_n]$,
the Tate algebra of $k$ and so our constructions will reduce to Berkovich theory. So, in this work, we have nothing new to add to the theory in the case of trivially valued field.
\end{rmk}

\begin{prop}
$W_k^n(\rho)$ is a complete bornological algebra whose bornology has a countable base.
\end{prop}
\begin{proof}
$W_k^n(\rho)$ is by definition a filtrant direct limit of a family $k$-Banach algebras with monomorphic system morphisms. By choosing a monotonic strictly decreasing
sequence of polyradii $\{\rho_i\}_{i \in \N}$ such that $\rho_i \to \rho$, then 
\[ \limind_{i \in \N} T_k^n(\rho_i) \cong W_k^n(\rho). \]
Hence, the bornology of $W_k^n(\rho)$ has the following explicit description. A subset $B \subset W_k^n(\rho)$ is bounded if and only if there exists a $\rho_i$ such that $B \subset T_k^n(\rho_i)$ and $B$ is bounded in $T_k^n(\rho_i)$. This readily implies that the bornology of $W_k^n(\rho)$ has a countable base.
\end{proof}

\begin{rmk}
It is easy to show that the functor ${}^b$ commutes with \emph{strict} monomorphic direct limits with closed image, cf. Theorem 1, page 105 of \cite{H1}. In our case the maps of the system which defines $W_k^n$
are not strict, since the image of $T_k^n(\rho)$ is dense in $T_k^n(\rho')$, for $\rho > \rho'$, and if the inclusion was a strict monomorphism then it would be
surjective, because both algebras are complete. The issue of understanding if the functor ${}^b$ commutes with a countable direct limits
\[ E = \limind_{i \in \N} E_i \]
of Fr\'echet spaces is a complicated problem, for which there is a vast literature. It turns out that, even if one considers a direct limit
of Banach spaces, the colimit may not commute with ${}^b$ (see Kh\"ote \cite{KHO}, §31.6). The condition that ${}^b$ commutes with the direct limit 
which defines the LF-algebra is called \emph{regularity} in literature, and since we already reserved this word in the previous chapter (definition
\ref{defn_born_regular}) we adopt the term \emph{LF-regularity}. Finally, we call \emph{LF-dagger structure} of \emph{LF-dagger topology}
or simply \emph{dagger topology} the locally convex LF topological vector structure given to $W_k^n$ by the system
\[ (W_k^n(\rho))^t \cong (\limind_{r > \rho} T_k^n(r))^t \cong \limind_{r > \rho} T_k^n(r)^t. \]
\end{rmk}

\begin{thm} \label{thm_W_regular}
The system 
\[ (W_k^n(\rho))^t \cong (\limind_{r > \rho} T_k^n(r))^t \]
is compactoid (\ie any morphism $\limind T_k^n(r) \rhook \limind T_k^n(r')$ is a compactoid map) hence $W_k^n(\rho)$ is a LF-regular
bornological algebra with the dagger LF-structure.
\end{thm}

\begin{rmk}
The concept of compactoid subset $X \subset E$ of locally convex topological vector space is defined in \cite{PGS}, definition 3.8.1.
This concept is essential when one is dealing with functional analysis over a non-locally compact complete field. If the base field is locally compact the concept
of compactoid subset reduces to the concept of pre-compact subset. Namely, the concept of compactoid subset is precisely the generalization of the notion of pre-compact subset when the base field is not locally compact, in which case pre-compactness lose its usefulness. The word compactoid is usually used in non-archimedean functorial analysis, and since we want to have a unified treatment for archimedean and non-archimedean base fields we use the word compactoid also to mean pre-compact subsets of locally convex spaces over archimedean base fields. This choice cannot cause any confusion because the definition of compactoid subset, as the one cited in \cite{PGS}, is equivalent to the definition of pre-compact subset when it is considered over $\R$ or $\C$.
\end{rmk}

\begin{proof}
This is a known result, but scattered in literature specialized in the complex case or in the non-archimedean case.
The non-archimedean case is solved in \cite{PGS}, remark 11.4.3, and the archimedean case can be deduced applying the results of 
Vogt in \cite{VOGT} or one can find an explicit proof in \cite{HOVE}. We give here a proof for the sake of clarity.

We only deal with the case when $k$ is locally compact just to explain the main ideas of the proof. It is enough to show that for any $\rho > 1$ the canonical map $T_k^n(\rho) \to T_k^n$ is compact, \ie it maps the unit ball of $T_k^n(\rho)$ to a compact subset of $T_k^n$. We consider the unit ball
\[ B_\rho = \{ f \in T_k^n(\rho) | \| f \|_\rho \le 1 \}. \]
For any $f \in B_\rho$ we can write
\[ f = \sum_{i = 0}^\infty a_i \rho^{-i} X^i \]
for a summable sequence $\underset{i = 0}{\overset{\infty}\sum} |a_i| \le 1$ . 
Consider a sequence
\[ \{ f_j \}_{j \in \N} = \l \{ \sum_{i = 0}^\infty a_{i,j} \rho^{-i} X^i \r \}_{j \in \N} \subset B_\rho. \]
We have to check that $\{f_j\}_{j \in \N}$ has a convergent subsequence in $T_k^n$. The sequences $a_{i,j}$, for fixed $i$ and $j \to \infty$, are bounded
sequences (bounded by $1$), so we can extract a convergent subsequence $\tilde{a}_{i, j}$ for any $i$. Therefore, we can define a subsequence
\[ \{ \tilde{f}_j \}_{j \in \N} = \l \{ \sum_{i = 0}^\infty \tilde{a}_{i,j} \rho^{-i} X^i \r \}_{j \in \N} \]
because we can always choose indexes in the appropriate way (continuing to extract convergent subsequence from bounded ones). If we denote
$a_i = \underset{j \to \infty}\lim \tilde{a}_{i , j}$, then $\tilde{f}_j \to f = \underset{i = 0}{\overset{\infty}\sum} a_i \rho^{-i} X^i$ in $T_k^n$ and
$f$ is a well-defined element of $T_k^n$ since $|a_i| \le 1$ for all $i$. Hence, for the norm of $T_k^n$ we have
\[ \|  \sum_{i = 0}^\infty a_i \rho^{-i} X^i \| = \sum_{i = 0}^\infty |a_i|\rho^{-i} \le \sum_{i = 0}^\infty \rho^{-i} = 
   \frac{\rho}{1 - \rho} < \infty, \]
and this shows that the image of $B_\rho$ is compact in $T_k^n$, so the system $(\underset{r > \rho}\limind T_k^n(r))^t$ is compact.
\end{proof}

The previous theorem asserts the following relations of commutation of colimits with functors: if we equip $W_k^n(\rho)$ with the dagger 
bornology then
\[ ((W_k^n(\rho))^t)^b = ((\limind_{r > \rho} T_k^n(r))^t)^b \cong (\limind_{r > \rho} T_k^n(r)^t)^b \cong \limind_{r > \rho} (T_k^n(r)^t)^b \cong W_k^n. \]
The preceding theorem has the following remarkable consequences.

\begin{cor} \label{cor:W_regular}
$W_k^n(\rho)^t$ is a Hausdorff, nuclear, reflexive, complete, locally complete LF-algebra.
\end{cor}
\begin{proof}
These properties are common to all compactiod inductive limits: for the non-archimedean version see \cite{PGS} 11.3.5;
for the archimedean version see the paragraph 8.5 in \cite{PCB}. Notice that in the archimedean case for $W_k^n(\rho)^t$ to be nuclear is a stronger requirement with respect to require that $(\underset{r > \rho}\limind T_k^n(r))^t$ is a compact system. One can easily show that the maps of this system are in fact nuclear, see for example \cite{HOVE}.
\end{proof}

We notice that $W_k^n(\rho)$ is a regular bornological algebra in the sense of definition \ref{defn_born_regular}, as a consequence of theorem \ref{thm:dagger_spectrum}. So, we can apply proposition \ref{prop_dagger_spectrum} to deduce the homeomorphism
\[ \cM(W_k^n(\rho)) \cong \cM(T_k^n(\rho)). \]
Moreover, if $k$ is algebraically closed we also have
\[ \Max(W_k^n(\rho)) \cong \Max(T_k^n(\rho)) \cong \D_k(\rho^+). \]

In the following lemmas and theorems we will use results established in the appendix \ref{app_coeff_wise} to which we refer for the notation and terminology we use. We now start to prove results for $W_k^n$ and strict dagger affinoid algebras and we will see later how to generalize them to $W_k^n(\rho)$ and non-strict dagger affinoid algebras.

\begin{lemma} \label{lemma:closed_born_1}
Let $\ol{k}$ be an algebraic closure of $k$, then for any $M \in \Max(W_k^n)$ there exists an injection
\[ (W_k^n)_M \rhook \ol{K}_n(\ol{M}) \]
where $\ol{K}_n(\ol{M})$ is the $\ol{k}$ stellen-algebra centred in $\ol{M} \in \D_{\ol{k}}((1, ..., 1)^+)$ which is a point over $M$, and $(W_k^n)_M$ is the localization of $W_k^n$ at $M$. 
\end{lemma}
\begin{proof}
We saw so far that
\[ \Max(W_k^n) \cong \limpro_{\rho > 1} \Max(T_k^n(\rho)) \cong \Max(T_k^n( (1,...,1) )). \]
If $\ol{k}$ is an algebraic closure of $k$ then 
\[ \Max(T_k^n( (1,...,1) )) \cong \frac{\D_{\ol{k}}((1, ..., 1)^+)}{\Gal(\ol{k}/k)} \]
giving to $\frac{\D_{\ol{k}}((1, ..., 1)^+)}{\Gal(\ol{k}/k)}$ the quotient topology. 
It is clear that $\ol{K}_n(\ol{M}) \cong \cO_{\ol{k}^n, \ol{M}}$, the ring of germs of analytic functions at $\ol{M}$, 
and the identity theorem implies that
\[ (W_k^n)_M \rhook K_n(M) \rhook \ol{K}_n(\ol{M}) \]
since any element of $(W_k^n)_M$ can be identified with an analytic function in a neighborhood of $\ol{M}$.
\end{proof}

In the next lemma and in next theorem when we consider $M \in \Max(W_k^n)$ we also suppose to fix a $\ol{M} \in \D_{\ol{K}}^n(1^+)$ which is mapped on $M$ by the canonical projection obtained by quotienting with the Galois action on $\D_{\ol{K}}^n(1^+)$.

\begin{lemma} \label{lemma:closed_born_2}
Let $M \in \Max(W_k^n)$ and $I \subset (W_k^n)_M$ be an ideal, then there exists an ideal $J \subset \ol{K}_n(\ol{M})$ such that
\[ I = (W_k^n)_M \cap J. \]
\end{lemma}
\begin{proof}
The ring homomorphism $(W_k^n)_M \rhook \ol{K}_n(\ol{M})$ is a local homomorphism of Noetherian rings both of which have dimension $n$. This means that there exists
a set of generators of the maximal ideal of $\ol{K}_n(\ol{M})$, $\fm = (f_1, ..., f_n)$ such that their preimages generate the maximal ideal
of $(W_k^n)_M$. Since the morphism is injective, this implies that $(W_k^n)_M$ contains a set of generators of $\fm$ (we call it again $f_1, \ldots, f_n$) and that the 
maximal ideal of $(W_k^n)_M$ is generated by $(f_1, ..., f_n)$.

Any ideal of $(W_k^n)_M$ and of $\ol{K}_n(\ol{M})$ is of the form $(P_1, ..., P_n)\ol{K}_n(\ol{M})$ (resp. $(P_1, ..., P_n) (W_k^n)_M)$ for some polynomials  $P_i \in \ol{K}_n(\ol{M})[f_1, ..., f_n]$ (resp. $P_i \in (W_k^n)_M[f_1, ..., f_n]$) since the local rings are regular, (cf. \cite{DIU}, Theorem 4.1), and therefore
\[ (P_1, ..., P_n) \ol{K}_n(\ol{M}) \cap (W_k^n)_M = (P_1, ..., P_n) (W_k^n)_M. \]
\end{proof}

The next theorem is our main result about $W_k^n$.

\begin{thm} \label{thm_W}
$W_k^n$ is Noetherian, factorial and every ideal of $W_k^n$ is closed for the dagger bornology.
\begin{proof} 
That $W_k^n$ is Noetherian and factorial is well-known: in \cite{FRI} (cf. theorem I.9) there is the first proof that $W_k^n$ is Noetherian and in \cite{DALES} a proof that $W_k^n$ is factorial (cf. theorem 1), for the archimedean case; for the non-archimedean case one can see the first section of \cite{GK}.  So, we have only to show that all the ideals are bornologically closed.

For any $M \in \Max(W_k^n)$ we have that $W_k^n \rhook \cO_{\ol{k}^n, \ol{M}} \cong \ol{K}_n(\ol{M})$ and by lemma \ref{lemma:closed_born_1} we have the factorization
\[ W_k^n \rhook (W_k^n)_M \rhook \ol{K}_n(\ol{M}). \]
If we put on $W_k^n$ the dagger bornology and on $\ol{K}_n(\ol{M})$ the inductive limit bornology given by the isomorphism
\[ \ol{K}_n(\ol{M}) \cong \limind_{\rho > 0} T_{\ol{k}}^n(\rho)[\ol{M}] \]
(where $T_{\ol{k}}^n(\rho)[\ol{M}]$ is the Tate algebra over $\ol{k}$ with polyradius $\rho$ centred in $\ol{M}$), we obtain that $W_k^n \to \ol{K}_n(\ol{M})$ is bounded because of the maximum modulus principle. Proposition \ref{prop_coeff_bound_univ_prop} implies that the 
identity $\ol{K}_n(\ol{M}) \to \ol{K}_n(\ol{M})$ is bounded if we put on the domain the just mentioned bornology and the bornology of 
coefficientwise boundedness on the codomain, see appendix \ref{app_coeff_wise} for the definition of the latter bornology. 
Hence, we get a bounded map $W_k^n \to \ol{K}_n(\ol{M})$, which factors through the inclusion
$(W_k^n)_M \rhook \ol{K}_n(\ol{M})$. Therefore, $W_k^n \rhook (W_k^n)_M$ is bounded if we put on $(W_k^n)_M$ the bornolgy induced by the coefficientwise
boundedness bornology on $\ol{K}_n(\ol{M})$.
Moreover, all the ideals of $(W_k^n)_M$ are closed for the bornology of coefficientwise boundedness since by lemma \ref{lemma:closed_born_2} all the 
ideals of $(W_k^n)_M$ are of the form $(W_k^n)_M \cap J$ for $J \subset \ol{K}_n(\ol{M})$ and all the ideals of $\ol{K}_n(\ol{M})$ are bornologically closed, by proposition \ref{prop:coeff_closed_ideals}.

Therefore, the preimage of any ideal of $(W_k^n)_M$ in $W_k^n$ is bornologically closed.
Since $W_k^n$ is Noetherian then any $I \subset W_k^n$ satisfies the relation
\[ I = \bigcap_{M \in \Max(W_k^n)} (I_M \cap W_k^n) = \bigcap_{M \in \Max(W_k^n)} \phi_M^{-1}(I_M) \]
where $\phi_M: W_k^n \to (W_k^n)_M$ is the canonical injection and $I_M$ is the extension of $I$ to $(W_k^n)_M$. So $I$ is an intersection of bornologically closed subsets of $W_k^n$, hence it is a bornologically closed subset.

\end{proof}
\end{thm}

\begin{rmk}
Notice that the same argument does not apply to $T_k^n(\rho)$ because the ideals of $T_k^n(\rho)$ are not finitely generated (for $k = \R, \C$),
and so not all localization of $T_k^n(\rho)$ are Noetherian.
Furthermore, it is known that not all ideals of $T_k^n(\rho)$ are closed, neither the finitely generated nor the principal ones.
In fact, one can show that a $\R$-Banach algebra on which every principal ideal is closed is a division algebra.
\end{rmk}

\begin{rmk}
For $k$ non-archimedean one could deduce previous theorem by noticing that all the ideals of $T_k^n(\rho)$ are closed for the topology induced by the norm, and this readily implies that all ideals of $W_k^n$ are closed for the direct limit bornology. As explained in the previous remark, there is no hope to find an analogous argumentation for archimedean base fields, and the advantage of theorem \ref{thm_W} is that it applies for any base field.
\end{rmk}

\begin{prop} \label{lattice_ideal}
Let $\cI_{T_k^n(\rho)}$ denotes the lattice of closed ideals of $T_k^n(\rho)$ and $\cI_{W_k^n}$ the lattice of ideals of $W_k^n$, then
\[ \cI_{W_k^n} \cong \limpro_{\rho > 1} \cI_{T_k^n(\rho)} \]
where the maps of the system are given by
\[ I_\rho \mapsto I_\rho \cap T_k^n(\rho') \]
for $\rho' > \rho$.
\end{prop}
\begin{proof} 
We can define a map $\sP(W_k^n) \to \sP(T_k^n(\rho))$ by mapping a subset $S \subset W_k^n$ to $S \cap T_k^n(\rho)$. 
If we consider an ideal $I \subset W_k^n$ we have for each 
$\rho$ a subset $I \cap T_k^n(\rho) \subset T_k^n(\rho)$ which is an ideal of $T_k^n(\rho)$.
Since any element of $I$ must belong to $T_k^n(\rho)$, for some $\rho > 1$, we see that 
\[ I = \bigcup_{\rho > 1} I \cap T_k^n(\rho) = \limind_{\rho > 1} I \cap T_k^n(\rho). \]
Given an ideal $I \subset W_k^n$, by theorem \ref{thm_W} we know that $I$ is closed and so $I \cap T_k^n(\rho)$ is a closed ideal of $T_k^n(\rho)$ (because $T_k^n(\rho) \to W_k^n$ is bounded).

On the other hand, take a filtrant system of closed ideals $I_\rho \subset T_k^n(\rho)$. Their limit
\[ I = \limind_{\rho > 1} I_\rho = \bigcup_{\rho > 1} I \cap T_k^n(\rho) \]
is an ideal of $W_k^n$. In fact, $I$ is clearly a subspace of $W_k^n$ and for each $a \in W_k^n$ and $x \in I$ we have that there exists a $\rho$
such that both $a, x \in T_k^n(\rho)$ and so $a x \in I_\rho \then a x \in I$. Moreover, $I$ is bornologically closed because $W_k^n$ carries
the direct limit bornology. So a sequence $\{ x_n \}_{n \in \N} \subset W_k^n$ converges in $W_k^n$ if and only if there exists a bounded subset $B$ such that for all $\la \in k^\times$, 
$\{ x_n \}_{n \in \N} \subset \la B$ for all $n > N$ for some $N = N(\la) \in \N$. But $B \subset T_k^n(\rho)$ for some $\rho > 1$, hence $\{x_n\}_{n \in \N}$ converges
to an element of $T_k^n(\rho)$ with respect to the norm of $T_k^n(\rho)$.

This shows that we have two maps $\sigma: \cI_{W_k^n} \to \underset{\rho > 1}\limpro \cI_{T_k^n(\rho)}$ and 
$\tau: \underset{\rho > 1}\limpro \cI_{T_k^n(\rho)} \to \cI_{W_k^n}$,
since elements of $\cI_{T_k^n(\rho)}$ are precisely systems of closed ideals. These maps are clearly inverse of each other and preserve the orderings
hence they are isomorphisms of partially ordered sets. 
\end{proof}

\begin{prop} \label{prop:quotient_W}
Let $I \subset W_k^n$ be an ideal then $W_k^n/I$ is a complete bornological algebra.
Moreover, if $I = (f_1, ..., f_r)$ then we can write
\[ \frac{W_n}{I} \cong \limind_{\rho > 1} \l ( \frac{T_k^{n}(\rho)}{\ol{(f_1, ..., f_r)_\rho}} \r ), \]
for $\rho$ small enough.
\end{prop}
\begin{proof} 
Since $I = (f_1, ..., f_r)$ is finitely generated there must exist a $\rho' > 1$ such that $f_i \in T_k^n(\rho)$ for each $\rho < \rho'$ and since
\[ \limind_{\rho' > \rho > 1} T_k^{n}(\rho) \cong \limind_{\rho > 1} T_k^{n}(\rho) \]
we can suppose that $f_i \in T_k^n(\rho)$ for all $\rho$. To conclude the proof it is enough to notice that calculating a quotient is a colimit operation, hence it commutes with all colimits by abstract non-sense. This commutation property is precisely a restatement of the claimed isomorphism.
\end{proof}

\section{The category of dagger affinoid algebras}

\index{strict dagger affinoid algebra}
\begin{defn}
A \emph{strict $k$-dagger affinoid algebra} is a complete bornological algebra which is isomorphic to a quotient $W_k^n/I$, with $I$ an ideal.
We denote by $\bS \bAff_k^\dagger$ the full subcategory of $\bBorn(\bAlg_k)$ whose objects are strict $k$-dagger affinoid algebras
and whose morphisms are bounded algebra morphisms.
\end{defn}

\begin{prop} \label{prop:strict_affinoid_closed_ideals}
If $A$ is a strict $k$-dagger affinoid algebra then $A$ is a complete bornological algebra for which all ideals are bornologically closed.
\end{prop}
\begin{proof} 
$A$ inherits the quotient bornology from $W_k^n$. With this bornology $A$ is a complete bornological 
m-algebra by proposition \ref{prop:quotient_W}. Suppose that $I \subset A$ is a non-closed ideal of $A$ then $A/I$ 
is a non-separated bornological algebra and the composition $W_k^n \to A \to A/I$ is a surjective morphism of $W_k^n$ to a 
non-separated bornological algebra, hence the preimage of $I$ is a non-closed ideal of $W_k^n$. This contradicts theorem \ref{thm_W}, therefore all ideals of $A$ are bornologically closed.
\end{proof}

\begin{prop}
Every strict $k$-dagger affinoid algebra is canonically a LF-regular algebra.
\end{prop}
\begin{proof} 
We showed that $W_k^n$ is defined by a compact (or compactoid if $k$ is not locally compact) direct limit and that $(W_k^n)^t$ is an LF-regular algebra. If $ A \cong \frac{W_k^n}{I}$ equipped with the quotient topology then we have that
\begin{equation} \label{eqn:system_A}
A \cong \limind_{\rho > 1} \frac{T_k^n(\rho)}{T_k^n(\rho) \cap I},
\end{equation}
by the open mapping theorem for LF-spaces.
Let $X \subset T_k^n(\rho)$ be an open subset such that $\varphi_{\rho, \rho'}(X) \subset Y \subset T_k^n(\rho')$ for $Y$ compact in $T_k^n(\rho')$
then $\pi_\rho(X) \subset \frac{T_k^n(\rho)}{I \cap T_k^n(\rho)}$ is an open subset in $\frac{T_k^n(\rho)}{I \cap T_k^n(\rho)}$ because $\pi_\rho$
is a strict epimorphism and $\pi_{\rho'}(Y) \subset \frac{T_k^n(\rho')}{I \cap T_k^n(\rho')}$ is compact(oid) (where $\varphi_{\rho, \rho'}$ denotes the system morphisms of (\ref{eqn:system_A}) and $\pi_\rho: T_k^n(\rho) \to \frac{T_k^n(\rho')}{I \cap T_k^n(\rho')}$ the quotient map). Moreover, we have that
\[ \varphi_{\rho, \rho'}(\pi_\rho(X)) \subset \pi_{\rho'}(Y) \]
hence the system of definition of $A$ is compact(oid) which implies that $A$ is LF-regular.
\end{proof}

We have the following corollary, which sum up the analogous corollary to the theorem \ref{thm_W_regular}.

\begin{cor} \label{cor_dagger_LF_regular}
Let $A$ be a strict $k$-dagger affinoid algebra then $A^t$ is Hausdorff, nuclear, reflexive, complete, locally complete LF-algebra.
\end{cor}
\begin{proof}
Same proof of corollary \ref{cor:W_regular}.
\end{proof}

\begin{prop} \label{prop:regular_dagger_affinoid_algebra}
Every strict $k$-dagger affinoid algebra is regular as bornological algebra.
\end{prop}
\begin{proof} 
Given a dagger affinoid algebra $A$ the proposition is deduced from the presentation
\[ A \cong \limind_{\rho > 1} \frac{T_k^{n}(\rho)}{\ol{(f_1, ..., f_r)_\rho}}  \]
for some $n \in \N$ and $f_1, \dots, f_r \in W_k^n$ given by proposition \ref{prop:quotient_W}, which permits to apply the same reasoning of theorem \ref{thm:dagger_spectrum} to deduce that $A$ is regular.
\end{proof}

For the next result we need to prove a lemma which is a generalization of the proposition 3.7.5 of \cite{BGR}. We state our lemma at the maximal
generality that we can prove it, and in doing so we recall a couple of definitions from \cite{BA}, which generalize to any base field the work done in \cite{GACH}.

\index{net on a bornological space}
\begin{defn} \label{defn:net}
Let $F$ be a $k$-vector space. A \emph{net} on $F$ is a map $\cN: \underset{j \in \N}\bigcup \N^j \to \sP(F)$ such that
\begin{enumerate}
	\item each $\cN(n_1, \ldots, n_j)$ is a disk;
	\item $\cN(\void) = F$;
	\item for every finite sequence $(n_0, \ldots, n_j)$ one has
	\[  \cN(n_0, \ldots, n_j) = \bigcup_{n \in \N} \cN(n_0, \ldots, n_j, n). \]
\end{enumerate}
\end{defn}

Notice that condition $(2)$ of previous definition is used to give sense to the formula
\[ F = \cN(\void) = \bigcup_{n \in \N} \cN(n). \]
If $s: \N \to \N$ is a sequence we will use the notation
\[ \cN_{s, j} = \cN(s(0), \ldots, s(j)). \]
 
\begin{defn} \label{defn:net_compatible}
Let $F$ be a separated bornological $k$-vector space of convex type. Then, we say that a net $\cN$ on $F$ is compatible with its bornology if
\begin{enumerate}
	\item for every sequence $s: \N \to \N$ there is a sequence of positive real numbers $b(s): \N \to \R_{> 0}$ such that for all $x_j \in \cN_{s,j}$ and $a_j \in K$ with $|a_j| \le b(s)_j$ the series
	\[ \sum_{j \in \N} a_j x_j \]
	converges bornologically in $F$ and $\underset{k \ge n}\sum a_j x_j \in \cN_{s,n}$ for every $n \in \N$.
	\item For every sequence $\{ \la_j\}_{j \in \N}$ of elements of $k$ and $s: \N \to \N$ the set
	\[ \bigcap_{j \in \N} \la_k \cN_{s, j} \]
	is bounded in $F$.
\end{enumerate}
We say that a separated bornological vector space \emph{has a net} if there exists a net on it which is compatible with its bornology. 
\end{defn}

The concept of bornological vector space with net is quite general and it is easy to show that the underlying bornological vector space of a dagger affinoid algebra has a net (e.g. example 2.3 (3) of \cite{BA}). The following lemma is a less general version of lemma 4.23 of \cite{BA}, which is enough for our scope. We reproduce here the full proof for the sake of clarity.

\begin{lemma} \label{lemma:closed_graph}
Let $A, B$ be bornological algebras for which the underlying bornological vector space of $A$ is complete and the one of $B$ has a net. Let $\phi: A \to B$ be an algebra morphism. Suppose that in 
$B$ there is a family of ideals $\cI$ such that
\ben
\item each $I \in \cI$ is bornologically closed in $B$ and each $\phi^{-1}(I)$ is bornologically closed in $A$;
\item for each $I \in \cI$ one has $\dim_k B/I < \infty$;
\item $\underset{I \in \cI}\bigcap I = (0)$.
\een
Then, $\phi$ is bounded.
\end{lemma}
\begin{proof}
Let $I \in \cI$ and let $\be: B \to B/I$ denotes the residue epimorphism and $\psi = \be \circ \phi$. Let $\ol{\psi}: A/\Ker(\psi) \to B/I$ denotes the canonical injection, which gives the following commutative diagram
\[
\begin{tikzpicture}
\matrix(m)[matrix of math nodes,
row sep=2.6em, column sep=2.8em,
text height=1.5ex, text depth=0.25ex]
{ A            & B   \\
  A/\Ker(\psi) & B/I \\};
\path[->,font=\scriptsize]
(m-1-1) edge node[auto] {$\phi$} (m-1-2);
\path[->,font=\scriptsize]
(m-1-1) edge node[auto] {} (m-2-1);
\path[->,font=\scriptsize]
(m-1-2) edge node[auto] {$\be$} (m-2-2);
\path[->,font=\scriptsize]
(m-2-1) edge node[auto] {$\ol{\psi}$} (m-2-2);
\path[->,font=\scriptsize]
(m-1-1) edge node[auto] {$\psi$} (m-2-2);
\end{tikzpicture}.
\]
We have that $\Ker(\psi) = \phi^{-1}(I)$, therefore, since by hypothesis $B/I$ is finite dimensional, also $A/\Ker(\psi)$ is. Thus, both $B/I$ and $A/\Ker(\psi)$ are finite dimensional separated bornological algebras, when they are equipped with the quotient bornology. Therefore, their underlying bornological vector spaces are isomorphic to the direct product of a finite number of copies of $k$. So, $\ol{\psi}$ is bounded and this implies the boundedness of $\psi$.

Let $\{a_n\}_{n \in \N} \subset A$ be a sequence such that $\lim a_n = 0$, bornologically. Then
\[ \be(\lim_{n \to \infty} \phi(a_n)) = \lim_{n \to \infty} (\be \circ \phi)(a_n) \]
because $\be$ is bounded, hence
\[ \lim_{n \to \infty} (\be \circ \phi)(a_n) = \lim_{n \to \infty} \psi(a_n) = \psi(\lim a_n) = 0 \]
which implies that $\phi(\underset{n \to \infty}\lim a_n) \in I$. 
Since this must be true for any $I \in \cI$ and $\underset{I \in \cI}\bigcap I = (0)$ we deduce that $\phi(\lim a_n) = 0$. This implies that the graph of $\phi$ is bornologically closed, because then for any sequence $\{ a_n \}_{n \in \N}$ in $A$ such that $\underset{n \to \infty}\lim a_n = a$ one has that
\[ \underset{n \to \infty}\lim ( a_n, \phi(a_n)) = (a, \phi(a)) \in \Gamma(\phi). \]
Now we can apply theorem 2.7 of \cite{BA} to infer that $\phi$ is bounded.
\end{proof}

\begin{thm}  \label{boundedness_thm}
Let $A, B \in \ob(\bS \bAff_k^\dagger)$ and $\phi: A \to B$ an algebra morphism, then $\phi$ is bounded.
\end{thm}
\begin{proof} 
It is enough to show that every morphism between strict $k$-dagger affinoid algebras fulfils the conditions of lemma \ref{lemma:closed_graph}.
Let $\cI$ be the family of ideals of $B$ defined by 
\[ \cI = \{ \fm^e | e \ge 1, \fm \text{ is maximal in } B \}. \]
We know all ideals of $A$ and $B$ are closed and that $B / \fm^e$ is a finite dimensional $k$-vector space. We have to check
that 
\[ \bigcap_{I \in \cI} I = 0 \]
but this is a classical consequence of Krull intersection lemma, because $B$ is Noetherian.
\end{proof}

Dagger affinoid algebra morphisms have also nice relations with respect to the natural filtrations with which dagger affinoid algebras are equipped. We need a lemma to explain this.

\begin{lemma}
Let $\rho = (\rho_1, ..., \rho_n)$ with $\rho_i \in \sqrt{k^\times}$ then there exist an $m \in \N$ and a surjective bounded algebra morphism
\[ T_k^m = T_k^m(1) \to T_k^n(\rho). \]
\end{lemma}
\begin{proof} 
The proof for the non-archimedean case is in \cite{BGR} 6.1.5/4. The same proof in the archimedean case shows that there is an isomorphism
\[ T_k^n \to T_k^n(\rho) \]
for any $\rho$.
\end{proof}

In the following proposition we follow the convention of choosing a base of neighborhoods of the closed polydisk $(1, \ldots, 1)$ made of polydisks of the form
$(\rho, \ldots, \rho)$ for $\rho > 1$.

\begin{prop} \label{prop_filtration_dagger}
Let 
\[ \phi: A = \frac{W_k^n}{(f_1, ..., f_r)} \to B = \frac{W_k^m}{(g_1, ..., g_s)} \]
be a morphism of strict $k$-dagger affinoid algebras. Then, there exists $\rho_A$ and $\rho_B$ such that
\[ A \cong \limind_{1 < \rho < \rho_A} \frac{T_k^n(\rho)}{(f_1, ..., f_r)_\rho} = \limind_{1 < \rho < \rho_A} A_\rho \]
\[ B \cong \limind_{1 < \rho < \rho_B} \frac{T_k^m(\rho)}{(g_1, ..., g_s)_\rho} = \limind_{1 < \rho < \rho_B} B_\rho, \]
for each $\rho \in \sqrt{|k^\times|}$ with $\rho < \rho_A$ and $\rho < \rho_B$ there exists a $\rho'$ such that the canonical morphism
\[ A_\rho \to A \to B \]
factorize through
\[ B_{\rho'} \to B. \]
\end{prop}
\begin{proof} 
For the non-archimedean case the result is proved by Grosse-Kl\"onne in \cite{GK}, lemma 1.8. The archimedean case can be solved by the same reasoning. Indeed, we have the morphisms
\[ T_k^n \stackrel{\sim}{\to} T_k^n(\rho) \to A_\rho \to A \stackrel{\phi}{\to} B. \]
We denote the composition $T_k^n = k \lt X_1, ..., X_n \gt \to B$ with the symbol $\a$. 
Since it is a composition of bounded morphisms we have that $\a(X_i) \in B$ is power-bounded for each $i$ and so there exists a $\rho'$ with 
$\a(X_i) \in B_{\rho'}$, for all $ 1 \le i \le n$ (because they are a finite number of elements). We have that
\[ \phi(A_\rho) = \a(T_k^n) \subset B_{\rho'} \to B \]
which proves the proposition.
\end{proof}

\begin{cor}
Let $\phi: A \to B$ be an algebra homomorphism between dagger affinoid algebras then $\phi$ ``is compatible with the filtrations''.
So, if we consider dagger algebras as ind-objects of the category of seminormed algebras, to give an algebra
morphism is equivalent to give a morphism of $A$ and $B$ as objects of $\bInd(\bSn(\bAlg_k))$.
\end{cor}
\begin{proof}
This is simply a restatement of last proposition.
\end{proof}

We can restate the previous corollary in the following way.

\begin{cor}
The following diagram
\[
\begin{tikzpicture}
\matrix(m)[matrix of math nodes,
row sep=2.6em, column sep=2.8em,
text height=1.5ex, text depth=0.25ex]
{ \bS \bAff_k^\dagger & \bInd(\bSn(\bAlg_k)) \\
                      & \bBorn(\bAlg_k) \\};
\path[->,font=\scriptsize]
(m-1-1) edge node[auto] {$\kappa$} (m-2-2);
\path[->,font=\scriptsize]
(m-1-1) edge node[above] {$\iota$} (m-1-2);
\path[->,font=\scriptsize]
(m-1-2) edge node[auto] {$$} (m-2-2);
\end{tikzpicture}
\]
is commutative and the immersions $\iota$ and $\kappa$ are fully faithful.
\end{cor}
\begin{proof} 
\end{proof}

The definition of $W_k^n$ is a particular case of the definition of the ring of bornological overconvegent power series, applied when $A = k$ 
(see definition \ref{defn:over_conv}). Moreover, the notation that we are using in the theory of dagger affinoid algebra is consistent with the notation introduced in the previous chapter. Therefore, we can state the universal property characterizing $W_k^n$:

\begin{prop}
Given any complete bornological algebra $A$ and $a_1, ..., a_n \in A^\ov$ weakly power-bounded elements, there exists a unique morphism
\[ \Phi: W_k^n = k \lt X_1, ..., X_n \gt^\dagger \to A \]
such that $\Phi(X_i) = a_i$.
\end{prop}

\index{affinoid generators}
In particular, if $\Phi: k \lt X_1, ..., X_n \gt^\dagger \to A$ is a strict epimorphism (and hence $A$ is a strict $k$-dagger affinoid algebra) 
we say that $f_1 = \Phi(X_1), ..., f_n = \Phi(X_n)$ form a \emph{system of affinoid generators} of $A$ and we write 
\[ A = k \lt f_1, ..., f_n \gt^\dagger. \]

The category of strict $k$-dagger affinoid algebras has tensor products, which  represent coproducts and push-outs. 
Given two morphisms of $k$-dagger affinoid algebras $A \to B$ and $A \to C$ we define
\[ B \otimes^\dagger_A C \doteq B \widehat{\otimes}_A C \]
where $\widehat{\otimes}$ denotes the completion of the projective bornological tensor product\footnote{The projective tensor product for bornological vector spaces, and its completed version, are discussed in chapter 4 of \cite{H1}. Also in section 3.3 of \cite{BABE} one can find an account of the basic properties of the completed projective tensor product with a particular emphasis on the closed symmetric monoidal structure that $\bCBorn_k$ gets from the functor $\cdot \widehat{\otimes} \cdot$.} of $B$ and $C$, over $k$, modulo the ideal generated by element of the form
\[ (a b) \otimes c - b \otimes (a c) \]
for $a \in A, b \in B, c \in C$. Notice that this ideal is always closed because $B \otimes^\dagger_k C$ is a strict $k$-dagger affinoid algebra.

\begin{prop} \label{prop:affinoid_A}
Let $A \in \ob(\bS \bAff_k^\dagger)$ then there is an isomorphism
\[ A \lt X_1, ..., X_n \gt^\dagger \to A \otimes_k^\dagger W_k^n. \]
\end{prop}
\begin{proof} 
It is enough to show that both algebras satisfy the same universal property in the category of complete bornological algebras. 
The universal property of $A \lt X_1, ..., X_n \gt^\dagger$ is given by theorem \ref{thm:over_univ_prop}. Let $\phi: A \to B$ be a morphism of $A$ to a complete bornological algebra and let $b_1, ..., b_n \in B^\ov$.  By the universal property of $W_k^n$ there exists a unique  morphism $\psi: W_k^n \to B$ such that $\psi(X_i) = b_i$.
And finally, by the universal property of the tensor product there exists a unique morphism $\gamma$ which makes the following diagram commutative
\[
\begin{tikzpicture}
\matrix(m)[matrix of math nodes,
row sep=2.6em, column sep=2.8em,
text height=1.5ex, text depth=0.25ex]
{A & & W_k^n \\
 & A \otimes_k^\dagger W_k^n \\
 & B \\};
\path[->,font=\scriptsize]
(m-1-1) edge node[auto] {$\a_1$} (m-2-2);
\path[->,font=\scriptsize]
(m-1-3) edge node[above] {$\a_2$} (m-2-2);
\path[->,font=\scriptsize]
(m-2-2) edge node[auto] {$\gamma$} (m-3-2);
\path[->,font=\scriptsize]
(m-1-1) edge node[below] {$\phi$} (m-3-2);
\path[->,font=\scriptsize]
(m-1-3) edge node[auto] {$\psi$} (m-3-2);
\end{tikzpicture}
\]
where $\a_i$ are the canonical maps of the tensor product. Therefore $A \lt X_1, ..., X_n \gt^\dagger$ and $A \otimes_k^\dagger W_k^n$ satisfy the same
universal property whence they are isomorphic.
\end{proof}

\begin{cor}
\[ W_k^n \otimes_k^\dagger W_k^m \cong W_k^{n + m} \]
and if $k'/k$ is an extension of complete valued fields then
\[ W_{k'}^n \cong W_k^n \otimes_k^\dagger k'. \]
\end{cor}
\begin{proof} 
It is easy to show the isomorphisms by checking universal properties as we did in last proposition.
\end{proof}

\begin{prop} \label{prop:affinoid_tensor_product}
Let $A, B, C$ be strict $k$-dagger affinoid algebras and $A \to B$, $A \to C$ two morphisms, then
\[ B \otimes^\dagger_A C \]
is a strict $k$-dagger affinoid algebra.
\end{prop}
\begin{proof} 
By hypothesis there are strict epimorphisms $W_k^n \to A$, $W_k^m \to B$, $W_k^p \to C$ and $B \otimes_k^\dagger C \to B \otimes_A^\dagger C$.
We then get strict epimorphisms
\[ W_k^m \otimes_k^\dagger W_k^p \to B \otimes_k^\dagger C \to B \otimes_A^\dagger C. \]
and by the isomorphism $W_k^m \otimes_k^\dagger W_k^p \to W_k^{m + p}$ showed in last corollary, we get that the strict epimorphisms
\[ W_k^{m + p} \to B \otimes_k^\dagger C \to B \otimes_A^\dagger C \]
shows that $B \otimes_A^\dagger C$ is strictly affinoid.
\end{proof}

\begin{prop} \label{prop:affinoid_finite_product}
The category $\bS \bAff_k^\dagger$ has finite direct products.
\end{prop}
\begin{proof} 
We check that the ring-theoretic direct product has the desired property. It is enough to check that $W_k^n \oplus W_k^m$ 
is an object of $\bS \bAff_k^\dagger$ and we notice that we can suppose $n = m$ (since we have to find a surjective homomorphism). The map 
$W_k^n \lt X \gt^\dagger \to W_k^n \oplus W_k^n$ defined
\[ \sum_{i = 0}^\infty a_i X^i \mapsto (a_0, \sum_{i = 0}^\infty a_i) \]
is a bounded and strict epimorphism, which proves the proposition.
\end{proof}

We gave the definition
\[ A \lt X \gt^\dagger = \limind_{\rho > 1} A \lt \rho^{-1} X_i \gt \]
for a bornological m-algebra. It is now clear that we can work out other similar constructions. For example, we can define
\[A \lt X^{-1} \gt^\dagger \doteq A \lt X, X^{-1} \gt^\dagger \doteq \limind_{\rho > 1} A \lt \rho^{-1} X_i, \rho X_i^{-1} \gt \]
The algebra $A \lt X^{-1} \gt^\dagger$ satisfies the following universal property: given any complete bornological algebra $B$,
a unit $f \in B$ such that $f^{-1} \in B^\ov$ and a morphism $\phi: A \to B$, then there exists a unique morphism
\[ \Phi: A \lt X^{-1} \gt^\dagger \to B \]
such that $\Phi|_A = \phi$ and $\Phi(X) = f$. We can also define
\[ A \lt X \gt^\dagger \lt Y^{-1} \gt^\dagger \] 
where now $X = (X_1, ..., X_n)$ and $Y = (Y_1, ..., Y_m)$ are vectors of variables, and verify that this algebra satisfies the following universal property: 
given any complete bornological algebra $B$, and $f_1,...,f_n \in B$ and units $g_1,...,g_m \in B$ such that $f_i, g_j^{-1} \in B^\ov$ and a morphism $\phi: A \to B$,
 then there exists a unique morphism
\[ \Phi: A \lt X \gt^\dagger \lt Y^{-1} \gt^\dagger \to B \]
such that $\Phi|_A = \phi$ and $\Phi(X_i) = f_i$ and $\Phi(Y_i) = g_i$. We will use the notation
\[ A \lt X \gt^\dagger \lt Y^{-1} \gt^\dagger = A \lt X, Y^{-1} \gt^\dagger. \]

If $A$ is an object of $\bS \bAff_k^\dagger$, we showed that $A \lt X_1, ..., X_n \gt^\dagger \in \ob(\bS \bAff_k^\dagger)$. Hence, we see that there is 
a canonical surjection
\[ A \lt X_1, X_2 \gt^\dagger \to A \lt X, X^{-1} \gt^\dagger \]
which shows that also $A \lt X, X^{-1} \gt^\dagger$ is a strict dagger affinoid algebra and more generally
\[ A \lt X, Y^{-1} \gt^\dagger \] 
for $X = (X_1, ..., X_n)$ and $Y = (Y_1, ..., Y_m)$ is a strict dagger affinoid algebra.

Now, let $A \in \ob(\bS \bAff_k^\dagger)$ and let $f = (f_1, ..., f_n), g = (g_1 ,..., g_m)$ be elements of $A$. 
Following classical rigid geometry, we define the following algebra
\[ A \lt f, g^{-1} \gt^\dagger \doteq \frac{A \lt X, Y \gt^\dagger}{(X - f, g Y - 1)}. \]
It is clear that $A \lt f, g^{-1} \gt^\dagger$ is a strict dagger affinoid algebra, which satisfies the following universal property.

\begin{prop} \label{prop:univ_property_laurent}
Let $\varphi: A \to B$ be a homomorphism between strict $k$-dagger affinoid algebras, such that
$\varphi(g_j)$ are units and $\varphi(f_i), \varphi(g_j)^{-1} \in B^\ov$, 
then there exists a unique morphism $\psi: A \lt f, g^{-1} \gt^\dagger \to B$ for which the diagram
\[
\begin{tikzpicture}
\matrix(m)[matrix of math nodes,
row sep=2.6em, column sep=2.8em,
text height=1.5ex, text depth=0.25ex]
{A & & A \lt f, g^{-1} \gt^\dagger \\
 & B & \\};
\path[->,font=\scriptsize]
(m-1-1) edge node[auto] {$\varphi$} (m-2-2);
\path[->,font=\scriptsize]
(m-1-1) edge node[auto] {$$} (m-1-3);
\path[->,font=\scriptsize]
(m-1-3) edge node[auto] {$\psi$} (m-2-2);
\end{tikzpicture}
\]
commutes. 
\end{prop}
\begin{proof} 
Immediate consequence of the universal property of $A \lt X \gt^\dagger$ and that of the quotient.
\end{proof}

We can mimic another classical construction. Given $g, f_1, ..., f_n \in A$ with no common zeros, we define
\[ A \l \lt \frac{f}{g} \r \gt^\dagger \doteq \frac{A \lt X \gt^\dagger}{(g X - f)}. \]
It is clear that $A \lt \frac{f}{g} \gt^\dagger$ is a strict dagger affinoid algebra, which satisfies the following universal property.

\begin{prop} \label{prop:univ_property_rational}
Let $\varphi: A \to B$ be a homomorphism between strict $k$-dagger affinoid algebras, such that
the elements $\frac{\varphi(f_i)}{ \varphi(g_j)} \in B^\ov$, 
then there exists a unique morphism $\psi: A \l \lt \frac{f}{g} \r \gt^\dagger \to B$ for which the diagram
\[
\begin{tikzpicture}
\matrix(m)[matrix of math nodes,
row sep=2.6em, column sep=2.8em,
text height=1.5ex, text depth=0.25ex]
{A & & A \l \lt \frac{f}{g} \r \gt^\dagger \\
 & B & \\};
\path[->,font=\scriptsize]
(m-1-1) edge node[auto] {$\varphi$} (m-2-2);
\path[->,font=\scriptsize]
(m-1-1) edge node[auto] {} (m-1-3);
\path[->,font=\scriptsize]
(m-1-3) edge node[auto] {$\psi$} (m-2-2);
\end{tikzpicture}
\]
commutes. 
\end{prop}
\begin{proof}
Again easy verifications.
\end{proof}

In order to understand Berkovich geometry from a dagger perspective, we introduce the notion of non-strict dagger affinoid algebra.

\index{dagger affinoid algebra}
\begin{defn}
A (non-strict) \emph{$k$-dagger affinoid algebra} is a bornological algebra which is isomorphic to a quotient 
\[ \frac{k \lt \rho_1^{-1} X_1, ..., \rho_n^{-1} X_n \gt^\dagger}{I} = \frac{W_k^n(\rho)}{I} \]
where $\rho = (\rho_i) \in \R_+^n$ is any given polyradius. We consider on $W_k^n(\rho) / I$ the quotient bornology and call it the dagger bornology.

We denote the category of all $k$-dagger affinoid algebras with bounded algebra morphisms with the symbol $\bAff^\dagger_k$.
\end{defn}

\begin{rmk} \label{rmk:non_strict_dagger_affinoids}
It is clear that the results of proposition \ref{prop_filtration_dagger}, \ref{prop:affinoid_A}, \ref{prop:affinoid_tensor_product}, \ref{prop:affinoid_finite_product}, \ref{prop:univ_property_laurent} and \ref{prop:univ_property_rational}, with appropriate modifications, can be stated also for non-strict $k$-dagger affinoid algebras. We omit the details which are pretty straightforward. Notice also that every dagger affinoid algebra $A$ can be written as a direct limit 
\[ A \cong \limind_{\rho > r_A} A_\rho \]
for some polyradius $r_A$, where $A_\rho$ are strict affinoid algebras (of course is enough to let $\rho$ varies only on polyradii whose coordinates values are in $\sqrt{|k^\times|}$).
\end{rmk}

\begin{lemma}
Let $K/k$ be an extension of complete valued fields and $A \in \ob(\bS  \bAff_k^\dagger)$, then $A \otimes_k^\dagger K \in \ob(\bS \bAff_K^\dagger)$.
\end{lemma}
\begin{proof} 
Let $W_k^n \to A$ be a strict epimorphism, then by the functoriality of the dagger tensor product we have a morphism
$W_k^n \otimes_k^\dagger K \to A \otimes_k^\dagger K$, and we showed that $W_k^n \otimes_k^\dagger K \cong W_K^n$. 
The lemma is proved noticing that $\otimes_k^\dagger$ is right exact, \ie it preserves surjective morphisms (cf. \cite{H1}, chapter 4, \S 1, $n^\circ$ 10).
\end{proof}

\begin{prop}
Let $A$ be a $k$-dagger affinoid algebra, then there exists a finitely generated extension $K/k$ of complete valued fields such that
\[ A \otimes_k^\dagger K \]
is a strict $K$-dagger affinoid algebra.
\end{prop}
\begin{proof} 
We note that only in the non-archimedean case one can meet non-strict algebras, so we suppose that $k$ is non-archimedean, without lost of generality. By definition
\[ A \cong \frac{k \lt \rho_1^{-1} X_1, ..., \rho_n^{-1} X_n \gt^\dagger}{I}. \]
and let $K/k$ be such that $\rho_i \in \sqrt{|K^\times|}$, for all $i$. We can always find such a field, as explained in section 2.1 of \cite{BER2}. It is enough to show that
\[ k \lt \rho_1^{-1} X_1, ..., \rho_n^{-1} X_n \gt^\dagger \otimes_k^\dagger K \cong W_K^n(\rho) \cong W_k^n \]
because then we get a surjection (because tensor product is right exact)
\[ k \lt \rho_1^{-1} X_1, ..., \rho_n^{-1} X_n \gt^\dagger \otimes_k^\dagger K \to \frac{k \lt \rho_1^{-1} X_1, ..., \rho_n^{-1} X_n \gt^\dagger}{I} \otimes_k^\dagger K. \]
So
\[ k \lt \rho_1^{-1} X_1, ..., \rho_n^{-1} X_n \gt^\dagger \otimes_k^\dagger K 
 = (\limind_{\la > \rho_1} T_k(\la)) \otimes_k^\dagger ... \otimes_k^\dagger (\limind_{\la > \rho_n} T_k(\la)) \otimes_k^\dagger K \cong \]
\[ \cong \limind_{\la_1 > \rho_1, ..., \la_n > \rho_n} \l (  T_k(\la_1)  \widehat{\otimes}_k ... \widehat{\otimes}_k T_k(\la_n) \widehat{\otimes}_k K  \r ) \cong \]
\[ \limind_{\la_1 > \rho_1, ..., \la_n > \rho_n} \l (  T_K(\la_1)  \widehat{\otimes}_K ... \widehat{\otimes}_K T_K(\la_n) \r ) 
 \cong (\limind_{\la > \rho_1} T_K(\la)) \otimes_K^\dagger ... \otimes_K^\dagger (\limind_{\la > \rho_n} T_K(\la)) \cong W_K^n(\rho) \]
where with $\widehat{\otimes}$ we denote the complete projective tensor product and where we use the fact that the complete projective tensor product commutes with arbitrary inductive limits because it is left adjoint to the internal hom functor of $\bBorn_k$ (again cf. \cite{H1}, chapter 4, \S 1, $n^\circ$ 10).
\end{proof}

Hence, we can often reduce problems about dagger affinoid algebras to problems about strict dagger affinoid algebras by tensoring with a suitable finitely generated extension of $k$, as Berkovich did in his foundational work on non-strict affinoid algebras, cf. chapter 2 of \cite{BER2}.

\begin{prop} \label{prop:non_strict_closed_ideals}
Let $A$ be a $k$-dagger affinoid algebra. All the ideals of $A$ are closed for the dagger bornology.
\end{prop}
\begin{proof}
Notice again that only in the non-archimedean case we can find non-strict dagger affinoid algebras. In this case the result follows easily from \cite{BER2} proposition 2.1.3.
\end{proof}

\begin{prop} \label{prop:non_strict_bounded_morphism}
All algebra morphisms between dagger affinoid algebras are bounded.
\end{prop}
\begin{proof}
Noticing that any non-strict dagger affinoid algebra can be written as a direct limit of strict ones, by the same reasoning of \ref{prop_filtration_dagger} (also cf. remark \ref{rmk:non_strict_dagger_affinoids}). Also, every algebra morphism can be written as a morphism of direct systems of algebras, in the same way we did in proposition \ref{prop_filtration_dagger}. Therefore, every morphism between non-strict dagger affinoid algebras can be written as a direct limit of bounded ones, hence it is bounded.
\end{proof} 

We conclude this section with a useful proposition.

\begin{prop} \label{prop_finite_dagger_algebra}
Let $B \in \ob(\bAff_k^\dagger)$ and $\varphi: B \to A$ a bounded finite morphism where $A$ is a bornological algebra whose underlying bornological vector space is complete. Then, $A \in \ob(\bAff_k^\dagger)$.
\end{prop}
\begin{proof} 
We may assume $B = W_k^n(\rho)$, since $B$ is a quotient of some $W_k^n(\rho)$ and we have to find a strict epimorphism of underlying bornological vector spaces from some $W_k^n$ to $A$. The assumption that $\varphi: B \to A$ is finite means that $A$ is a finite $B$-module with respect to the module structure given by $\varphi$, hence 
\[ A = \sum_{i = 1}^m \varphi(B) a_i = \sum_{i = 1}^m \varphi(W_k^n(\rho)) a_i \]
for some $a_i \in A$, and we may assume $a_i \in A^\ov$.
We can consider $W_k^n(\rho) \lt X_1, ..., X_m \gt^\dagger \cong W_k^{m + n}(\rho')$ and by its universal property, we can find 
\[ \Phi: W_k^{m + n}(\rho') \to A \]
such that $\Phi(X_i) = a_i$ for each $i$. Therefore, $\Phi$ is surjective and bounded and for theorem 4.9 of \cite{BA} it also a strict morphism, showing that $A$ is an object of $\bAff_k^\dagger$.
\end{proof}

\begin{rmk}
Proposition \ref{prop_finite_dagger_algebra} is a generalization to dagger context  of proposition 6.1.1/6 of \cite{BGR} valid for classical affinoid algebras. We will need this proposition in the second section of chapter 6.
\end{rmk}

\section{Complex dagger algebras as $\bInd$-Stein algebras} \label{sec_ind_stein}

We now define three different filtrations in $W_k^n$ and study the relations between them. 
These filtrations will be studied only in the archimedean case because they give the link between our work and the theory of Stein algebras, whereas in the non-archimedean case there is no such a well-established theory to compare with, or at least not so well-known. For more information about the non-archimedean case of the results of this section one can refer to section 6 of \cite{BABE}.
Thus, here we bound our discussion to the case when $k$ is archimedean.

We denote the closed polycylinder of polyradius $\rho = (\rho_1, ... ,\rho_n)$ with $\D_k(\rho^+)$. 
It has a neighborhood basis formed by the open polycylinders of polyradii $\rho' = (\rho_1', ..., \rho_n')$ with $\rho'_i > \rho_i$ for all $i$. 
On any open polycylinder, denoted $\D_k(\rho^-)$, we can consider three natural function algebras: the algebra of summable power-series
\[ T_k^n(\rho) \]
the algebra of holomorphic functions
\[ \cO(\D_k(\rho^-)) \]
and the disk algebra
\[ D(\D_k(\rho^-)) = \cO(\D_k(\rho^-)) \cap \sC(\D_k(\rho^+)) \]
where $\sC(\D_k(\rho^+))$ is the ring of continuous functions on $\D_k(\rho^+)$.
These algebras are in the following set-theoretic relations
\[ T_k^n(\rho) \subset D(\D_k(\rho^-)) \subset \cO(\D_k(\rho^-)) \]
and every inclusion is strict in general (\ie not every holomorphic function on $\D_k(\rho^-)$ has a continuous extension to the border of the polycylinder and
not every function which is holomorphic inside the polycylinder and continuous on the border can be represented by a power-series in $T_k^n(\rho)$).

As a set, the ring of germs of analytic functions on $\D_k(\rho^+)$ is defined by
\[ \cO(\D_k(\rho^+)) = \limind_{\stackrel{V \supset \D_k(\rho^+)}{\text{open}}} \cO(V) \]
and since the open polycylinders $\D_k(\rho'^-)$ of polyradius $\rho' > \rho$ form a neighborhood basis for $\D_k(\rho^+)$, we can write
\[ \cO(\D_k(\rho^+)) = \limind_{\rho' > \rho} \cO(\D_k(\rho'^-)). \]
This means that given any $f \in \cO(\D_k(\rho^+))$ we can always find a
$\rho' > \rho$ such that $f \in \cO(\D_k(\rho'^-))$ and then a $\rho' > \rho'' > \rho$ such that 
$f \in T_k^n(\rho'') \subset \cO(\D_k(\rho''^-))$, \ie $\cO(\D_k(\rho'^-)) \subset T_k^n(\rho'')$.
So, set-theoretically we have the identification
\[ \cO(\D_k(\rho^+)) \cong \limind_{\rho' > \rho} T_k^n(\rho'). \]
The inclusions 
\[ T_k^n(\rho') \subset D(\D_k(\rho'^-)) \subset \cO(\D_k(\rho'^-)) \]
give also the bijections of sets
\[ \cO(\D_k(\rho^+)) = \limind_{\rho' > \rho} \cO(\D_k(\rho'^-)) \cong \limind_{\rho' > \rho} T_k^n(\rho') \cong \limind_{\rho' > \rho} D(\D_k({\rho'}^-)). \]
Moreover, the bijection $\cO(\D_k(\rho^+)) \cong \underset{\rho' > \rho}\limind T_k^n(\rho')$
shows that $\cO(\D_k(\rho^+))$ coincides with what we called $W_k^n(\rho)$ in the previous sections.

\begin{prop} \label{prop:complex_S_O}
Consider on $\underset{\rho' > \rho}\limind \cO(\D_k(\rho'^-))$ and on $\underset{\rho' > \rho}\limind D(\D_k(\rho'^-))$ the direct limit bornologies. Then, we have the following isomorphism of bornological algebras
\[ \limind_{\rho' > \rho} \cO(\D_k(\rho'^-)) \cong \limind_{\rho' > \rho} D(\D_k(\rho'^-)). \]
\end{prop}

We give two proofs of this proposition. One using general results for the theory of bornological vector spaces, and a more elementary one in order to explain  the meaning of the result more explicitly.

\begin{proof} 
Since the bijections described so far are algebra morphisms, we just need to check that the underlying bornological vector spaces are isomorphic. Also, both bornological spaces are complete, because they are monomorphic direct limits of complete ones, and they both have a net compatible with their bornology, by theorem 2.8 (1) of \cite{BA}. We can therefore apply theorem 2.7 of \cite{BA} to deduce that the mentioned bijection is bounded in both directions from the fact that it is clearly bounded in one direction.
\end{proof}

Then the more elementary proof.

\begin{proof} 
A subset
\[ B \subset \limind_{\rho' > \rho} \cO(\D_k(\rho'^-)) \]
is bounded if and only if it is bounded in $\cO(\D_k(\rho'^-))$ for some $\rho' > \rho$. As a bornological $k$-algebra 
\[ \cO(\D_k(\rho'^-)) \cong \limpro_{\rho'' < \rho'} (\cO(\D_k(\rho'^-)), \|\cdot\|_{\D_k(\rho''^+)}) \]
where $(\cO(\D_k(\rho'^-)), \|\cdot\|_{\D_k(\rho''^+)})$ denotes the underlying algebra of $\cO(\D_k(\rho'^-))$ endowed with the norm
\[ \| f \|_{\D_k(\rho''^+)} \doteq \sup_{x \in \D_k(\rho''^+)} |f(x)|. \]
It is clear that if $\D_k(\rho'''^+) \subset \D_k(\rho''^+)$ then
\[ \| f \|_{\D_k(\rho'''^+)} \le \| f \|_{\D_k(\rho''^+)}, \]
hence the identity $(\cO(\D_k(\rho'^-)), \|\cdot\|_{\D_k(\rho''^+)}) \to (\cO(\D_k(\rho'^-)), \|\cdot\|_{\D_k(\rho'''^+)})$ is a bounded map, 
and this projective system gives precisely the Fr\'echet structure on $\cO(\D_k(\rho'^-))$. Therefore, a subset $B' \subset \cO(\D_k(\rho'^-))$ is bounded if and only if it is bounded in $\cO(\D_k(\rho''^-))$ for all $\rho'' < \rho'$. Hence, $B \subset \underset{\rho' > \rho}\limind \cO(\D_k(\rho'^-))$ is bounded if and only if there exists $\rho' > \rho$ such that $B$ is bounded in $D(\D_k(\rho'^-))$, 
which is precisely the property which characterizes the bornology of 
\[ \limind_{\rho' > \rho} D(\D_k(\rho'^-)). \]
\end{proof}

\begin{prop} \label{prop:complex_T_O}
We have also the bornological identification
\[ \limind_{\rho' > \rho} D(\D_k(\rho'^-)) \cong \limind_{\rho' > \rho} T_k^n(\rho'). \]
\end{prop}

The arguments of the first proof of proposition \ref{prop:complex_S_O} applies also for proving this proposition. Also in this case we find convenient to write an elementary proof, for the sake of a deeper understanding.

\begin{proof} 
Since $D(\D_k(\rho'^-))$ is the completion of $T_k^n(\rho')$ with respect to its spectral norm, there exists a canonical bounded morphism 
$T_k^n(\rho') \rhook D(\D_k(\rho'^-))$. This implies that the bijection 
\[ \limind_{\rho' > \rho} T_k^n(\rho') \to \limind_{\rho' > \rho} D(\D_k(\rho'^-)) \]
is bounded. It is harder to check that the inverse map is bounded too. This is equivalent to the existence of a canonical bounded map $D(\D_k(\rho'^-)) \to T_k^n(\rho'')$ (the restriction map) for any $\rho'' < \rho'$. One can easily reduce the general case to the case when $\rho' = 1$, so we have to check that for any $\rho < (1, \dots, 1)$ there is a bounded restriction map
\[ D^n = D(\D_k((1, \dots, 1)^-)) \to T_k^n(\rho). \]
The boundedness of this map can be deduced as a consequence of Cauchy formula for several complex variables. Indeed, if $f \in D^n$ then one can write
\[ f = \sum_{i \in \N^n} a_i X^i \]
with 
\[ a_i = \int_{\T^n} f(z) \ol{z}^i d m(z) \]
where $d m$ is the normalized Haar measure of the $n$-torus $\T^n$ (see for example \cite{RUD}, chapter 1 for a precise definition of $d m$). To check that a map between Banach spaces is bounded is enough to check that 
the unit ball of the domain space is mapped on a bounded subset of the codomain. The unit ball of $D^n$ consists of the elements $f \in D^n$ such that
\[ \|f\|_{\sup} = \sup_{x \in \cM(D^n)} |f(x)| \le 1 \]
hence for these elements we have the estimate $|a_i| \le 1$ for all $i$. It follows that the series $\underset{i \in \N^n}\sum a_i X^i$
is dominated by $\underset{i \in \N^n}\sum X^i$ and if $\rho < (1, \dots, 1)$
\[ \| f \|_\rho \le \| \sum_{i \in \N^n} X^i \|_\rho = \sum_{i \in \N^n} \rho^i < \infty. \]
Therefore the injection $D^n \rhook T_k^n(\rho)$ is bounded.
\end{proof}

\begin{rmk} \label{rmk_universal}
When we described the ring of overconvergent power-series we characterized its universal property by the identification
\[ k \lt X \gt^\dagger \cong \limind_{\rho > 1} k \lt \rho^{-1} X \gt \]
and we noticed that the $k$-algebras $k \lt \rho^{-1} X \gt$ are not complete with respect to their spectral norm if $k$ is archimedean, and 
that the universal property which characterizes the completion of $k \lt \rho^{-1} X \gt$ with respect to the spectral norm is (necessarily)
different from that of $k \lt \rho^{-1} X \gt$ itself. This fact emphasizes once more that if $k$ is archimedean there is no theory of affinoid algebras analogous to the one of rigid geometry. 

We saw in the previous section that $W_k^n$ has all the good properties of Tate algebras (interpreted with the right language, \ie one has to replace
the concept of closed ideal with bornologically closed, power-bounded with weakly power-bounded, etc..) and the previous proposition also shows that $W_k^n$ allows to have a uniform
theory \ie in all cases (archimedean and non-archimedean) and with all possible (reasonable) definitions we always get a unique algebra: $W_k^n$.
\end{rmk}

These identifications of algebras imply also the identification of spectra
\[ \cM(W_k^n(\rho)) \cong \limpro_{\rho' > \rho} \cM(\cO(\D_k(\rho'^-))) \cong \frac{\D_{\ol{k}}(\rho^+)}{\Gal(\ol{k}/k)}, \]
as a consequence of proposition \ref{prop_spectra_LB}.

\begin{thm} \label{thm:T_filtration}
Let 
\[ \phi: A = \frac{W_k^n}{(f_1, ..., f_r)} \to B = \frac{W_k^m}{(g_1, ..., g_s)} \]
be a morphism of $k$-dagger affinoid algebras. Then, there exist $\rho_A$ and $\rho_B$ such that
\[ A \cong \limind_{1 < \rho < \rho_A} \frac{D_k^n(\D_k(\rho^-))}{(f_1, ..., f_r)_\rho} = \limind A_\rho \]
\[ B \cong \limind_{1 < \rho < \rho_B} \frac{D_k^m(\D_k(\rho^-))}{(g_1, ..., g_s)_\rho} = \limind B_\rho, \]
so for each $\rho > 1$ with $\rho < \rho_A$ and $\rho < \rho_B$ there exists a $\rho'$ such that the canonical morphism
\[ A_\rho \to A \to B \]
factorizes through
\[ B_{\rho'} \to B. \]
\end{thm}
\begin{proof}
The isomorphism proved in proposition \ref{prop:complex_T_O} shows that all the arguments used in proposition \ref{prop_filtration_dagger} with the presentation
$W_k^n \cong \underset{\rho > 1}\limind T_k^n(\rho)$ applies also to the presentation $W_k^n \cong \underset{\rho > 1}\limind D(\D_k(\rho)^-)$.
So, we can derive the claimed isomorphisms and the required factorization with the same methods.
\end{proof}

\begin{thm} \label{thm:stein_filtration}
Let 
\[ \phi: A = \frac{W_k^n}{(f_1, ..., f_r)} \to B = \frac{W_k^m}{(g_1, ..., g_s)} \]
be a morphism of $k$-dagger affinoid algebras. Then, there exist $\rho_A$ and $\rho_B$ such that
\[ A \cong \limind_{1 < \rho < \rho_A} \frac{\cO_k^n(\D_k(\rho^-))}{(f_1, ..., f_r)_\rho} = \limind A_\rho \]
\[ B \cong \limind_{1 < \rho < \rho_B} \frac{\cO_k^m(\D_k(\rho^-))}{(g_1, ..., g_s)_\rho} = \limind B_\rho, \]
so for each $\rho > 1$ with $\rho < \rho_A$ and $\rho < \rho_B$ there exists a $\rho'$ such that the canonical morphism
\[ A_\rho \to A \to B \]
factorizes through
\[ B_{\rho'} \to B. \]
\end{thm}
\begin{proof}
The proof uses the same kind of argumentations of theorem \ref{thm:T_filtration}. 
\end{proof}

\begin{cor}  \label{cor_pro_stein}
If $k = \C$ then the category of $\C$-dagger affinoid algebras embeds fully faithfully in the category $\bInd(\bStein)$, where $\bStein$ denotes the category of
Stein algebras.
\end{cor}
\begin{proof}
Immediate consequence of last theorem.
\end{proof}

\begin{rmk}
We have already remarked that our closedness theorem for the ideals of dagger affinoid algebras is coherent with respect to the results on classical affinoid
algebras. With last theorem we can underline that our theorem \ref{thm_W} is coherent with respect to the theory of Stein spaces and Stein algebras, as explained for example in chapter 5 of \cite{REM}. Furthermore, our result about the boundedness of any algebra morphism between dagger algebras 
is coherent not only with respect to the theory of affinoid algebras but also with the theory of Stein algebras: see page 185 of \cite{REM}
(and the discussion in the preceding pages) where we show that every algebra morphism between Stein algebras is continuous.  Notice that since
Stein algebras are Fr\'echet algebras, the boundedness of morphism is equivalent to continuity.
\end{rmk}

\section{Dagger affinoid spaces}

In this section we study the dual category of the category of dagger affinoid algebras: The category of dagger affinoid spaces.

\index{dagger affiond space}
\begin{defn}
The category of strict $k$-dagger affinoid spaces has as objects the topological spaces $\cM(A)$, for $A$ strict $k$-dagger affinoid algebras,
and as morphisms the continuous maps which are duals of strict $k$-dagger algebra morphisms.
The category of (non-strict) $k$-dagger affinoid spaces has as objects the topological spaces $\cM(A)$, for $A$ (non-strict) $k$-dagger affinoid algebras, 
and as morphisms the continuous maps which are duals of $k$-dagger algebra morphisms.

We do not reserve particular symbols for the categories of dagger affinoid spaces, we denote them simply by $(\bS \bAff_k^\dagger)^\circ$ 
and $(\bAff_k^\dagger)^\circ$.
\end{defn}

We use the bornological spectra of dagger affinoid algebras as a mean to give a concrete representation of the objects of $(\bS \bAff_k^\dagger)^\circ$ 
and $(\bAff_k^\dagger)^\circ$.

\begin{prop} \label{prop_strict_germs}
If $k$ is non-archimedean then there is an equivalence of categories
\[ (\bS \bAff_k^\dagger)^\circ \cong \bS \bGerms_k \]
between the category of strict $k$-dagger affinoid spaces and the category of wide strict $k$-affinoid germs of analytic spaces in the sense of Berkovich. Moreover, the functor which associates to a strict $k$-dagger affinoid space its corresponding analytic germ induces a homeomorphism of underlying topological spaces.
\end{prop}
\begin{proof} 
Berkovich has shown that the category of strict wide analytic germs is anti-equivalent to the category of 
strict $k$-dagger affinoid algebras, defined by Grosse-Kl\"onne in \cite{GK} (see lemma 5.1.1, page 73, in \cite{BER}). Clearly, our category of strict dagger affinoid algebras is equivalent to that of Grosse-Kl\"onne, if $k$ is non-archimedean. Therefore, we get the desired equivalence of categories. 
Given a strict $k$-dagger affinoid space $X = \cM(A)$, then $X$ is homeomorphic to the underlying space of the associated wide strict $k$-affinoid germs of analytic spaces because, by the same reasoning used in proposition \ref{prop_dagger_spectrum}, $X$ is homeomorphic to the (Berkovich) spectrum of the affinoid algebra $A'$ associated to $A$ as in theorem 1.7 of \cite{GK}. We will discuss further the association $A \mapsto A'$ after theorem \ref{thm:dagger_ber_to_dagger_rigid}.
\end{proof}

We can generalize last proposition in the following way.

\begin{thm} \label{prop_germs}
Let $k$ be non-archimedean, then there is an equivalence
\[ (\bAff_k^\dagger)^\circ \cong \bGerms_k \]
between the category of $k$-dagger affinoid spaces and the category of wide $k$-affinoid germs of analytic spaces of Berkovich. Moreover, the functor which associates to a $k$-dagger affinoid space its corresponding analytic germ induces a homeomorphism of underlying topological spaces.
\end{thm}
\begin{proof} 
Any dagger affinoid algebra $A$ can be presented as a direct limit
\[ A \cong \limind_{\rho' > \rho} A_{\rho'} \]
with $A_{\rho'}$ strict $k$-dagger affinoid algebras. So, one can associate a wide $k$-affinoid germ of Berkovich analytic space and vice versa to a dagger affinoid algebra. Any morphism of dagger affinoid algebras 
\[ A = \limind_{\rho' > \rho} A_{\rho'} \to B = \limind_{r' > r} B_{r'} \]
is given by a system of morphism, as we have already discussed in the proof of proposition \ref{prop:non_strict_bounded_morphism}. By definition, a morphism of germs of Berkovich analytic spaces is precisely a morphism of the systems of its neighborhoods, or, using language of appendix \ref{pro_appendix}, is a morphism of pro-analytic spaces. Therefore, we obtain the required equivalence of categories. The claim about the homeomorphism of the underlying topological spaces follows directly from the analogous claim of \ref{prop_strict_germs}.
\end{proof}

\begin{prop} \label{prop:max_arch}
Let $A$ be a dagger affinoid algebra over $\C$ or $\R$, then there is a bijection of sets
\[ \cM(A) \cong \Max(A). \]
\end{prop}
\begin{proof} 
Let $|\cdot| \in \cM(A)$. $\Ker(|\cdot|)$ is a closed prime ideal of $A$ and we can extend the multiplicative seminorm $|\cdot|$ to
\[ K = \Frac\l ( \frac{A}{\Ker(|\cdot|)} \r ). \]
Then, we can take the completion $\widehat{K}$ which is therefore a field complete with respect to an archimedean absolute value, hence $\widehat{K} \cong \C$ or $\widehat{K} \cong \R$.
This means that we have a character $\chi: A \to \C$ and to give such a character is the same as giving a system of characters $\chi_\rho: A_\rho \to \C$, where each of the $A_\rho$ can be chosen to be a Stein algebra, by theorem \ref{thm:stein_filtration}. By the theory of Stein spaces, cf. chapter 5 of \cite{REM}, $\cM(A_\rho) \cong \Max(A_\rho)$ for Stein algebras (where $\Max$ denotes finitely generated maximal ideals). Finally, as a consequence of proposition \ref{lattice_ideal}, we have that
\[ \Max(A) = \limpro_{\rho > 1} \Max(A_\rho). \]
So, we showed that each bounded character of $A$ defines a point of $\Max(A)$, proving the proposition.
\end{proof}

Finally, one can easily prove the following proposition.

\begin{prop} \label{prop:compac_stein}
Let $A$ be a dagger affinoid algebra, then
\[ \cM(A) \cong \{ x \in \cM(W_k^n(\rho)) | f_1(x) = 0, \dots, f_r(x) = 0, f_i \in W_k^n(\rho) \} \]
for some $n \in \N$ and some polyradius $\rho > 0$.
\end{prop}
\begin{proof} 
The statement of the proposition is simply a geometric restatement of the definition of dagger affinoid algebras. Indeed, for any such algebra it must exist an isomorphism
\[ A \cong \frac{W_k^n(\rho)}{I}, \ \ n \in \N, \ \ I \subset W_k^n(\rho) \]
which correspond to the claimed closed immersion (in the sense of definition \ref{defn:closed_immersion}) of $\cM(A)$ in $\cM(W_k^n(\rho))$.
\end{proof}

\begin{cor} \label{prop_loc_arcwise_connectedness}
Every $k$-dagger affinoid space $\cM(A)$ is locally arcwise connected.
\end{cor}
\begin{proof} 
If $k$ is non-archimedean the claim has been proved by Berkovich in theorem 3.2.1 of  \cite{BER2}, for affinoid spaces. Thanks to proposition \ref{prop_strict_germs} and theorem \ref{prop_germs}. We know that the underlying topological space of a dagger affinoid space is homeomorphic to the associated affinoid space. 

In the archimedean case, we can apply proposition \ref{prop:compac_stein} to describe $\cM(A)$ as a Zariski closed subset of $\cM(W_k^n)$, for some $n \in \N$ and we can deduce the corollary using the fact that analytic sets (in complex geometry) are well-known to be locally arcwise connected.
\end{proof}

\section{Dagger affinoid subdomains}

\index{dagger affiond subdomain}
\index{dagger affiond localization}
\begin{defn}
Let $A$ be a $k$-dagger affinoid algebra. A closed subset $V \subset \cM(A)$ is called a \emph{$k$-dagger affinoid subdomain} of $\cM(A)$ if there exists a $k$-dagger affinoid algebra
$A_V$ and a morphism $\phi^*: A \to A_V$ whose opposite morphism of spaces $\phi: \cM(A_V) \to \cM(A)$ satisfies $\Im(\phi) \subset V$ and for any given $\psi^*: A \to B$ with
$\Im(\psi) \subset V$ there exists a unique bounded homomorphism $A_V \to B$ such that the following diagram
\[
\begin{tikzpicture}
\matrix(m)[matrix of math nodes,
row sep=2.6em, column sep=2.8em,
text height=1.5ex, text depth=0.25ex]
{A & A_V  \\
& B \\};
\path[->,font=\scriptsize]
(m-1-1) edge node[auto] {$\phi^*$} (m-1-2);
\path[->,font=\scriptsize]
(m-1-1) edge node[auto] {$\psi^*$} (m-2-2);
\path[->,font=\scriptsize]
(m-1-2) edge node[auto] {}  (m-2-2);
\end{tikzpicture}
\]
is commutative, for any $K$-dagger affinoid algebra $B$, where $K/k$ is an extension of valued fields.

We say that $V \subset \cM(A)$ is a \emph{strict $k$-dagger affinoid subdomain} if the associated algebra $A_V$ is a strict $k$-dagger algebra. 
The map on spectra
$\phi: \cM(A_V) \to \cM(A)$ opposite to a $\phi^*$ as above is called \emph{dagger affinoid subdomain embedding}. Furthermore, we will often refer to a map $A \to A_V$ whose map on spectra $\cM(A_V) \to \cM(A)$ is a dagger affinoid subdomain embedding as a \emph{dagger affinoid localization} or simply a localization. 
\end{defn}

\begin{prop}
Let $U \subset V \subset X$ be two inclusions of $k$-dagger affinoid subdomains, then $U$ is a subdomain in $X$.
\begin{proof} 
Easy consequence of the universal properties involved.
\end{proof}
\end{prop}

\begin{prop} \label{prop:fiber_product_subdomain}
Let $\phi: X = \cM(A) \to Y = \cM(B)$ be a morphism of $k$-dagger affinoid spaces and $U \subset Y$ a $k$-dagger affinoid subdomain. Then, $\phi^{-1}(U)$ is a subdomain of $X$.
\end{prop}
\begin{proof} 
Consider the following cartesian square
\[
\begin{tikzpicture}
\matrix(m)[matrix of math nodes,
row sep=2.6em, column sep=2.8em,
text height=1.5ex, text depth=0.25ex]
{X \times_Y U & U  \\
X & Y \\};
\path[->,font=\scriptsize]
(m-1-1) edge node[auto] {$p_1$} (m-1-2);
\path[->,font=\scriptsize]
(m-1-1) edge node[auto] {$p_2$} (m-2-1);
\path[->,font=\scriptsize]
(m-1-2) edge node[auto] {$i$}  (m-2-2);
\path[->,font=\scriptsize]
(m-2-1) edge node[auto] {$\phi$}  (m-2-2);
\end{tikzpicture}.
\]
Then $p_2$ is injective and $\Im(p_2) \subset \phi^{-1}(U)$, by the definition of fibre product. 
Suppose now that $\psi: Z \to X$ is another affinoid map such that $\Im(\psi) \subset \phi^{-1}(U)$. This implies that $\Im(\phi \circ \psi) \subset U$ so $\phi \circ \psi: Z \to Y$ factors through $U$. By the fact that $U$ is a dagger affinoid subdomain in $Y$ it follows that there exists a unique map $\psi': Z \to U$ which makes 
the triangle
\[
\begin{tikzpicture}
\matrix(m)[matrix of math nodes,
row sep=2.6em, column sep=2.8em,
text height=1.5ex, text depth=0.25ex]
{Z & U  \\
& Y \\};
\path[->,font=\scriptsize]
(m-1-1) edge node[auto] {$\psi'$} (m-1-2);
\path[->,font=\scriptsize]
(m-1-1) edge node[auto] {$\phi \circ \psi$} (m-2-2);
\path[->,font=\scriptsize]
(m-1-2) edge node[auto] {$i$}  (m-2-2);
\end{tikzpicture}
\]
commutative. Thus, we have maps $Z \to U$ and $Z \to X$ and by the universal property that characterizes $X \times_Y U$ we deduce that there exists a unique map $\psi'': Z \to X \times_Y U$ of objects over $Y$. Hence, $\psi = p_2 \circ \psi'$, with $\phi'$ determined in a unique way, so $p_2$ represents every map with $\Im(\psi) \subset \phi^{-1}(U)$.
Therefore $\phi^{-1}(U)$ is a $k$-dagger affinoid subdomain of $B$, whose associated dagger affinoid algebra is $B \otimes^\dagger_{A} A_U$.
\end{proof}

\begin{cor}
Let $U, V \subset X = \cM(A)$ be two $k$-dagger affinoid subdomains of $X$ then $U \cap V$ is a $k$-dagger affinoid subdomain of $X$.
\end{cor}
\begin{proof} 
Particular case of the previous proposition.
\end{proof}

\begin{prop} \label{prop_subdomain_max_ideals}
Let $A$ be a strict $k$-dagger affinoid algebra and $\phi: \cM(A_U) \rhook \cM(A)$ a strict $k$-dagger affinoid subdomain embedding, for a subdomain $U \subset \cM(A)$. Then
\ben
\item $\phi$ is injective and satisfies $\phi(\cM(A_U)) = U$;
\item for any $x \in \Max(A_U)$ and $n \in \N$, the map $\phi^*: A \to A_U$, opposite of $\phi$, induces an isomorphism 
      \[ \frac{A}{\fm_{\phi(x)}^n} \to \frac{A_U}{\fm_x^n}; \]
\item for any $x \in \Max(A_U)$, $\fm_x = \phi^*(\fm_{\phi(x)}) A_U$.
\een
\begin{proof} 
The proof is very similar to the proof of proposition 7.2.2/1 of \cite{BGR}. We write the full proof for the sake of clarity.
Let $y \in \Max(A) \cap U$ and consider, for any fixed $n \in \N$, the following diagram
\begin{equation} \label{eqn:max_ideal}
\begin{tikzpicture}
\matrix(m)[matrix of math nodes,
row sep=2.6em, column sep=2.8em,
text height=1.5ex, text depth=0.25ex]
{A                         & A_U  \\
 \frac{A}{\fm_y^n} & \frac{A_U}{\phi^*(\fm_y^n) A_U} \\};
 \path[->,font=\scriptsize]
(m-1-1) edge node[auto] {$\phi^*$} (m-1-2);
\path[->,font=\scriptsize]
(m-1-1) edge node[auto] {$\pi$} (m-2-1);
\path[->,font=\scriptsize]
(m-1-2) edge node[auto] {$\pi_U$} (m-2-2);
\path[->,font=\scriptsize]
(m-2-1) edge node[auto] {$\sigma$}  (m-2-2);
\path[->,font=\scriptsize, dashed]
(m-1-2) edge node[auto] {$\a$}  (m-2-1);
\end{tikzpicture}
\end{equation}
where $\pi$ and $\pi_U$ are the canonical strict epimorphisms and $\sigma$ is induced by $\phi^*$. Since $\phi$ represent all the dagger affinoid maps into
$U$ there exists a unique map $\a: A_U \to \frac{A}{\fm_y^n}$ making the upper triangle commutative. Both the maps $\pi_U$ and 
$\sigma \circ \a$ make the diagram 
\[
\begin{tikzpicture}
\matrix(m)[matrix of math nodes,
row sep=2.6em, column sep=2.8em,
text height=1.5ex, text depth=0.25ex]
{A                         & A_U  \\
  & \frac{A_U}{\phi^*(\fm_y^n) A_U} \\};
 \path[->,font=\scriptsize]
(m-1-1) edge node[auto] {$\phi^*$} (m-1-2);
\path[->,font=\scriptsize]
(m-1-1) edge node[auto] {$\sigma \circ \pi$} (m-2-2);
\path[->,font=\scriptsize, dashed]
(m-1-2) edge node[auto] {$\pi_U \text{ or } \sigma \circ \a$ }  (m-2-2);
\end{tikzpicture}
\]
commutative. Therefore, for the universal property of $A_U$, $\pi_U = \sigma \circ \a$, which implies that also the lower triangle of the
diagram (\ref{eqn:max_ideal}) commutes. This commutativity relation easily implies that $\sigma$ and $\a$ are surjective and also
\[ \ker(\pi_U) = \phi^*(\fm_y^n) A_U \subset \ker \a, \]
proving that $\sigma$ is bijective. This settles claims (2) and (3). For claim (1) we can consider $n = 1$ and we see that $\phi^*(\fm_y) A_U$ is a maximal ideal, showing the injectivity of $\phi$ on the subsets of maximal ideals
\[
\begin{tikzpicture}
\matrix(m)[matrix of math nodes,
row sep=2.6em, column sep=2.8em,
text height=1.5ex, text depth=0.25ex]
{ \Max(A_U) & \Max(A)  \\
  \cM(A_U)  & \cM(A) \\};
 \path[right hook->,font=\scriptsize]
(m-1-1) edge node[auto] {$$} (m-1-2);
\path[right hook->,font=\scriptsize]
(m-1-1) edge node[auto] {$$} (m-2-1);
\path[right hook->,font=\scriptsize]
(m-1-2) edge node[auto] {$$ }  (m-2-2);
\path[->,font=\scriptsize]
(m-2-1) edge node[auto] {$$ }  (m-2-2);
\end{tikzpicture}.
\]
This settles the archimedean case, because in this case $\Max(A_U) = \cM(A_U)$ and $\Max(A) = \cM(A)$, by proposition \ref{prop:max_arch}. The non-archimedean case follows by the application of proposition 2.1.15 of \cite{BER2} to the injection $\Max(A_U) \rhook \Max(A)$.
\end{proof}
\end{prop}

\begin{cor} \label{cor:subdomain_homeo}
Let $A$ be a $k$-dagger affinoid algebra and $\phi: \cM(A_U) \rhook \cM(A)$ a $k$-dagger affinoid embedding into
$U \subset \cM(A)$. Then, $\cM(A_U) = U$.
\end{cor}
\begin{proof}
Notice that the non-strict case occurs only for non-archimedean base fields, and in this case the same reasoning of Proposition 2.2.4 of \cite{BER2} applies.
\end{proof}

We conclude this section by introducing some particular classes of dagger affinoid subdomains. Let $A$ be a $k$-dagger affinoid algebra and $X = \cM(A)$ its spectrum. We use the following notations for $f = (f_1, \ldots, f_n) \in A^n$ and
$g = (g_1, \ldots, g_m) \in A^m$, $h \in A$, $r = (r_i) \in \R_{> 0}^n$, $s = (s_i) \in \R_{> 0}^m$
\[ X (r_1^{-1} f_1, \ldots, r_n^{-1} f_n) = X(r^{-1} f) = \l \{ x \in X | |f_i(x)| \le r_i  \r \} \]
and 
\[ X (r_1^{-1} f_1, \ldots, r_n^{-1} f_n, s_1^{-1} g_1^{-1}, \ldots, s_m^{-1} g_m^{-1}) = \]
\[ = X(r^{-1} f, s^{-1} g^{-1}) = \l \{ x \in X | |f_i(x)| \le r_i, |g_j(x)| \ge s_j \r \}. \]
While if $f_1, \ldots, f_n, h$ do not have any common zero
\[ X \l (\frac{r_1^{-1} f_1}{h}, \ldots, \frac{r_n^{-1}f_n}{h} \r ) = X \l( \frac{r^{-1} f}{h} \r ) = \l \{ x \in X | \frac{|f_i(x)|}{|h(x)|}  \le r_i \r \} \]

\index{Weierstrass subdomain}
\index{Laurent subdomain}
\index{rational subdomain}
\begin{prop}
Let $A$ be a $k$-dagger affinoid algebra and $X = \cM(A)$ its spectrum then
\ben
\item for any $f_1, \ldots, f_n \in A$ and $r = (r_i) \in \R_{> 0}^n$, $X (r_1^{-1} f_1, \ldots, r_n^{-1} f_n)$ is a $k$-dagger affinoid subdomain of $X$, called a \emph{Weierstrass subdomain};
\item for any $f_1, \ldots, f_n, g_1, \ldots, g_m \in A$ and $r = (r_i) \in \R_{> 0}^n$, $s = (s_i) \in \R_{> 0}^m$, $X ( r_1^{-1} f_1, \ldots,r_n^{-1} f_n,  s_1^{-1} g_1^{-1}, \ldots, s_m^{-1} g_m^{-1})$ is a $k$-dagger affinoid subdomain of $X$, called a \emph{Laurent subdomain};
\item for any $f_1, \ldots , f_n, g \in A$ and $r = (r_i) \in \R_{> 0}^n$, such that $f_1, \ldots, f_n, g$ do not have any common zero, $X( \frac{r^{-1} f}{g})$ is a $k$-dagger affinoid subdomain of $X$, called a \emph{rational subdomain}.
\een
\end{prop}
\begin{proof} 
\ben
\item Let $\phi^*: A \to B$ be a homomorphism of dagger affinoid algebras whose opposite map $\phi: \cM(B) \to \cM(A)$ is such that $\Im(\phi) \subset X(r_1^{-1} f_1, \ldots, r_n^{-1} f_n)$. For all $x \in \cM(B)$ one has
      \[ |\phi(f_i)(x)| \le r_i \]
      so $\phi(f_i) \in r_i B^\ov$. Then, by the universal property of $A \lt r^{-1} f \gt^\dagger$ (cf. proposition \ref{prop:univ_property_laurent}) we have a unique map
      \[ A \to A \lt r_i^{-1} f \gt^\dagger = \frac{ A \lt r_i^{-1} X \gt^\dagger}{(X - f)} \]
      which represents all dagger affinoid maps with image contained in $X(r_1^{-1} f_1, \ldots, r_n^{-1} f_n)$. So, the association $X(r^{-1} f) \mapsto A \lt r^{-1} f \gt^\dagger$ gives 
      the required $k$-dagger affinoid subdomain.
\item We can use the same reasoning on the universal property characterizing $A \lt r^{-1} f, s^{-1} g^{-1} \gt^\dagger$, given in proposition \ref{prop:univ_property_laurent}.
\item Again the same reasoning on the universal property characterizing $A \lt \frac{r^{-1} f}{g} \gt^\dagger$, showed in proposition \ref{prop:univ_property_rational}.
\een
\end{proof}

\begin{prop} \label{prop:prop_subdomains}
Let $A$ be a $k$-dagger affinoid algebra and $X = \cM(A)$ its spectrum;
\ben
\item let $U, V \subset X$ be dagger affinoid subdomains. If they are both Weierstrass then also $U \cap V$ is Weierstrass. The same is true for Laurent and rational;
\item every Weierstrass subdomain is a Laurent subdomain and any Laurent subdomain is a rational subdomain;
\item let $\phi: X \to Y$ be a map between dagger affinoid spaces and $U \subset Y$ a dagger affinoid subdomain. If $U$ is Weierstrass (resp. Laurent, resp. rational) then
     $\phi^{-1}(U)$ is Weierstrass (resp. Laurent, resp. rational);
\een
\end{prop}
\begin{proof} 
We use the same reasoning of \cite{BGR}, 7.2.3/7, 7.2.3/8 and 7.2.3/6.
\ben
\item The Weierstrass and Laurent cases are very easily to deduce from definitions, so we only check that the intersection of two rational subdomains is a rational subdomain. Let 
\[ U = X \l ( \frac{r_1^{-1} f_1}{h}, \ldots, \frac{r_n^{-1} f_n}{h} \r ) \]
\[ V = X \l ( \frac{s_1^{-1} g_1}{l}, \ldots, \frac{s_m^{-1} g_m}{l} \r ) \]
then $(r_1^{-1} f_1, \ldots, r_n^{-1} f_n, h) = A$ and $(s_1^{-1} g_1, \ldots, s_m^{-1} g_m, l) = A$. Thus, also the products $f_i g_j, f_i l, g_j l, h l$ generates the unit ideal. So 
\[ Y = X \l ( \frac{r_1^{-1} s_1^{-1} f_1 g_1}{h l}, \ldots, \frac{r_n^{-1}  s_m^{-1} f_n g_m}{h l} \r ) \]
is a well-defined rational subdomain of $X$. We have that $U \cap V \subset Y$ since if $x \in U \cap V$ then
\[ |f_i(x)| \le r_i |h(x)|, |g_j(x)| \le s_j |l(x)| \]
for all $i,j$, so by the multiplicativity of $|\cdot|$ we have that for all $i$ and $j$
\[ |(f_i g_j) (x)| \le r_i s_j |(h l) (x)|. \]
For the other inclusion, let $x \in Y$ then it satisfies the relation
\[ |f_i(x)| |l (x)| \le r_i |h(x)| |l (x)| \]
for all $i$ . Since $|h(x)| |l (x)| \ne 0$, otherwise there is a common zero between the functions defining $Y$, one can cancel $l(x)$ from the above relation obtaining
\[ |f_i(x)| \le r_i |h (x)|. \]
The same can be done for $g_j$ and $h$ in place of $f_i$ and $l$, proving the claim.
\item Any Weierstrass subdomain is Laurent by definition. Furthermore, by the previous point we see that for any Laurent subdomain we have the relation
      \[ X(r_1^{-1} f_1, \ldots, r_n^{-1} f_n, s_1^{-1} g_1^{-1}, \ldots, s_m^{-1} g_m^{-1}) = X(\frac{r_1^{-1} f_1}{1})  \cap \ldots \cap  X(\frac{1}{s_m^{-1}  g_m}) \]
      whose right hand side is a finite intersection of rational subdomains, therefore it is a rational subdomain.
\item Let $U = Y(r_1^{-1} f_1, \ldots, r_n^{-1} f_n)$ be a Weierstrass subdomain of $Y$. To the map of dagger affinoid spaces $\phi: X \to Y$ corresponds a map of dagger affinoid algebras 
      $\phi^*: A_Y \to A_X$. So, we have that
      \[ \phi^{-1}(Y(r_1^{-1} f_1, \ldots, r_n^{-1} f_n)) = \{ x \in X | |\phi(f_i)(x)| \le r_i \} = \]
      \[ = X(r_1^{-1} \phi^*(f_1), \ldots, r_n^{-1} \phi^*(f_n)). \]
      The same argument works for Laurent and rational subdomains.
\een
\end{proof}

The last proposition we discuss is needed in Chapter 5.

\begin{prop} \label{prop:epi_weie}
Suppose that $f: A \to B$ and $g: B \to C$ are maps of dagger affinoid spaces. Suppose that $f$ and $g \circ f$ are Weierstrass (resp. Laurent, resp. rational) localizations, then $g$ is a Weierstrass (resp. Laurent, resp. rational) localization.
\end{prop}
\begin{proof}
The proof can be easily worked out by checking the universal properties that characterize $f$, $g \circ f$ and $g$. The three cases are similar. We work out the details only of the (strict) Weierstrass one for the sake of brevity. 

Under the hypothesis of the proposition we can write
\[ B \cong A \lt f_1, \ldots, f_n \gt^\dagger, \ \  C \cong A \lt g_1, \ldots, g_m \gt^\dagger  \]
for some $f_1, \ldots, f_n, g_1, \ldots, g_m \in A$. Notice that we also have a commutative diagram of morphism of dagger affinoid spaces
\[
\begin{tikzpicture}
\matrix(m)[matrix of math nodes,
row sep=2.6em, column sep=2.8em,
text height=1.5ex, text depth=0.25ex]
{ \cM(C) & \cM(B)  \\
         & \cM(A) \\};
 \path[->,font=\scriptsize]
(m-1-1) edge node[auto] {$g^*$} (m-1-2);
\path[right hook->,font=\scriptsize]
(m-1-1) edge node[auto] {$(g \circ f)^*$} (m-2-2);
\path[right hook->,font=\scriptsize]
(m-1-2) edge node[auto] {$f^*$ }  (m-2-2);
\end{tikzpicture}.
\]
Therefore, $g^*$ is injective and $\Im(g^*)$ is contained in the dagger affinoid subdomain of $\cM(A)$ identified by $f^*$. So, the map $g$ is unique by the definition of dagger affinoid subdomain. This shows that there is an isomorphism
\[ A \lt g_1, \ldots, g_m \gt^\dagger \cong A \lt f_1, \ldots, f_n \gt^\dagger \lt f(g_1), \ldots, f(g_m) \gt^\dagger \cong C \]
which proves that $g$ is a Weierstrass localization.
\end{proof}

\chapter{Tate's acyclicity for dagger rational subdomains} \label{chap:tate}

This chapter is devoted to the proof of the Tate's acyclity theorem in the context of dagger affinoid algebras. For the moment we will prove this theorem for finite rational coverings of strict affinoid algebras. Our aim is to show that Tate's argument can be easily adapted in the dagger settings providing a proof that works uniformly over any base field. We will deal in chapter \ref{chp:global} with the non-strict case, which occurs only for non-archimedean base field, where we use the same reduction argument to the strict case discussed in \cite{BER2}. In the next chapter we will prove the Gerritzen-Grauert theorem for dagger affinoid spaces, whose main consequence is that the Tate's acyclicity theorem we prove in this chapter extends to all finite dagger affinoid coverings.

\section{Tate's acyclicity theorem}

In this chapter we use the notation $\cO_X(U)$ to denote the $k$-dagger affinoid algebra associated to a $k$-dagger affinoid subdomain $U \subset X$ of a $k$-dagger affinoid space $X = \cM(A)$. This algebra is simply what we denoted by $A_U$ in the previous chapter. We suppose also that $A$ is a strict $k$-dagger affinoid algebra.

\index{dagger affinoid covering}
\begin{defn}
We say that a family of dagger affinoid subdomains $U_i \subset X = \cM(A)$ is a \emph{dagger affinoid covering} if $\underset{i \in I}\bigcup U_i = X$.
\end{defn}

By the properties of the subdomains of $\cM(A)$ it is clear that the family of finite affinoid coverings give to $\cM(A)$ the structure of a $G$-topological space.
Moreover, the association $U \mapsto \cO_X(U)$ is a presheaf for this $G$-topology. In this chapter we check that this presheaf is an acyclic sheaf for rational coverings. To do this we use the same reduction argument to Laurent coverings used in the classical Tate's proof, that can be found for example in \cite{BGR}, section 8.2. We need some lemmata.

\begin{lemma} \label{inj_lemma}
Let $X = \cM(A)$ be a strict affinoid space and $X = \underset{i \in I}\bigcup U_i$ a strict dagger affinoid covering, then $\cO_X(X) \to \underset{i \in I}\prod \cO_X(U_i)$ is injective.
\end{lemma}
\begin{proof}
We can prove this lemma reasoning in the same way of proposition 7.3.2/3 of \cite{BGR}, but we can give a shorter argument. We note
that the non-archimedean case is an easy consequence of 2.4 of \cite{GK} in combination with proposition 7.3.2/3 of \cite{BGR}. The archimedean case follows directly from \ref{thm_steinness} (where we prove that in fact the algebras
$\cO_X(U_i)$ coincides with the algebras of analytic functions on $U_i$ given by the immersion $\cM(A) \rhook \C^n$, when both $\cM(A)$ and $\cM(A_{U_i})$ are thought as compact Stein subsets of $\C^n$, for some $n \in \N$). It is easy to check that there are no circular arguments in this forward reference.
\end{proof}

\begin{lemma}
Let $A$ be a strict $k$-dagger affinoid algebra. Then, for any $f \in A$ there exists a point $x \in \cM(A)$ such that
\[ |f(x)| = \sup_{x \in \cM(A)} |f_i(x)|. \] 
\end{lemma}
\begin{proof}
The spectral semi-norm is well-defined because $\cM(A)$ is compact as consequence of the fact that $A$ is a bornological m-algebra. So, $\underset{x \in \cM(A)}\sup |f_i(x)|$ must be attained at some $x \in \cM(A)$.
\end{proof}

\begin{lemma} \label{lemma:max_covering}
Let $X = \cM(A)$ with $A$ a strict $k$-dagger affinoid algebra and $f_1, \ldots, f_n \in A$. Then, the function
\[ \a(x) \doteq \max_{1 \le i \le n} |f(x)| \]
assume its minimum in $X$.
\end{lemma}
\begin{proof} 
The proof is identical to \cite{BGR} lemma 7.3.4/7.
If the $f_i$'s have a common zero then the statement holds trivially. Otherwise, they generate the unit ideal in $A$. Consider the covering of $X$
given by
\[ X_i = X \l ( \frac{f_1}{f_i}, \ldots, \frac{f_r}{f_i} \r ) = \l\{ x \in X | \a(x) = |f_i(x)| \r\}. \]
By the previous lemma $|f_i(x)|$ assumes its minimum on $X_i$, because $f_i$  is invertible on $X_i$ and the minimum of $|f_i(x)|$ on $X_i$ is realised as the maximum of $f_i^{-1}$. So, $\a$ assumes its minimum in $X$ which is the least of the minimum of the $f_i$'s over the $X_i$'s.
\end{proof}

\index{Laurent dagger affinoid covering}
\begin{defn}
Let $A$ be a strict $k$-dagger affinoid algebra, $f_1, \ldots, f_n \in A$ and $X = \cM(A)$. Then, each 
\[ \sU_i = \l \{ X(f_i), X(f_i^{-1}) \r \} \]
is a $k$-dagger affinoid covering of $X$. We denote by $\sU_1 \times \ldots \times \sU_n$ the covering consisting of all the intersections of the form 
$U_1 \cap \ldots \cap U_n$ where $U_i \in \sU_i$. We call $\sU_1 \times \ldots \times \sU_n$ the \emph{Laurent covering of $X$ generated by $f_1, \ldots, f_n$}.
\end{defn}

More explicitly, the elements of the Laurent covering generated by $f_1, \ldots, f_n$ are Laurent subdomains of the form
\[ X(f_1^{\mu_1} , \ldots, f_n^{\mu_n}) \]
with $\mu_i = \pm 1$.

\index{rational dagger affinoid covering}
\begin{defn}
Let $A$ be a strict $k$-dagger affinoid algebra, $f_1, \ldots, f_n \in A$ with no common zeros and $X = \cM(A)$. Then
\[ \sU \doteq \l \{ X \l ( \frac{f_1}{f_i} , \ldots, \frac{f_n}{f_i} \r ) \r \}_{i = 1, \ldots, n} \]
is a $k$-dagger affinoid covering called \emph{the rational covering of $X$ generated by $f_1, \ldots, f_n$}.
\end{defn}

We call the $G$-topology on $X$ induced by finite covering by strict rational domains the rational $G$-topology on $X$.

\begin{lemma}
Let $\sU$ be a rational covering of $X$, then there exists a Laurent covering $\sV$ of $X$ such that for any $V \in \sV$ the covering $\sU|_V$\footnote{$\sU|_V$ denotes the elements of $\sU$ which are subsets of $V$.} is a rational
covering of $V$, which is generated by units in $\cO_X(V)$.
\end{lemma}
\begin{proof}
The proof is very similar to \cite{BGR} lemma 8.2.2/3.
Let $f_1, \ldots, f_n \in \cO_X(X)$ be the elements which generate the rational covering $\sU$. Choose a constant $c \in k^\times$ such that
\[ |c|^{-1} < \inf_{x \in X} ( \max_{1 \le i \le n} |f_i(x)|). \] 
The element $c$ can always be found and $|c| > 0$ by lemma \ref{lemma:max_covering}.

Let $\cV$ be the Laurent covering generated by $c f_1, \ldots, c f_n$. Consider the sets 
\[ V = X((c f_1)^{\mu_1}, \ldots, (c f_n)^{\mu_n}) \in \sV \]
where $\mu_i = \pm 1$. We can assume that there exists a $s \in \{0 , 1, \ldots, n\}$ such that $\mu_1 = \ldots = \mu_s = 1$ and $\mu_{s + 1} = \ldots = \mu_n = -1$. Then
\[ X \l (\frac{f_1}{f_i}, \ldots, \frac{f_n}{f_i} \r ) \cap V = \void \]
for $i = 1, \ldots, s$, because 
\[ \max_{1 \le i \le s} |f_i(x)| \le |c|^{-1} < \max_{1 \le i \le n} |f_i(x)| \]
for all $x \in V$. In particular for all $x \in V$
\[ \max_{1 \le i \le n } |f_i(x)| = \max_{s + 1 \le i \le n} |f_i(x)| \]
hence $\sU|_V$ is the rational covering of $V$ generated by $f_i|_V$ for $s+ 1 \le i \le n$, which are units in $V$.
\end{proof}

\begin{lemma}
Let $\sU$ be a rational covering of $X$ which is generated by units $f_1, \ldots,f_n \in \cO_X(X)$. Then, there exists a Laurent covering $\sV$ of $X$ which is a refinement of $\sU$.
\end{lemma}
\begin{proof}
Yet the proof of this lemma follows the one of a lemma of \cite{BGR}, cf. lemma 8.2.2/4. We take as $\sV$ the Laurent covering generated by all the products $f_i f_j^{-1}$ with $ 1 \le i < j \le n$. Consider $V \in \sV$. 
Defining $I = \{(i, j) \in \N^2 | 1 \le i < j \le n\}$, we can find a partition $I_1 \coprod I_2 = I$ such that
\[ V = \bigcap_{(i, j) \in I_1} X(f_i f_j^{-1}) \cap \bigcap_{(i, j) \in I_2} X(f_j f_i^{-1}). \]
We define a partial order $\tilde{<}$ on $\{1, \ldots, n\}$ with the following rule: If $(i, j) \in I_1$ then $i \tilde{<} j$, otherwise if $(i, j) \in I_2$ then $j \tilde{<} i$.
For each $i, j \in \{1, \ldots, n\}$ with $i \ne j$ then $i \tilde{<} j$ or $j \tilde{<} i$. Consider a maximal chain (which always exists because we are dealing with finite posets)
$i_1 \tilde{<} \ldots \tilde{<} i_r$ of elements of $\{1, \ldots, n\}$. Since $i_r$ is maximal we have that for any $i \in \{1, \ldots, n\}$ the relation $i \tilde{<} i_r$ must holds, which implies $|f_i(x)| \le |f_{i_r}(x)|$ for all $x \in V$, \ie
\[ V \subset X \l ( \frac{f_1}{f_{i_r}}, \ldots, \frac{f_n}{f_{i_r}} \r ). \]
\end{proof}

\begin{lemma}
Let $\cF$ be a presheaf for the rational $G$-topology of $X = \cM(A)$. If the Laurent coverings are acyclic then all rational coverings are acyclic. 
\end{lemma}
\begin{proof}
Last lemma shows that any rational covering can be refined by a Laurent covering, directly proving the proposition.
\end{proof}

\begin{thm} \label{thm:tate}
Let $X = \cM(A)$ be a strict $k$-dagger affinoid space. The presheaf $\cO_X$ defined so far is acyclic for the rational $G$-topology.
\end{thm}
\begin{proof}

By previous lemmata we can reduce the theorem to the case of Laurent coverings. Then, by induction we can reduce it to the case of the Laurent covering $\sU = \{ X(f), X(f^{-1}) \}$ generated by one element $f \in \cO_X(X)$. Writing $\cO_X(X) = A$ we have to verify that
\begin{equation} \label{eqn:tate_laurent}
0 \to A \stackrel{\a}{\to} A \lt f \gt^\dagger \times A \lt f^{-1} \gt^\dagger \stackrel{\be}{\to} A \lt f, f^{-1} \gt^\dagger \to 0
\end{equation} 
is an exact sequence, where $\a$ is the canonical injection and where we have identified $A \lt f, f^{-1} \gt^\dagger = \cO_X(X(f, f^{-1}))$ with the module of alternating 2-cochains of the \u{C}ech complex. Also, the second map is defined $\be(f_1, f_2) = f_1 - f_2$.

We can construct the following commutative diagram of $A$-modules
\[
\begin{tikzpicture}
\matrix(m)[matrix of math nodes,
row sep=2.6em, column sep=2.8em,
text height=1.5ex, text depth=0.25ex]
{    &    &  0                                                &  0 \\
    &    &  (X - f) A \lt X \gt^\dagger \times (1 - f Y) A \lt Y \gt^\dagger &  (X - f) A \lt X, X^{-1} \gt^\dagger  & 0  \\
  0 & A  &          A \lt X \gt^\dagger \times A \lt Y \gt^\dagger           &  A \lt X, X^{-1} \gt^\dagger          & 0 \\
  0 & A  &          A \lt f \gt^\dagger \times A \lt f^{-1} \gt^\dagger      &  A \lt f, f^{-1} \gt^\dagger          & 0 \\
    &    &  0                                                &  0 \\};
\path[->,font=\scriptsize]
(m-1-3) edge node[auto] {} (m-2-3);
\path[->,font=\scriptsize]
(m-2-3) edge node[auto] {} (m-3-3);
\path[->,font=\scriptsize]
(m-3-3) edge node[auto] {$\zeta$} (m-4-3);
\path[->,font=\scriptsize]
(m-4-3) edge node[auto] {} (m-5-3);
\path[->,font=\scriptsize]
(m-1-4) edge node[auto] {} (m-2-4);
\path[->,font=\scriptsize]
(m-2-4) edge node[auto] {} (m-3-4);
\path[->,font=\scriptsize]
(m-3-4) edge node[auto] {$\xi$} (m-4-4);
\path[->,font=\scriptsize]
(m-4-4) edge node[auto] {} (m-5-4);
\path[->,font=\scriptsize]
(m-2-3) edge node[auto] {$\gamma$} (m-2-4);
\path[->,font=\scriptsize]
(m-2-4) edge node[auto] {} (m-2-5);
\path[->,font=\scriptsize]
(m-3-1) edge node[auto] {} (m-3-2);
\path[->,font=\scriptsize]
(m-3-2) edge node[auto] {$\cong$} (m-4-2);
\path[->,font=\scriptsize]
(m-3-2) edge node[auto] {$i$} (m-3-3);
\path[->,font=\scriptsize]
(m-3-3) edge node[auto] {$\delta$} (m-3-4);
\path[->,font=\scriptsize]
(m-3-4) edge node[auto] {} (m-3-5);
\path[->,font=\scriptsize]
(m-4-1) edge node[auto] {} (m-4-2);
\path[->,font=\scriptsize]
(m-4-2) edge node[auto] {$\a$} (m-4-3);
\path[->,font=\scriptsize]
(m-4-3) edge node[auto] {$\be$} (m-4-4);
\path[->,font=\scriptsize]
(m-4-4) edge node[auto] {} (m-4-5);
\end{tikzpicture}
\]
where $X$ and $Y$ are variables, $\a$ and $\be$ are as in (\ref{eqn:tate_laurent}), $i$ is the canonical injection, $\delta(h_1(X), h_2(Y)) = h_1(X) - h_2(X^{-1})$ and
$\gamma$ is induced by $\delta$. The map $\gamma$ is well-defined, because if 
\[ (h_1(X), h_2(Y)) \in (X - f) A \lt X \gt^\dagger \times (1 - f Y) A \lt Y \gt^\dagger \subset A \lt X \gt^\dagger \times A \lt Y \gt^\dagger \]
then $\delta(h_1(X), h_2(Y)) \in (X - f) A \lt X, X^{-1} \gt^\dagger$,
because we can write
\[ h_1(X) = (X - f) r_1(X), \ \ h_2(Y) = (1 - f Y) r_2(Y)  \]
for some $r_1(X) \in A \lt X \gt^\dagger, r_2(Y) \in A \lt Y \gt^\dagger$, so
\[ h_1(X) - h_2(X^{-1}) = (X - f) r_1(X) - (1 - f X^{-1}) r_2(X^{-1}) = \]
\[ = (X - f) r_1(X) - X X^{-1}(1 - f X^{-1}) r_2(X^{-1}) =  \]
\[ = (X - f) r_1(X) - X^{-1} (X - f) r_2(X^{-1}) \in (X - f) A \lt X, X^{-1} \gt^\dagger. \]
The vertical maps are defined by $X \mapsto f$ and $Y \mapsto f^{-1}$. 

We check the exactness of the maps in the diagram. The column on the left is exact by the definition of $A \lt f \gt^\dagger$ and of
$A \lt f^{-1} \gt^\dagger$. Now consider the column on the right: By definition 
\[ A \lt f, f^{-1} \gt^\dagger = \frac{A \lt X, Y \gt^\dagger}{(X - f, 1 - f Y)} \]
and 
\[ A \lt X, X^{-1} \gt^\dagger = \frac{A \lt X, Y \gt^\dagger}{(1 - X Y)}. \]
The map induced by $X \mapsto f$ from $A \lt X, X^{-1} \gt^\dagger$ to $A \lt f, f^{-1} \gt^\dagger$ is surjective and its kernel is cleary equal $(X - f) A \lt X, X^{-1} \gt^\dagger$. Hence also the second column is exact.

Now we consider the rows. We have that as an $A$-module
\[ A \lt X, X^{-1} \gt^\dagger \cong A \lt X \gt^\dagger \oplus X^{-1} A \lt X^{-1} \gt^\dagger \]
and so
\[ (X - f) A \lt X, X^{-1} \gt^\dagger = (X - f)A \lt X \gt^\dagger \oplus (1 - f X^{-1}) A \lt X^{-1} \gt^\dagger \]
hence $\gamma: (X - f)A \lt X \gt^\dagger \times (1 - f X^{-1}) A \lt X^{-1} \gt^\dagger \to (X - f) A \lt X, X^{-1} \gt^\dagger$ is surjective, and consequently also $\delta$ is surjective. Now, suppose that 
\[ \delta(h_1(X), h_2(Y)) = 0. \]
This means that
\begin{equation} \label{eqn:h1_h2}
h_1(X) - h_2(X^{-1}) = 0,
\end{equation}
so that $h_1$ is a power-series in $X$ and $h_2$ is a pwoer-series in $X^{-1}$. Therefore, the equation (\ref{eqn:h1_h2}) implies that all coefficients of the two series must be zero but the ones on degree $0$. This is equivalent to say that $h_1$ and $h_2$ must be constant proving that the second row is exact. 

Consider the last row. The injectivity of $\a$ follows from lemma \ref{inj_lemma}. $\be$ is surjective because it sits in a commutative where three maps are surjective, so $\be \circ \zeta$ and $\zeta$ are surjective which implies that $\be$ is surjective too.
Finally, by the definition of $\be$ it follows easily that $\Im(\a) \subset \Ker(\be)$. To prove the other inclusion, consider $h \in \Ker(\be)$. This means that $\be(h) = 0$, so
\[ \xi^{-1}(\be(h)) \subset (X - f) A \lt X, X^{-1} \gt^\dagger \]
\[ \delta^{-1}(\xi^{-1}(\be(h))) \subset A + (X - f) A \lt X \gt^\dagger \times (1 - f Y) A \lt  Y \gt^\dagger \]
therefore
\[ h \in \zeta(\delta^{-1}(\xi^{-1}(\be(h)))) \subset A = \Im(\a). \]
This proves the theorem.
\end{proof}

\begin{rmk}
The non-archimedean case of theorem \ref{thm:tate} was already proved by Grosse-Kl\"onne in \cite{GK}, proposition 2.6. In this section we saw that the theory developed until here permits to extend the reasoning used rigid geometry to encompasse also dagger affinoid spaces defined over archimedean fields.
\end{rmk}

\chapter{The Gerritzen-Grauert theorem}  \label{chap:gg}

In this chapter we prove the main structural theorem of the theory of dagger affinoid spaces, \ie the dagger version of the celebrated Gerritzen-Grauert theorem. We have two different strategies for the proof: One works for non-archimedean base fields and the other for the archimedean ones. In particular, in the first part of this chapter we study in more details the archimedean side of the theory we developed so far and we prove some results needed in the subsequent section to deduce the dagger version of the Gerritzen-Grauert theorem.

\section{Archimedean dagger affinoid subdomains}

In this section we study more in details dagger affinoid subdomains when $k$ is an archimedean complete valued field. We see some particular features which are different from the non-archimedean case.

\index{roomy dagger affinoid subdomain}
\begin{defn}
We say that a $k$-dagger affinoid subdomain $U \subset \cM(A)$ is \emph{roomy} if has non-empty topological interior with respect to the topology of $\cM(A)$ (\ie Berkovich's topology of the spectrum).
\end{defn}

Using the same arguments of proposition 2.2.3 (iii) of \cite{BER2}, one can show that if $k$ is non-archimedean then every strict $k$-dagger affinoid subdomain is roomy, although the non-strict ones may be not. If $k$ is archimedean only the strict case exists but anyway non-roomy subdomains exist, as the following example shows.

\begin{exa}
Consider $f_1 = 2 z_1, ..., f_n = 2 z_n \in W_k^n$, with $k$ archimedean. Then 
\[ \cM(W_k^n \lt f_1, ..., f_n \gt^\dagger) = \{ x \in \D_k((1, ..., 1)^+) | |f_1| \le 1, ..., |f_n| \le 1 \} = \]
\[ =   \{ x \in \D_k((1, ..., 1)^+) | |z_1| \le \frac{1}{2}, ..., |z_n| \le \frac{1}{2} \}. \]
Moreover,
\[ \cM(W_k^n \lt f_1^{-1}, ..., f_n^{-1} \gt^\dagger) = \{ x \in \D_k((1, ..., 1)^+) | |f_1| \ge 1, ..., |f_n| \ge 1 \} = \]
\[ = \{ x \in \D_k((1, ..., 1)^+) | |z_1| \ge \frac{1}{2}, ..., |z_n| \ge \frac{1}{2} \} \]
hence
\[ \cM(W_k^n \lt f_1, ..., f_n \gt^\dagger) \cap \cM(W_k^n \lt f_1^{-1}, ..., f_n^{-1} \gt^\dagger) = 
\{ x \in \D_k((1, ..., 1)^+) | |z_1| = \frac{1}{2}, ..., |z_n| = \frac{1}{2} \} \]
which has empty topological interior as subset of $\cM(W_k^n)$.
\end{exa}

\index{Stein dagger affinoid subdomain}
\begin{defn}
Let $X$ be a $k$-dagger affinoid space and $U \subset X$ a $k$-dagger affinoid subdomain. We say that $U$ is \emph{Stein} in the following case:
\begin{itemize}
\item if $k$ is non-archimedean for every $U$;
\item if $k = \R, \C$, if $U$ is a compact Stein subset of $X$ in the sense explained below.
\end{itemize}
\end{defn}

We will focus in the case $k = \C$ whence the case $k = \R$ can be easily deduced. We proved in proposition \ref{prop:compac_stein} that for any dagger affinoid algebra there is an embedding $\cM(A) \rhook \C^n$ for a suitable $n$ ($\C^n$ modulo complex conjugation if $k = \R$). This embedding gives to $\cM(A)$ a the natural analytic structure of compact Stein set,
or equivalently a pro-analytic structure in the sense of appendix \ref{pro_appendix}. In fact, this structure agrees with the one given by the representation
\begin{equation} \label{eqn:M_stein}
\cM(A) = \limpro_{\rho > 1} \cM(A_\rho)
\end{equation}
as pro-Stein space obtained in theorem \ref{thm:stein_filtration}. 
Let $U \subset \cM(A)$ be a pro-open subset. We denote $\wtilde{\cO}_{\cM(A)}(U)$ the following algebra
\begin{equation} \label{eqn:o_tilde}
 \wtilde{\cO}_{\cM(A)}(U) = \limind_{V \supset U} \wtilde{\cO}_{\cM(A)}(V) 
\end{equation}
where $V$ ranges over all open neighborhoods of $U$ in $\C^n$. This pre-sheaf is in fact the structural sheaf of the pro-analytic site structure that $\C^n$ induces on $\cM(A)$, following the definitions we give in appendix \ref{pro_appendix}. We say that a subset $U \subset \cM(A)$ is \emph{compact Stein} if it is a compact Stein subset of $\C^n$ when $\cM(A)$ is thought as a compact Stein subset of $\C^n$. Next theorem shows that this is a weel-defined notion.


\begin{thm} \label{thm_steinness}
Let $U \subset X = \cM(A)$ be a $\C$-dagger rational subdomain of $X$, then
\[ \wtilde{\cO}_X(U) \cong A_U \]
as algebras.
\end{thm}
\begin{proof}
Let $U = X \l ( \frac{f}{g}  \r )$, with $f =(f_1, \ldots, f_n), g \in A$. It is easy to show that there is an homeomorphism of $U$ with the set
\[ Y = \cM(A \lt \frac{f}{g} \gt^\dagger) = \l \{ (x, y) \in \cM(A) \times \cM(W_k^n) | g(x) y_1 - f_1(x) = 0, \ldots, g(x) y_n - f_n(x) = 0 \r \}, \]
given by the projection on $\cM(A)$, because $f_1, \ldots, f_n, g$ does not have a common zero. We write $\pi: Y \to X$ for this projection. $\pi$ is an analytic map and so we have the induced map
\[ \pi^*: \wtilde{\cO}_X(U) \to \wtilde{\cO}_Y(Y) \cong A \lt \frac{f}{g} \gt^\dagger. \]

$\pi^*$ is injective because $\pi$ is induced by the system of maps of Stein spaces  
\[ Y_\rho = \cM(A_\rho \lt \frac{f}{g} \gt^\dagger) =  \]
\[ = \l \{ (x, y) \in \cM(A_\rho) \times \cM(\cO(\D_\C(\rho^-))) | g(x) y_1 - f_1(x) = 0, ..., g(x) y_n - f_n(x) = 0 \r \} \]
to
\[ X_\rho = \{ x \in A_\rho | |f_1(x)| < \rho |g(x)|, \ldots, |f_n(x)| < \rho |g(x)| \} \]
for all $1 < \rho < \rho'$ for some $\rho' > 1$, where $A_\rho$ is as in equation (\ref{eqn:M_stein}). So, writing $\pi_\rho: Y_\rho \to X_\rho$ for these projections we have that $\pi = \underset{\rho > 1}\limpro \pi_\rho$, which implies that $\pi^* = \underset{\rho > 1}\limind \pi_\rho^*$. The maps $\pi_\rho^*$ are all injective for $1 < \rho < \rho'$ because $\pi_\rho: Y_\rho \to X_\rho$ is a homeomorphism and therefore we can apply the identity theorem. This shows that $\pi^*$ is injective because $\limind$ is an exact functor.

To prove surjectivity, note that the $\frac{f_i}{g}$ define holomorphic functions on a neighborhood $U$, therefore they are elements of $\wtilde{\cO}_X(U)$. Their pullback to $\wtilde{\cO}_Y(Y)$ is precisely mapped on $\iota \l( \frac{f_i}{g} \r )$, when we identify $\iota: A \lt \frac{f}{g} \gt^\dagger \cong \wtilde{\cO}_Y(Y)$. Notice that $U$ is a multiplicatively convex bornological algebra such that
\[ U \cong \cM(\wtilde{\cO}_X(U)) \]
just because it is the algebra of germs of analytic functions on a compact Stein and $\frac{f_i}{g}$ is spectrally power-bounded because
\[ \max_{x \in U} | \frac{f_i}{g}(x)| \le 1.\]
Since $\wtilde{\cO}_X(U)$ we can apply remark \ref{rmk:LB_spectral_power_bounded} to deduce that $\wtilde{\cO}_X(U)$ spectrally power-bounded elements agree with weakly power-bounded elements which implies that $\frac{f_i}{g} \in \wtilde{\cO}_X(U)^\ov$. Since $\wtilde{\cO}_X(U)$ is obviously an $A$-algebra, by the universal property of $A \lt X_1, \ldots, X_n \gt^\dagger$, there is a unique map $A \lt X_1, \ldots, X_n \gt^\dagger \to \wtilde{\cO}_X(U)$ that sends $X_i$ to $\frac{f_i}{g}$. It is also obvious to check that this maps makes commutative the following diagram
\[
\begin{tikzpicture}
\matrix(m)[matrix of math nodes,
row sep=2.6em, column sep=2.8em,
text height=1.5ex, text depth=0.25ex]
{                    & A \lt X_1, \ldots, X_n \gt^\dagger  \\
  \wtilde{\cO}_X(U)  & A \lt \frac{f}{g} \gt^\dagger  \\};
\path[->,font=\scriptsize]
(m-1-2) edge node[auto] {} (m-2-1);
\path[->,font=\scriptsize]
(m-1-2) edge node[auto] {} (m-2-2);
\path[->,font=\scriptsize]
(m-2-1) edge node[auto] {$\iota^{-1} \circ \pi^*$} (m-2-2);
\end{tikzpicture}
\]
where the vertical map is defined by $X_i \mapsto \frac{f_i}{g}$ and it is surjective by the definition of $A \lt \frac{f}{g} \gt^\dagger$. Therefore, $\pi^*$ is surjective.
\end{proof}

\begin{rmk}
In a more concise way, one can restate the previous proof by saying that the morphism
\[ \pi_\rho: Y_\rho \to X_\rho \]
is \'etale and bijective, for all $1 < \rho < \rho'$, for some $\rho' > 1$, therefore it is an isomorphism of Stein spaces. So $\pi = \underset{\rho > 1}\limpro \pi_\rho$ (and dually $\pi^* = \underset{\rho > 1}\limind \pi^*_\rho$) are limits of isomorphisms which implies that $\pi$ is an isomorphism of pro-analytic spaces, as defined in appendix \ref{pro_appendix}.
\end{rmk}

\begin{cor}
Let $U \subset X = \cM(A)$ be a $k$-dagger rational subdomain of $X$, then $U$ is a compact Stein subset of $X$.
\end{cor}
\begin{proof}
$\tilde{\cO}_X(U) \cong A \lt \frac{f}{g} \gt^\dagger$ is an $\bInd$-Stein algebra.
\end{proof}

\begin{cor}
	Let $U \subset X = \cM(A)$ be a $k$-dagger rational subdomain of $X$, then $\tilde{\cO}_X(U)$ does not depend on the embedding $\cM(A) \to \C^n$.
\end{cor}
\begin{proof}
	Immediate consequence of the isomorphism $\tilde{\cO}_X(U) \cong A \lt \frac{f}{g} \gt^\dagger$.
\end{proof}

\begin{cor}
Let $U \subset X = \cM(A)$ be a $k$-dagger rational subdomain of $X$, then
\[ \wtilde{\cO}_X(U) \cong A_U \]
as bornological algebras when on $\wtilde{\cO}_X(U)$ is given the direct limit bornology given by (\ref{eqn:o_tilde}).
\end{cor}
\begin{proof}
This corollary can be proven with the same reasoning used to prove propositions \ref{prop:complex_S_O} and \ref{prop:complex_T_O}. In particular, following the first proof we gave of proposition \ref{prop:complex_S_O} we remark that is easy to show that the bornologies of $A_U$ and $\wtilde{\cO}_X(U)$ have nets. Therefore, we can apply the closed graph theorem for bornological vector spaces with nets, see theorem 2.7 of \cite{BA} or theorem 3.2 of \cite{GACH}, to deduce the bornological isomorphism.
\end{proof}

Theorem \ref{thm_steinness} can be interpreted by saying that the ringed site defined by the weak G-topology of rational dagger subdomains is equivalent to the ringed site defined by rational compact Stein subsets of $\cM(A)$ endowed with the analytic structure induced by $\C^n$. Thus, from now, on we will use the same symbol, $\cO_{\cM(A)}$, for both sheaves.

\index{good dagger affinoid subdomain}
\begin{defn}
We say that a $k$-dagger affinoid subdomain $U \subset \cM(A)$ is \emph{good} in the following cases: if $k$ is non-archimedean, and if $k$ is archimedean if $\cO_X(V) \cong A_V$.
\end{defn}

Theorem \ref{thm_steinness} shows that for $k$ archimedean all rational dagger affinoid subdomains are good in this sense. We will see that the Gerritzen-Grauert theorem implies that all dagger affinoid subdomains are good.

\begin{cor}
Any $x \in \cM(A)$ has a neighborhood basis, for the topology of the spectrum, made of good roomy $k$-dagger subdomains.
\end{cor}
\begin{proof}
Since $k$ is archimedean, Weierstrass roomy subdomains which contains $x$ in their interior form a basis of neighborhoods for the topology of $x$ in $\cM(A)$, and they are good by theorem \ref{thm_steinness}. This is a consequence of the fact that the spectrum consists only of the maximal ideals of $A$.

In the non-archimedean case the proposition is proved in the chapter $2$ of \cite{BER2} (because by theorem \ref{prop_germs} the topological space $\cM(A)$ corresponds to the topological space of the germ, in the sense of Berkovich, associated to $\cM(A)$ whose topological space agrees with the associated affinoid space).
\end{proof}

\index{inner morphism of dagger affinoid algebras}
\begin{defn}
Let $\phi: B \to C$ be a morphism of dagger affinoid algebras both of which are $A$-algebras for a $A \in \ob(\bAff_k^\dagger)$. $\phi$ is called \emph{inner with 
respect to $A$} if there exists a strict epimorphism $\pi: A \lt r_1^{-1} X_1, \dots, r_n^{-1} X_n \gt^\dagger \to B$ such that
\[ \rho_C(\phi(\pi(X_i))) < r_i \]
for all $1 \le i \le n$.
\end{defn}

\index{relative interior of a morphism of dagger affinoid spaces}
\begin{defn}
Let $\phi: \cM(A) = X \to \cM(B) = Y$ be a morphism of $k$-dagger affinoid spaces. The \emph{relative interior of $\phi$} is
the set 
\[ \Int(X/Y) \doteq \{ x \in X | A \to \cH_x \text{ is inner w.r.t. $B$} \}. \]
The complement of $\Int(X/Y)$ is called the \emph{relative boundary of $\phi$} and it is denoted by $\partial(X/Y)$. If $B = k$, the 
sets $\Int(X/Y)$ and $\partial(X/Y)$ are denoted by $\Int(X)$ and $\partial(Y)$ and they are called the \emph{interior} and the \emph{boundary}
of $X$.
\end{defn}

\begin{prop}
Let $\phi: \cM(A) = X \to \cM(B) = Y$ be a morphism of $k$-dagger affinoid spaces.
The $\Int(X/Y) \subset \cM(A)$ is an open subset with respect to the topology of the spectrum.
\end{prop}
\begin{proof}
If $\Int(X/Y) = \void$ the statement of the proposition holds trivially, so suppose that $\Int(X/Y) \ne \void$. Consider $x \in X$ and denote by $\chi_x: A \to \cH_x$ a character that belongs to the equivalence class identified by $x$. $\chi_x$ is by hypothesis an inner morphism with respect to $B$, so we can find a strict epimorphism
\[ \pi: B \lt r_1^{-1} X_1, \ldots, r_n^{-1} X_n \gt^\dagger \to A \]
such that 
\[ |\chi_x(\pi(X_i))| < r_i \]
for all $1 \le i \le n$. We define the functions $\zeta_i: X \to \R_{\ge 0}$, given by the formula
\[ \zeta_i(x') = |\chi_{x'}(\pi(X_i))| \]
for $x' \in X$. Since the topology of $X$ is by definition the weakest topology making all functions of the form $|\cdot| \mapsto |f|$ continuous, for $f \in A$, then the functions $\zeta_i$ are continuous. Therefore, the preimages of $[0, r_i) \subset \R_{\ge 0}$ by the functions $\zeta_i$ are open subsets of $X$, hence 
\[ U = \zeta_1([0, r_1)) \cap \ldots \cap \zeta_n([0, r_n)) \]
is an open subset of $X$. It is clear that $U \subset \Int(X/Y)$ because all points of $U$ are inner with respect to the same $\pi$ we chose above, proving that $\Int(X/Y)$ is an open subset of $\cM(X)$.
\end{proof}

Although we just shown that the interior of a morphism is open it might be the empty set, as the next examples show.

\begin{exa}  \label{exa:non_roomy}
\ben
\item Let $k$ be a non-Archimedean field and let $r \notin \sqrt{|k^\times|}$. The non-strict (and not dagger) affinoid algebra
\[ K_r \doteq \frac{k \lt r^{-1} X, r Y^{-1} \gt}{(X Y - 1)} \]
is a valued field extending $k$. This can be shown by noticing that every element $f \in K_r$ can be written as a power-series
\[ f = \sum_{i = - \infty}^{\infty} |a_i| X^i, \ \ \text{with } |a_i| r^i \to 0, \text{ for } |i| \to \infty. \]
The condition $r \notin \sqrt{|k^\times|}$ implies that there exists a unique $i \in \N$ such that $\max \{ |a_i| r^i \}$ is realized and without loss of generality we can suppose that $i = 0$ and $a_0 = 1$. So, we can write
\[ f = 1 - h \]
with $|h| < 1$, which implies that the series
\[ \sum_{n = 0}^\infty h^n \]
converges and is an inverse of $f$. There is also a dagger version of this example. Indeed, one can define
\[ K_r^\dagger \doteq \frac{k \lt r^{-1} X, r Y^{-1} \gt^\dagger}{(X Y - 1)} \]
which is a bornological field over $k$. To see that this algebra is a field one can notice that $K_r^\dagger$ is a one-point subdomain of the one dimensional disk $X = \cM(k \lt X \gt^\dagger)$. If we denote $x \in X$ the image of the dagger affinoid embedding $\cM(K_r^\dagger) \to X$, then 
\[ K_r^\dagger \cong \cO_{X,x} \]
and the fact that $\cO_{X,x}$ is a field is discussed in \cite{BER5} exercise 1.3.6 (ii). Geometrically one can explain the fact that $K_r^\dagger$ is field by the absence of type 1 points on the spectrum of $\cM(K_r^\dagger)$. Since analytic functions on $X$ can vanish only on points of type 1 (as a consequence of the identity theorem), for any given analytic function on $X$ it is always possible to find a neighborhood of $x$ such that $f$ does not vanish on it.

Notice that it is also possible to consider algebras like 
\[ A \otimes_k^\dagger K_r^\dagger \]
for any $A \in \ob(\bAff_k^\dagger)$, but these algebras are not $K_r^\dagger$-dagger affinoid, because $K_r^\dagger$ is not a valued field. The last thing to notice is that the completion of $K_r^\dagger$ with respect to the norm $\max \{ |a_i| r^i \}$ is $K_r$, as already noticed in \cite{BER5} exercise 1.3.6 (ii).

The interior of $\cM(K_r)$ is empty, cf. \cite{BER2} page 39, and this is also true for $\cM(K_r^\dagger)$ as a consequence of the same reasoning or as a consequence of proposition \ref{prop:interior_dimension}.

\item An analogous example can be given also for archimedean base fields. In this case, only strict dagger affinoid algebras exists and we can consider 
\[ \frac{k \lt X, Y^{-1} \gt^\dagger}{(X Y - 1)}. \]
The interior of the dagger affinoid space associated to this dagger affinoid algebra is empty. This can be deduced using proposition \ref{prop:interior_dimension} which says that the spectrum of dagger affinoid spaces with non-empty interior over archimedean base fields must have topological dimension equal two times the dimension of the affinoid space. Since the affinoid space associated to $\frac{k \lt X, Y^{-1} \gt^\dagger}{(X Y - 1)}$ is one dimensional and its spectrum is also one dimensional as topological space (it is homeomorphic to $S^1$), we deduce that $\cM(\frac{k \lt X, Y^{-1} \gt^\dagger}{(X Y - 1)})$ has empty interior.
\een
\end{exa}

The important fact to notice about the examples given so far is that they are some sort of manifestations of the same phenomena. In both cases the underlying topological spaces have the ''wrong" topological dimension: $\cM(K_r)$ is a point, hence it is zero dimensional, but $\cM(K_r)$ as affinoid space has dimension one. This is also true for $\cM(K_r^\dagger)$, defining the dimension of a dagger affinoid space as in \cite{BER2}, page 34.
Then, $\cM( \C  \lt X, X^{-1} \gt^\dagger)$ is a $1$-dimensional dagger affinoid subdomain of the unital disk but it has topological dimension $1$ and not $2$ which is the expected topological dimension in complex analytic geometry. We now make this observation more precise.

\index{dimension dagger affinoid space}
\begin{defn}
Let $A \in \ob(\bAff_k^\dagger)$. We define the \emph{dimension} of $\cM(A)$ as the Krull dimension of the ring $A \otimes_k^\dagger K$, for any $K$ such that $A \otimes_k^\dagger K$ is strictly $K$-affinoid, and it is denoted $\dim(X)$.
\end{defn}

\begin{lemma} \label{lemma:topologial_dimension}
	Let $X = \cM(A)$ be a dagger affinoid space. 
	\ben
	\item If $k$ is non-archimedean then the topological dimension of $X$ is less or equal $\dim(X)$.
	\item If $k$ is archimedean then the topological dimension of $X$ is less or equal $2 \cdot \dim(X)$.
	\een
\end{lemma}
\begin{proof}
The (non-dagger) non-archimedean case of the proposition was shown by Berkovich in theorem 3.2.6 of \cite{BER2}, therefore the dagger affinoid case follows by a direct application of theorem \ref{prop_germs}.

Consider the archimedean case. It is enough to notice that $\cM(A)$ is a compact Stein subset of $\C^n$ endowed with its ring of germs of analytic functions, as shown in theorem \ref{thm_steinness}, and the result is well-known\footnote{One can show that the topological dimension of $\cM(A)$ is less or equal to its dimension as follow. Let $n = \dim(\cM(A))$. Since $\cM(A)$ is a compact Stein it can be embedded as a closed subset of a Stein space of dimension $n$. Therefore, it is enough to show that the topological dimension of a Stein space of dimension $n$ is precisely $2 n$. It is well known that a Stein space $X$, of dimension $n$, can be stratified as a disjoint union $X = X^{(0)} \cup \ldots \cup X^{(n)}$, where $X^{(i)}$ is a Stein manifold of dimension $n - i$. This stratification can be obtained by defining $X^{(0)}$ as the non-singular locus of $X$ and, inductively, $X^{(i)}$ as the non-singular locus of $X^{(i - 1)}$. Clearly, the topological dimension of a Stein manifold of dimension $n$ is $2 n$ because it is a complex analytic manifold, proving the claim.}.
\end{proof}

Last lemma can be strengthened under the hypothesis that $X$ has non-empty interior.

\begin{prop} \label{prop:interior_dimension}
Let $X = \cM(A)$ be a dagger affinoid space which is connected. 
\ben
\item If $k$ is non-archimedean and $\Int(X) \ne \void$ then the topological dimension of $X$ is equal $\dim(X)$.
\item If $k$ is archimedean and $\Int(X) \ne \void$ then the topological dimension of $X$ is equal $2 \cdot \dim(X)$.
\een
\end{prop}
\begin{proof}
We only need to check that under the hypothesis of the proposition, then inequality of lemma \ref{lemma:topologial_dimension} is in fact an equality.

Let $x \in \Int(X)$ and let $\chi_x: A \to \cH_x$ be a character that belongs to the equivalence class identified by $x$. $\chi_x$ is inner with respect to the morphism $k \to A$, so we can find an strict epimorphism
\[ \pi: k \lt r_1^{-1} X_1, \ldots, r_n^{-1} X_n \gt^\dagger \to A \]
such that $|\chi_x(\pi(X_i))| < r_i$, for all $1 \le i \le n$. We denote $r = (r_i)$ and $r' = (r_i')$ two $n$-tuples of real numbers such that with $|\chi_x(\pi(X_i))| < r_i'< r_i$. We can calculate the following pushout diagram
\[
\begin{tikzpicture}
\matrix(m)[matrix of math nodes,
row sep=2.6em, column sep=2.8em,
text height=1.5ex, text depth=0.25ex]
{ W_k^n(r) & A  \\
  W_k^n(r') & A \otimes_{W_k^n(r)}^\dagger W_k^n(r')  \\};
\path[->,font=\scriptsize]
(m-1-1) edge node[auto] {} (m-1-2);
\path[->,font=\scriptsize]
(m-1-2) edge node[auto] {} (m-2-2);
\path[->,font=\scriptsize]
(m-1-1) edge node[auto] {} (m-2-1);
\path[->,font=\scriptsize]
(m-2-1) edge node[auto] {} (m-2-2);
\end{tikzpicture}.
\]
By proposition \ref{prop:fiber_product_subdomain} the map $A \to A \otimes_{W_k^n(r)}^\dagger W_k^n(r')$ induces a dagger affinoid subdomain embedding $\cM( A \otimes_{W_k^n(r)}^\dagger W_k^n(r')) \rhook \cM(A)$. The conditions of $\chi_x$ of being inner and the choice of $r'$, readily implies that $\chi_x$ also defines a point of $\cM(A \otimes_{W_k^n(r)}^\dagger W_k^n(r'))$ which is then non-empty and therefore $A \otimes_{W_k^n(r)}^\dagger W_k^n(r') \ne 0$. For the non-archimedean case, we can conclude by noticing that we can consider a polyradius $r''$ with $r' < r'' < r$ and with $r''_i \in \sqrt{|k|}$. In this way, $A \otimes_{W_k^n(r)}^\dagger W_k^n(r'')$ is a strict dagger affinoid algebra with a dagger affinoid embedding $\cM( A \otimes_{W_k^n(r)}^\dagger W_k^n(r'')) \rhook \cM(A)$. Clearly $\dim (A \otimes_{W_k^n(r)}^\dagger W_k^n(r'')) = \dim(A)$ (by the fact that $X$ is supposed to be connected) and the topological dimension of strict dagger affinoid spaces is equal to their dimension (cf. theorem 3.2.6 of \cite{BER2}), therefore the non-archimedean case of the proposition is proved. 

In the archimedean case, we can work out a similar reasoning using the theory of Stein spaces. Fixing again a polyradius $r' < r'' < r$ we can define a Stein space by the push-out square
\[
\begin{tikzpicture}
\matrix(m)[matrix of math nodes,
row sep=2.6em, column sep=2.8em,
text height=1.5ex, text depth=0.25ex]
{ W_k^n(r) & A  \\
  \cO(\D((r'')^-)) & A \otimes_{W_k^n(r)}^\dagger \cO(\D((r'')^-))  \\};
\path[->,font=\scriptsize]
(m-1-1) edge node[auto] {} (m-1-2);
\path[->,font=\scriptsize]
(m-1-2) edge node[auto] {} (m-2-2);
\path[->,font=\scriptsize]
(m-1-1) edge node[auto] {} (m-2-1);
\path[->,font=\scriptsize]
(m-2-1) edge node[auto] {} (m-2-2);
\end{tikzpicture}.
\]
$A \otimes_{W_k^n(r)}^\dagger \cO(\D((r'')^-))$ is a Stein algebra (because the bottom horizontal morphism is surjetive) canonically endowed with an open embedding (of pro-analytic spaces) $\cM( A \otimes_{W_k^n(r)}^\dagger \cO(\D((r'')^-)) \rhook \cM(A)$. This shows that $\dim(A \otimes_{W_k^n(r)}^\dagger \cO(\D((r'')^-))) = \dim(A)$ (by the fact that $X$ is supposed to be connected) as analytic spaces because the open embedding $\cM( A \otimes_{W_k^n(r)}^\dagger \cO(\D((r'')^-)) \rhook \cM(A)$ is an isomorphism on stalks. By the theory of Stein spaces, we know that the topological dimension of $\cM( A \otimes_{W_k^n(r)}^\dagger \cO(\D((r'')^-))$ is $2 \cdot \dim(\cM( A \otimes_{W_k^n(r)}^\dagger \cO(\D((r'')^-)))$ which implies that the topological dimension of $\cM(A)$ is $2 \cdot \dim(A)$.
\end{proof}

\begin{rmk}
The connectedness hypothesis in proposition \ref{prop:interior_dimension} is used to ensure avoiding the case when a space has several connected components some with non-empty interior and some with empty interior. One can replace it with the hypothesis that $\Int(X)$ has non-empty intersection with all connected components of $X$.
\end{rmk}

\begin{rmk}
For examples of dagger affinoid spaces $X$ for which the topological dimension is less than the expected value we refer to examples \ref{exa:non_roomy}. We also notice that considering, for instance, $K_r^\dagger$ one can show that for every $r' < r$
\[ k \lt (r')^{-1} X, r' Y \gt^\dagger \otimes_{k \lt r^{-1} X, r Y \gt^\dagger}^\dagger K_r = 0. \]
This proves that the affinoid embedding $\cM(k \lt (r')^{-1} X, r' Y \gt^\dagger \otimes_{k \lt r^{-1} X, r Y \gt^\dagger}^\dagger K_r) \to \cM(A)$ is trivial, because $\cM(k \lt (r')^{-1} X, r' Y \gt^\dagger \otimes_{k \lt r^{-1} X, r Y \gt^\dagger}^\dagger K_r) = \void$. This is where the resoning of proposition \ref{prop:interior_dimension} fails for dagger affinoid spaces with empty interior.
\end{rmk}

\begin{conj}
We expect that the converse of proposition \ref{prop:interior_dimension} holds: if the topological space $\cM(A)$ has the "right" topological dimension, then $\Int(\cM(A)) \ne \void$. 
\end{conj}

This conjecture seems harder to prove with respect to proposition \ref{prop:interior_dimension} and less interesting for our scopes.

\index{good open set}
\begin{defn}
Let $A \in \ob(\bAff_k^\dagger)$ and $G \subset \cM(A)$ be an open subset for the topology of the spectrum. We say that $G$ is \emph{good} if $G$ admits an exhaustion by good dagger affinoid subdomains, 
\ie if there exists a sequence of good subdomains $U_i \subset \cM(A)$ such that $U_i \subset \Int(U_{i + 1})$ and
\[ \bigcup_{i \in \N} U_i = G. \]
\end{defn}

For any good open subset we can define the algebra
\[ \cO_X(G) = A_G \doteq \limpro_{i \in \N} A_{U_i}, \]
where $\{ U_i \}_{i \in \N}$ is an exhaustion of $G$. In the next chapter we will study more systematically how to extend the structural sheaf we defined so far for the weak G-topology to the topology of the spectrum. Here we record some basic properties we will need later on.

\begin{prop} \label{prop_good_open_subsets}
Let $G \subset \cM(A) = X$ be a good open subset. Then,
\ben
\item the algebra $A_G$ can be equipped canonically with a structure of a bornological Fr\'echet algebra;
\item there exists a homeomorphism $G \cong \cM(G)$;
\item the space $\cM(G)$ does not depend on the choice of the exhaustion.
\een
\end{prop}
\begin{proof}
\ben
\item Let $\{U_i\}_{i \in \N}$ be an exhaustion of $G$. Then there is a canonical sequence of maps
\[ \cdots \to A_{i + 1} \to A_i \to A_{i - 1} \to \cdots \]
given by the dagger subdomain embeddings $U_i \rhook U_{i + 1}$, where $U_i = \cM(A_i)$. The $A_i$ are dagger affinoid algebras, hence we can write $A_i = \underset{\rho_i > r_i}\limind A_i^{\rho_i}$ and we get a presentation
\[ A_G = \limpro_{i \in \N} \limind_{\rho_i > r_i} A_i^{\rho_i}. \]
The request that $U_i \subset \Int(U_{i + 1})$ implies that the $U_i$ must have non-empty interior, at least for a final part of the colimit. Therefore, we can suppose that all the $U_i$ have non-empty interior. Also, the condition $U_i \subset \Int(U_{i + 1})$ implies that for any $i$ there exists a $\rho_i > r_i$ such the map $A_{i + 1} \to A_i$ factors through a map $ A_{i + 1} \to A_i^{\rho_i}$. So, reasoning by cofinality, we get an isomorphism of bornological vector spaces
\[ A_G =  \limpro_{i \in \N}  A_i^{\rho_i}. \]
Thus, we can write $A_G$ as a projective limit of Banach algebras therefore this gives to $A_G$ the structure of a Fr\'echet algebra (It is easy to check that the normality property of metrizable topological vector spaces implies that the bornology of $\underset{i \in \N}\limpro A_i^{\rho_i}$ coincides with the bornology of $\underset{i \in \N}\limpro (A_i^{\rho_i})^t$ which is a Frech\'et space by definition). 

\item We notice that the limit
\[ A_G^t \cong \l ( \limpro_{i \in \N}    A_i^{\rho_i} \r )^t \cong \limpro_{i \in I}  \l ( A_i^{\rho_i} \r )^t \]
satisfies the condition of theorem \ref{spectrum_loc_conv}, hence
\[ \cM(A_G) \cong \cM(A_G^t) \cong \bigcup_{i \in \N} \cM(A_i^t) \cong \bigcup_{i \in \N} U_i = G. \]
\item $A_G$ is a Frech\'et algebra, thus the definition of the spectrum is intrinsic.
\een
\end{proof}

\begin{exa}
The main examples of good open subsets are given by open subsets of dagger affinoid spaces exhausted by rational subdomains. We can see it with a simple example. In the disk $X = \cM(W_k^1)$, consider
\[ G = \{ x \in X | |x| < \frac{1}{2} \}. \]
$G$ is exhausted by the family of subdomains
\[ U_i = \{ x \in X | |x| \le \frac{1}{2} - \frac{1}{i}  \}, \ \ i \in \N, \ \ i > 1. \]
It is clear that if $k = \C$
\[ \limpro_{i \in \N} A_{U_i} = \cO_X(G) \]
the usual set of analytic functions on the open disk $G$ and the projective limit bornology we are considering on it is the von Neumann bornology of the Frech\'et structure classically considered on $\cO_X(G)$.
\end{exa}


We give now another proof of Tate's acyclicity for rational coverings, in the case $k$ is archimedean, exploiting the isomorphism given in \ref{thm_steinness}.

\begin{thm} \label{thm_arch_tate}
(Tate's acyciclity for $k$ archimedean) \\
Let $\cM(A) = X$ be a dagger affinoid space over $\C$ and consider on it the $G$-topology given by finite coverings by rational subdomains. 
Then, the presheaf $U \mapsto A_U$ is acyclic.
\end{thm}
\begin{proof}
Consider a covering of $X$, $\{ U_1, ..., U_n \}_{1 \le i \le n}$, made by rational $k$-dagger subdomains. We can suppose that each $U_i$ is roomy, otherwise
it has zero measure and is necessarily contained in some other (roomy) element of the covers. 
$X$ can be seen as a pro-Stein space that can be presented with the formula 
\[ X = \bigcap_{\rho > 1} X^\rho \]
where $X^\rho$ are Stein spaces, as a consequence of theorem \ref{thm:stein_filtration}. By theorem \ref{thm_steinness} we know that $A_{U_i} \cong \cO_X(U_i)$, where $\cO_X$ is the canonical sheaf of the pro-analytic site of $X$ given by (one of) its embedding in $\C^n$. So there exists a sequence of open Stein subspaces $U_i^\rho \subset X^\rho$ such that 
\[ U_i = \bigcap_{\rho > 1} U_i^\rho \]
and
\[ \cO_X (U_i) \cong \limind_{\rho > 1} \cO_{X^\rho}(U_i^\rho). \]
Since the $U_i$ cover $X$, we can find a $\rho > 1$ small enough (and suitale $U_i^{\rho'}$) such that for any $\rho' < \rho$ the $U_i^{\rho'}$ cover $X^{\rho'}$. Hence, we get
that the \u{C}ech complex
\[ 0 \to \cO_{X^{\rho'}} \to \prod_{1 \le i \le n} \cO_{U_i^{\rho'}} \to \prod_{1 \le i, j \le n} \cO_{U_i^{\rho'} \cap U_j^{\rho'}} \to \cdots \]
is acyclic for $\rho' < \rho$. Taking the direct limit on the $\rho$ we get the diagram
\[ 0 \to \limind_{\rho > 1} \cO_{X^\rho}(X_\rho) \to \limind_{\rho > 1} \prod_{1 \le i \le n} \cO_{U_i^\rho} (U_i^\rho) \to \limind_{\rho > 1} \prod_{1 \le i, j \le n} \cO_{U_i^\rho \cap U_j^\rho} (U_i^\rho \cap U_j^\rho) \to \cdots \]
which calculate
\begin{equation} \label{eqn:cech}
0 \to \cO_{X}(X) \to \prod_{1 \le i \le n} \cO_X(U^i) \to \prod_{1 \le i, j \le n} \cO_X(U_i \cap U_j) \to \cdots 
\end{equation} 
because the part of diagram with $\rho' < \rho$ is final and the functor $\limind$ commutes with finite limits. By the fact that $\limind$ is an exact functor we deduce that the complex (\ref{eqn:cech}) is acyclic.
\end{proof}

\begin{rmk} \label{rmk_conrad}
The fact that affinoid subdomains of rigid analytic spaces can be seen as the non-archimedean counterpart of the compact Stein subsets of classical complex spaces is a classical observation that goes back at the origin of rigid geometry. The fact that this analogy can also reversed was suggested
by Brian Conrad in \cite{CON}. If $X$ is a separated complex analytic space
we can define a site $X_s$, which we can call the \emph{compact Stein site} of $X$, whose objects are the compact Stein
subsets of $X$ and the coverings are the locally finite coverings. By the properties characterizing the compact Stein subsets of a separated analytic space it is easy to deduce that this site is well-defined and that the \u{C}ech cohomology of coherent sheaves of $X_s$ coincides with the usual sheaf cohomology of $X$.
\end{rmk}

Up to now we saw a lot of similarities between analytic geometry over $\C$ and over non-archimedean base fields. Now we discuss a property that distinguish the archimedean and the non-archimedean base field cases. 
If $k$ is a non-archimedean base field, then the (dagger) Weierstrass subdomains satisfy the transitivity property, \ie if $V \subset U$ is a Weierstrass subdomain and $U \subset X$ is a Weierstrass subdomain, then
$V \subset X$ is a Weierstrass subdomain. The proof of this result is based on the ``rigidity'' of Weierstrass subdomains, \ie on the fact that 
if $X(f_1, ..., f_n)$ and $X(g_1, ..., g_n)$ are Weierstrass subdomains and $|f_i - g_i|_{sup} \le 1$ for any $i$ then $X(f_1, ..., f_n) = X(g_1, ..., g_n)$ (cf. theorem 7.4.2/2 of \cite{BGR}).
This is consequence of the non-archimedean nature of $|\cdot|_{sup}$ which fails if the base field is archimedean, making the ``rigidity'' of Weierstrass subdomains no more true, as next counter-example shows.

\begin{exa} \label{exa:non-transitive}
Consider the functions $z, \frac{3 z}{2} \in W_\C^1$. Then notice that $|z - \frac{3 z}{2}|_{sup} = |- \frac{z}{2}|_{sup} \le 1$ but $X(z) \ne X(\frac{3 z}{2})$ as dagger affinoid subdomains of $\cM(W_\C^1)$. 
We can also produce an explicit counter-example to the transitivity of Weierstrass dagger subdomains in the archimedean case. Let
$f(z) = 4 z - \frac{1}{2}$ and consider the subdomain $X(f) \subset \cM(W_\C^1)$. By theorem \ref{thm_steinness} we know 
that $A_{X(f)} \cong T_\C^1(\frac{1}{4})$ in a non-canonical way, because the center of the disk is placed at $z = \frac{1}{2}$.
Now, the function $\frac{1}{z}$ is analytic and well-defined on $X(f)$ and so the subset
\[ U = \{ x \in X(f) | |\frac{1}{z}| \le \frac{8}{3} \} \subset X(f) \]
is a Weierstrass subdomain of $X(f)$ but it is not a Weierstrass subdomain of $\cM(W_\C^1)$, for which it is a Laurent subdomain.
\end{exa}

Examples as \ref{exa:non-transitive} do not exist for non-archimedean base fields because it is impossible to ``move'' the centres of disks because 
all the $k$-rational points belongings in disks are centres for them.

\index{local ring at a point of an affinoid space}
\begin{defn}  \label{defn:local_rings}
Let $X = \cM(A)$ be a $k$-dagger affinoid space. Then, we define the \emph{ring of germs of analytic functions} at $x \in X$ as
\[ \cO_{X, x} = \limind_{x \in V} A_V, \]
where $V$ varies over the family of $k$-dagger affinoid subdomains of $X$ that contains $x$.
\end{defn}

\begin{rmk}
The previous definition is different from the definition of $\cO_{X, x}$ given by Berkovich at the beginning of section 2.3 of \cite{BER2}, for the fact that we do not ask for $V$ to vary over the set of dagger affinoid neighborhoods (for the topology of the spectrum) of $x$. We can give this different definition, obtaining the right result (\ie $\cO_{X, x}$ coincides with the local ring of $x$ of Berkovich geometry for $k$ non-archimedean when we associate to $X$ its classical affinoid spaces, if $x$ is not in the Shilov boundary) because we are working with dagger affinoid subdomains, which are by definition equipped with germs of analytic functions.
\end{rmk}

\index{open immersion of affinoid spaces}
\begin{defn} \label{defn:open_immersion}
Let $\phi: X \to Y$ be a map of $k$-dagger affinoid spaces, we say that $\phi$ is an \emph{open immersion} if it is injective and the induced map
\[ \phi^*_x: \cO_{Y, \phi(x)} \to \cO_{X, x} \]
is a bijection for any $x \in X$.
\end{defn}

We notice that with respect to this definition, dagger analytic geometry improves some bad features of Berkovich geometry. In fact, if a definition like the one we just gave of open immersions is given in Berkovich geometry then the subdoman embedding $\cM(K_r) \to \cM(T_k^1)$ of example \ref{exa:non_roomy}, is not an open immersion, whereas the dagger version $\cM(K_r^\dagger) \to \cM(W_k^1)$ is. So, the use of dagger analytic geometry allows to restore compatibility of the notion of affinoid subdomain with the notion of open immersion, as in classical rigid geometry. 

\begin{prop} \label{prop:stalk_weie}
Let $X$ be a $k$-dagger affinoid space and $z \in X$. If $k$ is archimedean then $\cO_{X, z}$ is a $k$-dagger affinoid algebra and the open immersion $\{z\} \rhook X$ is a Weierstrass
subdomain embedding.
\end{prop}
\begin{proof}
Consider the case $X = \cM(W_k^n)$ and $k = \C$. Then, given any $z = (z_1, \dots, z_n)$ with $|z_i| \le 1$ we can choose other $2 n$ $n$-tuples of complex numbers, for any $r > 0$ defined
\[ \a_i = (\a_{1, i} = z_1, \dots, \a_{i, i} =  z_i + r, \dots, \a_{n, i} = z_n) \]
\[ \be_i = (\be_{1, i} = z_1, \dots, \be_{i, i} =  z_i - r, \dots, \be_{n, i} = z_n) \]
Then, the $\C$-dagger affinoid subdomain
\[ U = \{ y \in \cM(W_k^n) | |g_1(y)| \le r, \dots, |g_n(y)| \le r, |h_1(y)| \le r, \dots, |h_n(y)| \le r \}  = \{ z \} \]
where
\[ g_i(X_1, \ldots, X_n) = X_i + \a_i, h_i(X_1, \ldots, X_n) = X_i + \be_i. \]
This is the required representation of $\{z\}$ as a Weierstrass subdomain. If $k = \R$ and $z$ is $\R$-rational then we can reason as before. Otherwise,
if $z = (z_1, \dots, z_n)$ is not an $\R$-rational point, then we can consider $\ol{z} = (\ol{z}_1, \dots, \ol{z}_n)$ and the real number
$z' = (z_1 \ol{z}_1, \dots, z_2 \ol{z}_2) = (z_1', \ldots, z_n')$. Defining $\a$ and $\be$ 
\[ \a_i = (\a_{1, i} = z_1', \dots, \a_{i, i} =  z_i' + r, \dots, \a_{n, i} = z_n') \]
\[ \be_i = (\be_{1, i} = z_1', \dots, \be_{i, i} =  z_i' - r, \dots, \be_{n, i} = z_n') \]
the subdomain
\[ U = \{ y \in \cM(W_k^n) | |g_1(y)| \le r, \dots, |g_n(y)| \le r, |h_1(y)| \le r, \dots, |h_n(y)| \le r \}=\{ z \} \]
where
\[ g_i(X_1, \ldots, X_n) = X_i^2 - \Im(z_i), h_i(X_1, \ldots, X_n) = X_i^2 - \Im(z_i) \]
is the required subdomain.

Finally, if $X$ is any $k$-dagger affinoid space then there exists a closed immersion $\phi: X \to W_k^n$ and the pullback of 
$\{ \phi(x) \} \subset W_k^n$ is a subdomain of $X$ of Weierstrass type whose underlying subset is $\{ x\}$, $\forall x \in X$, because $\phi$ is injective and because of proposition \ref{prop:prop_subdomains}.
\end{proof}

\begin{cor}
In the hypothesis of last proposition, for any $x \in X$, one has that $A_{\{x\}} \cong \cO_{X, x}$.
\end{cor}
\begin{proof}
Particular case of theorem \ref{thm_steinness}.
\end{proof}

\begin{cor}
$\cO_{X, x}$ is a local Noetherian ring.
\end{cor}
\begin{proof} 
Well-known both in archimedean and non-archimedean case, cf. \cite{REM} and \cite{BER2}.
\end{proof}

\begin{prop} \label{prop:embed_open}
Let $\phi: U \rhook X = \cM(A)$ be a $k$-dagger affinoid subdomain embedding. Then, $\phi$ is an open immersion.
\end{prop}
\begin{proof} 
	This is an immediate consequence of definition \ref{defn:local_rings}, because we define local rings using the weak G-topology induced by dagger affinoid subdomains (compare our definition with the one given at the beginning of section 2.3 of \cite{BER2}).
\end{proof}

\begin{prop} \label{prop_base_countable_infinity}
Every $k$-dagger affinoid space $\cM(A)$ has a basis for the weak topology formed by open connected subsets countable at infinity.
\end{prop}
\begin{proof} 
We can reason in the same way as \cite{BER2}, proposition 2.2.8. Let $U$ be an open neighborhood of $x \in \cM(A)$ and choose a roomy 
Laurent subdomain $V (f_1, \ldots, f_n, g_1^{-1}, \ldots, g_m^{-1})$ such that $V \subset U$ (we can always do it because they are a base of closed neighborhood of the topology of $\cM(A)$). We choose a monotonic increasing sequence of real numbers $\{\epsilon_i\}$ such that $0 < \epsilon_i < 1$ for all 
$i \in \N$. For any $i \in \N$ define the sequence of subdomains
\[ V_i = X(\epsilon_i f_1, \ldots, \epsilon_i f_n, \epsilon_i^{-1} g_1^{-1}, \ldots, \epsilon_i^{-1} g_m ). \]
We choose for any $V_i$ the connected component $W_i \subset V_i$ which contains $x$. So, $W_i$ is a Laurent roomy connected neighborhood 
of $x$ and each $W_i$ lies in the topological interior of $W_{i + 1}$, then
\[ W = \bigcup_{i = 1}^\infty W_i \]
is an open connected neighborhood which is countable at infinity.
\end{proof}

\section{The Gerritzen-Grauert theorem}

\index{closed immersion of affinoid spaces}
\begin{defn} \label{defn:closed_immersion}
A morphism of $k$-dagger affinoid spaces $f: X = \cM(A) \to Y = \cM(B)$, is called \emph{closed immersion} if $f^*: B \to A$ (the correspondent map of $k$-dagger affinoid algebras) is a surjective homomorphism of algebras.
\end{defn}

\index{locally closed immersion of affinoid spaces}
\begin{defn}
A morphism of $k$-dagger affinoid spaces $f: X \to Y$ is called \emph{locally closed immersion} if the induced homomorphism of local rings $\cO_{Y, f(x)} \to \cO_{X, x}$ is surjective for all
$x \in X$.
\end{defn}

As Temkin noticed in \cite{TEM} and as we already pointed our in last section, the notion of locally closed immersion in Berkovich geometry is not meaningful for non-strict affinoid spaces. This issue forced Temkin to state his version of the Gerritzen-Grauert theorem in terms of monomorphisms of affinoid spaces. Although in our context we can avoid to discuss the relation between locally closed immersions and monomorphisms of (dagger) affinoid spaces, we include such a study for the sake of completeness and for further applications like the ones in \cite{BABE}.

Recall that a morphism $\phi: X \to Y$ in a category $\cC$ is called \emph{monomorphism} if for all $Z \in \ob(\cC)$ the induced maps 
$\Hom(Z, Y) \to \Hom(Z, X)$ are injective. If $\cC$ has fiber products then $\phi: X \to Y$ is a monomorphism if and only if the diagonal morphism $\Delta: X \to X \times_Y X$ is an isomorphism (cf. \cite{EGA1}, 5.3.8). The dual statement holds: if $\cC$ has pushouts then $\phi: X \to Y$ is an epimorphism if and only if the codiagonal morphism $Y \coprod_X Y \to Y$ is an isomorphism.

\begin{lemma}\label{lem:inj2stalk}
Let $A$ be a dagger affinoid algebra. Then, the canonical morphism $A \to \underset{x \in \cM(A)}\prod \cO_{\cM(A), x}$ is injective.
\end{lemma}
\begin{proof} 
If $A$ is a strict dagger affinoid algebra the map $A \to \underset{x \in \cM(A)}\prod \cO_{\cM(A), x}$ is injective by lemma \ref{inj_lemma}. If $A$ is not strict, we can find a valued field extension $K/k$ such that $A \otimes_k^\dagger K$ is a strict dagger affinoid algebra. In the commutative diagram
\[
\begin{tikzpicture}
\matrix(m)[matrix of math nodes,
row sep=2.6em, column sep=2.8em,
text height=1.5ex, text depth=0.25ex]
{ A & \underset{x \in \cM(A)} \prod \cO_{\cM(A), x}  \\
  A \otimes_k^\dagger K & \underset{x \in \cM(A \otimes_k^\dagger K)}\prod \cO_{\cM(A \otimes_k^\dagger K), x} \\};
\path[->,font=\scriptsize]
(m-1-1) edge node[auto] {$$} (m-1-2);
\path[->,font=\scriptsize]
(m-1-1) edge node[auto] {$$} (m-2-1);
\path[->,font=\scriptsize]
(m-2-1) edge node[auto] {$$}  (m-2-2);
\path[->,font=\scriptsize]
(m-1-2) edge node[auto] {$$} (m-2-2);
\end{tikzpicture}
\]
the vertical arrows are injective, the bottom horizontal arrow is injective and hence also the top arrow is injective.
\end{proof}

\begin{lemma}\label{lem:DaggerLocIsEpi}
	If a morphism of dagger affinoid spaces $\phi: X = \cM(B) \to Y = \cM(A)$ is a locally closed immersion. Then, $\phi^*: A \to B$ is an epimorphism of dagger affinoid algebras.
\end{lemma}
\begin{proof} 
Let $\phi: X \to Y$ be a locally closed immersion. Then for each $x \in X$ the induced morphism on stalks $\cO_{Y, \phi(x)} \to \cO_{X, x}$ is surjective. Let $f^*,g^*: C \to B$ be morphisms of dagger affinoid algebras such that $\phi^* \circ f^* = \phi^* \circ g^*$. Let $Z = \cM(C)$ and $f, g: Z \to X$ be the maps of dagger affinoid spaces that correspond to $f^*$ and $g^*$. Then, we have the commutative diagram
\begin{equation}\label{eqn:StalksDiag}
\begin{tikzpicture}
\node (a) at (0,0) {$A$};
\node (b) at (2,0) {$\underset{z \in Z}\prod \cO_{Y, \phi(f(z))}$};
\node (c) at (0,-2) {$B$};
\node (d) at (2,-2) {$\underset{z \in Z}\prod \cO_{X, f(z)}$};
\node (e) at (0,-4) {$C$};
\node (f) at (2,-4) {$\underset{z \in Z}\prod \cO_{Z, z}$};
\path[->,font=\scriptsize]
(a) edge (b)
(a) edge node[auto] {$\phi^*$} (c)
(b) edge node[auto] {$\underset{z \in Z}\prod \phi_z^*$} (d)
(c) edge (d)
(e) edge node[auto] {$\iota$} (f)
([xshift= 2pt]c.south) edge node[right] {$g^{*}$} ([xshift= 2pt]e.north)
([xshift= -2pt]c.south) edge node[left] {$f^{*}$} ([xshift= -2pt]e.north)
([xshift= 2pt]d.south) edge node[right] {$\underset{z \in Z}\prod g^{*}_z$} ([xshift= 2pt]f.north)
([xshift= -2pt]d.south) edge node[left] {$\underset{z \in Z}\prod f^{*}_z$} ([xshift= -2pt]f.north);
\end{tikzpicture}
\end{equation}
where $\underset{z \in Z}\prod \phi_{z}^*$ is surjective. So, from $\phi^* \circ f^* = \phi^* \circ g^*$ we can deduce $(\underset{z \in Z}\prod f_{z}^*) \circ (\underset{z \in Z}\prod \phi_{z}^*) = (\underset{z \in Z}\prod g_{z}^*) \circ (\underset{z \in Z}\prod \phi_{z}^*)$ which implies $\underset{z \in Z}\prod f^{*}_z = \underset{z \in Z}\prod g^{*}_{z}$. Since the bottom square in diagram \ref{eqn:StalksDiag} is commutative it follows that $\iota \circ f^{*} = \iota \circ g^{*}$.  By Lemma \ref{lem:inj2stalk} $\iota: C \to \underset{z \in Z}\prod \cO_{Z, z}$ is injective and so $f^{*} = g^{*}$. 
\end{proof}

To prove the converse of lemma \ref{lem:DaggerLocIsEpi} we need some preparatory lemmas.

\begin{lemma} \label{lemma:ringed_spaces}
	Let $\{ f_i : (X_i, \cO_{X_i}) \to (Y_i, \cO_{Y_i}) \}_{i \in I}$ be a cofiltered projective system of morphisms of locally ringed spaces such that all $f_i$ are injective onto their image and they induces isomorphisms (resp. surjection, resp. injection) on all stalks. Then, $\limpro f_i: \underset{i \in I}\limpro (X_i, \cO_{X_i}) \to \underset{i \in I}\limpro (Y_i, \cO_{Y_i})$ has the same property.
\end{lemma}
\begin{proof} 
	First notice that if $\{(X_i, \cO_{X_i}) \}_{i \in I}$ is such an inverse system of locally ringed spaces and $\pi_i: X \to X_i$ denote the canonical morphisms of the projective limit of the underlying topological spaces, then $\underset{i \in I}\limpro (X_i, \cO_{X_i}) \cong (\underset{i \in I} \limpro X_i, \underset{i \in I} \limind \pi_i^{-1} \cO_{X_i})$, where $\underset{i \in I} \limind \pi_i^{-1}\cO_{X_i}$ is the direct limit of the direct system of sheaves of rings obtained by pulling back to $X$ the structural sheaves of $X_i$ by $\pi_i$ (for a proof of the fact that this space is indeed the projective limit in the category of ringed spaces, and even locally ringed spaces, see theorem 4 of \cite{GIL}). Thus, if we have a morphism $\{ f_i: (X_i, \cO_{X_i}) \to (Y_i, \cO_{Y_i}) \}_{i \in I}$ of two such systems $\{(X_i, \cO_{X_i}) \}_{i \in I}$ and $\{(Y_i, \cO_{Y_i}) \}_{i \in I}$ such that all $f_i$ are injective on the underlying topological space, then $f = \limpro f_i$ is injective at the level of spaces because projective limits of injective maps are injective. The maps on stalks are therefore obtained as a filtered direct limit of homomorphisms, hence $\cO_{Y_i, f(x)} \to \cO_{X_i, x}$ are isomorphisms (resp. surjections, resp. injections) because filtered colimits are functorial (resp. colimits, resp. exact).
\end{proof}

\begin{lemma} \label{lem:Stein_mono}
	Let $f: X \to Y$ be a morphism of Stein spaces over $k$, where $k = \R$ or $k = \C$. Then, $f$ is a monomorphism if and only if it is a locally closed immersion.
\end{lemma}
\begin{proof} 
The if part of the lemma can be proved in the same way of lemma \ref{lem:DaggerLocIsEpi}. To prove the converse consider a morphism $f: X \to Y$ which is not a locally closed immersion. Suppose that there exists a point $y \in Y$ such that $x, x' \in f^{-1}(y)$ with $x \ne x'$. Then, for the two immersions $\iota:\{ x \} \to X$ and $\iota': \{ x' \} \to X$ we have that $f \circ \iota = f \circ \iota'$, thus $f$ is not a monomorphism. Suppose now that $f$ is injective and that there is a $y \in Y$ such that $\cO_{Y, y} \to \cO_{X, f^{-1}(y)}$ is not surjective. 
Choosing a representative for an element of $h \in \cO_{X, f^{-1}(y)}$ which is not in the image of $\cO_{Y, y} \to \cO_{X, f^{-1}(y)}$ one can see that on any (small enough) neighborhood of $f^{-1}(y)$ the map $f$ does not restrict to a monomorphism. In fact, one can find small enough neighborhoods $U_X$ of $f^{-1}(y)$ and $U_Y$ of $y$ for which $h$ can be lifted to an analytic function on $U_X$ such that the identity morphism $\cO_{X}(U_X) \to \cO_{X}(U_X)$ is not distinguishable from the morphism $\cO_{X}(U_X) \to \cO_{X}(U_X)$ which sends $h$ to $0$, by the restriction of $f$ given by $\cO_{Y}(U_Y) \to \cO_{X}(U_X)$, proving that $f$ is not a monomorphism. Thus, monomorphisms of Stein spaces are necessarily locally closed immersions.
\end{proof}

\begin{lemma} \label{lem:non_arch_strict_surj}
	Let $k$ be non-archimedean. Let $f: \cM(B) \to \cM(A)$ be a closed immersion of classical strictly affinoid spaces. Then, for any $x \in \cM(B)$ the morphism $\cO_{\cM(A), f(x)} \to \cO_{\cM(B), x}$ is surjective.
\end{lemma}
\begin{proof} 
Let $\{U_i\}_{i \in I}$ be a base of neighborhoods of $f(x)$ in $\cM(A)$. Each $U_i$ can be chosen to be a strictly affinoid neighborhood of $f(x)$. The family $\{ M(B) \cap U_i\}_{i \in I}$ is therefore a base of neighborhoods of $x$ made of strictly affinoid subdomains (because $f$ is an injective morphism between compact Hausdorff spaces, therefore it is a homeomorphism onto its image). The following diagram
\[
\begin{tikzpicture}
\matrix(m)[matrix of math nodes,
row sep=2.6em, column sep=2.8em,
text height=1.5ex, text depth=0.25ex]
{ A       & B  \\
  A_{U_i} & A_{U_i} \what{\otimes}_A B  \\};
\path[->,font=\scriptsize]
(m-1-1) edge node[auto] {$$} (m-1-2);
\path[->,font=\scriptsize]
(m-1-1) edge node[auto] {$$} (m-2-1);
\path[->,font=\scriptsize]
(m-2-1) edge node[auto] {$$}  (m-2-2);
\path[->,font=\scriptsize]
(m-1-2) edge node[auto] {$$} (m-2-2);
\end{tikzpicture}
\]  
shows that the maps $A_{U_i} \to A_{U_i \cap \cM(B)}$ are surjective, because $A_{U_i \cap \cM(B)} \cong A \what{\otimes}_A B$ and because $\what{\otimes}_A$ preserve surjective maps (since it is a left adjoint functor). To conclude, we notice that $\cO_{\cM(A), f(x)} \cong \underset{i \in I} \limind A_{U_i}$ and $\cO_{\cM(B), x} \cong \underset{i \in I} \limind A_{U_i \cap \cM(B)}$, and since $\underset{i \in I} \limind$ is exact, because $I$ is filtered, we get that the map on stalks is surjective.
\end{proof}

\begin{lemma} \label{lem:EpiIsLocClosed}
	Let $A, B \in \bAff_{k}^\dagger$ and let $f: \cM(B) \to \cM(A)$ be a morphism in the category of $k$-dagger affinoid spaces. Assume that $f$ is a monomorphism, then $f$ is a locally closed immersion.
\end{lemma}
\begin{proof} 
If $k$ is non-Archimedean we can write $\cM(A) = \underset{\rho > r_A}\limpro \cM(A_\rho)$, $\cM(B) = \underset{\rho > r_B}\limpro \cM(B_\rho)$ and $f = \limpro f_\rho$ a monomorphism in $(\bAff_k^\dagger)^\circ$, where $A_\rho$ and $B_\rho$ are strictly affinoid algebras and $r_A, r_B$ are polyradii obtained from presentations of $A$ and $B$ as in remark \ref{rmk:non_strict_dagger_affinoids}. The morphisms $f_\rho: \cM(B_\rho) \to \cM(A_\rho)$ can be chosen to be monomorphisms as a consequence of proposition \ref{prop_filtration_dagger} and remark \ref{rmk:non_strict_dagger_affinoids}. Indeed, $f$ by hypothesis is a monomorphism which implies that the morphism of dagger affinoid algebras $B \otimes_A^\dagger B \to B$ is an isomorphism, therefore it is easy to check that $B_\rho \otimes_{A_\rho}^\dagger B_\rho \to B_\rho$ must be isomorphism of strictly affinoid algebras for $\rho$ small enough, which implies that $f_\rho$ are monomorphisms. 
Then, proposition 1.2 of \cite{TEM} yields that $f_\rho$ are locally closed immersions. Therefore, by the Gerritzen-Grauert Theorem for strictly affinoid spaces there exists a finite covering $\{ X_{i, \rho} \}_{i \in I}$ of $\cM(A_\rho)$ such that $f_\rho$ restricted to $f^{-1}(X_{i, \rho})$ is a Runge immersion. Recall that closed immersions are finite morphisms, so their relative interior is the whole domain (cf. corollary 2.5.13 (i) of \cite{BER2}) and the relative interior of affinoid domain embeddings coincides with the topological interior of the underlying topological spaces. Therefore, using lemma \ref{lem:non_arch_strict_surj} and lemma 2.5.8 (iii) of \cite{BER2} we see that for all $x \in \Int(\cM(B_\rho)/\cM(A_\rho))$ the morphisms induced on stalks by $f_\rho$ are surjective. As locally ringed spaces we have the isomorphism
\[ \limpro_{\rho > r_B} \Int(\cM(B_\rho)/\cM(A_\rho)) \cong \cM(B) \]
Therefore, applying lemma \ref{lemma:ringed_spaces} we deduce the claim.

If $k$ is archimedean, we notice that the representations $\cM(A) = \underset{\rho > 1}\limpro \cM(A_\rho)$, $\cM(B) = \underset{\rho > 1}\limpro \cM(B_\rho)$ now can be chosen with $\cM(A_\rho)$ and $\cM(B_\rho)$ Stein spaces by applying theorem \ref{thm:stein_filtration}. Also in this case one can easily check that $f_\rho$ are monomorphisms using the same argument used in the non-archimedean case. Since the underlying locally ringed spaces of $\cM(A)$ and $\cM(B)$ coincide with the projective limits $\underset{\rho > r_A}\limpro \cM(A_\rho)$ and $\underset{\rho > r_B}\limpro \cM(B_\rho)$ of the underlying locally ringed spaces, it is enough to show that monomorphisms of Stein spaces are locally closed immersions and then again the lemma follows directly applying lemma \ref{lemma:ringed_spaces}. The fact that monomorphisms of Stein spaces are locally closed immersions is proved in lemma \ref{lem:Stein_mono}.
\end{proof}

Putting together lemma \ref{lem:EpiIsLocClosed} and lemma \ref{lem:DaggerLocIsEpi} we have the following direct corollary.

\begin{cor} \label{cor:mono_loc_closed}
	Let $A, B \in \bAff_{k}^\dagger$ and let $f: \cM(B) \to \cM(A)$ be a morphism in the category of $k$-dagger affinoid spaces. $f$ is a monomorphism if and only if $f$ is a locally closed immersion.
\end{cor}
\begin{proof} 
\end{proof}

\index{Runge immersion}
\begin{defn}
A morphism of $k$-dagger affinoid spaces $f: X \to Y$, is called \emph{Runge immersion} if it factors in a diagram
\[
\begin{tikzpicture}
\matrix(m)[matrix of math nodes,
row sep=2.6em, column sep=2.8em,
text height=1.5ex, text depth=0.25ex]
{X & & Y  \\
& Y' \\};
\path[->,font=\scriptsize]
(m-1-1) edge node[auto] {$f$} (m-1-3);
\path[->,font=\scriptsize]
(m-1-1) edge node[auto] {$g$} (m-2-2);
\path[->,font=\scriptsize]
(m-2-2) edge node[auto] {$h$}  (m-1-3);
\end{tikzpicture}
\]
where $g$ is a closed immersion and $h: Y' \to Y$ is a Weierstrass domain embedding. The map of $k$-dagger affinoid algebras that corresponds to $f$ is called \emph{Runge localization}.
\end{defn}

We will need the following two lemmas about Runge localizations.

\begin{lemma} \label{lemma:runge1}
Let $f: A \to B$ be a morphism of dagger affinoid algebras that can be factored as $h \circ g$ where $h$ is a Weierstrass localization and $g$ is a closed immersion, then $f$ is a Runge localization.
\end{lemma}
\begin{proof}
We notice that if $k$ is non-archimedean we can reason like proposition 7.4.3/2 of \cite{BGR} and characterize Runge localizations as morphisms of $k$-dagger affinoid algebras with dense set-theoretic image. But this argument does not work for archimedean base fields, therefore we give another argument.

consider morphisms $A \stackrel{g}{\to} C \stackrel{h}{\to} B$ as in the statement. Then $C \cong \frac{A}{(g_1, \ldots, g_n)}$ for some $g_1, \ldots, g_n \in A$ and $B = C \lt h_1, \ldots, h_m \gt^\dagger$ for some $h_1, \ldots, h_m \in C$. Consider any representative $\tilde{h}_i$ of $h_i$ in $A$. The algebra
\[ B' = \frac{A \lt \tilde{h}_1, \ldots, \tilde{h}_m \gt^\dagger}{(g_1, \ldots, g_n)} \] 
is isomorphic to $B$ because it can be characterized by the same universal property. By definition $A \to B'$ is a Runge localization.
\end{proof}

\begin{lemma} \label{lemma:runge2}
Let $f: A \to B$ be a Runge localization and let $h \circ g$ be a factorization of $f$ such that $g$ is a Weierstrass localization. Then, $h$ is a Runge localization.
\end{lemma}
\begin{proof}
We write $g: A \to C$ and $h: C \to B$. Since $f$ is a Runge localization it admits a factorization $f = h' \circ g'$ where $g': A \to D$ is a Weierstrass localization and $h': D \to B$ is a closed immersion. Consider the diagram
\[
\begin{tikzpicture}
\matrix(m)[matrix of math nodes,
row sep=2.6em, column sep=2.8em,
text height=1.5ex, text depth=0.25ex]
{ A & C  \\
  D & C \otimes_A^\dagger D \\};
\path[->,font=\scriptsize]
(m-1-1) edge node[auto] {$$} (m-1-2);
\path[->,font=\scriptsize]
(m-1-1) edge node[auto] {$$} (m-2-1);
\path[->,font=\scriptsize]
(m-2-1) edge node[auto] {$$}  (m-2-2);
\path[->,font=\scriptsize]
(m-1-2) edge node[auto] {$$}  (m-2-2);
\end{tikzpicture}.
\]
The canonical maps $C \to C \otimes_A^\dagger D$ and $D \to C \otimes_A^\dagger D$ are Weierstrass localizations by the stability of Weierstrass localizations by base change (cf. proposition \ref{prop:prop_subdomains}). Therefore, by the universal property characterizing $C \otimes_A^\dagger D$, we get a commutative diagram
\[
\begin{tikzpicture}
\matrix(m)[matrix of math nodes,
row sep=2.6em, column sep=2.8em,
text height=1.5ex, text depth=0.25ex]
{ C &   \\
  C \otimes_A^\dagger D & B \\
  D &   \\};
\path[->,font=\scriptsize]
(m-1-1) edge node[auto] {$$} (m-2-1);
\path[->,font=\scriptsize]
(m-1-1) edge node[auto] {$h$} (m-2-2);
\path[->,font=\scriptsize]
(m-2-1) edge node[auto] {$$}  (m-2-2);
\path[->,font=\scriptsize]
(m-3-1) edge node[auto] {$$}  (m-2-1);
\path[->,font=\scriptsize]
(m-3-1) edge node[below] {$h'$}  (m-2-2);
\end{tikzpicture}.
\]
The fact that $h'$ is surjective implies that $C \otimes_A^\dagger D \to B$ is surjective. Therefore, we have factored $h$ as a Weierstrass localization followed by a surjective map proving that it is a Runge localization.
\end{proof}

We can state our generalization of the Gerritzen-Grauert theorem.

\index{Gerritzen-Grauert theorem}
\begin{thm}
Let $\phi: X \to Y$ be a locally closed immersion of $k$-dagger affinoid spaces. Then, there exists a finite covering of $\phi(X)$ by rational subdomains $Y_i$ of $Y$ such that all morphisms $\phi_i: \phi^{-1}(Y_i) \to Y_i$ are Runge immersions.
\end{thm}

We divide the proof in two separated case. The case when the base field is archimedean and the one when it is non-archimedean.

Non-archimedean case:

\begin{proof}

One possible strategy to prove the non-achimedean case of the theorem is to write a dagger version of the arguments used by Temkin in \cite{TEM}. But it is more easy to deduce the dagger version of the Gerritzen-Grauert theorem relying on the classical one. Hence, we use the second approach.

Let $\phi: X = \cM(A) \to Y = \cM(B)$ be a locally closed immersion of dagger affinoid spaces. Consider, as in the first part of the proof of lemma \ref{lem:EpiIsLocClosed} a presentation $\cM(A) = \underset{\rho > r_A}\limpro \cM(A_\rho)$, $\cM(B) = \underset{\rho > r_B}\limpro \cM(B_\rho)$ and $\phi = \limpro \phi_\rho$ such that $A_\rho$ and $B_\rho$ are strictly affinoid algebras and $r_A, r_B$ are polyradii and the morphisms $f_\rho: \cM(B_\rho) \to \cM(A_\rho)$ can be chosen to be monomorphisms for all $\rho$ small enough (see ibid. for an explanation of why such a choice is always possible). We will also suppose that $r_A = r_B$, because in the case $r_A \ne r_B$ the proof needs only an easy adaptation of the indexes. 

Therefore, fixing a $\rho$ small enough to have that $\phi_\rho$ is a monomorphism, we can apply the classical Gerritzen-Grauert theorem to get a finite covering $\{ Y_{\rho, i} = \cM(B_{\rho, i}') \}_{i \in I}$ of $Y_\rho$ such that 
\[ \phi_{\rho, i}: \phi_\rho^{-1}(Y_{\rho, i}) \to Y_{\rho, i} \]
are Runge immersions of classical affinoid spaces. Thus, for such a fixed $\rho$ and any $i \in I$ we have a diagram of the form
\begin{equation} \label{eqn:GG}
\begin{tikzpicture}
\matrix(m)[matrix of math nodes,
row sep=2.6em, column sep=2.8em,
text height=1.5ex, text depth=0.25ex]
{ B_\rho & A_\rho \widehat{\otimes}_{B_\rho} B_{\rho,i}'  \\
  B_{\rho,i}' & C_{\rho,i}  \\};
\path[->,font=\scriptsize]
(m-1-1) edge node[auto] {} (m-1-2);
\path[->,font=\scriptsize]
(m-1-1) edge node[auto] {} (m-2-1);
\path[->,font=\scriptsize]
(m-2-2) edge node[auto] {} (m-1-2);
\path[->,font=\scriptsize]
(m-2-1) edge node[auto] {} (m-2-2);
\end{tikzpicture}
\end{equation}
where the left vertical arrow corresponds to a rational subdomain embedding, the bottom horizontal to a Weierstrass embedding and the right vertical arrow to a closed immersion. We fix presentations $T_k^{n_B}(\rho) \to B_\rho$, $T_k^{n_{B',i}}(\rho) \to B_{\rho,i}'$ and $T_k^{n_{C,i}}(\rho) \to C_{\rho,i}$. Of course we can always put $n_B = n_{B',i} = n_{C,i} = n$ by taking $n = \max \{ n_B,  n_{B',i}, n_{C,i} \}$.
For any $r_B < \rho' < \rho$ we can calculate the complete tensor product along the morphism $T_k^n(\rho) \to T_k^n(\rho')$ obtaining a diagram of strictly affinoid algebras
\[
\begin{tikzpicture}
\matrix(m)[matrix of math nodes,
row sep=2.6em, column sep=2.8em,
text height=1.5ex, text depth=0.25ex]
{ B_{\rho'} & A_{\rho'} \widehat{\otimes}_{B_{\rho'}} B_{\rho',i}'  \\
	B_{\rho', i}' & C_{\rho',i}  \\};
\path[->,font=\scriptsize]
(m-1-1) edge node[auto] {} (m-1-2);
\path[->,font=\scriptsize]
(m-1-1) edge node[auto] {} (m-2-1);
\path[->,font=\scriptsize]
(m-2-2) edge node[auto] {} (m-1-2);
\path[->,font=\scriptsize]
(m-2-1) edge node[auto] {} (m-2-2);
\end{tikzpicture}
\]
whenever the coordinates of $\rho'$ lie in $\sqrt{|k^\times|}$ (where $B_{\rho', i}' = B_{\rho, i}' \widehat{\otimes}_{T_k^n(\rho)} T_k^n(\rho')$ and $C_{\rho', i} = C_{\rho, i} \widehat{\otimes}_{T_k^n(\rho)} T_k^n(\rho')$). It is straightforward to check that $B_{\rho'} \to B_{\rho', i}'$ is a rational subdomain localization, $B_{\rho', i}' \to C_{\rho',i}$ is a Weierstrass subdomain localization and that $C_{\rho',i} \to A_{\rho'} \widehat{\otimes}_{B_{\rho'}} B_{\rho',i}'$ is surjecutve. Indeed, since $T_k^n(\rho) \to T_k^n(\rho')$ is a Weierstrass subdomain localization then applying (3) of proposition \ref{prop:prop_subdomains} we obtain that the morphisms $B_\rho \to B_{\rho'}$, $B_{\rho,i}' \to B_{\rho',i}'$, $C_{\rho,i} \to C_{\rho',i}$ and $A_\rho \widehat{\otimes}_{B_\rho} B_{\rho,i}' \to A_{\rho'} \widehat{\otimes}_{B_{\rho'}} B_{\rho',i}'$ are Weierstrass subdomain localization. Therefore, the commutative diagram
\[
\begin{tikzpicture}
\matrix(m)[matrix of math nodes,
row sep=2.6em, column sep=2.8em,
text height=1.5ex, text depth=0.25ex]
{ B_\rho & B_{\rho'}  \\
  B_{\rho,i}' & B_{\rho',i}'  \\};
\path[->,font=\scriptsize]
(m-1-1) edge node[auto] {} (m-1-2);
\path[->,font=\scriptsize]
(m-1-1) edge node[auto] {} (m-2-1);
\path[->,font=\scriptsize]
(m-1-2) edge node[auto] {} (m-2-2);
\path[->,font=\scriptsize]
(m-2-1) edge node[auto] {} (m-2-2);
\end{tikzpicture}
\]
directly yields that $B_{\rho'} \to B_{\rho',i}'$ is rational subdomain embedding (since that other three maps in the diagram are and we can apply proposition \ref{prop:epi_weie}) and the commutative diagram
\[
\begin{tikzpicture}
\matrix(m)[matrix of math nodes,
row sep=2.6em, column sep=2.8em,
text height=1.5ex, text depth=0.25ex]
{ B_{\rho, i}' & B_{\rho', i}'  \\
  C_{\rho, i} & C_{\rho', i}  \\};
\path[->,font=\scriptsize]
(m-1-1) edge node[auto] {} (m-1-2);
\path[->,font=\scriptsize]
(m-1-1) edge node[auto] {} (m-2-1);
\path[->,font=\scriptsize]
(m-1-2) edge node[auto] {} (m-2-2);
\path[->,font=\scriptsize]
(m-2-1) edge node[auto] {} (m-2-2);
\end{tikzpicture}
\]
yields that $B_{\rho', i} \to C_{\rho'}'$ is a Weierstrass subdomain
localization for the same reason (again applying proposition \ref{prop:epi_weie}). Whereas, to prove that $C_{\rho', i}' \to A_{\rho'} \widehat{\otimes}_{B_{\rho'}} B_{\rho',i}'$ is a surjective one has to notice that the diagram 
\[
\begin{tikzpicture}
\matrix(m)[matrix of math nodes,
row sep=2.6em, column sep=2.8em,
text height=1.5ex, text depth=0.25ex]
{C_{\rho,i} & C_{\rho',i}  \\
 A_{\rho} \widehat{\otimes}_{B_\rho} B_{\rho,i}' & A_{\rho'} \widehat{\otimes}_{B_{\rho'}} B_{\rho',i}'  \\};
\path[->,font=\scriptsize]
(m-1-1) edge node[auto] {} (m-1-2);
\path[->,font=\scriptsize]
(m-1-1) edge node[auto] {} (m-2-1);
\path[->,font=\scriptsize]
(m-1-2) edge node[auto] {} (m-2-2);
\path[->,font=\scriptsize]
(m-2-1) edge node[auto] {} (m-2-2);
\end{tikzpicture}
\]
is cocartesian, in the category of complete bornological algebras, because (writing $D_\rho = A_{\rho} \widehat{\otimes}_{B_\rho} B_{\rho,i}'$ and $D_{\rho'} = A_{\rho'} \widehat{\otimes}_{B_{\rho'}} B_{\rho',i}'$)
\[ D_{\rho'} = D_\rho \widehat{\otimes}_{T_k^n(\rho)} T_k^n(\rho') \cong (D_\rho \widehat{\otimes}_{C_{\rho, i}} C_{\rho, i} ) \widehat{\otimes}_{T_k^n(\rho)} T_k^n(\rho') \cong \]
\[ \cong D_\rho \widehat{\otimes}_{C_{\rho, i}} (C_{\rho, i}  \widehat{\otimes}_{T_k^n(\rho)} T_k^n(\rho')) \cong D_\rho \widehat{\otimes}_{C_{\rho, i}} C_{\rho', i}  \]
and therefore $C_{\rho', i} \to D_{\rho'}$ is surjective because complete tensor products preserve the surjectivity of algebra morphisms. 

Hence, we constructed dagger affinoid algebras $B_i = \underset{\rho > r_B} \limind B_{\rho, i}'$, $C_i = \underset{\rho > r_B} \limind C_{\rho, i}$ such that $\{ \cM(B_i) \}_{i \in I}$ is a rational covering of $Y$ and the morphisms
\[
\begin{tikzpicture}
\matrix(m)[matrix of math nodes,
row sep=2.6em, column sep=2.8em,
text height=1.5ex, text depth=0.25ex]
{ B & A \otimes_B^\dagger B_i \\
  B_i & C_i  \\};
\path[->,font=\scriptsize]
(m-1-1) edge node[auto] {} (m-1-2);
\path[->,font=\scriptsize]
(m-1-1) edge node[auto] {} (m-2-1);
\path[->,font=\scriptsize]
(m-2-2) edge node[auto] {} (m-1-2);
\path[->,font=\scriptsize]
(m-2-1) edge node[auto] {} (m-2-2);
\end{tikzpicture}
\]
are suitable to prove the theorem, \ie $B_i \to C_i$ is a dagger Weierstrass subdomain localization and $C_i \to A \otimes_B^\dagger B_i$ is surjective, because all properties are clearly stable by countable filtered colimits. The theorem is therefore proved. 
\end{proof}

We now settle the archimedean case of the theorem (the idea of the proof is inspired by the proof of the Gerritzen-Grauert theorem for affinoid spaces of  Temkin in \cite{TEM}).
For this case, we need some further lemmata. Notice that now $k = \R, \C$ and so only strict dagger affinoid algebras are considered.

\begin{lemma} \label{lemma:surjective_1}
	Let $\phi: A \to B$ be a morphism of $k$-dagger affinoid algebras then $\phi$ is surjective if and only if
	there exists a system of dagger affinoid generators $h_1, ..., h_t \in B^\ov$ and some generators $l_1, ..., l_t \in A^\ov$ such that $\phi(l_i) = h_i$
	(the $\{ l_i \}_{1 \le i \le t}$ might not be a full system of generators of $A$).
\end{lemma}
\begin{proof}
	If $B = k \lt h_1, ..., h_t \gt^\dagger$ and there exists $l_1, ..., l_t \in A^\ov$ such that $\phi(l_i) = h_i$ then we can compose the canonical 
	morphism $\pi: W_k^n \to A$ for which $\pi(X_1) = l_1, \ldots, \pi_t(X_t) = l_t$ and $n \ge t$ with $\phi$ obtaining a map $W_k^n \to B$, which is surjective by the definition of generators. Hence, $\phi$ is surjective because $\pi$ and $\phi \circ \pi$ are.
	
	On the other hand, let $\phi: A \to B$ be surjective and let $\psi: W_k^n \to A$ be the canonical map defining $A$. The composition
	$\phi \circ \psi: W_k^n = k \lt X_1, ..., X_n \gt^\dagger \to B$ is surjective and so $(\phi \circ \psi)(X_i)$ is a system of
	affinoid generators of $B$ and $\psi(X_i) \in A^\ov$ and the image of $\psi(X_i)$ is a system of affinoid generators of $B$.
\end{proof}

\begin{lemma} \label{lemma_surjectivity}
	Let $\phi: A \to B$ be a surjective morphism of $k$-dagger affinoid algebras. We apply proposition \ref{prop_filtration_dagger} to write 
	\[ A \cong \limind_{1 < \rho} \frac{T_k^n(\rho)}{\ol{(f_1, ..., f_r)_\rho}} = \limind_{1 < \rho} A_\rho \]
	\[ B \cong \limind_{1 < \rho} \frac{T_k^m(\rho)}{\ol{(g_1, ..., g_s)_\rho}} = \limind_{1 < \rho} B_\rho, \]
	for some $f_1, \ldots, f_r \in A$ and $g_1, \ldots, g_s \in B$.
	Then, we can find $\rho > 1$ such that for any $1 < \rho' < \rho$ the induced maps
	$\phi_{\rho'}: A_{\rho'} \to B_{\rho'}$ are surjective.
\end{lemma}
\begin{proof}
	We can write 
	\[ B = k \lt h_1, ..., h_t \gt^\dagger \]
	for a finite system of affinoid generators of $B$. By lemma \ref{lemma:surjective_1} the surjectivity of the map $\phi: A \to B$ implies that there is a set of element $l_1, ..., l_t \in A^\ov$
	such that $\phi(l_i) = h_i$. So, since $\{ l_i \}_{1 \le i \le t}$ is a finite number of elements there exists a $\rho'$ such that $l_i \in A_{\rho'}$ for any $i$. Hence, $\phi_\rho: A_\rho \to B_\rho$ is surjective for any $\rho < \rho'$ by the same reasoning of lemma \ref{lemma:surjective_1} applied to the map $\phi_\rho$.
\end{proof}

\begin{lemma} \label{lemma_weie_domain}
Let $A = \underset{\rho > 1}\limind A_\rho$ be a $k$-dagger affinoid algebra and let $\rho' > 1$ be such that $A_{\rho'}$ contains a set of affinoid generators of $A$, $f_1, \ldots, f_n \in A_{\rho'}$. Then, the algebra
\[ A' = \limind_{\rho > \rho'} A_\rho \]
is a $k$-dagger affinoid algebra and the canonical map $A' \to A$ is a dagger affinoid subdomain localization, \ie $\cM(A) \rhook \cM(A')$ is a subdomain. Moreover
\[ A \cong A' \lt f_1, ..., f_n \gt^\dagger. \]
\end{lemma}
\begin{proof}
Let $\pi: W_k^n \to A$ be a presentation of $A$ such that $X_i \mapsto f_i$ and let $\rho' > 1$ be a polyradius such that $f_i \in A_{\rho'}$ for all $i$. It is clear that $W_k^n(\rho')$ is a dagger affinoid algebra. Then, the map
\[ W_k^n(\rho') \to A' \]
obtained by restricting $\pi$ is clearly a strict epimorphism, because
\[ A' \cong k \lt (\rho'_{1})^{-1} X_1, \ldots, (\rho'_{n})^{-1} X_n \gt^\dagger. \] 
We know that the map
$\cM(A) \rhook \cM(A')$ is injective. It is easy to check that the commutative diagram
\begin{equation} \label{eqn:comm_weie}
\begin{tikzpicture}
\matrix(m)[matrix of math nodes,
row sep=2.6em, column sep=2.8em,
text height=1.5ex, text depth=0.25ex]
{ W_k^n(\rho') & W_k^n  \\
	A'           & A  \\};
\path[->,font=\scriptsize]
(m-1-1) edge node[auto] {} (m-1-2);
\path[->,font=\scriptsize]
(m-1-2) edge node[auto] {} (m-2-2);
\path[->,font=\scriptsize]
(m-1-1) edge node[auto] {} (m-2-1);
\path[->,font=\scriptsize]
(m-2-1) edge node[auto] {} (m-2-2);
\end{tikzpicture}
\end{equation}
is cocartesian, where horizontal maps are canonical injections. The injection $W_k^n(\rho') \rhook W_k^n$ is a subdomain embedding and is easy to check the isomorphism
\[ W_k^n(\rho') \lt X_1, \dots, X_n \gt^\dagger \cong W_k^n, \]
where $X_i \in W_k^n(\rho')$ and $W_k^n(\rho') \lt X_1, \dots, X_n \gt^\dagger$ is the Weierstrass localization of $W_k^n(\rho')$ corresponding to the subset defined by the inequalities $|X_i| \le 1$, for all $i$. Proposition \ref{prop:fiber_product_subdomain} applied to the cocartesian square (\ref{eqn:comm_weie}) implies that $A' \to A$ is an affinoid localization and by the commutativity of the previous diagram we see that for all $i$ the element $f_i \in A'$ is mapped to a generator of $A$ (\ie over itself) so that
\[ A \cong A' \lt f_1, ..., f_n \gt^\dagger. \]
Hence, $\cM(A)$ is a Weierstrass subdomain of $\cM(A')$.
\end{proof}

\begin{lemma} \label{lemma:penultimo}
Let $A = \underset{\rho > 1}\limind A_\rho$ be a $k$-dagger affinoid algebra and let $\phi: A \to B = \underset{\rho > 1}\limind B_\rho$ be a dagger subdomain localization, where we represented $A$ and $B$ as inductive limits of Stein algebras using theorem \ref{thm:stein_filtration}. Let $\rho' > 1$ be such that $B_{\rho'}$ contains a set of affinoid generators of $B$ and $A_{\rho'}$ a set of affinoid generators of $A$.
Then, $\cM(B_\rho)$ is an open Stein subspace of $\cM(A_\rho)$ for all $1 < \rho < \rho'$.
\end{lemma}
\begin{proof}
$\phi$ can be written as a direct limit $\underset{\rho} \limind \phi_\rho = \underset{\rho > 1} \limind A_\rho \to \underset{\rho > 1} \limind B_\rho$, for $\rho$ small enough, where $A_\rho$ and $B_\rho$ are Stein algebras. We claim that for $\rho > 1$ small enough $\phi_\rho$ is an open embedding of Stein spaces. To prove that, consider presentations $A \cong \frac{W_k^n}{I}$ and $B \cong \frac{W_k^m}{J}$ where $I = (f_1, \ldots, f_{r_A})$ and $J = (g_1, \ldots, g_{r_B})$. The data of the dagger affinoid subdomain localization $\phi$ is therefore equivalent to the inclusion
\[ \{ x \in \A_k^m | |X_i(x)| \le 1,  g_1(x) = \ldots = g_{r_B}(x) = 0 \} \subset \{ x \in \A_k^n | |X_i(x)| \le 1, f_1(x) = \ldots = f_{r_A}(x) = 0 \} \]
of the associated dagger affinoid spaces. Since the series $f_1, \ldots, f_{r_A}$ and $g_1, \ldots, g_{r_B}$ are over-convergent it makes sense to consider them in a neighborhood of $\cM(W_k^n)$ and $\cM(W_k^m)$ respectively. Therefore, it makes sense to ask if the inclusion 
\[ \{ x \in \A_k^m | |X_i(x)| \le 1 + \epsilon,  g_1(x) = \ldots = g_{r_B}(x) = 0 \} \subset \]
\[ \subset \{ x \in \A_k^n | |X_i'(x)| \le 1 + \epsilon,  f_1(x) = \ldots = f_{r_A}(x) = 0 \} \]
holds for $\epsilon > 0$ small enough. Indeed, the simple change of variables $X_i \mapsto \la Y_i$ and $X_i' \mapsto \la Y_i'$, with $|\la| \ge 1 + \epsilon$ shows that the last inclusion holds. Using strict or non-strict inequalities does not change the validity of the inclusions. Therefore, for dimension reasons (at each point of the spectra), the inclusions $\phi_\rho: \cM(A_\rho) \to \cM(B_\rho)$ are open embeddings.
\end{proof}

\begin{lemma} \label{lemma:ultimo}
Let $\phi: A = \underset{\rho > 1}\limind A_\rho \to B = \underset{\rho > 1}\limind B_\rho$ be a dagger affinoid subdomain localization and $\rho' > 1$ be such that 
$B_{\rho'}$ contains a set of generators of $B$ and $A_{\rho'}$ contains a set of generators of $A$. Defining
\[ A' = \limind_{\rho > \rho'} A_\rho, B' = \limind_{\rho > \rho'} B_\rho \]
then
\[ A \to A \otimes_{A'}^\dagger B' \]
is a roomy dagger subdomain of $\cM(A)$ containing $\phi^*(\cM(B))$.
\end{lemma}
\begin{proof}
First we notice that, reasoning like in the proof of lemma \ref{lemma:penultimo} we can show that, in the commutative square
\[
\begin{tikzpicture}
\matrix(m)[matrix of math nodes,
row sep=2.6em, column sep=2.8em,
text height=1.5ex, text depth=0.25ex]
{ A' & B'  \\
  A  & B  \\};
\path[->,font=\scriptsize]
(m-1-1) edge node[auto] {} (m-1-2);
\path[->,font=\scriptsize]
(m-1-2) edge node[auto] {} (m-2-2);
\path[->,font=\scriptsize]
(m-1-1) edge node[auto] {} (m-2-1);
\path[->,font=\scriptsize]
(m-2-1) edge node[auto] {} (m-2-2);
\end{tikzpicture}
\]
all the maps are subdomain embeddings. Moreover, by lemma \ref{lemma_weie_domain} the vertical maps are Weierstrass subdomain embeddings.
Now, considering the diagram 
\[
\begin{tikzpicture}
\matrix(m)[matrix of math nodes,
row sep=2.6em, column sep=2.8em,
text height=1.5ex, text depth=0.25ex]
{ A' & B' \\
 A &    \\};
\path[->,font=\scriptsize]
(m-1-1) edge node[auto] {$\phi'$} (m-1-2);
\path[->,font=\scriptsize]
(m-1-1) edge node[auto] {$\psi$} (m-2-1);
\end{tikzpicture}
\]
we can calculate the pushout $A \otimes_{A'}^\dagger B'$. Moreover, by the definition of the pushout there is a unique map 
$A \otimes_{A'}^\dagger B' \to B$ which makes the following diagram commutative
\[
\begin{tikzpicture}
\matrix(m)[matrix of math nodes,
row sep=2.6em, column sep=2.8em,
text height=1.5ex, text depth=0.25ex]
{    & B'    \\
 A & A \otimes_{A'}^\dagger B'   \\
    &                  & B   \\};
\path[->,font=\scriptsize]
(m-1-2) edge node[auto] {} (m-2-2);
\path[->,font=\scriptsize]
(m-2-1) edge node[auto] {} (m-2-2);
\path[->,font=\scriptsize]
(m-2-2) edge node[auto] {} (m-3-3);
\path[->,font=\scriptsize]
(m-2-1) edge node[below] {$\phi$} (m-3-3);
\path[->,font=\scriptsize]
(m-1-2) edge node[auto] {} (m-3-3);
\end{tikzpicture}.
\]
The maps $A \to A \otimes_{A'}^\dagger B'$ and $A \otimes_{A'}^\dagger B' \to B$ are dagger subdomain embeddings because dagger affinoid embeddings are stable by base change. $A \otimes_{A'}^\dagger B'$ is roomy in $A$ because $\cM(A \otimes_{A'}^\dagger B') = \cM(A) \times_{\cM(A')} \cM(B')$, which is indeed an intersection of subsets of $\cM(A')$. So $\cM(A) \cap \cM(B_\rho) \subset \cM(A \otimes_{A'}^\dagger B')$ for any $1 < \rho < \rho'$, and $\cM(B_\rho)$ is a Stein space by lemma \ref{lemma:penultimo}.
\end{proof}

\begin{rmk}
In the previous lemma we used dagger affinoid algebras instead of dagger affinoid spaces. This can be confusing if one wants to grasp the geometrical aspects of the discussion. For the convenience of the reader we put the dual diagram of maps that we discussed so far to help to visualize what is happening
\[
\begin{tikzpicture}
\matrix(m)[matrix of math nodes,
row sep=2.6em, column sep=2.8em,
text height=1.5ex, text depth=0.25ex]
{ \cM(B)  &  \\
          & \cM(A \otimes_{A'}^\dagger B') & \cM(A) \\
          &             \cM(B')           & \cM(A') \\};
\path[->,font=\scriptsize]
(m-2-2) edge node[auto] {$$} (m-2-3);
\path[->,font=\scriptsize]
(m-2-2) edge node[auto] {$$} (m-3-2);
\path[->,font=\scriptsize]
(m-3-2) edge node[auto] {$$} (m-3-3);
\path[->,font=\scriptsize]
(m-2-3) edge node[auto] {$$} (m-3-3);
\path[->,font=\scriptsize]
(m-1-1) edge node[auto] {$$} (m-2-2);
\path[->,font=\scriptsize]
(m-1-1) edge node[auto] {$$} (m-2-3);
\path[->,font=\scriptsize]
(m-1-1) edge node[auto] {$$} (m-3-2);
\end{tikzpicture},
\]
keeping the notation of lemma \ref{lemma:ultimo}, where all maps are dagger affinoid embeddings.
\end{rmk}

So, we can now prove the archimedean case of the Gerritzen-Grauert theorem.

\begin{proof}
Let $\phi: X = \cM(A) \to Y = \cM(B)$ be a locally closed immersion. This means that $\phi$ is injective and for each $x \in X$
there is a surjection
\[ \phi_x^*: \cO_{Y, \phi(x)} \to \cO_{X, x}. \]
By proposition \ref{prop:stalk_weie}, in the archimedean case $\cO_{Y, \phi(x)}$ and $\cO_{X, x}$ are $k$-dagger affinoid algebras and $\{x\}, \{\phi(x)\}$ are dagger affinoid subdomains of Weierstrass type of $X$ and $Y$ respectively. Thus, we can write
\[ \cO_{Y, \phi(x)} = B \lt f_1, ..., f_n \gt^\dagger = D, \ \ \cO_{X, x} = A \lt g_1, ..., g_m \gt^\dagger = C \]
for $f_1, ..., f_n \in B$ and $g_1, ..., g_m \in A$.  We can write presentations
\[ C = A \lt g_1, ..., g_m \gt^\dagger = \limind_{\rho > 1} (A \lt g_1, ..., g_m \gt )_\rho \to B \lt f_1, ..., f_n \gt^\dagger = \limind_{\rho > 1} (B \lt f_1, ..., f_n \gt )_\rho = D \]
and since morphisms of dagger affinoid algebras can always be written as morphisms of inductive systems, there exists a $\rho' > 1$ such that for every $1 < \rho < \rho'$, the morphism $\phi_x^*$ can be written as the inductive limit of the morphisms
\[ \phi_\rho^*: (B \lt f_1, ..., f_n \gt )_{\rho} \to (A \lt g_1, ..., g_m \gt )_{\rho}. \]
Moreover, since $\phi_x^*$ is surjective also the morphisms $\phi_\rho^*$ can be chosen to be surjective, applying lemma \ref{lemma_surjectivity}.

We define $C_{\rho} = (A \lt g_1, ..., g_m \gt )_{\rho}$ and $D_{\rho} = (B \lt f_1, ..., f_n \gt )_{\rho}$ and also
\[ A' = \limind_{\rho > \rho'} A_\rho, B' = \limind_{\rho > \rho'} B_\rho \]
and
\[ C' = \limind_{\rho > \rho'} C_\rho, D' = \limind_{\rho > \rho'} D_\rho \]
for a fixed $\rho' > 1$. These algebras fits into the commutative diagram
\[
\begin{tikzpicture}
\matrix(m)[matrix of math nodes,
row sep=2.6em, column sep=2.8em,
text height=1.5ex, text depth=0.25ex]
{ B' & A' \\
  D' & C' \\};
\path[->,font=\scriptsize]
(m-1-1) edge node[auto] {$$} (m-1-2);
\path[->,font=\scriptsize]
(m-1-1) edge node[auto] {$$} (m-2-1);
\path[->,font=\scriptsize]
(m-2-1) edge node[auto] {$$}  (m-2-2);
\path[->,font=\scriptsize]
(m-1-2) edge node[auto] {$$}  (m-2-2);
\end{tikzpicture}
\]
where the bottom map is surjective. We can define the pushouts $C'' = A \otimes_{A'}^\dagger C'$ and $D'' = B \otimes_{B'}^\dagger D'$. We claim that these data define two roomy subdomains of $X$ and $Y$ containing $x$ and $\phi(x)$, respectively. We will focus on the inclusion $\cM(C'') \rhook X$; the same reasoning can be done for $\cM(D'') \rhook Y$ changing $A$ with $B$ and $C$ with $D$.
In the diagram
\[
\begin{tikzpicture}
\matrix(m)[matrix of math nodes,
row sep=2.6em, column sep=2.8em,
text height=1.5ex, text depth=0.25ex]
{ A' & A \\
  C' & C \\};
\path[->,font=\scriptsize]
(m-1-1) edge node[auto] {$$} (m-1-2);
\path[->,font=\scriptsize]
(m-1-1) edge node[auto] {$$} (m-2-1);
\path[->,font=\scriptsize]
(m-2-1) edge node[auto] {$$}  (m-2-2);
\path[->,font=\scriptsize]
(m-1-2) edge node[auto] {$$}  (m-2-2);
\end{tikzpicture}
\]
all maps are dagger affinoid localizations and it is easy to check that, by its construction also $A' \to C'$ is a Weierstrass subdomain localization. This fact, in combination with lemma \ref{lemma:ultimo}, implies that the inclusions $\cM(D'') \rhook Y$ and $\cM(C'') \rhook X$ define roomy Weirstrass localizations. 
Therefore, we get a commutative diagram
\[
\begin{tikzpicture}
\matrix(m)[matrix of math nodes,
row sep=2.6em, column sep=2.8em,
text height=1.5ex, text depth=0.25ex]
{ D' & C' \\
  D'' & C'' \\};
\path[->,font=\scriptsize]
(m-1-1) edge node[auto] {$$} (m-1-2);
\path[->,font=\scriptsize]
(m-1-1) edge node[auto] {$$} (m-2-1);
\path[->,font=\scriptsize]
(m-2-1) edge node[auto] {$$}  (m-2-2);
\path[->,font=\scriptsize]
(m-1-2) edge node[auto] {$$}  (m-2-2);
\end{tikzpicture}
\]
where the upper horizontal map is a surjective and the vertical maps are Weierstrass localizations (the map $D'' \to C''$ is obtained as the composition of the canonical maps $B \otimes_{B'}^\dagger D' \to A \otimes_{A'}^\dagger D' \to A \otimes_{A'}^\dagger C'$ and the commutativity of the diagram is clear). Applying lemma \ref{lemma:runge1} we deduce that the map $D' \to C''$ is a Runge localization and applying lemma \ref{lemma:runge2} we deduce that $D'' \to C''$ is a Runge localization. It is easy to check that
\[ C'' \cong D'' \otimes_B^\dagger A \]
because both algebras define the same subdomain of $\cM(A)$. Therefore, we get a commutative diagram
\[
\begin{tikzpicture}
\matrix(m)[matrix of math nodes,
row sep=2.6em, column sep=2.8em,
text height=1.5ex, text depth=0.25ex]
{ B & D'' \otimes_B^\dagger A \\
  D'' & \\};
\path[->,font=\scriptsize]
(m-1-1) edge node[auto] {$$} (m-1-2);
\path[->,font=\scriptsize]
(m-1-1) edge node[auto] {$$} (m-2-1);
\path[->,font=\scriptsize]
(m-2-1) edge node[auto] {$$}  (m-1-2);
\end{tikzpicture}
\]
where the vertical map is a Weierstrass localization and the diagonal map is a Runge localization (compare this diagram with the diagram (\ref{eqn:GG})).

Finally, since our choice of $x \in \cM(A) = X$ was arbitrary and $\cM(D'') \subset Y$ is roomy we can always cover $\phi(X)$ made of a finite number of this subdomains $D_1'', \dots, D_r''$ (since $\cM(B) = Y$ is compact) which define roomy subdomains of a finite number of points $\phi(x_1), \ldots, \phi(x_r) \in \phi(X)$. Therefore, we showed that there exists a finite covering $D_1'', \dots, D_n''$ of $\phi(X)$ as required to prove the theorem.
\end{proof}

The main corollaries to the Gerritzen-Grauert theorem are the following ones.

\begin{cor} \label{cor:main_gg}
Let $X$ be a $k$-dagger affinoid space and $U \subset X$ a $k$-dagger affinoid subdomain. Then, there exists a finite number of rational subdomains $U_i \subset X$ such that $\bigcup U_i = U$.
\end{cor}
\begin{proof}
As for the classical case, it is enough to apply the Gerritzen-Grauert theorem to the open immersion $U \rhook X$ which is a particular locally closed immersion.
\end{proof}

\begin{cor}
Let $X$ be a dagger affinoid space. Then the $G$-topology generated by finite coverings made of dagger subdomains and the one generated by finite coverings by rational subdomains coincide.
\end{cor}
\begin{proof}
As corollary \ref{cor:main_gg} shows, every affinoid subdomain can be written as a union of finitely many rational subdomains. Therefore, every dagger affinoid subdomain is an admissible open for the G-topology generated by finite coverings made of rational subdomains.
\end{proof}

\begin{rmk} \label{rmk_gerritzen_weierstrass}
We remark a major difference between the archimedean and the non-archimedean case of the Gerritzen-Grauert theorem. In the proof of the archimedean version we showed that we can factor a locally closed immersion in a Runge immersion followed by a Weierstrass domain embedding, choosing a suitable dagger affinoid covering of the codomain. 
Using the intuition developed in non-archimedean rigid geometry one would deduce from this fact that any dagger affinoid subdomain is a finite union of Weierstrass subdomains. But this is false! Indeed, we saw in the previous section (see the discussion after the remark \ref{rmk_conrad} and example \ref{exa:non-transitive}) that in the archimedean case,
unlike the non-archimedean case, the transitivity property of Weierstrass subdomains fails. So the composition of two Weierstrass embeddings is, in general, a rational embedding. With this difference in mind, we see that in the archimedean case we cannot obtain a stronger statement with respect to corollary \ref{cor:main_gg}, which is the main consequence of the Gerritzen-Grauert theorem. Therefore, the the last two corollaries are the precise generalization of the known ones
in the classical theory of affinoid algebras, valid in the same fashion uniformly in the archimedean case and in the non-archimedean case.
\end{rmk}

\begin{rmk} \label{rmk:gg_end}
In this section we discussed two very different strategies to prove the Gerritzen-Grauert theorem, one for the case when the base field is archimedean and one for the case when base field in archimedean. Both these proofs fail to be extended to encompass all cases at once, at least if one uses only obvious adaptations. Hence, a natural question to ask is if there is a way to uniformize the proof of the Gerritzen-Grauert theorem. We have no answer to this question.
\end{rmk}

We end our study of the dagger affinoid subdomains by showing that the dagger site is strictly contained in compact Stein site, when $k$ is archimedean,
\ie that in general there exist compact Stein subsets of $\cM(W_k^n)$ that are not dagger subdomains. One easy example is the following: by \cite{FRI} (page 124), we see that $\cO_{\cM(W_\R^2)}(K)$ is non-Noetherian if $K \subset \cM(W_\R^2)$ is the holomorphic convex hull of
\[ K' \doteq \l \{ (x, y) \in \R^2 | \exp \l ( - \frac{1}{x^2} \r ) \sin \l ( \frac{1}{x} \r ) \le y \le 1, 0 \le x \le 1 \r \}. \]
$K$ is a compact Stein subset of the two dimensional dagger unitary polydisk whose algebra of holomorphic germs is non-Noetherian,
therefore it cannot be written as a quotient of $W_\R^n$ for any $n$. $\cO_{\cM(W_\R^2)}(K)$ is non-Noetherian because the function 
\[ \exp \l ( - \frac{1}{x^2} \r ) \sin \l ( \frac{1}{x} \r ) \]
oscillates wildly near $(x, y) = (0, 0)$ and the following theorem of Siu.

\begin{thm} (Siu) \\ \label{thm:siu}
Let $K$ be a compact Stein subset of an analytic space $(X, \cO)$. Then $\cO_X(K)$ is Noetherian if and only if $V \cap K$ has finitely many topological components for each complex-analytic subvariety $V$
defined in an open neighbourhood of $K$.
\end{thm}
\begin{proof}
See \cite{SIU}.
\end{proof}

Indeed, $K' \cap \{ x = 0 \}$ has infinitely many connecetd components and $\cO_{\cM(W_\R^2)}(K)$ is therefore non-Noetherian. Siu's characterization of Noetherianity of $\cO_X(K)$ permits to deduce the following consquence.

\begin{prop}
Let $K \subset X$ be a compact Stein sub-space of a dagger affinoid space such that $\cO_X(K)$ is non-Noetherian, then $K$ cannot be written as a finite union of dagger affinoid subdomains of $X$.
\end{prop}
\begin{proof}
Let $U_1, \ldots, U_n$ be dagger affinoid subdomains of $X$ such that 
\[ \bigcup_{1 \le i \le n} U_i = K. \]
Let $K \rhook V$ be an open embedding of $K$, as a pro-analytic space, in a Stein space $V$, \ie an injective map such that $\cO_{V, x} \to \cO_{K, x}$ is an isomorphism of stalks for every $x \in K$. For every $i$ and any closed analytic variety $Y$ of $V$, the intersection $U_i \cap Y$ has finitely many topological connected components by Siu's theorem. Clearly, the number of connected components of $K \cap Y$ is less or equal to the sum of the connected components of $U_i \cap Y$ over $i$. Hence, applying again Siu's theorem we deduce that $\cO_X(K)$ is Noetherian, in contraddiction with our hypothesis.
\end{proof}

Even more interesting is to describe a complex compact Stein space such that its algebra of germs of analytic functions is Noetherian which is not dagger affinoid. 

\index{complex $n$-dimensional ball}
\begin{defn}
We denote
\[ B^n \doteq \{ z = (z_1, \ldots, z_n) \in \C^n | |z| = \sqrt{|z_1|^2 + \ldots + |z_n|^2} \le 1 \} \]
and we call it the \emph{$n$-dimensional ball} in $\C^n$.
\end{defn}

One can show that $B^n$ is a compact Stein subset of $\C^n$ writing it as the intersection of the open balls of radius bigger than $1$. The open balls are known to be Stein, cf. \cite{RUD2}.
It is an easy consequence of Siu's theorem that $\cO_{\C^n}(B^n)$ is Noetherian.

\begin{thm}
$B^n$ is not a dagger affinoid space for $n \ge 2$.
\end{thm}

We give two proofs of this result.

\begin{proof}
If $B^n$ is a dagger affinod space then there exists an $m$ and a closed embedding $\iota: B^n \rhook \cM(W_\C^m)$. It is easy to check that $\partial B^n = \iota^{-1}(\iota(B^n) \cap \partial (\cM(W_\C^m)))$ (where now $\partial$ means topological boundary with respect to $\C^m$). Therefore, $\iota$ restricts to a map from the interior of $B^n$ to the interior of $\cM(W_\C^m)$, \ie from the open unit ball to the open polydisk. This restriction is clearly a proper map, but no such map exists thanks to theorem 15.2.4 of \cite{RUD2}.
\end{proof}

The second proof.

\begin{proof}
Again, if $B^n$ is a dagger affinod space then there exists an $m$ and a closed embedding $\iota: B^n \rhook \cM(W_\C^m)$. This implies that there exists a morphism of algebras
\[ \iota^*: W_\C^m \to \cO_{\C^n}(B^n) \]
which is surjective. We can write
\[ \cO_{\C^n}(B^n) = \limind_{\rho > 1} \cO_{\C^n}(B^n(\rho^-)) \]
where $\cO_{\C^n}(B^n(\rho^-))$ are the Stein algebras of analytic functions on open balls of radius $\rho$. Reasoning like proposition \ref{prop_filtration_dagger} one can show that $\iota^*$ is induced by a map of inductive systems and hence it is induced by (\ie it is the restriction of) a proper map
\[ \iota(\rho): B^n(\rho^-) \rhook \D^m(\rho^-) \]
where $\D^m(\rho^-)$ is the $m$-dimensional open polydisk of radius $\rho > 1$. But again by theorem 15.2.4 of \cite{RUD2} $\iota(\rho)$ cannot exist.
\end{proof}

We end this section by recalling the work of Liu \cite{Liu}. In ibid. the notion of compact Stein space in classical non-archimedean geometry is introduced. Liu performs a very precise analysis of these spaces showing their main properties: The category they form is anti-equivalent to the category of their algebras of analytic functions, their algebras of analytic functions are Noetherian and he also gives an example of a compact Stein space which is not an affinoid space. Therefore, these spaces looks like a generalization of affinoid spaces akin to compact Stein spaces of complex geometry, of which dagger affinoid spaces form a subclass. 
Indeed, these analogies can be made precise but it is not the aim of this work to deal with it. Notice also that Liu is able to show that the algebras of analytic functions on non-archimedean compact Stein spaces are always Noetherian, in contrast with what can happen over $\C$, as we saw so far. But this is only a by-product of the definition he uses.
We refer the reader to \cite{BA2}, where we discuss the details of these issues.

\chapter{Dagger analytic spaces} \label{chp:global}

This last chapter is devoted to the construction of the category of $k$-dagger analytic spaces following the methods that Berkovich developed in \cite{BER4}. 
In the first section we prove some auxiliary results on finite modules on a dagger affinoid algebra and on coherent sheaves over a dagger affinoid space. In particular, we show that the category of finite modules over a dagger affinoid algebra $A$ is equivalent to the category of finite dagger modules, \ie the category of finite $A$-modules endowed with a canonical structure of complete bornological module from the bornology of $A$. Then, we deduce the general version of the Tate's acyclicity theorem, extending the results of chapter 4 to the weak $G$-topology induced by dagger affinoid subdomains and encompassing also to the non-strict case. We end the first section with a proof of the Kiehl's theorem on coherent sheaves adapting the classical proof of rigid geometry.
In the subsequent section we introduce the notion of $k$-dagger analytic space. We use the methods of Berkovich nets, so a $k$-dagger analytic space is a triple $(X, A, \tau)$ where $X$ is a locally Hausdorff topological space, $\tau$ a Berkovich net on $X$ and $A$ an atlas of dagger affinoid algebras for the net $\tau$. We see that we can perform all the basic constructions of \cite{BER4}, underlining the differences of the two approaches and the relations of dagger analytic spaces with the pro-analytic structure on the pro-site of $X$ (as defined in appendix \ref{pro_appendix}). The next section discusses the relations between our dagger analytic spaces and the other definitions of analytic spaces already present in literature. First we deal with the non-archimedean case: We are mainly interested to compare our dagger analytic spaces with the Grosse-Kl\"onne ones and with Berkovich's analytic spaces. Then, we show that in the complex case, the category of classical complex analytic spaces embeds in a fully faithful way in the category of dagger analytic spaces. In the last section we see how, for non-archimedean base fields, the notion of flat morphism in dagger analytic geometry behave better than the respective notions in Berkovich geometry. More precisely, the naive definition of flatness for dagger analytic spaces is stable by base change and is equivalent to the non-naive one given by Ducros in the context of Berkovich geometry in order to fix the non-stability by base change of naive flatness for Berkovich spaces.

\section{Coherent sheaves over a dagger affinoid space}

\begin{notation}
In this section $A$ will always denote a $k$-dagger affinoid algebra. We will also suppose, when needed, to have fixed a representation $A \cong \underset{\rho > r}\limind A_\rho$, as the one given in proposition \ref{prop_filtration_dagger} and theorems \ref{thm:T_filtration} and \ref{thm:stein_filtration}, when needed.
\end{notation}

For any ring $A$ we denote with $\bMod_A$ the category of modules over $A$ with $A$-linear morphisms and with $\bMod^f_A$ the full subcategory of 
$\bMod_A$ made of finite $A$-modules.

\index{finite dagger module}
\begin{defn}
Let $A$ be a $k$-dagger affinoid algebra, we say that a bornological $A$-module is a \emph{finite dagger module} (or simply a \emph{dagger module})
if there exist a strict epimorphism $A^n \to M$.  
We denote with $\bMod_A^\dagger$ the category of finite dagger $A$-modules with bounded $A$-linear morphisms.
\end{defn}

For any $k$-dagger affinoid algebra, the free module $A^\mu$, for any cardinal $\mu$, has a canonical filtration induced by the filtration of $A$ that we denote
\[ A^\mu = \limind_{\rho > r} (A^\mu)_\rho. \]

\index{canonical filtration of a dagger module}
\begin{defn}
Let $A$ be a $k$-dagger affinoid algebra and $M$ an $A$-module. Then, any surjection
\[ \psi: A^\mu \to M \]
induces a filtration on $M$ and we call it the \emph{filtration} on $M$ with respect to $\psi$ and denotes its elements by
\[ M = \limind_{\rho > r} (M_\psi)_\rho \]
where $(M_\psi)_\rho = \psi((A^\mu)_\rho)$.
\end{defn}

In general this filtration is not preserved by bounded $A$-module morphisms.

\begin{prop} \label{prop:filtration}
Let $A$ be a $k$-dagger affinoid algebra, $M$ a finite $A$-module and $\psi: A^n \to M$, $\psi': A^m \to M$ two presentation of $M$. Then, there exists a 
$r' > r$ such that for any $r < \rho < r'$
\[ (M_\psi)_\rho = (M_{\psi'})_\rho \]
and $(M_\psi)_\rho$ are all finitely generated $A_\rho$-modules.
\end{prop}
\begin{proof}
First of all, we notice that the fact that $A \cong \underset{\rho > r}\limind A_\rho$ is Noetherian implies that $M$ is of finite presentation. Hence, we can write the following exact sequences
\[ A^{n'} \stackrel{\phi}{\to} A^n \stackrel{\psi}{\to} M \to 0 \]
\[ A^{m'} \stackrel{\phi'}{\to} A^m \stackrel{\psi'}{\to} M \to 0 \]
and $M \cong \Coker(\phi) \cong \Coker(\phi')$. For any $r < \rho < r'$ we have the inclusion morphisms $(A_\rho)^{n'} \rhook A^{n'}$ and $(A_\rho)^n \rhook A^n$.
We write $\{ f_1, \ldots, f_{n'} \}$ for a base for $A^{n'}$ and $\{ e_1, \ldots, e_n \}$ for a base for $A^n$, in order to write
\[ \phi(f_i) = \sum_{j = 1}^n a_{i,j} e_j \]
for some $a_{i, j} \in A$. Since $\{ a_{i, j} \}$ is a finite subset of $A$, we can find a $\rho > r$ such that $a_{i,j} \in A_{\rho'}$ for any 
$r < \rho' < \rho$. This shows that we can find a commutative diagram
\[
\begin{tikzpicture}
\matrix(m)[matrix of math nodes,
row sep=2.6em, column sep=2.8em,
text height=1.5ex, text depth=0.25ex]
{ A^{n'} & A^n & M & 0 \\
  (A_\rho)^{n'} & (A_\rho)^n & \Coker(\phi_\rho) & 0 \\};
\path[->,font=\scriptsize]
(m-1-1) edge node[auto] {$\phi$} (m-1-2);
\path[->,font=\scriptsize]
(m-1-2) edge node[auto] {$\psi$} (m-1-3);
\path[->,font=\scriptsize]
(m-2-1) edge node[auto] {$\phi_\rho$}  (m-2-2);
\path[->,font=\scriptsize]
(m-2-2) edge node[auto] {$\psi_\rho$}  (m-2-3);
\path[right hook->,font=\scriptsize]
(m-2-1) edge node[auto] {}  (m-1-1);
\path[right hook->,font=\scriptsize]
(m-2-2) edge node[auto] {}  (m-1-2);
\path[right hook->,font=\scriptsize]
(m-2-3) edge node[auto] {}  (m-1-3);
\path[->,font=\scriptsize]
(m-1-3) edge node[auto] {}  (m-1-4);
\path[->,font=\scriptsize]
(m-2-3) edge node[auto] {}  (m-2-4);
\end{tikzpicture}
\]
where $\Coker(\phi_\rho) = (M_\phi)_\rho$ and the map $\phi_\rho$ is simply defined by mapping $f_i \mapsto \underset{j = 1}{\overset{n}\sum} a_{i,j} e_j$. The same reasoning can be applied to the presentation
\[ A^{m'} \stackrel{\phi'}{\to} A^m \stackrel{\psi'}{\to} M \to 0 \]
for a base $\{e_1', \ldots, e_m' \}$ of $A^m$.

Now consider the two systems of generators $\{\psi(e_1), \ldots, \psi(e_n) \}$ and $\{\psi'(e_1'), \ldots, \psi'(e_m') \}$ of $M$. Since they generate
$M$ we can write
\[ \psi(e_i) = \sum_{j = 1}^m b_{i,j} \psi'(e_j') \]
and 
\[ \psi'(e_i') = \sum_{j = 1}^n c_{i,j} \psi(e_j) \]
for some $c_{i, j}, b_{i,j} \in A$. Taking $\rho > r$ small enough to have $c_{i, j}, b_{i,j} \in A_\rho$ for every $i$ and $j$, the $A_\rho$-submodule of $M$ generated by $\{\psi(e_1), \ldots, \psi(e_n) \}$ and by $\{\psi'(e_1'), \ldots, \psi'(e_m') \}$ is the same, hence $(M_\psi)_\rho = (M_{\psi'})_\rho$, which concludes the proof.
\end{proof}

By the previous proposition, when we deal with finite $A$-module we suppress the dependence of the filtration with respect to $\psi$, and therefore we can deduce the following immediate corollary.

\begin{cor} \label{cor:1}
Let $M$ be a finitely generated $A$-module. Then, on $M$ there is a uniquely determined, up to isomorphism, bornology which makes $M$ into a finite dagger $A$-module.
\end{cor}
\begin{proof}
Since $M = \underset{\rho > r} \limind M_\rho$ and $M_\rho$ are $A_\rho$-Banach modules, then $M$ with the direct limit bornology induced by this direct limit is a finite dagger $A$-module. The cofinality result of proposition \ref{prop:filtration} immediately implies that any two such bornologies on $M$ are isomorphic.
\end{proof}

\begin{cor} \label{prop_dagger_modules}
Let $A$ be a $k$-dagger affinoid algebra, then the category of finite $A$-modules is equivalent to the category of finite dagger modules.
\end{cor}
\begin{proof}
It is enough to show that every submodule $N \subset M$ of a finitely generated $A$-module is bornologically closed because then we can apply the bornological closed graph theorem (cf. theorem 3.3 of \cite{BA}) to deduce that every homomorphism between finitely generated $A$-modules is bounded.

Consider an inclusion $N \subset M$ of finite $A$-modules where $M$ is equipped with its canonical dagger bornology described in corollary \ref{cor:1}. It is easy to check that the bornology on $M$ induces a bornology on $N$ which makes $N$ into a finite dagger $A$-module. Therefore, the bornology induced by $M$ on $N$ is equivalent to the quotient bornology induced by a presentation $A^n \to N$. Hence, $N$ is a complete bornological vector space and by lemma 3.14 of \cite{BABEKR} is a closed subspace of $M$ (see ibid. for a discussion of why LB bornological spaces are proper and therefore lemma 3.14 can be applied to $M$).
\end{proof}

\index{$A$-dagger affinoid algebra}
\begin{defn}
Let $A$ be a commutative, unital complete bornological $k$-algebra.  We say that a complete bornological $k$-algebra $B$ is a \emph{$A$-dagger affinoid algebra} if there exists a surjective morphism $A \lt \rho_1 X_1, \ldots, \rho_n X_n \gt^\dagger \to B$ of algebras which is a strict epimorphism of underlying bornological vector spaces.
\end{defn}

\begin{rmk}
When $A$ is a dagger affinoid algebra then also $A \lt \rho_1 X_1, \ldots, \rho_n X_n \gt^\dagger$ is a dagger affinoid algebra and therefore it has a countable base for the bornology. Therefore, Buchwalter's theorem (cf. theorem 4.9 of \cite{BA}) implies that the condition on strict epimorphism of underlying bornological vector spaces is equivalent to bounded surjectivity (provided that $B$ is complete as we specified).
\end{rmk}

\begin{lemma}
Let $A$ be a $k$-dagger affinoid algebra, $B$ an $A$-dagger affinoid algebra and $M, N \in \bMod_A^\dagger$. Then
\ben
\item $M \otimes_A^\dagger N \cong M \otimes_A N \in \bMod_A^\dagger$;
\item $M \otimes_A^\dagger B \cong M \otimes_A B \in \bMod_B^\dagger$;
\item any $A$-linear morphism $M \to N$ is strict.
\een
\end{lemma}
\begin{proof}
All the three assertion are formal consequence of the equivalence showed in corollary \ref{prop_dagger_modules}.
\end{proof}

\begin{lemma} \label{lemma_BGR_7_3_2_3}
Let $X = \cM(A)$ be a strict $k$-dagger affinoid subdomain, $x \in \Max(A) \subset X$ and $\fm_x$ the maximal ideal corresponding to $x$. The canonical
map $A \to \cO_{X, x}$ factors through the localization $A_{\fm_x}$, which embeds into $\cO_{X, x}$. We also get isomorphisms
\begin{equation} \label{eqn:quot_iso}
\frac{A}{\fm_x^n A} \to \frac{A_{\fm_x}}{\fm_x^n A_{\fm_x}} \to \frac{\cO_{X, x}}{\fm_x^n \cO_{X, x}}
\end{equation} 
for any $n \in \N$ and
\[ \what{A} \to \what{A_{\fm_x}} \to \what{\cO_{X, x}} \]
between the $\fm_x$-adic completions.
\end{lemma}
\begin{proof}
Consider the closed immersion $\cM(A/\fm_x^n) \rhook \cM(A)$ given by the quotient map $A \to \frac{A}{\fm_x^n}$. The space $\cM(A/\fm_x^n)$
has only one point and it is a direct consequence of proposition \ref{prop_subdomain_max_ideals} that 
\[ \frac{A}{\fm_x^n A} \cong \frac{\cO_{X, x}}{\fm_x^n \cO_{X, x}}. \]
Since 
\[ \fm_x \cO_{X, x} = \{ h \in \cO_{X, x} | h(x) = 0 \}, \]
then every element of $A$ not in $\fm_x$ is mapped to a unit of $\cO_{X, x}$. This proves that the map $A \to \cO_{X, x}$ factors through $A_{\fm_x}$. The isomorphisms (\ref{eqn:quot_iso}) just showed, easily imply the isomorphisms on $\fm_x$-adic completions. 

It remains to show that the map $A_{\fm_x} \to \cO_{X, x}$ is injective. This follows from Krull's intersection theorem which implies that $\what{A_{\fm_x}}$ is Hausdorff and the maps to the $\fm_x$-adic completion is injective. Indeed, we can then write the commutative diagram
\[
\begin{tikzpicture}
\matrix(m)[matrix of math nodes,
row sep=2.6em, column sep=2.8em,
text height=1.5ex, text depth=0.25ex]
{ A_{\fm_x}        & \cO_{X, x} \\
  \what{A_{\fm_x}} & \what{\cO_{X, x}} \\};
\path[->,font=\scriptsize]
(m-1-1) edge node[auto] {$$} (m-1-2);
\path[->,font=\scriptsize]
(m-1-1) edge node[auto] {$$} (m-2-1);
\path[->,font=\scriptsize]
(m-2-1) edge node[auto] {$\cong$}  (m-2-2);
\path[->,font=\scriptsize]
(m-1-2) edge node[auto] {$$}  (m-2-2);
\end{tikzpicture}
\]
where all maps but the top horizontal are known to be injective. Hence also $A_{\fm_x} \to \cO_{X, x}$ is injective.
\end{proof}

The next lemma is an important generalization of proposition 2.1.2 of \cite{BER2} to the case of bornological vector spaces. In analogy with \cite{BER2}, we use this lemma to reduce proofs for modules over non-strict $k$-dagger affinoid algebras to the case of modules over strict $k$-dagger affinoid algberas.

\begin{lemma} \label{lemma_BER2_2_1_2}
Let $r_1, \ldots, r_n \in \R_{\ge 0}$ with $r_1, \ldots, r_n \notin \sqrt{|k^\times|}$ and $K$ a valued extension of $k$ such that $r_1, \ldots, r_n \in \sqrt{|K^\times|}$. Then, a sequence $M \to N \to P$ of $k$-dagger modules morphisms is exact if and only if the sequence
\begin{equation} \label{eqn:ber_seq}
M \otimes_k^\dagger K \to N \otimes_k^\dagger K \to P \otimes_k^\dagger K.
\end{equation} 
\end{lemma}
\begin{proof}
It is an easy consequence of proposition \ref{prop:filtration} that the data of an exact sequence $M \to N \to P$ is equivalent to the data of a system of exact sequences
\[ M_\rho \to N_\rho \to P_\rho \]
of $k$-Banach spaces for some sufficiently small $\rho > r$. Then, we can apply proposition 2.1.2 of \cite{BER2} to deduce that the exactness of $ M_\rho \to N_\rho \to P_\rho$ is equivalent to the exactness of the sequences
\[ M_\rho \what{\otimes}_k K \to N_\rho \what{\otimes}_k K \to P_\rho \what{\otimes}_k K \]
which, by another application of proposition \ref{prop:filtration}, is equivalent to the exactness of the sequence (\ref{eqn:ber_seq}).
\end{proof}

\begin{rmk}
We underline that, also if not explicit mentioned, the previous lemma applies only for $k$ non-archimedean because the non-strict case does not exist for archimedean base fields.
\end{rmk}

\begin{prop} \label{prop_dagger_flat}
Let $A$ be a $k$-dagger affinoid algebra and $V \subset \cM(A)$ a dagger affinoid subdomain. Then, $A_V$ is a flat $A$-algebra.
\end{prop}
\begin{proof}
Suppose first that $A$ is a strict $k$-dagger affinoid algebra and that $V$ is a strict dagger affinoid subdomain of $\cM(A)$.
Let $x \in \Max(A_V) \subset \cM(A)$, $\fm \subset A$ and $\fm_V \subset A_V$ be the maximal ideals that correspond to $x$. By lemma
\ref{lemma_BGR_7_3_2_3} the map $A \to A_V$ induces a bijection $\what{A} \to \what{A_V}$ between the $\fm$-adic completion of $A$ and the
$\fm_V$-adic completion of $A_V$. It is a well-known result of commutative algebra that this implies that $A_V$ is a flat $A$-algebra, see \cite{BourComm}, Alg\'ebre Commutative Chap. III, \S 5.2.

For the non-strict case, we can reason in the same way of \cite{BER2}, proposition 2.2.4. 
Consider a complete valued field $K/k$ such that both $A \otimes_k^\dagger K$ and $A_V \otimes_k^\dagger K$ are strict dagger affinoid algebras, and an injective morphism $\phi: M \to N$
of finite $A$-modules. We can endow $M$ and $N$ with their canonical $k$-dagger bornology and so the injection $\phi: M \to N$ becomes a strict
monomorphism of bornological modules. By lemma \ref{lemma_BER2_2_1_2} also the morphism $M \otimes_k^\dagger K \to N \otimes_k^\dagger K$ induced by $\phi$ is injective, and it is also necessarily strict because $M \otimes_k^\dagger K$ and $N \otimes_k^\dagger K$ are finite modules over $A \otimes_k^\dagger K$. Then, $A_V \otimes_k^\dagger K$ is a flat $A \otimes_k^\dagger K$ algebra, hence the morphism
\[ \phi': (M \otimes_k^\dagger K) \otimes_{A \otimes_k^\dagger K}^\dagger (A_V \otimes_k^\dagger K) \to 
   (N \otimes_k^\dagger K) \otimes_{A \otimes_k^\dagger K}^\dagger (A_V \otimes_k^\dagger K) \]
is injective. Moreover, we have
\[ (M \otimes_k^\dagger K) \otimes_{A \otimes_k^\dagger K}^\dagger A_V \otimes_k^\dagger K \cong (M \otimes_A A_V) \otimes_k^\dagger K \]
because both objects satisfy the same universal property. The analogous statement holds for $N$. We obtain the exact sequence
 \[ 0 \to (M \otimes_A A_V) \otimes_k^\dagger K \to (N \otimes_A A_V) \otimes_k^\dagger K \]
which by lemma \ref{lemma_BER2_2_1_2} implies the exactness of the sequence
\[ 0 \to M \otimes_A A_V \to N \otimes_A A_V. \]
This implies that $A_V$ is a flat $A$-algebra because any $A$-module can be written as the filtered direct limit of its finitely generated submodules and the filtered direct limit functor of $k$-vector spaces is an exact functor.
\end{proof}

Now we can state the general version of the Tate's acyclicity theorem.

\begin{thm} \label{thm_tate_gen}
Let $A$ be a $k$-dagger affinoid algebra and $M$ an $A$-module. Then, for any finite dagger affinoid covering $\{V_i\}_{i \in I}$ of $\cM(A)$ the \u{C}ech complex
\[ 0 \to M \to \prod_{i \in I} M \otimes_A A_{V_i} \to \prod_{i,j \in I} M \otimes_A A_{V_i \cap V_j} \to \cdots \]
is exact. Moreover, if $M \in \ob(\bMod_A^\dagger)$ then the morphism of the complex are all strict.
\end{thm}
\begin{proof}
For the proof we can use the same reduction to the case of strict dagger affinoid algebras and to the case $M = A$ used in proposition 2.2.5 of  \cite{BER2}, for the classical affinoid case. This is possibile because we have all the required ingredients: we proved that the category of finite $A$-modules is equivalent to the category of finite $A$-dagger modules, that $A \to A_V$ is a flat morphism of rings and we proved the dagger version of Tate's acyclicity theorem for rational coverings and Gerritzen-Grauert theorem in chapters \ref{chap:tate} and \ref{chap:gg}.
\end{proof}

\begin{rmk}
In the case $k = \C$ we saw that the dagger subdomains of $\cM(A)$ are compact Stein subsets of $\cM(A)$ and the acyclicity of coherent
analytic sheaves for finite coverings of compact Stein subsets is a well-known property.
\end{rmk}

We recall the definition of coherent sheaf over a ringed $G$-topological space (or $G$-ringed space).

\index{coherent sheaf}
\begin{defn} \label{defn:coherent}
Let $(X, \cO_X)$ be a $G$-ringed space and $\sF$ an $\cO_X$-module.
\ben
\item $\sF$ is called \emph{of finite type}, if there exists an admissible covering $\{X_i\}_{i \in I}$ of $X$ together with exact sequences of type
      \[ \cO_X^{s_i}|_{X_i} \to \sF|_{X_i} \to 0, i \in I, \]
      with $s_i \in \N$.
\item $\sF$ is called \emph{of finite presentation}, if there exists an admissible covering $\{X_i\}_{i \in I}$ of $X$ together with exact 
      sequences of type
      \[ \cO_X^{r_i}|_{X_i} \to O_X^{s_i}|_{X_i} \to \sF|_{X_i} \to 0, i \in I, \]
      with $s_i, r_i \in \N$.
\item $\sF$ is called \emph{coherent} if $\sF$ is of finite type and if for any admissible open
      subspace $U \subset X$ the kernel of any morphism $\cO_X^s |_U \to \sF|_U$ is of finite type.
\een
\end{defn}

We now consider on dagger affinoid spaces $X$ the weak $G$-topology given by finite coverings by dagger affinoid subdomains.
The category of coherent $\cO_X$-modules for this $G$-topology is denoted by $\bCoh(X)$.

Consider a $k$-dagger affinoid space $X = \cM(A)$ and a $A$-module $M$. We can look at the functor $\sF$ from the category of $k$-dagger affinoid subdomains in $X$ 
which associates to any affinoid subdomain $\cM(A') \subset X$ the tensor product $M \otimes_A A'$. $\sF$ is a presheaf on $X$ with respect 
to the weak $G$-topology, and this presheaf is, in fact, a sheaf by Tate's acyclicity theorem, \ref{thm_tate_gen}. 

\begin{defn}
Let $M \in \bMod_A$. We call the $\cO_X$-module given by
\[ \sF(U) = M \otimes_A A' \]
for any dagger affinoid subdomain $U = \cM(A') \subset \cM(A)$, the \emph{$\cO_X$-module associated} to the $A$-module $M$, and denote it by $\sF = M \otimes_A \cO_X$.
\end{defn}

\begin{prop} \label{prop_coherent}
Let $X = \cM(A)$ be a $k$-dagger affinoid space.
\ben
\item The functor
      \[ \cdot \otimes_A \cO_X : M \mapsto M \otimes_A \cO_X \]
      from $A$-modules to $\cO_X$-modules is fully faithful:
\item The functor $\cdot \otimes_A \cO_X$ commutes with the formation of kernels, images, cokernels, and tensor products.
\item A sequence of $A$-modules $0 \to M' \to M \to M'' \to 0$ is exact if and only if the associated sequence
      \[ 0 \to M' \otimes_A \cO_X \to M \otimes_A \cO_X \to M'' \otimes_A \cO_X \to 0 \]
      is exact.
\een
\end{prop}
\begin{proof}
\ben
\item The canonical map
      \[ \Hom_A(M,M') \to \Hom_{\cO_X}( M \otimes_A \cO_X, M' \otimes_A \cO_X) \]
      is bijective, since any $\cO_X$-morphism $M \otimes_A \cO_X \to M' \otimes_A \cO_X$ is uniquely
      determined by the $A$-morphism between $M = M \otimes_A \cO_X(X)$ and
      $M' = M' \otimes_A \cO_X(X)$. 
\item By its definition the functor commutes with tensor products. If
      \[ 0 \to M' \to M \to M'' \to 0 \]
      is an exact sequence of $A$-modules, the induced sequence
      \[ 0 \to M' \otimes_A A' \to M \otimes_A A' \to M''\otimes_A A' \to 0 \]
      is exact for any affinoid subdomain $\cM(A') \subset X$, since the corresponding map $A \to A'$ is flat, by lemma \ref{prop_dagger_flat}. From this one easily concludes that the functor of
      taking associated $\cO_X$-modules is exact.
\item We already proved the exactness of the functor $\cdot \otimes_A \cO_X : M \mapsto M \otimes_A \cO_X$, so it remain only to prove the only if part of last claim. It is enough to notice that an $A$-module $M$ is trivial if and only if $M \otimes_A \cO_X$ is trivial.
\een
\end{proof}

\begin{defn}
Let $(X, \cO_X)$ be a $G$-ringed space. We say that a sheaf $\sF$ of $\cO_X$-modules, is \emph{$\sU$-coherent}, for an admissible covering $\sU$,
if $\sF|_{X_i}$ is associated to a finite $\cO_{X_i}(X_i)$-module for all $i \in I$.
\end{defn}

\begin{lemma} \label{lemma:U_coh}
Let $X = \cM(A)$ be a $k$-dagger affinoid space and $\sF$ a sheaf of $O_X$-modules over $X$. Then, $\sF$ is coherent if and only if there exists an admissible (hence finite) dagger affinoid convering $\{U_i\}_{i \in I}$ of $X$ such that $\sF|_{U_i}$ is associated to a finite $\cO_{X}(U_i)$-modules for each $i \in I$.
\end{lemma}
\begin{proof}
This is a simple restatement of the condition of definition \ref{defn:coherent} in the particular case of dagger affinoid spaces equipped with the weak G-topology.
\end{proof}

The following proof is an adaptation of the classical argument for proving Kiehl's theorem, as one can find, for example, in \cite{BOSCH}, theorem 1.14.4.

\begin{thm} \label{thm:kiehl}(Kiehl's theorem) \\
Let $A$ be a strict $k$-dagger affinoid algebra and $X = \cM(A)$. Then, a $\cO_X$-modules is coherent if and only if it is associated to a finite $A$-dagger module.
\end{thm}
\begin{proof} 
If $\sF$ is associated to a finite $A$-module, then it is coherent as a trivial consequence of lemma \ref{lemma:U_coh}. The other implication is the hard part of the theorem, and we proceed by steps following the classical proof done in rigid geometry.
\begin{itemize}
\item First step: By lemma \ref{lemma:U_coh} we can reduce
      the problem of checking coherence of $\sF$ to check the $\sU$-coherence for a suitable covering. We will use the strict Laurent coverings of the form $X(f) \cup X(f^{-1}) = X$ for $f \in A$. So now, we check that the classical proof of Kiehl's theorem needs only the properties of the category of dagger affinoid algebras we proved so far and the
      non-archimedean nature of the base field plays no role in this proof.
\item Second step: We consider the covering $\sU = \{ U_1 = X(f),  U_2 = X(f^{-1}) \}$ and let $\sF$ be a $\sU$-coherent $\cO_X$-module.
      We claim that $H^1(\sU, \sF) = 0$. By our assumptions, $M_1 = \sF(U_1), M_2 = \sF(U_2)$, $M_{12} = \sF(U_1 \cap U_2)$
      are finite modules over $A \lt f \gt^\dagger$, $A \lt f^{-1} \gt^\dagger$, and $A \lt f, f^{-1} \gt^\dagger$ respectively, and the \u{C}ech
      complex of alternating cochains $C_a^{\bullet}(\sU,\sF)$ becomes
      \[ 0 \to M_1 \times M_2 \stackrel{d_0}{\to}  M_{12} \to 0. \]
      Since $H^1(\sU,\sF)$ can be computed using alternating cochains it is only necessary to show that $d_0 : M_1 \times M_2 \to M_{12}$ is
      surjective. Recall that we already proved that the map $A \lt f \gt^\dagger \times A \lt f^{-1} \gt^\dagger \to A \lt f, f^{-1} \gt^\dagger$ is surjective in theorem \ref{thm:tate}.
      Choose elements $v_1, \ldots , v_m \in M_1$ and $w_1, \ldots ,w_n \in M_2$ which generate $M_1$ as a $A \lt f \gt^\dagger$-module and $M_2$ as a $ A \lt f^{-1} \gt^\dagger$-module. Using the fact that $\sF$ is
      $\sU$-coherent, the restrictions of $v_1, \ldots , v_m$ to $U_1 \cap U_2$, as well as the restrictions of
      $w_1, \ldots, w_n$ to $U_1 \cap U_2$ generate $M_{12}$ as $A \lt f, f^{-1} \gt^\dagger$-module.
      We consider two epimorphisms
      \[ (A \lt f \gt^\dagger)^m \to M_1 \]
      \[ (A \lt f^{-1} \gt^\dagger)^n \to M_2. \]
      Since $\sF$ is $\sU$-coherent, $\sF|_{U_1}$ and $\sF|_{U_2}$ are generated by global sections, so the above epimorphisms give two epimorphisms
      \[ (A \lt f, f^{-1} \gt^\dagger)^m \to M_{12} \]
      \[ (A \lt f, f^{-1} \gt^\dagger)^n \to M_{12}. \]
      Now denoting $B = A \lt f \gt^\dagger \times A \lt f^{-1} \gt^\dagger$ we get (by construction) the following commutative diagram
      \[
      \begin{tikzpicture}
      \matrix(m)[matrix of math nodes,
      row sep=2.6em, column sep=2.8em,
      text height=1.5ex, text depth=0.25ex]
      { B^m \times B^n & (A \lt f, f^{-1} \gt^\dagger)^m \times (A \lt f, f^{-1} \gt^\dagger)^n & M_{12} \times M_{12} \\
       M_1 \times M_2 & & M_{12} \\};
      \path[->>,font=\scriptsize]
      (m-1-1) edge node[auto] {} (m-1-2);
      \path[->>,font=\scriptsize]
      (m-1-2) edge node[auto] {} (m-1-3);
      \path[->>,font=\scriptsize]
      (m-1-1) edge node[auto] {}  (m-2-1);
      \path[->>,font=\scriptsize]
      (m-1-3) edge node[auto] {$\Delta$}  (m-2-3);
      \path[->,font=\scriptsize]
      (m-2-1) edge node[auto] {$d_0$}  (m-2-3);
      \end{tikzpicture}
      \]
      where $\Delta$ is the codiagonal map, which is well-defined because $M_{12} \times M_{12} \cong M_{12} \oplus M_{12}$. All maps of the diagram are surjective by construction except for $d_0$, which is then necessarily surjective. (We remark that the argument we used to settle the second step of the proof is much shorter than the one in theorem 1.14.4 of \cite{BOSCH}).
      
\item Third step: We claim that from the previous step we can deduce that $\sF$ is associated to a finite $A$-module. 
      The proof of this step is identical to the classical proof of Kiehl's theorem. We write here all details for the convenience of the reader. 

      We claim that the fact that for any $\sU$-coherent module one has $H^1(\sU, \sF) = 0$, implies that $\sF$ is associated to a finite
      $A$-module. For this last step it is not necessary to make a difference between Laurent and more general admissible strict dagger affinoid coverings. 
      Therefore, consider a covering $\sU =  \{ U_i \}_{i=1,...,n}$ of $X$ consisting of dagger affinoid subdomains $U_i = \cM(A_i) \subset X$.
      Since $\sF$ is $\sU$-coherent, $\sF|_{U_i}$ is associated to a finite $A_i$-module $M_i$.
      Consider a maximal ideal $\fm_x \subset A$ (where with $x$ we denote the corresponding point $x \in \cM(A)$)
      and denote by $\fm_x \cO_X$ the associated coherent ideal sheaf of the structural sheaf $\cO_X$. Its product
      with $\sF$ yields a submodule $\fm_x \sF \subset \sF$ which is $\sU$-coherent, since its restriction
      to each $U_i$ is associated to the submodule $\fm_x M_i \subset M_i$ which is finite, because
      $M_i$ is a finite module over a Noetherian ring. Then $\sF/ \fm_x \sF$ is $\sU$-coherent by proposition \ref{prop_coherent} and
      \[ 0 \to \fm_x \sF \to \sF \to \sF/\fm_x \sF \to 0 \]
      is a short exact sequence of $\sU$-coherent $\cO_X$-modules. If $U' = \cM(A')$ is a dagger affinoid subdomain of $X$ which is contained in $U_i$ for
      some index $i$, then the above short exact sequence restricts to a short exact
      sequence of coherent modules on $U'$. More precisely, as the modules $\fm_x \sF$, $\sF$,
      and $\sF/\fm_x \sF$ are $\sU$-coherent, their restrictions to $U_i$ are associated to finite
      $A_i$-modules, from which it follows that the same is true for restrictions to $U'$ in terms of $A'$-modules.
      By proposition \ref{prop_coherent}, the above short exact sequence leads to a short exact sequence of $A'$-modules
      \[ 0 \to \fm_x \sF(U') \to \sF(U') \to \sF/\fm_x \sF(U') \to 0. \]
      In particular, $U'$ can be any intersection of sets in $\sU$, and we thereby see that the canonical sequence of \u{C}ech complexes
      \[ 0 \to C^\bullet( \sU, \fm_x \sF) \to C^\bullet(\sU, \sF) \to C^\bullet(\sU,\sF/ \fm_x \sF) \to 0 \]
      is exact. Since we are assuming that $H^1(U, \fm_x \sF) = 0$, the associated long exact sequence in cohomology yields a short exact sequence
      \[ 0 \to \fm_x \sF(X) \to \sF(X) \to \sF/\fm_x \sF(X) \to 0. \]
      Next we claim that the restriction homomorphism $\sF/\fm_x \sF(X) \to \sF/ \fm_x \sF(U_j)$ is bijective for any index $j$ such that 
      $x \in U_j$. To justify the claim, consider a dagger affinoid subdomain $U' = \cM(A') \subset X$ such that $\sF|_{U'}$ is associated to a finite
      $A'$-module $M'$ and write $U' \cap U_j = \cM(A_j')$. Then $\sF/\fm_x|_{U'}$ is associated to the quotient $M'/\fm_x M'$, and the canonical map
      \[ M'/\fm_x M' \to (M'/ \fm_x M') \otimes_A A_j' \to (M'/\fm_x M') \otimes_{A' / \fm_x A'} (A_j'/\fm_x A_j') \]
      is bijective for $x \in U_j$. This follows from proposition \ref{prop_subdomain_max_ideals} if $x \in U' \cap U_j$, 
      since the restriction map $A'/\fm_x A' \to A_j'/\fm_x A_j'$ 
      is bijective in this case. However, the latter is also true for $x \notin U'$ since in this case the quotients $A'/\fm_x A'$ and $A_j'/\fm_x A_j'$ 
      are trivial. Now, since $\sF$ is $\sU$-coherent, we look at the canonical diagram
      \[
      \begin{tikzpicture}
      \matrix(m)[matrix of math nodes,
      row sep=2.6em, column sep=2.8em,
      text height=1.5ex, text depth=0.25ex]
      { (\sF/ \fm_x \sF)(X) & \underset{i = 1}{\overset{n}\prod} (\sF/\fm_x \sF)(U_i) & \underset{i,k}{\overset{n}\prod} (\sF/\fm_x \sF)(U_i \cap U_k)   \\
        (\sF/ \fm_x \sF)(U_j) & \underset{i = 1}{\overset{n}\prod} (\sF/\fm_x \sF)(U_i \cap U_j) & \underset{i,k}{\overset{n}\prod} (\sF/\fm_x \sF)(U_i \cap U_k \cap U_j)  \\};
      \path[->,font=\scriptsize]
      (m-1-1) edge node[auto] {} (m-1-2);
      \path[->,font=\scriptsize]
      (m-1-2) edge node[auto] {} (m-1-3);
      \path[->,font=\scriptsize]
      (m-1-2) edge node[auto] {} (m-2-2);
      \path[->,font=\scriptsize]
      (m-1-1) edge node[auto] {}  (m-2-1);
      \path[->,font=\scriptsize]
      (m-2-1) edge node[auto] {}  (m-2-2);
      \path[->,font=\scriptsize]
      (m-2-2) edge node[auto] {}  (m-2-3);
      \path[->,font=\scriptsize]
      (m-1-3) edge node[auto] {}  (m-2-3);
      \end{tikzpicture}
      \]
      with exact rows. By the consideration above, the middle and right restriction morphisms are bijective. Thus, the same will hold 
      for the left one. To conclude, we look at the commutative diagram
      \[
      \begin{tikzpicture}
      \matrix(m)[matrix of math nodes,
      row sep=2.6em, column sep=2.8em,
      text height=1.5ex, text depth=0.25ex]
      { \sF(X) & \sF/\fm_x(X) \\
        M_j = \sF(U_j) & \sF/\fm_x(U_j) = M_j/\fm_x M_j \\};
      \path[->,font=\scriptsize]
      (m-1-1) edge node[auto] {} (m-1-2);
      \path[->,font=\scriptsize]
      (m-1-2) edge node[auto] {} (m-2-2);
      \path[->,font=\scriptsize]
      (m-1-1) edge node[auto] {}  (m-2-1);
      \path[->,font=\scriptsize]
      (m-2-1) edge node[auto] {}  (m-2-2);
      \end{tikzpicture}
      \]
      for $x \in U_j$, the exact sequence $0 \to \fm_x \sF(X) \to \sF(X) \to \sF/\fm_x \sF(X) \to 0$ and the bijection 
      $\sF/\fm_x \sF(X) \to \sF/ \fm_x \sF(U_j)$ imply that $M_j/\fm_x M_j$, as an $A_j$-module, is generated by the image of $\sF(X)$. Hence, by
      the classical lemma of Nakayama, $\sF(X)$ generates $M_j$ locally at each point $x \in \Max(A_j)$. Then the submodule of $M_j$ 
      which is generated by the image of $\sF(X)$ must coincide with $M_j$. Therefore, we can choose elements $f_1, \ldots, f_s  \in \sF(X)$
      such that their images generate all modules $M_i = \sF(U_i)$ simultaneously for $i = 1, \ldots, n$. As a consequence, 
      the morphism of $\cO_X$-modules $\phi: \cO_X^s \to \sF$ given by $f_1, \ldots , f_s$ is an epimorphism of $\sU$-coherent $\cO_X$-modules,
      and its kernel $\ker \phi$ is a $\sU$-coherent submodule of $\cO_X^s$ by proposition \ref{prop_coherent}. 
      We can work now in the same way as before with $\ker \phi$ in place of $\sF$ and construct an epimorphism $\psi: \cO_X^r \to \ker \phi$, 
      thus obtaining an exact sequence
      \[ \cO_X^r \stackrel{\psi}{\to}  \cO_X^s \stackrel{\phi}{\to} \sF \to 0 \]
      of $\cO_X$-modules. Hence, we see that $\sF$ is isomorphic to the cokernel of $\psi$, and so $\sF$ is associated to the cokernel 
      of the $A$-module morphism $\wtilde{\psi} : A^r \to A^s$ by proposition \ref{prop_coherent}. The latter is finite and, hence, $\sF$ 
      is associated to a finite $A$-module. This finishes the proof of the third step and also the proof of theorem.
\end{itemize}
\end{proof}

\begin{cor}
Let $A$ be a $k$-dagger affinoid algebra and $X = \cM(A)$. Then a $\cO_X$-modules is coherent if and only if is associated to a finite $A$-dagger module.
\end{cor}
\begin{proof}
The non-strict case is easily reduced to the strict case base changing the base field with $K/k$ for which $A \otimes_k^\dagger K$ becomes a strict $K$-dagger affinoid algebra and using lemma \ref{lemma_BER2_2_1_2}.
\end{proof}

\begin{cor}
Let $A$ be a $k$-dagger affinoid algebra and $X = \cM(A)$. Then, the category of finite $A$-modules is equivalent to the category of coherent shaves on $X$.
\end{cor}
\begin{proof}
By proposition \ref{prop_coherent} the association of a finite $A$-module to its associated coherent sheaf is a fully faithful functor. To prove that this functor is an equivalence it remains only to show that it is essentially surjective, which is precisely what theorem \ref{thm:kiehl} asserts.
\end{proof}

\begin{rmk}
The main contribution of theorem \ref{thm:kiehl} is to unify the proofs of known results in both archimedean and non-archimedean analytic geometry. Indeed, the dagger version of Kiehl's theorem has been already discussed in Grosse-Kl\"onne work \cite{GK} whereas the analog of Kiehl's theorem for compact Steins in complex geometry has been discussed in theorem 11.9.2 of \cite{TAY}. 
\end{rmk}
%

\index{special dagger subset}
\begin{defn}
Let $A$ be a $k$-dagger affinoid algebra. A closed subset $U \subset \cM(A)$ is called \emph{special} if it is a finite union of $k$-dagger affinoid subdomains.
\end{defn}

On $X = \cM(A)$ we can consider the \emph{special} $G$-topology whose admissible open subsets are given by special subsets and whose
coverings coincide with the coverings of the weak $G$-topology of $X$.

\begin{defn}
Let $A$ be a $k$-dagger affinoid algebra and $V \subset \cM(A)$ a special subset. Given a finite covering $V_1, \ldots, V_n \subset V$ of $V$ by dagger affinoid subdomains
we define
\[ A_V \doteq \ker \l ( \prod_{1 \le i \le n} A_{V_i} \to \prod_{1 \le i,j \le n} A_{V_i \cap V_j} \r ). \]
\end{defn}

\begin{prop} \label{prop_special_dagger}
Let $V \subset \cM(A) = X$ be a special subset. Then
\ben
\item $A_V$ does not depend on the choice of the dagger affinoid covering of $V$;
\item the correspondence $V \mapsto A_V$ is a sheaf for the special $G$-topology;
\item $V$ is a dagger affinoid subdomain if and only if $A_V$ is a $k$-dagger affinoid algebra, and in this case there is a canonical homeomorphism $\cM(A_V) \cong V$.
\een
\end{prop}
\begin{proof}
\ben
\item Since the association $U \mapsto A_U$ is a sheaf for the $G$-topology whose coverings are finite coverings by dagger affinoid subdomains the value of
 \[ \ker \l ( \prod_{1 \le i \le n} A_{V_i} \to \prod_{1 \le i,j \le n} A_{V_i \cap V_j} \r ) \]
does not depend on the choice of the covering of $V$ used to calculate it.
\item The family of coverings for the special $G$-topology and for the weak $G$-topology of dagger affinoid subdomains (essentially) coincide.
\item The direct implication is proved in corollary \ref{cor:subdomain_homeo}. For the converse, let $V_i$ be a finite covering of $V$ by rational subdomains (we can always do it
      by our generalization of the Gerritzen-Grauert theorem), \ie subdomains of the form
      \[ V_i = X \l (\frac{\ol{f}_i}{g_i} \r ), i = 1, \ldots, n. \]
      We have a canonical morphism
      \[ A \to \prod_{1 \le i \le n} A_{V_i} \]
      whose image lands in the kernel of 
      \[ \prod_{1 \le i \le n} A_{V_i} \to \prod_{1 \le i,j \le n} A_{V_i \cap V_j}. \]
      We will deduce that there exists a unique map $\a: A \to A_V$ that identifies $\cM(A_V)$ with the subspace $V \subset \cM(A)$. Given any morphism of dagger affinoid algebras $\be: A \to B$ such that $\Im(\cM(B)) \subset V$, there is a unique continuous map $\gamma^*: \cM(B) \to \cM(A_V)$ for which the following diagram 
      \[
      \begin{tikzpicture}
      \matrix(m)[matrix of math nodes,
      row sep=2.6em, column sep=2.8em,
      text height=1.5ex, text depth=0.25ex]
      { \cM(A) & \cM(A_V) \\
               & \cM(B) \\};
      \path[->,font=\scriptsize]
      (m-1-2) edge node[auto] {$\a^*$} (m-1-1);
      \path[->,font=\scriptsize]
      (m-2-2) edge node[auto] {$\be^*$} (m-1-1);
      \path[->,font=\scriptsize]
      (m-2-2) edge node[auto] {$\gamma^*$}  (m-1-2);
      \end{tikzpicture}
      \]
      commutes. We have to show that this map is induced by a map of affinoid algebras $\gamma: A_V \to B$. Let $1 \le i_1, \ldots, i_m \le n$ be the indexes such that $g_i$ does not vanish in $V$. Then, it is clear that the morphism of $k$-algebras $A [g_{i_1}^{-1}, \ldots, g_{i_m}^{-1}] \rhook A_V$ has dense image. Moreover, $\be(g_i)$ is also a unit in $\cM(B)$ hence there is a unique morphism $A [g_{i_1}^{-1}, \ldots, g_{i_m}^{-1}] \to B$ extending $A \to B$. Putting all together we get the diagram
      \[
      \begin{tikzpicture}
      \matrix(m)[matrix of math nodes,
      row sep=2.6em, column sep=2.8em,
      text height=1.5ex, text depth=0.25ex]
      { A &  & A_V \\
          &  A[g_{i_1}^{-1}, \ldots, g_{i_m}^{-1}]   & \\
          &  & B  \\};
      \path[->,font=\scriptsize]
      (m-1-1) edge node[auto] {$$} (m-1-3);
      \path[->,font=\scriptsize]
      (m-1-1) edge node[auto] {$$} (m-2-2);
      \path[->,font=\scriptsize]
      (m-2-2) edge node[auto] {$$}  (m-1-3);
      \path[->,font=\scriptsize]
      (m-2-2) edge node[auto] {$$}  (m-3-3);
      \path[->,font=\scriptsize, dashed]
      (m-1-3) edge node[auto] {$\gamma$}  (m-3-3);
      \end{tikzpicture}
      \]
      where the dashed arrow $\gamma$ is uniquely defined by the uniqueness of $A [g_{i_1}^{-1}, \ldots, g_{i_m}^{-1}] \to B$ and the density of the image of $A [g_{i_1}^{-1}, \ldots, g_{i_m}^{-1}] \rhook A_V$. This is the required morphism of dagger affinoid algebras representing $\gamma^*$.
\een
\end{proof}

\section{Dagger analytic spaces}

In this section we address the problem of the construction of the category of dagger analytic spaces over $k$. We start by recalling the properties of Berkovich nets that we will use later on. We will follow very closely
the discussion of the first section of \cite{BER4}, checking that all the main results can be translated in our settings with minor changes in some proofs.

\begin{notation}
Let $X$ be a set, $\tau \subset \sP(X)$ a collection of subsets of $X$ and $Y \subset X$. We write
\[ \tau|_Y \doteq \{ V \in \tau | V \subset Y \} \]
and call it the \emph{restriction} of $\tau$ to $Y$.
\end{notation}

\index{Berkovich net}
\index{Berkovich quasi-net}
\begin{defn}
Let $X$ be a topological space and $\tau \subset \sP(X)$ we say that:
\ben
\item $\tau$ is \emph{dense} if for any $V \in \tau$ each point $x \in V$ has a neighborhood basis of elements of $\tau|_V$;
\item $\tau$ is a \emph{quasi-net on $X$} if for any $x \in X$ there exists $V_1, \cdots, V_n \in \tau$ such that 
      $x \in V_1 \cap \cdots \cap V_n$ and $V_1 \cup \cdots \cup V_n$ is a neighborhood of $x$;
\item $\tau$ is a \emph{net on $X$} (or better a \emph{Berkovich net}) if $\tau$ is a quasi-net and for any $U, V \in \tau$ the collection $\tau|_{U \cap V}$
      is a quasi-net on $U \cap V$.
\een 
\end{defn}

\begin{prop}
Let $\tau$ be a quasi-net on $X$ then
\ben
\item a subset $U \subset X$ is open if and only if for any $V \in \tau$ the subset $U \cap V$ is open in $V$;
\item if every element of $\tau$ is compact, then $X$ is Hausdorff if and only if for any $U, V \in \tau$ the intersection $U \cap V$ is compact.
\een
\end{prop}
\begin{proof} 
\cite{BER4}, lemma 1.1.1.
\end{proof}

\begin{prop} \label{prop_lemma_ber_nets}
Let $\tau$ be a net of compact subsets of $X$ then
\ben
\item for any $U, V \in \tau$ the intersection $U \cap V$ is locally closed in $U$ and $V$;
\item if $V \subset V_1 \cup \cdots \cup V_n$ for $V, V_1, \cdots, V_n \in \tau$ then there exist $U_1, ..., U_m \in \tau$ such that
      \[ V = U_1 \cup \cdots \cup U_m \]
      and $U_i \subset V_{j_i}$ for any $i \in \{1, \cdots, m \}$ and $j_i \in \{1, \cdots, n\}$.
\een
\end{prop}
\begin{proof} 
\cite{BER4}, lemma 1.1.2.
\end{proof}

\index{$k$-dagger analytic space}
\begin{defn}
Let $X$ be a locally Hausdorff topological space. A $k$-dagger affinoid atlas on $X$ is the data of a net $\tau$ on $X$ 
such that to each $V \in \tau$ there is assigned a $k$-dagger affinoid algebra $A_V$ and a homeomorphism $V \cong \cM(A_V)$ such that 
if $V' \in \tau$ and $V' \subset V$ we assign a map $\a_{V/V'}: A_V \to A_{V'}$ which is the dagger affinoid subdomain embedding of $V'$ in $V$.

The triple $(X, A, \tau)$ is called a \emph{$k$-dagger analytic space}. If all $A_V$ are strict $k$-dagger affinoid algebras then this triple is
called a strictly $k$-dagger analytic space.
\end{defn}

By propositions \ref{prop_loc_arcwise_connectedness} and \ref{prop_base_countable_infinity} each point of a dagger affinoid space has a fundamental 
system of open arcwise connected subsets that are countable at infinity. It follows that a $k$-dagger analytic space has a basis for the topology formed by open locally compact paracompact arcwise connected subsets.

\begin{lemma}
Let $(X, A, \tau)$ be a $k$-dagger analytic space. If $W$ is a $k$-dagger affinoid subdomain in some $U \in \tau$, then it is a $k$-dagger affinoid 
subdomain in any $V \in \tau$ that contains $W$. 
\end{lemma}
\begin{proof} 
Since $\tau|_{U \cap V}$ is a quasi-net and $W$ is compact, we can find $U_1, \ldots, U_n \in \tau|_{U \cap V}$ with 
$W \subset U_1 \cup \ldots \cup U_n$. Furthermore, since $W$ and $U_i$, are dagger affinoid subdomains in $U$, then $W_i \doteq W \cap U_i$ and
$W_i \cap W_j$ are dagger affinoid subdomains in $U_i$ and a fortiori also in $V$. By Tate's Acyclicity Theorem, applied to the dagger affinoid 
covering $\{ W_i \}$ of $W$ the bornological algebra 
\[ A_W = \Ker(\prod_{1 \le i \le n} A_{W_i} \to \prod_{1 \le i,j \le n} A_{W_i \cap W_j}) \]
is $k$-dagger affinoid and $W \cong \cM(A_W)$. So $W$ is a $k$-dagger affinoid subdomain in $V$ by proposition \ref{prop_special_dagger}.
\end{proof}

\index{completion of a Berkovich net}
\begin{defn}
Let $(X, A, \tau)$ be a $k$-dagger analytic space. With $\ol{\tau}$ we denote the family of all $W \subset X$ such that $W$ is a $k$-dagger affinoid 
subdomain in some $V \in \tau$. Then, $\ol{\tau}$ is a dense net in $X$ and we call $\ol{\tau}$ the \emph{completion} of the net $\tau$.
\end{defn}

\begin{prop}
Let $(X, A, \tau)$ be a $k$-dagger analytic space. The family $\ol{\tau}$ is a net on $X$, and there exists a unique (up to a 
canonical isomorphism) dagger affinoid atlas $\ol{A}$ for the net $\ol{\tau}$ which extends $\tau$. 
\end{prop}
\begin{proof} 
The proof is identical to \cite{BER4}, proposition 1.2.6.
\end{proof}

\index{strong morphism of $k$-dagger analytic spaces}
\begin{defn}
A \emph{strong morphism} of $k$-dagger analytic spaces $\phi: (X, A, \tau) \to (X', A', \tau')$ consists of a continuous map 
$\phi: X \to X'$, such that for each $V \in \tau$ there exists $V' \in \tau'$ with $\phi(V) \subset V'$, and for any such 
pair $(V, V')$ a system of compatible (with respect to the atlases $\tau$ and $\tau'$) morphisms of $k$-dagger affinoid spaces $\phi_{V/V'}: (V, A_V) \to (V', A_{V'})$. 
\end{defn}

\begin{prop}
Any strong morphism $\phi: (X, A, \tau) \to (X', A', \tau')$ induces in a unique way to a strong morphism 
$\ol{\phi}: (X, A, \ol{\tau}) \to (X', A', \ol{\tau}')$. 
\end{prop}
\begin{proof} 
Let $U$ and $U'$ be dagger affinoid subdomains in $V \in \tau$ and $V' \in \tau'$, respectively, and suppose that $\phi(U) \subset U'$.
Take $W' \in \tau'$ with $\phi(V) \subset W'$. Then $\phi(U) \subset W_1 \cup \ldots \cup W_n$ for some $W_1, \ldots, W_n \in \tau'|_{V' \cap W'}$, because $\phi(U)$ is compact. 
The morphism of dagger affinoid spaces $\phi_{V/W'}$ induces a morphism $V_i \doteq \phi_{V/W'}^{-1}(W_i) \to W_i$ that induces, in its turn, 
a morphism of dagger affinoid spaces $U_i \doteq U \cap V_i \to U'$. Thus, we have morphisms of dagger affinoid spaces $U_i \to U \to U'$ that are
compatible on intersections. So they give a system of morphisms of dagger affinoid spaces $\ol{\phi}_{U/U'}: (U, A_U) \to (U', A_{U'})$ which are clearly compatible, 
showing that we can extend $\phi$ to the morphism $\ol{\phi}$ in a unique way.
\end{proof}

The composition of two strong morphisms $\phi$ and $\psi$ is simply the composition of the underlying continuous map and of the system of dagger affinoid maps, 
which can always be arranged in a compatible way passing to the maps $\ol{\phi}$ and $\ol{\psi}$ canonically induced by $\phi$ and $\psi$ by previous proposition.
The category of $k$-dagger analytic spaces with strong morphisms is denoted by $\wtilde{\bAn}_k^\dagger$.

\index{quasi-isomorphism of $k$-dagger analytic spaces}
\begin{defn}
A strong morphism $\phi: (X, A, \tau) \to (X', A', \tau')$ is said to be a \emph{quasi-isomorphism} if $\phi$ induces a homeomorphism 
between $X$ and $X'$ and, for any pair $V \in \tau$ and $V' \in \tau'$ with $\phi(V) \subset V'$, $\phi$ identifies $V$ with a
dagger affinoid subdomain in $V'$. 
\end{defn}

It is easy to see that if $\phi$ is a quasi-isomorphism, then so is $\ol{\phi}$. The following easy lemma is needed to show that the system of quasi-isomorphisms admits
calculus of fractions.

\begin{lemma} \label{lemma:quasi_iso}
Let $\phi: (X, A, \tau) \to (X', A', \tau')$ be a strong morphism. Then, for any pair $V \in \ol{\tau}$ and $V' \in \ol{\tau}'$ the intersection 
$V \cap \phi^{-1}(V')$ is a finite union of $k$-dagger affinoid subdomains in $V$.
\end{lemma}
\begin{proof}
Take $U' \in \ol{\tau}'$ with $\phi(V) \subset U'$. Then, we can find $U_1', \ldots, U_n' \in \tau'|_{U' \cap V'}$ such that $\phi(V) \subset U'_1 \cup \dots \cup U_n'$, and 
\[ V \cap \phi^{-1}(V') = \bigcup_{i = 1}^n \phi^{-1}_{V/U'}(U_i'). \]
\end{proof}

\begin{prop} \label{prop:quasi_iso}
The system of quasi-isomorphisms in $\wtilde{\bAn}_k^\dagger$ admits calculus of right fractions.
\end{prop}
\begin{proof} 
We have to verify that the system of quasi-isomorphisms satisfies the following properties:
\ben
\item all the identity morphisms are quasi-isomorphisms;
\item the composition of two quasi-isomorphisms is a quasi-isomorphism; 
\item any diagram of the form $(X, A, \tau) \stackrel{\phi}{\to} (X', A', \tau') \stackrel{g}{\leftarrow} (X'', A'',\tau'')$, where $g$ is a quasi-isomorphism, can be completed to a commutative square 
      \[
      \begin{tikzpicture}
      \matrix(m)[matrix of math nodes,
      row sep=2.6em, column sep=2.8em,
      text height=1.5ex, text depth=0.25ex]
      { (X, A, \tau) & (X', A', \tau') \\
        (X''', A''', \tau''') & (X'', A'', \tau'') \\};
      \path[->,font=\scriptsize]
      (m-1-1) edge node[auto] {$\phi$} (m-1-2);
      \path[->,font=\scriptsize]
      (m-2-2) edge node[auto] {$g$} (m-1-2);
      \path[->,font=\scriptsize]
      (m-2-1) edge node[auto] {$f$}  (m-1-1);
      \path[->,font=\scriptsize]
      (m-2-1) edge node[auto] {}  (m-2-2);
      \end{tikzpicture}
      \]
      where $f$ is a quasi-isomorphism;
\item given two strong morphisms $\phi, \psi: (X, A, \tau) \to (X', A', \tau')$ and a quasi-isomorphism $g: (X', A', \tau') \to (X'', A'', \tau'')$ such that 
      $g \circ \phi = g \circ \psi$, then there exists a quasi-isomorphism $f: (X''', A''', \tau''') \to (X, A, \tau)$ with $\phi \circ f = \psi \circ f$.
      (We shall show, in fact, that in this situation $\phi = \psi$).
\een

The property (1) is clearly valid. To verify (2) it suffices to apply the definition of the composition and the third statement of proposition \ref{prop_special_dagger}. 
Suppose then, that we have a diagram as in (3). We can assume that $X' = X''$ and then $\tau'' \subset \tau'$. Let $\tau'''$ denote the family of all $V \in \ol{\tau}$ for 
which there exists $V'' \in \tau''$ with $\phi(V) \subset V''$. By lemma \ref{lemma:quasi_iso} it follows that $\tau'''$ is a net on $X$. We can attach to $\tau'''$ a $k$-affinoid atlas 
by restricting the dagger affinoid atlas $\ol{\tau}$ and the strong morphism $\ol{\phi}$ induces a strong morphism 
$\sigma: (X, A''', \tau''') \to (X', A', \tau')$. Then $\sigma$ factors through the canonical quasi-isomorphism $(X, A''', \tau''') \to (X, A, \tau)$, obtaining the required commutative diagram. 

Finally, we claim that in the situation (4) the morphisms $\phi$ and $\psi$ coincide. First of all, it is clear that the underlying maps of topological spaces coincide. 
Furthermore, let $V \in \tau$ and $V' \in \tau'$ be such that $\phi(V) \subset V'$. Take $V'' \in \tau''$ with $g(V') \subset V''$. Then we have two morphisms $k$-dagger affinoid 
spaces $\phi_{V/V'}, \psi_{V/V'}: V \to V'$ such that their compositions with $g_{V'/V''}$ coincide. Since $V'$ is a dagger affinoid domain in $V''$, it follows that $\phi_{V/V'} = \psi_{V/V'}$, because $g_{V'/V''}$ is a monomorphism by corollary \ref{cor:mono_loc_closed}. 
\end{proof}

\begin{defn}
Let $(X, A, \tau)$ be a $k$-dagger analytic space. If $\sigma$ is a net on $X$, we write $\sigma < \tau$ if $\sigma \subset \ol{\tau}$ and we say that $\sigma$ is \emph{finer}
then $\tau$.
\end{defn}

The system of nets $\{\sigma\} \subset \sP(X)$ such that $\sigma < \ol{\tau}$ is filtered and there is a canonical quasi-isomorphism
\[ (X, A_\sigma, \sigma) \to (X, \ol{A}, \ol{\tau}) \]
where $A_\sigma$ is the restriction of $\ol{A}$ to $\sigma$.

\index{category of $k$-dagger analytic spaces}
\begin{defn}
The category obtained by inverting quasi-isomorphisms of $\wtilde{\bAn}_k^\dagger$ is denoted $\bAn_k^\dagger$ and is called the \emph{category of $k$-dagger analytic spaces}.
\end{defn}

The category of $k$-dagger affinoid spaces embeds fully faithfully in $\bAn_k^\dagger$.
We can state the following corollary to proposition \ref{prop:quasi_iso}, which is simply a more explicit way to describe the localization we
just defined.

\begin{cor}
Let $(X, A, \tau), (X', A', \tau')$ be two $k$-dagger analytic spaces, then
\[ \Hom_{\bAn_k^\dagger}((X, A, \tau), (X', A', \tau')) = \limind_{\sigma < \tau} \Hom_{\wtilde{\bAn}_k^\dagger}((X, A_\sigma, \sigma), (X', A', \tau')).  \]
\end{cor}
\begin{proof}
See \cite{GABZIS} chapter 1 \S 2, where it is shown that previous formula holds when the system of morphisms to invert satisfies calculus of right fractions, as we showed in proposition \ref{prop:quasi_iso}.
\end{proof}

To continue our study of $k$-dagger analytic spaces we need to substitute the concept of $k$-dagger affinoid subdomain with a more general concept.

\index{$\tau$-special subset}
\begin{defn} \label{defn:t_special}
Let $(X, A, \tau)$ be a $k$-dagger analytic space. We say that a subset $W \subset X$ is \emph{$\tau$-special} if it is compact and there exists a covering $W = W_1 \cup \ldots \cup W_n$
such that $W_i, W_i \cap W_j \in \tau$, for all $1 \le i,j \le n$, and 
\[ A_{W_i} \otimes_k^\dagger A_{W_j} \to A_{W_i \cap W_j} \]
is a strict epimorphism, for all $1 \le i,j \le n$. A covering of $W$ of the above type is said to be a \emph{$\tau$-special covering} of $W$. 
\end{defn}

\begin{prop} \label{prop:t_special}
Let $W$ be a $\tau$-special subset of $(X, A, \tau)$. If $U, V \in \tau|_W$, then $U \cap V \in \ol{\tau}$ and
\[ A_U \otimes_k^\dagger A_V \to A_{U \cap V} \]
is a strict epimorphism.
\end{prop} 
\begin{proof} 
Let $\{W_i\}$ be a covering  of $W$ as in definition \ref{defn:t_special}, \ie a $\tau$-special covering. 
Since the subsets $U \cap W_i$, and $V \cap W_i$, are compact for any $i$, we can find finite coverings $\{ U_{i, k} \}$ of $U \cap W_i$, and 
$\{ V_{i,k} \}$ of $V \cap W_i$, by subsets from $\tau$. Furthermore, since the diagonal maps $W_i \cap W_j \to W_i \times W_j$, are closed immersions, it follows that 
$U_{i, k} \cap V_{j, l} \in \ol{\tau}$
and $U_{i, k} \cap V_{j, l} \to U_{i, k} \times V_{j, l}$ are closed immersions for all $i,j,k,l$. Consider now the finite dagger affinoid covering 
$\{ U_{i, k} \times V_{j, l} \}$ of the dagger affinoid space $U \times V$. For each quadruplet $i, j, k, l$, the algebra $A_{U_{i, k} \cap V_{j, l}}$ is a finite $A_{U_{i, k} \times V_{j, l}}$-algebra, and the system $\{ A_{U_{i, k} \cap V_{j, l}} \}$ satisfies the condition of Kiehl's theorem
applied to the $k$-dagger affinoid space $U \times V$, \ie they glue to give a finite $A_{U \times V}$-module. So this system defines a finite $A_{U \times V}$-algebra whose underlying module is isomorphic to 
\[ A_{U \cap V} \doteq \Ker( \prod A_{U_{i, k} \cap V_{j, l}} \to \prod A_{U_{i, k} \cap V_{j, l} \cap U_{i', k'} \cap V_{j', l'}} ), \]
and therefore $A_{U \cap V}$ is $k$-dagger affinoid (by proposition \ref{prop_finite_dagger_algebra}). Since $\cM(A_{U \cap V}) \cong U \cap V$ we see that 
$U \cap V$ is a dagger affinoid subdomain in U and V \ie, $U \cap V \in \ol{\tau}$, by proposition \ref{prop_special_dagger}.
\end{proof}

\begin{cor}
Let $W$ be a $\tau$-special subset of $X$, then any finite covering of $W$ by subsets from $\tau$ is a $\tau$-special covering.
\end{cor}

\index{specialization of a Berkovich net}
\begin{defn}
Let $(X, A, \tau)$ be a $k$-dagger analytic space. We denote by $\what{\tau}$ the set of $\ol{\tau}$-special subsets such that the algebra
\[ A_W \doteq \Ker \l ( \prod A_{W_i} \to \prod A_{W_i \cap W_j } \r ), \]
for a $\ol{\tau}$-special covering $\{W_i\}$, satisfies the following conditions:
\ben 
\item $A_W$ is a $k$-dagger affinoid algebra;
\item $\cM(A_W) \cong W$;
\item the $k$-dagger affinoid spaces $(W_i, A_{W_i})$ defines $k$-dagger affinoid subdomains in $(\cM(A_W), A_W)$.
\een
We call $\what{\tau}$ the \emph{specialization} of the net $\tau$.
\end{defn}

\begin{prop} \label{prop:max_atlas}
Let $(X, A, \tau)$ be a $k$-dagger analytic space, then
\ben
\item the collection $\what{\tau}$ is a net;
\item for any net $\sigma$ such that $\sigma < \tau$ one has $\what{\sigma} = \what{\tau}$;
\item there exists a unique (up to a canonical isomorphism) dagger analytic atlas $\what{A}$ on the net $\what{\tau}$ that extends the atlas $A$;
\item $\what{\what{\tau}} = \what{\tau}$.
\een
\end{prop}
\begin{proof} 
\ben
\item Let $U, V \in \tau$ and take $\ol{\tau}$-special coverings $\{ U_i \}$ of $U$ and $\{ V_j \}$ of $V$. Since $U \cap V = \underset{i,j}\bigcup (U_i \cap V_j)$ and $\ol{\tau}|_{U_i \cap V_j}$
      are quasinets, it follows that $\what{\tau}|_{U \cap V}$ is a quasinet.
\item Let $\sigma$ be a net with $\sigma < \tau$. By proposition \ref{prop:t_special}, to verify the equality $\what{\sigma} = \what{\tau}$, it suffices to show that for any 
      $V \in \ol{\tau}$ there exist $U_1, \ldots, U_n \in \ol{\sigma}$ with $V = U_1 \cup \ldots \cup U_n$. Since $\sigma$ is a net on $X$, we can find $W_1, \ldots, W_n \in \sigma$ with 
      $V \subset W_1 \cup \ldots \cup W_n$, because $V$ is compact. Since $V, W_i \in \ol{\tau}$ and $\ol{\tau}$ is a net of compact subsets we can find $U_1, \ldots, U_n \in \ol{\tau}$ such that $V = U_1 \cup \ldots \cup U_n$ and each $U_i$ is contained in some $W_j$ applying proposition \ref{prop_lemma_ber_nets}. Finally, since $W_j \in \sigma$, it follows that $U_i \in \ol{\sigma}$.
\item For each $V \in \what{\tau}$ we fix a $\ol{\tau}$-special covering $\{ V_i \}$ and assign to $V$ the algebra $A_V$ and the homeomorphism $\cM(A_V) \cong V$. Then, we have to construct for each $U \in \what{\tau}$ with $U \subset V$ a canonical bounded homomorphism $A_V \to A_U$  that identifies $(U, A_U)$ 
      with a dagger affinoid subdomain in $(V, A_V)$. Consider first the case when $U \in \tau$. By proposition \ref{prop:t_special}, $U \cap V_i$ is a dagger affinoid subdomain in $V_i$
      and therefore in $V$. It follows that $U$ is a dagger affinoid subdomain in $V$. If $U$ is arbitrary we can consider a $\ol{\tau}$-special covering $\{U_i\}$ of $U$ and then by what we said each $U_i$ is a dagger affinoid subdomain in $V$. It follows that $U$ is a dagger affinoid domain in $V$. 
\item From proposition \ref{prop:t_special} it follows that $\ol{\what{\tau}} = \what{\tau}$. Let $\{ V_i \}$ be a $\what{\tau}$-special covering of some $V \in \what{\what{\tau}}$. 
      For each $i$ we take a $\ol{\tau}$-special covering $\{V_{i,j}\}$, of $V_i$. Then $\{V_{i,j}\}$, is a $\ol{\tau}$-special covering of $V$, and therefore $V \in \what{\tau}$.
\een
\end{proof}

\begin{rmk}
(4) of proposition \ref{prop:max_atlas} tells that $\widehat{\tau}$ is the maximal atlas on $X$ in the isomorphism class of $k$-dagger analytic spaces in which $(X, \tau, A)$ belongs.
\end{rmk}

\index{dagger affinoid domain}
\begin{defn}
Let $(X, A, \tau)$ be a $k$-dagger analytic space. The subsets of $X$ which are in $\what{\tau}$ are said to be \emph{$k$-dagger affinoid domains} in $X$.
\end{defn}

\index{good $k$-dagger analytic space}
\begin{defn}
We say that a $k$-dagger analytic space $(X, A, \tau)$ is \emph{good} if for each point $x \in X$ there is a $k$-dagger affinoid domain which is neighborhood of $x$.
\end{defn}

Now we want to study how to glue $k$-dagger analytic spaces, how to put a structural sheaf on them and to study the relations between the weak $G$-topology induced by $k$-dagger affinoid domains and its 'true' topology of the underlying space. The key definition to do such a study is the next one. From now on we often suppress the explicit reference to the net $\tau$ and the atlas $A$ of a dagger analytic space and we always suppose to work always with its maximal atlas given by $\what{\tau}$.

\index{analytic domain}
\begin{defn}
A subset $Y$ of a $k$-dagger analytic space $(X, A, \tau)$ is said to be a \emph{$k$-dagger analytic domain} if, for any point $x \in Y$, there exist $k$-dagger affinoid domains $V_1, \ldots, V_n \in \tau$ that are contained in $Y$ and such that $x \in V_1 \cap \ldots \cap V_n$ 
and the set $V_1 \cup \ldots \cup V_n$ is a neighborhood of $x \in Y$ (\ie, the restriction of the net affinoid domains on $Y$ is a net on $Y$). 
\end{defn}

We remark that the intersection of two dagger analytic domains is a dagger analytic domain, and the preimage of a dagger analytic domain with respect to a morphism of a dagger analytic 
spaces is a dagger analytic domain. Furthermore, the family of dagger affinoid domains that are 
contained in a dagger analytic domain $Y \subset X$ defines a dagger affinoid atlas on $Y$ and there is a canonical morphism of dagger analytic spaces $\iota: Y \to X$ which satisfy the following universal property. For any morphism $\phi: Z \to X$ 
with $\phi(Z) \subset Y$ there exists a unique morphism $\psi: Z \to Y$ with $\phi = \iota \circ \psi$.
A morphism $\phi: Y \to X$ that induces an isomorphism of $Y$ with a dagger analytic domain in X is said to be an \emph{open immersion}
also if $Y$ may not be open in the topology of $X$. We remark also that then all open subsets of $X$ are dagger analytic domains. 
It is clear that a dagger analytic domain that is isomorphic to a $k$-dagger affinoid space is a $k$-dagger affinoid domain. 

\begin{prop} \label{prop_exact_sequence_analytic_subdomain}
Let $\{Y_i\}_{i \in I}$ be a covering of a $k$-dagger analytic space $(X, A, \tau)$ by dagger analytic domains such that each point of 
$x \in X$ has a neighborhood of the form $Y_1 \cup \ldots \cup Y_n$, with $x \in Y_1 \cap \ldots Y_n$ 
(\ie, $\{Y_i\}_{i \in I}$ is a quasinet on $X$). Then for any dagger analytic space $X'$ the following sequence of sets is exact 
\[ \Hom(X, X') \to \prod_{i \in I} \Hom(Y_i, X') \to \prod_{i,j \in I}\Hom(Y_i \cap Y_j, X'). \]
\end{prop}
\begin{proof} 
Let $\phi_i: Y_i \to X'$ be a family of morphisms such that, for all pairs $i,j \in I$, $\phi_i|_{Y_i \cap Y_j} =\phi_j|_{Y_i \cap Y_j}$. Then, these $\phi_i$ define a map of underlying topological space  
$X \to X'$ which is continuous. Furthermore, let $\sigma$ be the collection of dagger affinoid domains $V \subset X$ such that there exist $i \in I$ and a dagger affinoid domain 
$V' \subset X'$ with $V \subset Y_i$, and $\phi(V) \subset V'$. $\sigma$ is a net on $X$ and $\sigma < \tau$, therefore there is a morphism $\phi: X \to X'$ 
that gives rise to all the morphisms $\phi_i$, by the localization we used to define morphism of dagger analytic spaces.
\end{proof}

We now describe two ways of gluing $k$-dagger analytic spaces. Let $\{ X_i \}_{i \in I}$ be a family 
of $k$-dagger analytic spaces, and suppose that, for each pair $i,j \in I$, we are given dagger analytic 
domains $X_{i, j} \subset X_i$ and an isomorphism of $k$-dagger analytic spaces $\phi_{i, j}: X_{i,j} \to X_{j,i }$ so that 
$X_{i,i} = X_i$, $\phi_{i,j}(X_{i, j} \cap X_{i ,l}) = X_{j, i} \cap X_{j,l}$ and $\phi_{i,l} = \phi_{j,l} \circ \phi_{i,j}$ on 
$X_{i,j} \cap X_{i,l}$, for any $i,j,l \in I$. We are looking for a $k$-dagger analytic space $X$ with a family of morphisms 
$\mu_i : X_i \to X$ such that: 
\ben
\item $\mu_i$ is an isomorphism of $X_i$ with a $k$-dagger analytic domain in $X$; 
\item the family $\{ \mu_i(X_i) \}_{i \in I}$ covers $X$; 
\item $\mu_i(X_{i,j}) = \mu_i(X_i) \cap \mu_j(X_j)$; 
\item $\mu_i = \mu_j \circ \phi_{i,j}$ on $X_{i,j}$.
\een

If such $X$ exists, we say that it is obtained by \emph{gluing $X_i$ along $X_{i,j}$}.

\begin{prop} \label{prop_glueing}
The space $X$ obtained by gluing $X_i$ along $X_{i,j}$ exists and is unique (up to a canonical isomorphism) in each of the following cases: 
\ben
\item all $X_{i,j}$ are open in $X_i$; 
\item for any $i \in I$, all $X_{i,j}$ are closed in $X_i$ and the number of $j \in I$ with $X_{i, j} \ne \void$ is finite. 
\een
Furthermore, in the case 1) all $\mu_i(X_{i,j})$ are open in $X$. In the case 2) all $\mu_i(X_{i,j})$ are closed in $X$ 
and if all $X_i$ are Hausdorff (resp. paracompact), then $X$ is Hausdorff (resp. paracompact).
\end{prop}
\begin{proof} 
Let $\wtilde{X}$ be the disjoint union $\underset{i \in I}\coprod X_i$. The system $\{ \phi_{i, j} \}$ defines an equivalence relation $R$ on $\wtilde{X}$. We denote by $X$ the quotient space 
$\wtilde{X}/R$ and by $\mu_i: X_i \to X$ the induced maps. In the case 1. the equivalence relation $R$ is open, 
and therefore all $\mu_i(X_i)$ are open in $X$. In the case 2., the equivalence relation $R$ is closed, and therefore all $\mu_i(X_i)$ are closed in $X$ and $\mu_i$ induces 
a homeomorphism $X_i \to \mu_i(X_i)$. Moreover, if all $X_i$ are Hausdorff then $X$ is Hausdorff by the fact that $R$ is closed. If all $X_i$ are paracompact, then $X$ is paracompact because it has a 
locally finite covering by closed paracompact subsets.  Furthermore, let $\tau$ denote the collection of all subsets $V \subset X$ for which there 
exists $i \in I$ such that $V \subset \mu_i(X_i)$ and $\mu_i^{-1}(V)$ is a $k$-dagger affinoid domain in $X_i$ (in this case $\mu_j^{-1}(V)$ is a $k$-dagger affinoid domain in $X_j$, 
for any $j$ with $V \subset  \mu_j(X_j)$). It is easy to see that $\tau$ is a net, and there is an evident $k$-dagger affinoid atlas $A$ with this net induced by the atlases on the $X_i$. In this way, 
we get a $k$-dagger analytic space ($X, A, \tau)$ that satisfies the properties 1.-4. above. That $X$ is unique up to canonical isomorphism follows from proposition
\ref{prop_exact_sequence_analytic_subdomain}.
\end{proof}

\begin{defn}
Let $X$ be a $k$-dagger analytic space. The family of its dagger analytic domains can be considered as a category, which can be equipped with a Grothendieck topology 
generated by the pretopology for which the set of coverings of a dagger analytic domain $Y \subset X$ is formed by the families $\{ Y_i \}$ of analytic domains in $Y$ that 
are quasinets on $Y$. For brevity, the above Grothendieck topology is called the \emph{(Berkovich) $G$-topology on $X$} and the corresponding site is denoted by $X_G$. 
\end{defn}

From Proposition \ref{prop_exact_sequence_analytic_subdomain} it follows that any representable presheaf on $X_G$ is a sheaf.

\index{dagger affine space}
\begin{defn}
We define the \emph{$n$-dimensional affine space} $\A_k^n$ over $k$ as the spectrum of the Fr\'echet algebra
\begin{equation} \label{eqn:affine_spectrum}
\A_k^n = \cM(\limpro_{\rho \to \infty} W_k^n(\rho)) \cong \limind_{\rho \to \infty} \cM(W_k^n(\rho))
\end{equation}
equipped with its canonical $k$-dagger analytic space structure.
\end{defn}

To see that indeed $\underset{\rho \to \infty}\limpro W_k^n(\rho)$ is a Fr\'echet algebra we can use the same reasoning of proposition 
\ref{prop_good_open_subsets} and the homemorphism of (\ref{eqn:affine_spectrum}) is proved in corollary \ref{cor:frechet}. The family of all closed polydisks centered in zero in $\A_k^n$, with their canonical $k$-dagger structure, defines on $\A_k^n$ a $k$-dagger
affinoid atlas with respect to which $\A_k^n$ is a good analytic space. Moreover, $\A_k^1$ is a ring object in the category of $k$-dagger analytic spaces and we have the following proposition.

\begin{prop}
Let $(A, \cM(A))$ be a $k$-dagger affinoid space equipped with its canonical structure of $k$-dagger analytic space, then
\[ \Hom_{\bAn_k^\dagger}(\cM(A), \A_k^1) \cong A.  \]
\end{prop}
\begin{proof}
By the universal property that charcterize the algebras $W_k(\rho)$ we see that
\[ \Hom_{\bAn_k^\dagger}(\cM(A), \A_k^1) = \Hom_{\bAn_k^\dagger}(\cM(A), \limind_{\rho \to \infty} \cM(W_k(\rho))) \cong \]
\[ \cong \limind_{\rho \to \infty}\Hom_{\bAn_k^\dagger}(\cM(A), \cM(W_k(\rho))) = \limind_{\rho \to \infty} \rho A^\ov = A \]
because
\[ A = \limind_{\rho \to \infty} \rho A^\ov \]
follows by the spectral characterization of weakly power-boundedness. Notice that we can write the isomorphism
\[ \Hom_{\bAn_k^\dagger}(\cM(A), \limind_{\rho \to \infty} \cM(W_k(\rho))) \cong \limind_{\rho \to \infty}\Hom_{\bAn_k^\dagger}(\cM(A), \cM(W_k(\rho))) \]
because the spectrum of $A$ is compact and therefore its image in $\A_k^1$ by a continuous map must land in $\cM(W_k(\rho))$ for some $\rho$.
\end{proof}

So, it make sense to give the following definition of the structural sheaf of
a $k$-dagger analytic space.

\begin{defn} \label{defn:structural_sheaf}
Let $X \in \ob(\bAn_k^\dagger)$ we define $\cO_{X_G}$, \emph{the structural sheaf of $X$ for the $G$-topology}, as the sheaf
\[ U \mapsto \cO_{X_G}(U) = \Hom_{\bAn_k^\dagger}(U, \A_k^1) \]
for any dagger analytic domain of $X$.
\end{defn}

\begin{rmk}
In the complex analytic case this definition naturally generalize the classical definition of structural sheaf (we will see in next sections how to embed fully faithfully the category of classical complex analytic spaces in the category of $\C$-dagger analytic spaces). Thus, given any complex analytic space $X$, thought as a dagger analytic space, any open subset of $X$ is a dagger analytic domain. We have that the structural sheaf just defined gives the association
\[ U \mapsto \cO_{X_G}(U) = \Hom_{\bAn_\C^\dagger}(U, \C) \]
for any open set in $X$, \ie it precisely restricts to the classical structural sheaf of complex analytic geometry on open subsets of $X$.
\end{rmk}

\begin{rmk}
We also remark that if $(X, A, \tau)$ is a Berkovich analytic space over $k$, which we now suppose non-archimedean, and $X$ is such that every affinoid algebra $A_U$, for $U \in \tau$, can be written as
\[ A_U \cong \frac{T_k^{n_U}(\rho)}{(f_1^U, \ldots, f_n^U)} \]
with $n_U \in \N$ and $f_1^U, \ldots, f_n^U \in T_k^{n_U}$ overconvergent analytic functions on the disk $\cM(T_k^{n_U}(\rho))$, 
then we can associate to $A_U$ canonically a dagger affinoid algebra $A_U^\dagger$ and hence to $A$ canonically
an atlas $A^\dagger$ of dagger affinoid algebras. In this way, we can associate to $(X, A, \tau)$ the $k$-dagger analytic space $(X, A^\dagger, \tau)$. We will say more on this topic
later on. What we want to underline now is that if $X$ is without borders for any open (with respect to the topology of $X$) subset $V \subset X$ we get the canonical isomorphism
\[ \Hom_{\bAn_k^\dagger}(U, \A_k^1) \cong \cO_{X_G}^{\text{Ber}}(U) \]
of Fr\'echet algebras, where $\cO_{X_G}^{\text{Ber}}(U)$ is the structural sheaf of $(X, A, \tau)$ as defined in \cite{BER4}.
\end{rmk}

\section{Sheaves over dagger analytic spaces}

Let $X$ be a $k$-dagger analytic space. The category of $\cO_{X_G}$-modules over $X_G$ is denoted by $\bMod(X_G)$.

\index{coherent sheaf on a $k$-dagger analytic spaces}
\begin{defn}
An object $\sF \in \ob(\bMod(X_G))$ is said to be \emph{coherent} if there exists a quasinet $\tau$ of $k$-dagger affinoid domains in $X$ such that, for each $V \in \tau$, 
$\sF|_{V_G}$ is isomorphic to the cokernel of a homomorphism of free $\cO_{V_G}$-modules of finite rank.
The category of coherent $\cO_{X_G}$-modules over $X_G$ is denoted by $\bCoh(X_G)$.
\end{defn}

By Kiehl's theorem (cf. \ref{thm:kiehl}) this is equivalent to require that for the given quasinet $\tau$ of $k$-dagger affinoid domains in $X$ the restriction $\sF|_{V_G}$, for $V \in \tau$, is
associtated to a finite $\cO_{V_G}(V)$-module.

\index{support of a sheaf}
\begin{defn}
For any object $\sF \in \ob(\bCoh(X_G))$ we define the \emph{support} of $\sF$ as the set
\[ \Supp(\sF) \doteq \{ x \in X | \sF_x \ne 0 \}. \]
Moreover, for any ideal sheaf $\sI \subset \cO_{X_G}$ the \emph{variety} of $\sI$ is defined
\[ V(\sI) \doteq \Supp \l( \frac{\cO_{X_G}}{\sI} \r). \]
\end{defn}

\index{Picard group}
\begin{defn}
We define $\bPic(X_G)$ to be the group of isomorphisms classes of invertible sheaves on $X_G$ and we call it \emph{the Picard group} of $X_G$, equipping it with 
the multiplication given by the tensor product of sheaves.
\end{defn}

We notice that $\bPic(X_G) \cong H^1(X_G, \cO_{X_G}^\times)$, by a general result on $G$-ringed spaces, cf. \cite{HART} exercise III.4.5. 

Since every open subset of $X$ is a $k$-dagger analytic domain, then the identity map
\[ \pi: X_G \to X \]
is a morphism of $G$-topological spaces, from which we get a morphism of topoi that we now describe. The direct image functor is simply the restriction \ie for any presheaf $\sF$ on $X_G$ 
and any open subset $U$ we have that
\[ (\pi_* \sF)(U) = \sF(U). \]
In particular, this permits to equip canonically $X$ with the structural sheaf
\[ \cO_X \doteq \pi_*(\cO_{X_G}) \]
which gives to $(X, \cO_X)$ the structure of a locally ringed space. The inverse image functor is defined
\[ (\pi^{-1} \sF)(U) = \limind_{V \supset U} \sF(V). \]
for any presheaf $\sF$ on $X$, where the direct limit is taken on the open neighborhoods of $U$ for any $k$-dagger analytic domain $U$ of $X$. 

\begin{notation}
Given a site $X$, we denote with $\wtilde{X}$ its associated topos.
\end{notation}

\begin{prop}
The couple $\ol{\pi} = (\pi_*, \pi^{-1}): \wtilde{X}_G \to \wtilde{X}$ defines a morphism of topoi. Moreover, for any sheaf $\sF$ on $X$ we have that $\sF \cong \pi_* \pi^{-1} \sF$, 
so in particular $\pi^{-1}$ is fully faithful.
\end{prop}
\begin{proof}
First we notice that the relation $\sF \cong \pi_* \pi^{-1} \sF$ follows easily form the definitions because for any open set $U \subset X$
\[ (\pi_* \pi^{-1} \sF )(U) = \limind_{V \supset U} \sF (U) = \sF(U) \]
where $V \subset U$ runs through the open neighborhood of $U$, but $U$ is an open neighborhood of itself. This fact implies that $\pi^{-1}$ is fully faithful by a well-known result of abstract non-sense, after we show that $(\pi_*, \pi^{-1})$ is an adjoint pair of functors.

So, we have only to check that $(\pi_*, \pi^{-1})$ is an adjoint pair of functors.
Let $\phi: \pi^{-1} \sF \to \sG$ be a morphism of sheaves on $X_G$, then this induces a morphism
\[ \pi_* \phi: \pi_* \pi^{-1} \sF \cong \sF \to \pi_* \sG. \]
On the other hand, given a morphism $\phi: \sF \to \pi_* \sG$, this induces a morphism
\[ \pi^{-1} \phi: \pi^{-1} \sF \to \pi^{-1} \pi_* \sG. \]
Let now $U \subset X$ be an admissible open in $X_G$ then we have a canonical map
\[ (\pi^{-1} \pi_*) \sG(U) = \limind_{V \subset U} \sG(V) \to \sG(U)  \]
obtained taking the limit map of the restriction maps $\sG(V) \to \sG(U)$. This system of maps gives a canonical morphism of sheaves $\pi^{-1} \pi_* \sG \to \sG$ which
gives the required adjuntion. 
To end the proof we need to check that $\pi^{-1}$ is an exact functor, and this follows from the exactness of the functor $\limind$ for filtered index sets.
\end{proof}

\begin{exa}

 We use the same counter-example used in \cite{BER4}, remark 1.3.8, to see that $\pi_*$ is not fully faithful. 
      Let $X$ be the closed unit disk on $k$, with $k$ non-archimedean, and let $x_0$ be the maximal (Gauss) point of $X$. 
      We construct two sheaves $\sF$ and $\sF'$ on $X_G$ as follows. Let $Y$ be a dagger analytic domain in $X$. 
      Then define $\sF(Y) = \Z$ if $x_0 \in Y$, and $\sF(Y) = 0$ otherwise and $\sF'(Y) = \Z$ if $\{ x \in X | r < |T(x)| < 1 \} \cup \{x_0\} \subset Y$ for some $0 < r < 1$, 
      and $\sF'(Y) = 0$ otherwise. The sheaves $\sF$ and $\sF'$ are manifestly not isomorphic on $X_G$, but $\pi_* \sF = \pi_* \sF' = i_* \Z$, where $i$ is the embedding $\{ x_0 \} \to X$. 
\end{exa}

Let $X_\bPro$ denote the pro-site of $X$ relative to the Grothendieck topology of $X$ defined by the underlying topological space, as defined in definition \ref{defn_pro_site}.
There are canonical morphisms of topoi
\[ \ol{f} = (f_*, f^{-1}): \wtilde{X}_\bPro \to \wtilde{X} \]
and 
\[ \ol{g} = (g_*, g^{-1}): \wtilde{X}_\bPro \to \wtilde{X}_G \]

\begin{prop}
The following diagram of morphism of topoi 
\[
\begin{tikzpicture}
\matrix(m)[matrix of math nodes,
row sep=2.6em, column sep=2.8em,
text height=1.5ex, text depth=0.25ex]
{ \wtilde{X}_\bPro & \wtilde{X} \\
  \wtilde{X}_G     &  \\};
\path[->,font=\scriptsize]
(m-1-1) edge node[auto] {$\ol{f}$} (m-1-2);
\path[->,font=\scriptsize]
(m-1-1) edge node[auto] {$\ol{g}$} (m-2-1);
\path[->,font=\scriptsize]
(m-2-1) edge node[auto] {$\ol{\pi}$}  (m-1-2);
\end{tikzpicture}
\]
is commutative.
\end{prop}
\begin{proof}
The morphism $\ol{f}: \wtilde{X}_\bPro \to \wtilde{X}$ is the topos morphism as defined in the appendix A, after definition \ref{defn:pro_open}. The map $\ol{f}$ factors through $\ol{\pi}$ because the definitions of $\ol{f}$ and $\ol{\pi}$ agree on the elements of $X_G$.
\end{proof}

We recall, also this from appendix A, that we can pullback the structural sheaf of $X$ to $X_\bPro$ giving to $X_\bPro$ the structure of a ringed topos.
Hence, the morphisms $\ol{\pi}, \ol{f}, \ol{g}$ defines naturally adjoint pairs of functors 
\[ \pi_*: \bMod(X_G) \leftrightarrows \bMod(X): \pi^* \]
\[ f_*: \bMod(X_\bPro) \leftrightarrows \bMod(X): f^* \]
\[ g_*: \bMod(X_\bPro) \leftrightarrows \bMod(X_G): g^* \]
where the direct image functors are as above and the pullback functors are defined
\[ \pi^*(\sF) \doteq \pi^{-1} \sF \otimes_{\pi^{-1} \cO_X} \cO_{X_G} \]
\[ f^*(\sF) \doteq f^{-1} \sF \otimes_{f^{-1} \cO_X} \cO_{X_\bPro} \]
\[ g^*(\sF) \doteq g^{-1} \sF \otimes_{g^{-1} \cO_{X_G}} \cO_{X_\bPro}. \]

\begin{defn}
Let $\sF \in \bMod(X)$ we says that $\sF$ is \emph{coherent} if locally on $X$ is isomorphic to a cokernel of a morphism of free finite rank modules.
We denote with $\bCoh(X)$ the category of coherent $\cO_X$-modules over $X$ and with $\bPic(X) = H^1(X, \cO_X^\times)$ the Picard group of $X$.
\end{defn}

\begin{prop} \label{prop_good_cO_sheaves}
Let $X$ be a $k$-dagger analytic space, then 
\ben 
\item for any $\cO_X$-module $\sF$ one has 
\[ \sF \cong f_* f^* \sF; \]
so, in particular the functor $f^*$ is fully faithful;
\item for any $\cO_X$-module $\sF$ one has 
\[ \sF \cong \pi_* \pi^* \sF; \]
so, in particular the functor $\pi^*$ is fully faithful;
\een
and if $X$ is good
\ben 
\item the functor $f^*$ restricts to an equivalence of categories $\bCoh(X) \to \bCoh(X_\bPro)$; 
\item the functor $\pi^*$ restricts to an equivalence of categories $\bCoh(X) \to \bCoh(X_G)$.
\een
\end{prop}
\begin{proof}
The two parts of the proposition have very similar proofs that therefore are written together.

\ben
\item It is sufficient to verify that for any point $x \in X$ there is an isomorphism of stalks $\sF_x \cong (f_* f^* \sF)_x$ (resp. $\sF_x \cong (\pi_* \pi^* \sF)_x$). Since every open subset of $X$ is open in $X_\bPro$ and in $X_G$ then $f^* \sF (U) = \sF(U)$ and $\pi_* \pi^* \sF(U) = \sF(U)$. The fully faithfulness of $f^*$ and $\pi^*$ follows by abstract non-sense.

\item Since $f^*$ (resp. $\pi^*$) is fully faithful, it suffices to verify that for a coherent $\cO_{X_G}$-module $\sF$ the $\cO_X$-module $f_* \sF$ (resp. $\pi_* \sF$) is coherent and (resp. $\sF \cong \pi^* \pi_* \sF$). 
Since $X$ is good then $f_* \sF$ (resp. $\pi_* \sF$) is coherent because every point has an affinoid neighborhood where we can write it as a the sheaf associated to a finite module. By the coherence of $f_* \sF$ (resp. $\pi_* \sF$) it follows that on a net of $k$-dagger affinoid domains $\{U_i\}_{i \in I}$ of
$X_\bPro$ (resp. $X_G$) the sheaf $f^* f_* \sF$ (resp. $\pi^* \pi_* \sF$) is associated finite $\cO_{U_i}$-modules which agrees with $\sF(U_i)$, hence the assertion.
\een
\end{proof}

\begin{prop} \label{prop_good_cO_loc_free_sheaves}
If $X$ is a good $k$-dagger analytic space, then a coherent $\cO_X$-module $\sF$ is locally free if and only if $\pi^*(\sF)$ is locally free.
\end{prop}
\begin{proof}
We may assume that $X = \cM(A)$ is $k$-dagger affinoid since we are dealing with a local property on a good $k$-analytic space. Since $A$ is Noetherian it suffices to show that a 
      finite $A$-module $M$ is projective if and only if the $\cO_{X_G}$-module associated to $M$ is locally free. The direct implication follows by the characterization
      of projectivity on Noetherian modules. Indeed, a finitely presented module $M$ is projective if and only if it is locally free for the Zaraski topology of $\Spec(A)$, and
      this readily implies that $M$ is locally free on $\cM(A)$ because Zariski opens are open in $\cM(A)$.
      Conversely, suppose that for some finite $k$-dagger affinoid covering $\{ V_i \}_{i \in I}$ of $X$ the finite $A_{V_i}$-modules 
      $M \otimes_A^\dagger A_{V_i}$ are free. It suffices to verify that $M$ is flat over $A$, because for finitely presented modules over commutative rings flatness is equivalent to projectivity. 
      To show flatness, we take an injective homomorphism of finite $A$-modules $P \to Q$. Then the homomorphisms 
      \[ (M \otimes_A^\dagger P) \otimes_A^\dagger A_{V_i} \to (M \otimes_A^\dagger Q) \otimes_A^\dagger A_{V_i}, \]
      are also injective for all $i$, because we are supposing $M \otimes_A^\dagger A_{V_i}$ free. Applying Tate's Acyclicity Theorem to the finite $A$-modules 
      $M \otimes_A^\dagger P$ and $M \otimes_A^\dagger Q$, we obtain the injectivity of the homomorphism $M \otimes_A^\dagger P \to M \otimes_A^\dagger Q$.
\end{proof}

\begin{cor}
If $X$ is a good $k$-dagger analytic space, then there is an isomorphism $\bPic(X) \cong \bPic(X_G)$.
\end{cor}

\begin{rmk}
For proving the result of proposition \ref{prop_good_cO_loc_free_sheaves} for the pro-site of $X$ one needs a better local understanding of non good dagger analytic spaces.
\end{rmk}

\begin{prop}
Let $X$ be a good $k$-dagger analytic space:
\ben
\item for any abelian sheaf  $\sF$ on $X$, one has $H^q(X, \sF) \cong H^q(X_G, \pi^{-1} \sF)$, for any $q \ge 0$;
\item then $H^q(X, \sF) \cong H^q(X_G, \pi^* \sF)$, for any $q \ge 0$ and any $\sF \in \ob(\bCoh(X))$. 
\een
\end{prop}
\begin{proof}
\ben
\item Any open covering of $X$ is a covering both for the topology of $X$ and for the $G$-topology. Therefore it generates two Leray spectral sequences that are convergent to the 
      groups $H^q(X, \sF)$ and $H^q(X_G, \pi^{-1} \sF)$, respectively. Comparing them, we see that it suffices to verify the statement for sufficiently small $X$, 
      so we can assume that $X$ is paracompact. It suffices to verify that if $\sF$ is injective, then $H^q(X_G, \pi^{-1} \sF) = 0$ for $q \ge 1$. 
      Since $X$ is paracompact, it suffices to verify that the \u{C}ech cohomology groups of $\pi^{-1} \sF$ with respect to a locally finite covering by compact analytic domains are 
      trivial. And choosing such a covering $\{ U_i \}_{i \in I}$ for which the topological interiors of the $U_i$ cover $X$, which can always be done because $X$ is good,
      the \u{C}ech cohomology groups of $\pi^{-1} \sF$ with respect to this covering agree with the \u{C}ech cohomology groups of $\sF$ with respect to the covering of $X$ 
      obtained by the topological interior of $U_i$, and hence vanishes.
\item The same reasoning as before reduces the proof to the case when $X \subset \cM(A)$ is an open paracompact subset of a $k$-dagger affinoid space, and
      we can suppose that $X$ does not intersect the border of $\cM(A)$.
      In this case $H^q(X, \sF)$ is an inductive limit of the $q$-th cohomology groups of the \u{C}ech complexes associated with locally 
      finite open coverings $\{ U_i \}_{i \in I} $ of $X$. On other hand, since the cohomology groups of a coherent sheaf on a dagger affinoid space are trivial with respect
      to its weak $G$-topology, then 
      $H^q(X_G, \pi^* \sF)$ is the $q$-th cohomology group of the \u{C}ech complex associated with an arbitrary locally finite affinoid covering $\{ V_j \}_{j \in J}$ of $X$. 
      It remains to remark that for any open covering $\{ U_i \}_{i \in I}$ we can find a dagger affinoid covering $\{ V_j \}_{j \in J}$ such that each $V_j$ is contained 
      in some $U_i$ and $\bigcup_j \Int(V_j/X) = X$, because we are supposing $X \subset \Int(\cM(A))$ and $U_i$ are open in $X$.
\een
\end{proof}

We recall that we showed, in the previous section, that a morphism of $k$-dagger analytic spaces $\phi: Y \to X$ induces a morphism of $G$-ringed topological spaces $\phi_G: Y_G \to X_G$. 
If the spaces $X$ and $Y$ are good, then for any coherent $\cO_X$-module $\sF$ there is a canonical isomorphism of coherent $\cO_{Y_G}$-modules 
\[ \pi_Y^*(\phi^* \sF) \cong \phi_G^* (\pi_X^* \sF) \]
which can be readily verified by calculating the sections on any analytic domain of $Y$.

\begin{lemma}
The following properties of a morphism of $k$-dagger analytic spaces $\phi: Y \to X$ are equivalent
\ben 
\item for any point $x \in X$ there exist $k$-dagger affinoid domains $V_1, \ldots, V_n \subset X$ such that $x \in V_1 \cap \ldots \cap V_n$, $V_1 \cup \ldots \cup V_n$ is a neighborhood of $x$ and
      the induces maps $\phi^{-1}(V_i) \to V_i$ are finite morphisms (resp. closed immersions) of $k$-dagger affinoid spaces;
\item for any dagger affinoid domain $V \subset X$, $\phi^{-1}(V) \to V$ is a finite morphism (resp. a closed immersion) of $k$-dagger affinoid spaces. 
\een
\end{lemma}
\begin{proof}
The second condition implies the first obviously. So, suppose that the first condition is true. Then, the hypothesis implies that the collection $\tau$ of all dagger affinoid domains 
$V \subset X$ such that $\phi^{-1}(V) \to V$ is a finite morphism (resp. a closed immersion) of affinoid spaces is a net. So, since an arbitrary dagger affinoid domain $V$ is compact then $V \subset V_1 \cup \ldots \cup V_n$ for some $V \in \tau$.
By proposition \ref{prop_lemma_ber_nets}, we can find affinoid domains $U_1, \ldots, U_m \subset X$ such that $V = U_1 \cup \ldots \cup U_m$ and each 
$U_i$ is contained in some $V_i$. So this implies that $\phi$ induces a finite morphism (resp. closed immersion) on each $U_i$ (because they are subdomains of some $V_j$)
and so by Kiehl theorem $\phi$ induces the same kind of morphism on $V$.
\end{proof}

\index{finite morphism between dagger analytic spaces}
\index{closed immersion between dagger analytic spaces}
\begin{defn}
A morphism $\phi: Y \to X$ satisfying the equivalent properties of the previous lemma is said to be \emph{finite} (resp. a \emph{closed immersion}).
\end{defn}

A finite morphism $\phi: Y \to X$ induces a compact map with finite fibres on the underlying topological spaces, and $\phi_*(\cO_{Y_G})$ is a coherent 
$\cO_{X_G}$-module. If $\phi$ is a closed immersion, then it induces a homeomorphism of the topological space $Y$ with 
its image in $X$ and the homomorphism $\cO_{X_G} \to \phi_*(\cO_{Y_G})$ is surjective. 
Its kernel is a coherent sheaf of ideals in $\cO_{X_G}$. Furthermore, we say that a subset $Z \subset X$ is \emph{Zariski closed} 
if, for any dagger affinoid domain $V \subset X$, the intersection $Z \cap V$ is \emph{Zariski closed} in $V$. 
The complement to a Zariski closed subset is called \emph{Zariski open}. 
If $\phi: Y \to X$ is a closed immersion, then the image of $Y$ is Zariski closed in $X$. Conversely, if $Z$ is Zariski closed in $X$, then 
there is a closed immersion $Y \to X$ that identifies the underlying topological space of $Y$ with $Z$.

Furthermore, a morphism of $k$-dagger analytic spaces $\phi: Y \to X$ is said to be a $G$-\emph{locally} (resp. locally) \emph{closed immersion} if there exist a quasinet $\tau$
of analytic (resp. open analytic) domains in $Y$ and, for each $V \in \tau$, an analytic (resp. an open analytic) domain 
$U \subset X$ such that $\phi$ induces a closed immersion $V \to U$. 

\begin{prop}
The category $\bAn_k^\dagger$ admits fibre products. 
\end{prop}
\begin{proof}
The same argument of \cite{BER4}, proposition 1.4.1 applies. 
\end{proof}

\index{separated morphism between dagger analytic spaces}
\begin{defn}
A morphism of $k$-dagger analytic spaces $\phi: Y \to X$ is said to be \emph{separated} (resp. \emph{locally separated}) 
if the diagonal morphism $\Delta_{Y/X}: Y \to Y \times_X Y$ is a closed (resp. a locally closed) immersion. 
If the canonical morphism $X \to \cM(k)$ is separated (resp. locally separated), then $X$ is said to be \emph{separated} (resp. \emph{locally separated}).
\end{defn}

Recall that a continuous map of topological spaces $f: Y \to X$ is said to be \emph{Hausdorff} if for 
any pair of different points $x, y \in Y$ with $f(x) = f(y)$ there exist open neighborhoods 
$U_x$ of $x$ and $U_y$ of $y$ with $U_x \cap U_y = \void$. This is equivalent to require that the image of $Y$ in $Y \times_X Y$ is closed and so, applying this definition
to the identity map, one recovers the definition of Hausdorff topological space.
We remark that if $f: Y \to X$ is a Hausdorff map and $X$ is a Hausdorff topological space, then also $Y$ is a Hausdorff topological space. 

Furthermore, let $X$ and $Y$ be topological spaces and suppose that each point of $X$ has 
a compact neighborhood. In this situation a compact (or proper) map $f: Y \to X$ is Hausdorff, then it takes closed subsets of $Y$ to closed subsets of $X$ and each point of $Y$ has a compact neighborhood.

\begin{prop}
A locally separated morphism of $k$-dagger analytic spaces $\phi: Y \to X$ is separated if and only if the induced continuous map of underlying topological spaces is Hausdorff.
\end{prop}
\begin{proof}
If a morphism $\phi: Y \to X$ is separated, then the set-theoretic image of $Y$ is closed in the topological space $Y \times_X Y$, because it corresponds to $V(\sI)$
for some coherent ideal sheaf on the analytic space $Y \times_X Y$. By the construction of the fiber product, \cf \cite{BER4}, one can see that it exists a canonical map 
\[ \pi: |Y \times_X Y| \to |Y| \times_{|X|} |Y| \]
where $|\cdot|$ denotes the underlying topological spaces. $\pi$ is a compact map, then, as remarked so far, it is in particular a closed map, so $\pi(\phi(Y))$ is closed also in $|Y| \times_{|X|} |Y|$ and $\phi$ is a Hausdorff map. 

On the other hand, if the map $\phi: Y \to X $ induces a Hausdorff map of topological spaces, then $W = (Y \times_X Y) - \Delta_{Y/X}(Y)$ is an open subset.
Since by hypothesis the diagonal morphism $\Delta_{Y/X}$ is a composition of a closed immersion with an open immersion, it suffices to show that $\Delta_{Y/X}(Y)$ is closed 
in $| Y \times_X Y |$. To show this claim we consider the compact map $\pi$ as above and $z \in (Y \times_X Y) - \Delta_{Y/X}(Y)$, and let $\pi(z) = (y_1, y_2)$. 
If $y_1 \ne y_2$ then $\pi^{-1}(W)$ is an open neighborhood of $z$ that does not meet $\Delta_{Y/X}(Y)$. If $y_1 = y_2$, then we take an open neighborhood $V$ of $y_1 = y_2$ 
such that $\Delta_{V/X}: V \to V \times_X V$ is a closed immersion. Since $V \times_X V$ is an open subset of $Y \times_X V$ and $z \notin \Delta_{Y/X}(Y)$, then we can find an open neighborhood 
of $z$ that does not meet $\Delta_{Y/X}(Y)$. The required fact follows. 
\end{proof}

In particular it follows that if $Y$ is separated, then its underlying topological space is Hausdorff. 

\begin{cor}
A morphism of good $k$-dagger analytic spaces is separated if and only if the induces continuous map is Hausdorff. In particular a good $k$-dagger
analytic space is separated if and only if is Hausdorff.
\end{cor}
\begin{proof}
Good $k$-dagger spaces are locally separated because every point has a neighborhood isomorphic to a $k$-dagger affinoid space and $k$-dagger affinoid space are separated.
\end{proof}

The notion of boundary of a morphism between analytic spaces we gave so far for dagger affinoid spaces can be globalized. Our aim is to do it in the same fashion it has been done in definition 1.5.4 of \cite{BER4} but it seems to us that that definition has a circularity problem. Indeed, it asserts that "The relative interior of a morphism $\phi: Y \to X$ is the set $\Int(Y/X)$ of all points $y \in Y$ for which there exists an open neighborhood $\cV$ of $y$ such that the induced morphism $\cV \to X$ is closed". But, by definition (cf. page 49 of \cite{BER2}), a closed morphism is a morphism such that $\partial(Y/X) = \emptyset$. Therefore, the notion of closed morphism relies on the notion of interior and boundary of a morphism and cannot be used to define them! So, we propose an alternative, more explicit, definition.

\index{relative interior of a morphism of dagger analytic space}
\begin{defn} \label{defn:boundary_analytic_space}
The \emph{relative interior} of a morphism $\phi: (Y, A_Y, \tau_Y) \to (X, A_X, \tau_X)$ is the set $\Int(Y/X)$ 
of all points $y \in Y$ for which there exists $V_1, \ldots, V_n \in \tau_Y$ and $U_1, \ldots, U_m \in \tau_X$ such that 
\ben
\item $y \in V_1 \cap \ldots \cap V_n$ and $V_1 \cup \ldots \cup V_n$ is a neighbohood of $y$;
\item $\phi(y) \in U_1 \cap \ldots \cap U_m$ and $U_1 \cup \ldots \cup U_m$ is a neighbohood of $\phi(y)$;
\item for each $1 \le i \le n$ there exists a $1 \le j \le m$ such that $\phi(V_i) \subset U_j$ and $y \in \Int(Y/X)$.
\een
The complement of $\Int(Y/X)$ is called the \emph{relative boundary} of $\phi$ and is denoted by $\partial(Y/X)$. If $X = \cM(k)$, these sets are denoted by $\Int(Y)$ and $\partial(Y)$ and are called the \emph{interior} and the \emph{boundary} of $Y$, respectively.
\end{defn}

\index{closed morphism of dagger analytic spaces}
\begin{defn} \label{defn:closed_morphism}
A morphism $\phi: Y \to X$ of $k$-dagger analytic spaces is called \emph{closed} if $\partial(Y/X) = \emptyset$.
\end{defn}

\begin{prop}
Let $\phi: Y \to X$ be a morphism of $k$-dagger analytic spaces.
\ben
\item If $Y$ is a dagger analytic domain in $X$, then $\Int(Y/X)$ coincides 
with the topological interior of $Y$ in $X$. 
\item For a sequence of morphisms $Z \stackrel{\psi}{\to} Y \stackrel{\phi}{\to} X$, one has 
\[ \Int(Z/Y) \cap \psi^{-1}(\Int(Y/X)) \subset \Int(Z/X). \]
If $\phi$ is locally separated (resp. and good) then 
\[ \Int(Z/X) \subset \Int(Z/Y) \]
\[ (resp. \ \ \ \Int(Z/X) = \Int(Z/Y) \cap \psi^{-1}(\Int(Y/X) \ \ ). \]
\item For a morphism $f: X' \to X$, one has $(f')^{-1}(\Int(Y/X)) \subset \Int(Y'/X')$, with $f': Y \times_X X' \to X'$.
\item For a non-archimedean field $K$ over $k$, one has 
\[ \pi^{-1}(\Int(Y/X)) \subset \Int(Y \otimes_k^\dagger K /X \otimes_k^\dagger K), \] 
where $\pi: Y \otimes_k^\dagger K \to Y$. 
\een
\end{prop}

\section{Relations with classical analytic geometries}

In this section we study the relations between the analytic spaces we constructed so far and others present in the literature. 
In particular, we start from the non-archimedean case and we compare our $k$-dagger analytic spaces with the rigid analytic spaces of defined by Tate, with Berkovich's spaces and with Grosse-Kl\"onne's dagger spaces. After that, we discuss the archimedean complex case.

\subsection{The non-archimedean case}

Assume that $k$ is a non-archimedean complete valued field and that the valuation on $k$ is nontrivial. Let $X$ be a Hausdorff strictly 
$k$-dagger analytic space and consider the set
\[ X_0 = \{ x \in X | [\cH(x):k] < \infty \}. \]
From \cite{BER2}, proposition 2.1.15, it follows that the set $X_0$ is everywhere dense in $X$. If $X = \cM(A)$ is strictly $k$-dagger affinoid space
then $X_0 = \Max(A)$, so $X_0$ can be equipped with the structure of dagger analytic space in the sense of Gr\"osse-Klonne, \cite{GK}. Therefore, we consider on $X_0$ the weak (and strong) $G$-topology generated by the family of strictly dagger affinoid subdomains of $X_0$. We briefly recall the main definitions of dagger spaces in the sense of Gr\"osse-Klonne.

\index{rigid dagger analytic spaces}
\begin{defn} \label{defn:rigid_dagger}
A locally $G$-ringed space $(X, \cO_X)$ whose $G$-topology is saturated is said a \emph{rigid $k$-dagger space} if there exists a covering
\[ X = \bigcup_i U_i \]
by admissible open subsets such that $(U_i, \cO_X|_{U_i})$ is isomorphic to the locally $G$-ringed
space of a dagger affinoid space (equipped with the strong G-topology) associated to some strict $k$-dagger algebra $A_i$. A morphism of rigid $k$-dagger spaces is a morphism of locally $G$-ringed spaces.
\end{defn}

We can state for rigid $k$-dagger spaces the following standard definitions.

\begin{defn}
Let $f: X \to Y$ be a morphism of rigid $k$-dagger space. We say that $f$ is
\ben
\item an \emph{open immersion} if is a homeomorphism onto its image and if for any $x \in X$ the canonical morphism $f_x: \cO_{Y, f(x)} \to \cO_{X, x}$ is an isomorphism;
\item a \emph{closed immersion} if is injective and if the canonical map $\cO_Y \to f_*(\cO_X)$ is surjective;
\item a \emph{locally closed immersion} if is injective and if the canonical maps $f_x: \cO_{Y, f(x)} \to \cO_{X, x}$ are surjective for all $x \in X$;
\item \emph{quasi-compact} if for any quasi-compact open subspace $Y' \subset Y$ the inverse image $f^{-1}(Y')$ is a quasi-compact subset of $X$;
\item \emph{quasi-separated} if the diagonal morphism $\Delta_{X/Y}: X \to X \times_Y X$ is a closed immersion;
\item \emph{separated} if the diagonal morphism $\Delta_{X/Y}: X \to X \times_Y X$ is quasi-compact morphism.
\een
\end{defn}

\begin{defn}
A rigid $k$-dagger space $X$ is said
\ben
\item \emph{quasi-compact} if it has a finite admissible covering by dagger affinoid subdomains;
\item \emph{quasi-separated} (resp. proper, resp. \emph{separated}) if the canonical morphism $X \to \cM(k)$ is quasi-separated (resp. proper (see definition \ref{defn:proper}), resp. separated);
\item \emph{normal at a point $x \in X_0$} if $\cO_{X, x}$ is a normal ring; it is said \emph{normal} if is normal at every point in $X_0$;
\item \emph{reduced at a point $x \in X_0$} if $\cO_{X, x}$ is a reduced ring; it is said \emph{reduced} if is reduced at every point in $X_0$;
\item \emph{smooth at a point $x \in X_0$} if $\cO_{X, x}$ is a regular ring; it is said \emph{smooth} (or \emph{regular}) if is smooth at every point in $X_0$;
\item \emph{quasi-algebraic} if there exists an admissible covering $X = \bigcup X_i$ such that for all $i$ there exists an open immersion $X_i \rhook Y_i^\an$ into
      the analytification of a $k$-scheme $Y_i$ of finite type (see section 3.3 of \cite{GK} for the definition of the dagger analytification of a $k$-scheme of finite type).
\een
\end{defn}

Now, suppose that $X$ is an arbitrary Hausdorff strictly $k$-dagger analytic space. We want to associate to ita rigid $k$-daggger analytic space by defining on the points $X_0 \subset X$ a $G$-topology. Therefore, we say that a subset $U \subset X_0$ is \emph{admissible open} if, for any strictly dagger affinoid domain $V \subset X$, 
the intersection $U \cap V_0$ is an admissible open subset in the rigid $k$-dagger affinoid space $V_0$. In the same way, a covering $\{ V_i \}_{i \in I}$
of an admissible open subset $U \subset X_0$ by admissible open subsets is said \emph{admissible} if, for any 
strictly dagger affinoid domain $V \subset X$, $\{ V_i \cap V_0 \}_{i \in I}$ is an admissible open covering of $V_0 \cap U$. 
In this way we get a $G$-topology on the set $X_0$ and sheaves of rings $\cO_{V_0}$, with $V_0 = V \cap X_0$, where $V$ runs 
through the strictly dagger affinoid domains in $X$ and these sheaves are compatible on intersections. So these sheaves glue together to give a sheaf of rings $\cO_{X_0}$ on the 
$G$-topological space $X_0$ just defined. The locally $G$-ringed space $(X_0, \cO_{X_0})$ constructed in this way is thus a quasiseparated rigid $k$-dagger space.

\begin{defn}
A collection of subsets of a set is said to be \emph{of finite type} if each subset of the collection meets only a finite number of other subsets of the collection.
\end{defn}

The next theorem is the dagger analogous of theorem 1.6.1 of \cite{BER4}.

\begin{thm} \label{thm:dagger_ber_to_dagger_rigid}
The correspondence $(X, A, \tau) \mapsto (X_0, \cO_{X_0})$ is a fully faithful functor from the 
category of Hausdorff strictly $k$-dagger analytic spaces to the category of quasiseparated rigid $k$-dagger spaces. 
This functor induces an equivalence between, the category of paracompact strictly 
$k$-analytic spaces and the category of quasiseparated rigid $k$-dagger spaces that have an admissible dagger affinoid covering of finite type. 
\end{thm}
\begin{proof}
Let $X$ be a Hausdorff strictly $k$-dagger analytic space. We need the following dagger version lemma 1.6.2 of \cite{BER4}.

\begin{lemma}
\ben
\item Any open dagger affinoid domain in the rigid $k$-dagger space $X_0$ associated to $X$ is of the form $V_0$ for a $V \subset X$, a strictly dagger affinoid domain in $X$. 
\item A system $\{ V_i \}_{i \in I}$ of strictly dagger affinoid domains in $X$ induces an admissible covering
      of $X_0$ if and only if each point of $X$ has a neighborhood of the form $V_{i_1} \cup \ldots \cup V_{i_n}$ (\ie, if and only if $\{V_i\}_{i \in I}$ is a quasinet on $X$). 
\een
\end{lemma}
\begin{proof}
\ben 
\item Let $f: U_0 \to X_0$ be an open $k$-dagger affinoid domain in $X_0$ and $V \subset X$ a strictly $k$-dagger affinoid domain in $X$.
     Then, $f^{-1}(V_0)$ is a finite union of dagger affinoid domains in $U_0$ and the family $\{f^{-1}(V_0)\}$, where $V$ runs through strictly $k$-dagger affinoid domains in $X$, 
     is a dagger affinoid covering of $U_0$. It follows that we can find strictly $k$-dagger affinoid domains $U_1, \ldots, U_n \subset U$ and $V_1, \ldots, V_n \subset X$ 
     such that $U = U_1 \cup \ldots \cup U_n$ and $f|_{(U_i)_0}$ is induced by a morphism of strictly $k$-dagger affinoid spaces $\phi_i: U_i \to V_i$, that identifies 
     $U_i$ with a $k$-dagger affinoid domain in $V_i$. All morphisms $\phi_i$ are compatible on intersections, therefore, they induces a morphism of strictly $k$-dagger
     analytic spaces $\phi: U \to X$. Since $\phi$, as a map of topological spaces, is proper and induce an injection on the everywhere dense subset $U_0 \subset U$, 
     it follows that $\phi$ induces a homeomorphism of $U$ with its image in $X$. Finally, $\phi$ identifies $U_i$, with a strictly $k$-dagger affinoid
     domain in $V_i$ and therefore $\phi$ identifies $U$ with a strictly $k$-dagger affinoid domain in $X$.
\item Suppose first that $\{(V_i)_0\}_{i \in I}$ is an admissible covering of $X_0$ for some family of strictly $k$-dagger affinoid domain $V_i \subset X$. 
      This means that, for any strictly $k$-dagger affinoid domain $V \subset X$ the family
      \[ \{ V_0 \cap (V_i)_0 \}_{i \in I} \]
      is an admissible covering of $V_0$. It follows that $V$ is contained in a finite union of the form $V_{i_1} \cup \ldots \cup V_{i_n}$, and so it follows
      that each point of $X$ has a neighborhood made of a finite unions of elements of $\{V_i\}_{i \in I}$, because $V$ is strictly $k$-dagger affinoid. 
      
      Conversely, assume that the latter property is true. Then any strictly $k$-dagger affinoid domain is contained in a finite 
      union $V_{i_1} \cup \ldots \cup V_{i_n}$ and therefore $\{ (V_i)_0 \cap V_0 \}_{i \in I}$ is an admissible covering of $V_0$.
\een
\end{proof}

To check that the functor $(X, A, \tau) \mapsto (X_0, \cO_{X_0})$ is fully faithful one only needs to use the previous lemma and some properties 
of the rigid $G$-topology which are true also for rigid dagger $G$-topology we used to define $(X_0, \cO_{X_0})$. So, we omit this part of the proof that is identical to the corresponding part of theorem 1.6.1 of \cite{BER4}. Instead, we give the details of the proof of the equivalence between paracompact Hausdorff strictly $k$-dagger spaces and quasiseparated rigid $k$-dagger spaces.

If $X$ is a paracompact strictly $k$-dagger analytic space, then it has a strictly $k$-dagger affinoid 
atlas with a locally finite net, and therefore the rigid $k$-dagger space $X_0$ has an admissible affinoid covering of finite type. 
It is also quasiseparated because of the first point of previous lemma and lemma \ref{prop_lemma_ber_nets}. Conversely, 
let $\cX$ be a quasiseparated rigid $k$-dagger space that has an admissible affinoid covering $\{U_i\}$ of finite type. 
We fix $k$-dagger affinoid spaces $\cM(A_i) = V_i$ such that $(V_i)_0 \cong U_i$, where the $A_i$ are strictly $k$-dagger affinoid algebras.
Since $\cX$ is quasiseparated, for any pair $i,j \in I$ the intersection $U_i \cap U_j$ is a finite union of open affinoid domains in $\cX$. Thus, there are strictly dagger special domains $V_{i,j} \subset V_i$ and $V_{i,j} \subset V_j$ that correspond to $U_i \cap U_j$, under the identifications $U_i \mapsto (V_i)_0$. 
Let $\nu_{i,j}: V_i \to V_j$ denotes the induced isomorphism. It is clear that this system of maps satisfies the glueing properties needed for applying the second part of proposition \ref{prop_glueing} which gives a paracompact strictly $k$-dagger analytic space $X$ which, by its construction, is such that $X_0 \cong \cX$.
\end{proof}

Thanks to work of Grosse-Kl\"onne, it is not difficult to describe a functor from the category of rigid $k$-dagger spaces to the category of rigid $k$-spaces. This functor can be described starting from dagger affinoid algebras in the following way. If $A$ is a strict $k$-dagger affinoid algebra we can equip $A$ canonically with a norm which is intrinsic to $A$. This norm is the norm used by Grosse-Kl\"onne in \cite{GK} and can be canonically attached to $A$ from its dagger bornology in the following
way. Let $P = \{ p_{\rho} \}_{\rho > 1}$ be the family of norms on $A$ which define the dagger bornology of $A$. This family is a chain in the lattice of all seminorms
over $A$, so we can equip $A$ with the infimum of this family (which is not an element of $P$ because otherwise the bornology on $A$ would be given by a norm).

By the results in the first section of \cite{GK}, the completion of $A$ with respect to this norm gives a $k$-Banach algebra $\what{A}$ which is an affinoid algebra and if $f: A \to B$ is a morphism of 
strict $k$-dagger algebra, then $f$ induces a morphism $\what{f}: \what{A} \to \what{B}$ between the strict $k$-affinoid algebras. In \cite{GK} the association $A \mapsto \widehat{A}$ is studied only for the case of strictly (dagger) affinoid algebras, but this association extends obviously to the non-strict case. 

\begin{prop}
Consider the functor $\what{\cdot}: \bAff_k^\dagger \to \bAff_k$, from the category of $k$-dagger affinoid algebras to the category of $k$-affinoid algebras. Then, $\cM(A) \cong \cM(\what{A})$ and if $A$ is strictly $k$-dagger affinoid there is a canonical map of $G$-topological spaces $A_0 = \Max(\what{A}) \to \Max(A)$ is a isomorphism of $G$-topological spaces.
\end{prop}
\begin{proof}
The homeomorphism $\cM(A) \cong \cM(\what{A})$ is an easy consequence of proposition \ref{prop:seminorm_dagger_nondagger_spectrum}. Then, if $A$ is striclty $k$-affinoid, the isomorphism $\Max(\what{A}) \to \Max(A)$ one can show that every rational affinoid subdomain of $\what{A}$ can be defined by polynomials with coefficients in $\what{A}$ which are overconvergent analytic functions, therefore the $G$-topologies of $\Max(\what{A})$ and $\Max(A)$ coincide.
\end{proof}

We can extend the functor $\what{\cdot}: \bAff_k^\dagger \to \bAff_k$ to a global functor $\what{\cdot}:  \bRig_k^\dagger \to \bRig_k$ defining it locally. This functor has the following properties.

\begin{thm}
The functor $\what{\cdot}: \bRig_k^\dagger \to \bRig_k$ has the following properties
\ben
\item $\what{\cdot}$ induces a homeomorphism of underlying $G$-topological space $X \cong \widehat{X}$ and an isomorphisms between the local rings of the structure sheaves (notice that here only rigid points are considered);
\item $X$ is connected (resp. normal, resp. reduced, resp. smooth) if and only if $\what{X}$ is;
\item $f : X \to Y$ is a closed immersion (resp. an open immersion, resp. a locally closed immersion, resp. an isomorphism, resp. quasicompact, resp. separated) if and only if
      $\what{f} : \what{X} \to \what{Y}$ is;
\item if $f: X \to Y$ is finite, then $\what{f}$ is finite;
\item if $X$ is a quasialgebraic rigid space, there exists an admissible affinoid covering $X = \bigcup U_i$ and affinoid dagger spaces $V_i$ and isomorphisms $\what{V}_i \cong U_i$.
\een
\end{thm}
\begin{proof}
We refer to \cite{GK} for the proofs of the assertions of this theorem.
\end{proof}

\begin{defn}
Let $X$ be a rigid $k$-dagger space. We say that an $\cO_X$-module $\sF$ is \emph{coherent} if there exists an admissible affinoid covering $\{U_i\}$ of $X$ such that $\sF|_{U_i}$ is
associated to a finite $\cO_X(U_i)$-module.
\end{defn}

The functor $\what{\cdot}$ naturally extends to a functor $\what{\cdot}: \bCoh(\cO_{X}) \to \bCoh(\cO_{\what{X}})$ in the following way. Given a dagger affinoid space
$X = \Max(A)$ and an $A$-module $M$ we get functorially a $\what{A}$-module by considering the module
\[ \what{M} \doteq M \otimes_A \what{A}. \]
Clearly if $M$ is a finite $A$-module then $\what{M}$ is a finite $\what{A}$-module, so this gives a functor $\bCoh(\cO_{X}) \to \bCoh(\cO_{\what{X}})$ when $X$ is
dagger affinoid space because of the Kiehl theorem (cf. theorem \ref{thm:kiehl}). So, given a coherent $\cO_X$-module $\sF$, for a general rigid $k$-dagger space, then if $\{U_i\}$ is covering such that $\sF|_{U_i}$ is associated to a finite $\cO_X(U_i)$-module the data of the sheaves $\what{\sF}|_{\what{U}_i}$ associated to finite
$\cO_{\what{X}}(\what{U}_i)$-modules defines a coherent $\cO_{\what{X}}$-module.

\index{partially proper rigid dagger space}
\begin{defn} \label{defn:proper}
A morphism $f : X \to Y$ of rigid $k$-dagger spaces is called \emph{partially proper} if $f$ is separated and there exist admissible affinoid coverings $Y = \bigcup Y_i$ 
and $f^{-1}(Y_i) = \bigcup X_{i,j}$ (all $i$), such that for every $X_{i,j}$ there is an admissible open affinoid subset $\tilde{X}_{i,j} \subset f^{-1}(Y_i)$ with $X_{i,j} \rhook \tilde{X}_{i,j}$ inner  with respect to $\tilde{X}_{i,j} \rhook f^{-1}(Y_i)$.

A morphism $f: X \to Y$ is said \emph{proper} if is partially proper and quasi-compact.
\end{defn}

\begin{rmk}
In definition \ref{defn:proper} we used Berkovich's terminology of inner morphisms but our definition is equivalent to definition 2.24 of \cite{GK}.
\end{rmk}

\index{Stein rigid dagger space}
\begin{defn}
Let $X$ be a rigid $k$-dagger space. We says that $X$ is \emph{Stein} if it admits a dagger affinoid exhaustion \ie a covering $\{ U_i \}_{i \in \N}$ such that
\[ X = \bigcup_{i \in \N} U_i \]
and $U_i \subset \Int(U_{i + 1})$ for all $i$ such that each $U_i \rhook U_{i + 1}$ is a Weierstrass subdomain embedding.
\end{defn}

Clearly every Stein rigid $k$-dagger space is partially proper. 

\begin{thm}
If $A$ is a partially proper rigid $k$-dagger spaces then $\what{\cdot}: \bCoh(X) \to \bCoh(\what{X})$ is an equivalence and for every $\sF \in \bCoh(X)$ we have
\[ \sF(X) \cong \what{\sF}(\what{X}). \]
\end{thm}
\begin{proof}
We address first the case when $X = \underset{i \in \N}\bigcup U_i$ is Stein. Let $\sG$ be a coherent $\cO_{\what{X}}$-module. Since $U_i \subset \Int(U_{i + 1})$ then
$\sG$ defines a $\cO_X(U_i)$-module that we denote with $\sF(U_i)$, for any $i$. By Kiehl's theorem we get coherent $\cO_X|_{U_i}$-modules which by definition glues 
to give an $\cO_X$-module $\sF$. So, $\what{\sF}$ is a coherent $\cO_{\what{X}}$-module and $\what{\sF} = \sG$. Moreover $\sF(X) = \sG(\what{X})$ because
\[ \sF(X) = \limpro_{i \in \N} \cO_X(U_i) = \limpro_{i \in \N} \cO_{\what{X}}(\what{U}_i) = \sG(\what{X}) \]
which is due to the fact that the projective systems are cofinal one in each other, because of the condition $U_i \subset \Int(U_{i + 1})$. 

For the general case it is enough to show that any partially proper rigid $k$-dagger space has an admissible open covering of the form $X = \underset{j \in J}\bigcup S_j$ where $S_j$ are Stein spaces and their
finite intersection are Stein. And such a covering can be constructed from a dagger affinoid covering $X = \underset{i \in I}\bigcup V_i$ given by definition \ref{defn:proper}. So, for any admissible open dagger affinoid domain $Z_j \subset X$ we can find a $Z_j \subset \Int(V_i)$, as explained in the last part of the proof of theorem 2.26 of \cite{GK}, from which we can construct the desired Stein covering.
\end{proof}

\begin{thm} \label{thm:equiv_part_prop}
The functor $\what{\cdot}$ induces an equivalence between the category of partially proper rigid $k$-dagger spaces and the category partially proper rigid $k$-analytic spaces.
\end{thm}
\begin{proof}
We refer to theorem 2.27 of \cite{GK} for the proof of this result.
\end{proof}

\begin{exa}
The prototype of Stein rigid $k$-(dagger) space is the affine space, which is defined as the direct limit
\[ (\A^n_k)^\dagger \doteq \limind_{\rho \to \infty} \Max(W_k^n(\rho)) \]
in the dagger case and
\[ \A^n_k \doteq \limind_{\rho \to \infty} \Max(T_k^n(\rho)) \]
in the rigid case. And applying the functor $\what{\cdot}$ we get $\what{(\A^n_k)^\dagger} \cong \A^n_k$.
\end{exa}

\begin{prop}
Let $X$ be a partially proper rigid $k$-dagger space and $\sF$ a coherent $\cO_X$-module then
\[ H^i(X, \sF) \cong H^i(\what{X}, \what{\sF}). \]
\end{prop}
\begin{proof}
\cite{GK} theorem 3.2.
\end{proof}

Then, we see the relations between our $k$-dagger analytic spaces and Berkovich $k$-analytic spaces.

\index{wide affinoid subdomain}
\begin{defn}
Let $X$ be a $k$-analytic space in the sense of Berkovich. An affinoid subdomain $U \subset X$ is said \emph{wide} if $U \subset \Int(X)$.
\end{defn}

\index{wide Berkovich $k$-analytic space}
\begin{defn}
Let $X$ be a $k$-analytic space in the sense of Berkovich. We say that $X$ is \emph{wide} if there exists a net of wide affinoid subdomains $U_i \subset X$ such that
\[ \bigcup_{i \in I} U_i = X. \]
\end{defn}

\begin{rmk}
The wideness condition is equivalent to the closeness condition of definition \ref{defn:closed_morphism}. Namely, $X$ is wide if and only if
\[ \partial(X) = \emptyset. \]  
\end{rmk}

In particular, Stein spaces (defined in Berkovich geometry) are wide.

\begin{thm} \label{thm:dagger_comparison}
There is an equivalence of categories between the category of wide $k$-analytic spaces in the sense of Berkovich and the wide $k$-dagger analytic space.
\end{thm}
\begin{proof}
Working locally we can define a functor $\what{\cdot}: \bAn_k^\dagger \to \bAn_k$. Using the same ideas of the proof of theorem \ref{thm:equiv_part_prop}, it is easy to see that the functor $\what{\cdot}$ admits a quasi-inverse when the spaces are without boundary.
\end{proof}

\begin{cor}
There is an equivalence between, the category of wide paracompact strictly $k$-dagger analytic spaces and the category of partially proper quasi-separated 
rigid $k$-dagger spaces that have an admissible dagger affinoid covering of finite type. 
\end{cor}
\begin{proof}
It is enough to combine the equivalence given by theorem \ref{thm:dagger_comparison} with the one given by \ref{thm:dagger_ber_to_dagger_rigid}.
\end{proof}

\subsection{Complex analytic case}

In this section we show how the category of classical complex analytic spaces embeds fully faithfully in the category of $\C$-dagger analytic spaces.

\begin{lemma} \label{lemma:dagger_covering}
Let $\cX$ be a separated complex analytic space in the classical sense. Any $x \in \cX$ has a neighborhood basis made of compact Stein subsets whose algebra of germs of analytic functions is a $\C$-dagger affinoid algebra.
\end{lemma}
\begin{proof}
By definition each $x \in \cX$ has an open neighborhood $x \in U_x$ such that $(U_x, \cO_{\cX}|_{U_x})$ is isomorphic to an affine analytic set and by the separateness hypothesis $\cX$ is Hausdorff. Choosing an explicit isomorphism
\[ U_x \cong \{ z \in U \subset \C^n | f_1(z) = \ldots = f_r(z) = 0, f_1, \ldots,f_r \in \cO_{\C^n}(U) \} = \tilde{U}_x \]
and denoting by $\tilde{z}$ the image of $x$ in $\tilde{U}_x$ we can choose a basis of closed polydisks $\{ D_i \}_{i \in \N}$ centered in $\tilde{z}$ as base of neighborhoods of $\tilde{z}$ in $\C^n$. Then,
the subsets $\{D_i \cap \tilde{U}_x \}_{i \in \N}$ form a basis of neighborhoods of $\tilde{z}$ in $\tilde{U}_x$ and they are canonically $\C$-dagger affinoid spaces that form the required system of neighborhoods.
\end{proof}

Let $\cX$ be a separated complex analytic space in the classical sense. By the previous lemma it is clear that the family of $\C$-dagger compact Stein subsets of $\cX$ endow $\cX$ with a structure of a good $\C$-dagger analytic space.

\begin{defn}
Let $\cX$ be a separated complex analytic space, we denote by $\cX^\dagger$ its associated good $\C$-dagger analytic space.
\end{defn}

\begin{lemma} \label{lemma:complex_functorial}
The association $\cX \mapsto \cX^\dagger$ defines a functor from the category of complex analytic spaces to the category of $\C$-dagger analytic spaces.
\end{lemma}
\begin{proof}
Let $f: \cX \to \cY$ be a morphism of separated complex analytic spaces
and $U \subset \cX^\dagger$ a $\C$-dagger affinoid domain. Since $f(U)$ is compact then, by lemma \ref{lemma:dagger_covering}, it can be covered by a finite number of compact subsets $U_1, \ldots, U_n \subset \cY$ each of which is an analytic subset of some closed polydisk of $\C^{m_i}$. So, each $f(U) \cap U_i$ in a compact subset of an analytic subset of a closed polydisk of $\C^{m_i}$. Clearly, $U \cap f^{-1}(U_i)$ is a dagger affinoid subdomain of $\cX^\dagger$ such that $f(U \cap f^{-1}(U_i)) \subset U_i$. Moreover, the family of all such $U_i$'s for $U \in \tau_{\cX^\dagger}$ ( where $\tau_{\cX^\dagger}$ is the Berkovich net of $\C$-dagger affinoid subdomains of $\cX^\dagger$) is a net on $\cY$ which belongs in the equivalence class of the net of $\cY^\dagger$. Therefore, the collection of data $f|_{U \cap f^{-1}(U_i)}: U \cap f^{-1}(U_i) \to U_i$ for all $U \in \tau_{\cX^\dagger}$ defines a morphism of $\C$-dagger analytic spaces $f: \cX^\dagger \to \cY^\dagger$.
\end{proof}

\begin{lemma} \label{lemma:complex_open_subsets}
Let $\cX$ be a complex analytic space. The $\C$-dagger analytic and the complex analytic structural sheaves agree on open subsets, \ie if $U \subset \cX$ is an open set, then
\[ \cO_\cX(U) \cong \cO_{\cX^\dagger}(U). \] 
\end{lemma}
\begin{proof}
In \ref{defn:structural_sheaf} we defined 
\[ \cO_{\cX^\dagger}(U) = \Hom_{\bAn_\C^\dagger}(U, \A_\C^1) \]
for any analytic domain of $\cX^\dagger$. Without lost of generality we can suppose that $U$ is connected. Since $\cX$ has a countable base for the topology and it is locally compact, then $U$ is hemi-compact and it can be exhausted by compact subsets. Moreover, we can suppose that $U$ is Stein because complex analytic spaces can be covered by Stein subspaces. Therefore, we can suppose that
\[ U = \bigcup_{i \in \N} U_i \]
where $U_i$ are $\C$-dagger affinoid spaces and $U_i \rhook U_{i + 1}$ is a dagger affinoid embedding such that $U_i \subset \Int(U_{i + 1})$. Reasoning like proposition \ref{prop_good_open_subsets} one can directly verify that
\[ \Hom_{\bAn_\C^\dagger}(U, \A_\C^1) \cong \limpro_{i \in \N} A_{U_i} \]
which is isomorphic to the Stein algebra associated to $U$, therefore
\[ \Hom_{\bAn_\C^\dagger}(U, \A_\C^1) \cong \Hom_{\bAn_\C}(U, \A_\C^1) \cong \cO_\cX(U). \]
\end{proof}

\begin{thm} \label{thm:complex_embed}
The category of separated analytic spaces over $\C$ embeds fully faithfully in the category of good $\C$-dagger analytic spaces.
\end{thm}
\begin{proof}
The fact that $\cX^\dagger$ is a good $\C$-dagger analytic space for every complex analytic space $\cX$ is an immediate consequence of lemma \ref{lemma:dagger_covering}. So, by lemma \ref{lemma:complex_functorial}, only the fully faithfully part of the statement of the theorem remains to be proven. 

It is clear that a morphism $f: \cX^\dagger \to \cY^\dagger$ induces morphism of complex analytic spaces $\cX \to \cY$ just by restricting it to open subsets of $\cX$ and $\cY$, which are analytic domains. So, the association $\cX \to \cX^\dagger$ is a full functor, it remains to show that it is faithful. Let $f, g: \cX \to \cY$ be two morphisms which define the same morphisms of associated $\C$-dagger analytic spaces. Consider $U \in \tau_{\cX^\dagger}$ and $V \in \tau_{\cY^\dagger}$ such that $f(U) \subset V$. Since the canonical embeddings $\iota_U: U \rhook \cX^\dagger$ and $\iota_V: V \rhook \cY^\dagger$ are isomorphisms on their image, then
\[ A_U \cong \cO_{\cX^\dagger}(U), \ \ A_V \cong \cO_{\cY^\dagger}(V). \]
By lemma \ref{lemma:complex_open_subsets} for all open Stein neighborhoods $U \subset \tilde{U}$ and $V \subset \tilde{V}$ one has that
\[ \cO_{\cX^\dagger}(\tilde{U}) \cong \cO_{\cX}(\tilde{U}), \ \ \cO_{\cY^\dagger}(\tilde{V}) \cong \cO_{\cY}(\tilde{V}), \]
and since $U$ and $V$ are compact Stein subsets of $\cX$ and $\cY$ we have that
\[ A_U \cong \limind_{\tilde{U} \supset U} \cO_{\cX}(\tilde{U}), \ \
A_V \cong \limind_{\tilde{V} \supset V} \cO_{\cX}(\tilde{V}). \]
Hence, the fact that $f$ and $g$ agree on open Stein subsets of $\cX$ readily implies that $f$ and $g$ agree on compact Stein subsets of $\cX$ which is equivalent to say that $f$ and $g$ induces the same morphism of $\C$-dagger analytic spaces.
\end{proof}

We end by giving some examples of the spaces we have defined so far.

\begin{exa}
\ben
\item The dagger real affine line $\A_\R^1$ is homeomorphic to the upper half plane of $\C$ with its boundary which agrees with the real line $\R \subset \C$. More precisely, $\A_\R^1$ is homeomorphic to $\C$ quotiented by
      the action of the complex conjugation. Analogous statements are true for $\A_\R^n$.
\item In \cite{BER2}, remark 1.5.5, Berkovich suggested an application of its theory of analytic spaces for the archimedean base fields. But, using his theory he was able
      to show that one could obtain a topological manifold, that is something which does not belong to analytic geometry. We can now see how our approach fixes this problem.
      Berkovich considered for $k = \C$ and $n = 2$ the ``archimedean affinoid space''
      \[ \cM(D^n) \cong \{ (z_1, z_2) \in \C | |z_1| \le 1, |z_2| \le 1 \}  \]
      where $D^n$ is the disk algebra of the unitary polydisk in $\C^n$. The subset defined by the equation $X_1 X_2 - 1 = 0$, where $X_1$ and $X_2$ are coordinate is given by
      \[ X = \cM \l ( \frac{D^n}{(X_1 X_2 - 1)} \r) \cong \{ (z_1, z_2) \in \C | |z_1| = 1, |z_2| = 1 \}.  \]
      The elements of $D^n$ define continuous on the border and moreover any continuous function on the border of $\cM(D^n)$ ``induces'' a unique analytic function inside the polydisk which is an element of $D^n$. Hence, the structural sheaf by which $X$ is endowed in this way is the sheaf of continuous function from $X \cong S^1 \to \C$, which is not a space
      that one would like to study in analytic geometry. On the other hand, the same space considered as a $\C$-dagger affinoid space gives
      \[ X^\dagger = \cM \l ( \frac{W_\C^n}{(X_1 X_2 - 1)} \r) \cong \{ (z_1, z_2) \in \C | |z_1| = 1, |z_2| = 1 \}.  \]
      So for the underlying topological spaces we have the homeomorphism $|X^\dagger| \cong |X|$ but on $X^\dagger$ we consider functions that can be extended to analytic functions on a neighborhood of $X^\dagger$ in the hypersurface of 
      $\C^2$ given by the equation $X_1 X_2 - 1 = 0$. In this way, one can see easily that $X^\dagger$ is isomorphic to the circle of the complex plane equipped with the
      algebra of functions of the form
      \[  \frac{W_\C^n}{(X_1 X_2 - 1)} \cong \l \{ \sum_{i \in \Z} a_i X^i, \sum_{i \in \Z} |a| < \infty \r \}, \]
      which is a compact Stein subset of $\C$ and hence it is a meaningful object to study in analytic geometry.
\item Building on last example, one can easily see that (using proposition \ref{prop_glueing}) that it is possible to construct the projective line over $\C$ by glueing two copies of the unital polydisk of $\C$ on their boundaries. This allows to obtain projective geometry directly from compact Stein space and to use only Noetherian algebras in the process of constructing projective analytic space.
\een
\end{exa}

\section{Flatness for $k$-dagger analytic spaces}

In this last section, we want to discuss flatness for morphisms of $k$-dagger analytic spaces. We see how the natural concept of flatness behaves better and simpler in dagger analytic geometry and some pathologies of Berkovich geometry simply disappear. Indeed, these pathological behaviours are induced by non-overconvergent analytic functions, which, from the point of view discussed in this work, should not be considered analytic functions.

\index{flat morphism of dagger analytic spaces}
\begin{defn} \label{defn:flat_morphism}
A morphism $f: X \to Y$ of $k$-dagger analytic spaces is said to be \emph{flat at $x \in X$} if the morphism $f_x: \cO_{Y, f(x)} \to \cO_{X, x}$ is flat. We say that $f$ is \emph{flat} if it is flat at all $x \in X$.
\end{defn}

In Berkovich geometry there are two non-equivalent definitions of flat map between $k$-analytic spaces, that we borrow from \cite{DUC}.

\index{naively flat morphism Berkovich analytic spaces}
\begin{defn} \label{defn:naive_flat_morphism}
A morphism $f: X \to Y$ of $k$-analytic Berkovich spaces is said to be \emph{naively flat at $x \in X$} if the morphism $f_x: \cO_{Y, f(x)} \to \cO_{X, x}$ is flat. We say that $f$ is \emph{naively flat} if it is naively flat at all $x \in X$.
\end{defn}

The problem of this definition of flatness in Berkovich geometry is that it is not stable under base changes. But this failure can happen only at boundary points and they are given by the existence of non-overconvergent analytic functions, as explained in examples 2.19-2.23 of \cite{DUC}. Therefore, it is natural to expect that these problems do not arise in dagger analytic geometry. The next definition is the substitute to flatness proposed by Ducros for solving the problems of naive flatness in Berkovich geometry.

\index{universally flat morphism Berkovich analytic spaces}
\begin{defn} \label{defn:ducros_flat_morphism}
A morphism $f: X \to Y$ of $k$-analytic Berkovich spaces is said to be \emph{(universally) flat at $x \in X$} if:
\ben
\item when both $X$ and $Y$ are good, given any $g: Y' \to Y$ morphism of good $k$-anaytic spaces and $x' \in X \times_{Y} Y'$  is such that $p(x') = x$, where $p: X \times_{Y} Y' \to X$ is the canonical morphism, then $g^* \cO_{Y, f'(x')}$ is a flat $\cO_{X \times_{Y} Y', x'}$-module, for every $x'$;
\item for general $X$ and $Y$, if there exist $U \subset X$ and $Y \subset Y$ such that $U$ and $V$ are good analytic domains such that $f(U) \subset V$ and $f$ is (universally) flat at $x$. 
\een
\end{defn}

Clearly the definition of universally flat morphism can be given also for $k$-dagger analytic spaces. We do not give the details of such a definition because it is not important for us and it can be easily worked out from definition \ref{defn:ducros_flat_morphism}. The next lemma shows that the only problems that the notion of naively flat morphisms have in Berkovich geometry are placed on the boundary.

\begin{lemma} 
Let $f: X \to Y$ be a morphism of $k$-analytic Berkovich spaces and let $x \in \Int(X/Y)$. Then, $f$ is flat at $x$ if and only if it is naively flat.
\end{lemma}
\begin{proof}
See theorem 6.5.2.3 of \cite{DUC}.
\end{proof}

Finally, the dagger analytic improvement of the above results.

\begin{thm} \label{thm:flatness}
Let $f: X \to Y$ be a morphism of $k$-dagger analytic spaces. Then, $f$ is universally flat if and only if it is naively flat.
\end{thm}
\begin{proof}
We proved in proposition \ref{prop:embed_open} that dagger affinoid subdomain embeddings are open immersions of dagger affinoid spaces, in the sense of definition \ref{defn:open_immersion}. This means that for every point in the bornological spectrum of a dagger affinoid space a dagger affinoid subdomain embedding induces an isomorphism on the stalk of the structural sheaves of the spaces involved. This property readily implies the same property for $k$-dagger analytic spaces as they are built glueing dagger affinoid spaces. Notice that, as we already remarked so far, this property fails for boundary points of Berkovich spaces.

Therefore, since the local rings of points on the boundary of dagger affinoid spaces are not special, as it is the case in Berkovich geometry, they have the same properties of inner points. In particular, they always satisfy property 6.5.1 of \cite{DUC} which implies that universally flatness is equivalent to naively flatness.
\end{proof}

\chapter{Addendum: applications}

In this work we have not discussed very much the possible applications of the theory developed so far. This is due partially because we think that the theory is interesting per se and partially because this work is already long enough. We devote this small ending chapter to describe some applications that our results have found, which go beyond the scope of this work. 

In \cite{BEKR} Kreminzer and Ben-Bassat describe a formal framework in which one can interpret Berkovich geometry as a kind of relative algebraic geometry (in the sense of T\"oen-Vezzosi) on the quasi-abelian category of Banach spaces over $k$, for a complete non-archimedean field $k$. One of the key ideas of that work is to use the homological methods of the theory of quasi-abelian categories developed by Schneiders in \cite{SN} for giving an algebraic (more precisely derived algebraic) characterization of affinoid subdomain embeddings. This characterization is a key step toward the foundation of a theory of derived analytic spaces using the approach homotopical algebraic geometry of T\"oen-Vezzosi. 

In \cite{BABE} the author and Ben-Bassat show that the methods of dagger analytic geometry can be used for extending the results of \cite{BEKR} to encompass complex analytic geometry. The ideas of the present work are very suitable for this scope for the following reasons:
\ben
\item adapting the proofs of \cite{BEKR} to the classical complex analytic case is not a straightforward task because they heavily rely on the theory of affinoid spaces;
\item there is no direct parallel of the theory of affinoid spaces over $\C$, \ie using $\C$-Banach algebras, as we explained so far;
\item the category of complete bornological vector spaces is canonically a closed symmetric monoidal category whose monoidal functor is the completed projective tensor product we descrbed so far; this is in contrast with the category of locally convex spaces which has not such a structure;
\item in the main proofs of \cite{BEKR} the Gerritzen-Grauert theorem is used in a crucial way and the main result of the present work is the generalization of this theorem for $\C$-dagger affinoid spaces.
\een
These reasons render the extension of the results of \cite{BEKR} very natural in the context of dagger analytic geometry. Then, theorem \ref{thm:complex_embed} implies that classical complex spaces fit into this theory.

Another advantage of the dagger affinoid approach is that it seems be more naturally extendible over general Banach rings with respect to the affinoid approach. In \cite{Poi} Poineau proposed a notion of global analytic space with the main motivation of studying the case when the base Banach ring is $(\Z, |\cdot|_\infty)$. The proposed definition of global analytic space is the one suggested by Berkovich in the first chapter of \cite{BER2} and it uses local models that are akin to classical complex spaces instead of affinoid spaces. Also in this case one can ask for an affinoid approach in the aim of a deeper understanding of global analytic spaces. One of the main problems that one faces for developing such a theory is the same one that is met over $\C$: Using Banach algebras of convergent power-series does not lead to the desired results. We notice also that one would like to have a theory over $(\Z, |\cdot|_\infty)$ that gives a meaningful theory over $\C$ by base change, but since no theory of affinoid spaces exists over $\C$ it turns out that it cannot exist neither over $\Z$ (or if it exists it is bad behaving for base changes). It turns out that overconvergent power-series rings can be used for solving this problem, see for example \cite{Pau2}. We remark that over a general Banach ring the theory of bornological modules is not meaningful as the theory of bornological vector spaces is over a non-trivially valued field. Therefore, the right framework for dealing with overconvergent algebras of power-series over Banach rings is to use the ind-category of the category of Banach modules, see the first part of \cite{BABE} for a detailed study of this category. For a more detailed description of these issues we refer the reader to section 6 of \cite{BABE} where it is also explained how to simplify some technical issues of \cite{Pau2}.

Our last remark builds on the concluding remarks of chapter 5. There we have proved that there exist compact Stein subsets which are not isomorphic to $\C$-dagger affinoid spaces and we claimed that compact Stein subsets of $\C^n$ with an associated Noetherian algebra of germs of analytic functions are precisely the ones which are endowed with a structure of $\C$-dagger anlaytic spaces. This result follows from the classical theorem of Siu (cf. \ref{thm:siu}) and it will be studied in full details in the future work \cite{BA2}. We remarked also that this result unifies and simplifies the understanding of the notions of ``compact Steinness" present in the literature about complex geometry and non-archimedean geometry. Therefore, we have reasons to think that the methods proposed in this work can be very helpful in the study of the geometry of complex analytic and non-archimedean spaces.

\appendix

\chapter{The pro-site} \label{pro_appendix}

This appendix describes a construction that helps to study the relations between the classical analytic spaces and the ones we defined in this work. Since our approach uses overconvergent analytic functions, it is natural to think to dagger analytic spaces as pro-analytic spaces, \ie as objects in the pro-category of the category of analytic spaces, $\bAn_k$. Since $\bAn_k$ has a structure of a site one could ask if this structure induces a natural structure of site on its pro-category. We devote this appendix to work out such a construction in full generality, namely for any site.

\section{Pro-objects and Ind-objects}

Here we recall the definitions of pro and ind objects as given in the first expos\'e of \cite{SGA4}. 
In this section $\cC$ denotes a category and $\widehat{\cC}$ the category of presheaves of sets over $\cC$, \ie the objects of $\widehat{\cC}$ 
are functors $F: \cC^\circ \to \cS et$ and morphisms are natural transformations between functors. To avoid set-theoretical problems the reader can think that $\cC$ is a small category, but we neglect this issue that can always be avoided fixing suitable Grothendieck universes.

\index{filtered category}
\begin{defn}
A category $\cC$ is said \emph{pseudo-filtered} if
\ben
\item all diagram of the form
\[
\begin{tikzpicture}
\matrix(m)[matrix of math nodes,
row sep=2.6em, column sep=2.8em,
text height=1.5ex, text depth=0.25ex]
{  & j \\
 i  \\
   & k  \\};
\path[->,font=\scriptsize]
(m-2-1) edge node[auto] {} (m-1-2);
\path[->,font=\scriptsize]
(m-2-1) edge node[auto] {} (m-3-2);
\end{tikzpicture}
\]
can be inserted to a commutative diagram of the form
\[
\begin{tikzpicture}
\matrix(m)[matrix of math nodes,
row sep=2.6em, column sep=2.8em,
text height=1.5ex, text depth=0.25ex]
{  & j \\
 i &   & l\\
   & k  \\};
\path[->,font=\scriptsize]
(m-2-1) edge node[auto] {} (m-1-2);
\path[->,font=\scriptsize]
(m-2-1) edge node[auto] {} (m-3-2);
\path[->,font=\scriptsize]
(m-1-2) edge node[auto] {} (m-2-3);
\path[->,font=\scriptsize]
(m-3-2) edge node[auto] {} (m-2-3);
\end{tikzpicture};
\]
\item all diagram of the form $f,g: i \rightrightarrows j$ can be inserted in a diagram of the form $i \rightrightarrows j \stackrel{h}{\to} k$ with $h \circ f = h \circ g$.
\een
A category is said \emph{filtered} if is pseudo-filtered, non-empty and connected, \ie if any couple of objects can be connected by a sequence of arrows 
(without imposing any condition on the direction of the arrows). A category $\cC$ is said \emph{cofiltered} if $\cC^\circ$ is filtered.
\end{defn}

\index{cofinal functor}
\begin{defn}
Let $\phi: I \to I'$ be a functor. We say that $\phi$ is \emph{cofinal} if for any functor $u: I'^\circ \to \cC$ the canonical natural transformation 
$\limpro u \to \limpro u \circ \phi$ is an isomorphism.

If $\phi: I \to I'$ is an inclusion of categories we say that $I$ is a cofinal subcategory of $I'$.
\end{defn}

Given two objects $i, j \in \ob(I)$ we say that $i$ \emph{is bigger} than $j$ if $\Hom (i, j) \ne \void$.

\begin{prop} 
Let $\phi: I \to I'$ be a functor then:
\ben
\item if $\phi$ is cofinal then for any $j \in \ob(I')$ there exist an object $i \in \ob(I)$ such that $\phi(i)$ is bigger than $j$;
\item if $I$ is filtered, then $\phi$ is cofinal if and only if for any $i \in \ob(I)$ and couple of arrow $i' \stackrel{f,g}{\rightrightarrows} \phi(i)$ in $I'$ there 
     exists an arrow $h: i \to j$ in $I$ such that $\phi(h) \circ f = \phi(h) \circ g$;
\item if $I'$ is filtered and $\phi$ is fully faithful then $\phi$ is cofinal if and only if it satisfies condition (1).
\een
\end{prop}
\begin{proof}
Proposition 8.1.3 of the first expos\'e of \cite{SGA4}.
\end{proof}

\begin{cor} 
Let $\phi: I \to I'$ be a cofinal functor then:
\ben
\item if $I$ is filtered then $I'$ is filtered;
\item if $I'$ is filtered and $\phi$ fully faithful then $I$ is filtered.
\een
\end{cor}
\begin{proof}
Ibid.
\end{proof}

\begin{prop} (Deligne) \\
Let $I$ be a small filtered category, then there exists a partially ordered set $E$ and a cofinal functor $\phi: \cE \to I$, 
where $\cE$ is the category canonically associated to $E$.
\end{prop}
\begin{proof}
Proposition 8.1.6 of the first expos\'e of \cite{SGA4}.
\end{proof}

\index{ind-object}
\begin{defn}
An \emph{ind-object} (or \emph{inductive system}) of $\cC$ is a filtered functor $\phi: I \to \cC$. $I$ is called the \emph{index category} of the system. The category of all ind-objects of $\cC$ with system morphisms is denoted by $\bInd(\cC)$.
\end{defn}

\index{pro-object}
\begin{defn}
A \emph{pro-object} (or \emph{projective system}) of $\cC$ is a cofiltered functor $\phi: I^\circ \to \cC$. $I$ is called the \emph{index category} of the system. The category of all pro-objects of $\cC$ with system morphisms is denoted by $\bPro(\cC)$.
\end{defn}

\begin{prop}
Let $\cC$ be a category, then $\bPro(\cC)$ admits small filtered projective limits. If $\cC$ admits pullbacks then also $\bPro(\cC)$ admits pullbacks. Moreover, there is a canonical fully faithful functor $c: \cC \to \bPro(\cC)$ which commutes with pullbacks but not in general with small filtered projective limits.
\end{prop}
\begin{proof}
Proposition 8.5.1 of the first expos\'e of \cite{SGA4} and proposition 8.9.5 of the same expos\'e.
\end{proof}

\section{Sites}

We recall that a site is a category $\cC$ equipped with a Grothendieck topology and that to any site we can associate its category of sheaves of sets $\wtilde{\cC}$. This category is
defined to be the full subcategory of $\what{\cC}$ identified by objects $\sF \in \ob(\what{\cC})$ such that
\[ \Hom_{\what{\cC}}(X, \sF) \cong \Hom_{\what{\cC}}(R, \sF) \]
for any covering sieve of $\cC$ over any $X \in \ob(\cC)$. If the topology of $\cC$ is given by a pre-topology this is equivalent to the usual exact sequence condition in the
definition of sheaves. We recall also that the inclusion functor $i: \wtilde{\cC} \to \what{\cC}$ has a left adjoint functor, denoted $\ol{a}: \what{\cC} \to \wtilde{\cC}$, which associates to any presheaf a sheaf. 


We also need the following proposition from \cite{SGA4}, expos\'e II.

\begin{prop} \label{prop_presheaves_topology}
Let $\cC$ be a category and $\bF = \{\sF_i\}_{i \in I}$ a family of presheaves over $\cC$. We denote  with $J_\bF(X)$, for all object $X \in \ob(\cC)$, the class of sieves
$R \to X$ such that for any morphism $Y \to X$ of $\cC$ with codomain $X$, the sieve $R \times_X Y$ has the following property: The map
\[ \Hom_{\what{\cC}}(Y, \sF_i) \to \Hom_{\what{\cC}}(R \times_X Y, \sF_i) \]
is bijective (resp. injective) for all $i \in I$. Then, the collections $J_\bF(X)$ define a Grothendieck topology on $\cC$, which is the finest topology such that all $\sF_i$ are sheaves 
(resp. a separable presheaves).
\end{prop}
\begin{proof}
See ibid. proposition 2.2.
\end{proof}

\section{The pro-site and pro-analytic spaces}

Let $k$ be a complete valued field (in this appendix also the trivially valued case is allowed). It is possible to associate to $k$ a category $\bAn_k$ of analytic spaces over $k$.
If $k$ is archimedean the construction is classical and we are mainly interested in the case $k = \C$. If $k$ is non-archimedean there are
different possible choices for the definition of what $\bAn_k$ is. At least there are the following: 
Classical rigid analytic spaces, Berkovich analytic spaces and Huber adic spaces.
In our work we were mainly interested on Berkovich spaces so we will mainly refer to this category, but what we say in this appendix can be recasted in any settings with suitable adaptations. In particular, when we talk about open subsets of some analytic space we refer to the open sets of the topological spaces underlying the 
Berkovich space and not to the admissible opens for the $G$-topology of analytic domains over it.

\index{pro-analytic space}
\begin{defn}
We call an object of $\bPro(\bAn_k)$ a $k$-\emph{pro-analytic space}.
\end{defn}

$\bPro(\bAn_k)$ is a complete category because $\bAn_k$ admits fiber products. The inclusion $\bAn_k \rhook \bPro(\bAn_k)$ is fully faithful and so any analytic space over $k$ can be canonically considered as a pro-analytic space.
But $\bPro(\bAn_k)$ has much more objects, hence also objects of $\bAn_k$ have much more subobjects (\ie open immersions) when seen as objects of $\bPro(\bAn_k)$.

\index{pro-open subset}
\begin{defn} \label{defn:pro_open}
When we consider the small site induced by $\bPro(\bAn_k)$ on a pro-analytic space, we call its admissible open subset \emph{pro-open} subsets.
\end{defn}

So, in some sense the inclusion $\bAn_k \rhook \bPro(\bAn_k)$ enhance also the concept of open subspace. We define a topology on $\bPro(\bAn_k)$ in order to make the definition of pro-open subset rigorous.

Let $\cC$ be a site and $\bPro(\cC)$ the pro-category of $\cC$. For any sheaf $\sF$ on $\cC$ we can define a presheaf of $\bPro(\cC)$ in the following way
\[ \pi^{-1}(\sF)(U) = \pi^{-1}(\sF)(``\limpro_{i \in I}" U_i) \doteq \limind_{i \in I} \sF(U_i). \]
We want to put in a canonical way a topology on $\bPro(\cC)$ such that all the presheaves obtained in this way are sheaves. So, we give the following definition.

\index{pro-site}
\begin{defn} \label{defn_pro_site}
Let $\cC$ be a site. We denote by $\cC_\bPro$ the category $\bPro(\cC)$ endowed with the finest topology such that all presheaves of the form $\pi^{-1}(\sF)$, for
$\sF \in \ob(\wtilde{\cC})$, are sheaves and we call it the \emph{pro-site} of $\cC$ and its topology the \emph{pro-topology} of $\cC_\bPro$.
\end{defn}

\begin{rmk}
From proposition \ref{prop_presheaves_topology} it follows that the definition of the pro-site is well-posed.
\end{rmk}

\begin{prop}
There is a canonical morphism of topoi $(\pi_*, \pi^{-1}): \wtilde{\cC}_\bPro \to \wtilde{\cC}$ for any site $\cC$, where 
\[ \pi_* \sF (X) = \sF( ``\limpro_{i \in I}" X_i ) \]
where $X \in \cC$ and $``\underset{i \in I}\limpro" X_i$ is the trivial system with $X_i = X$ and morphisms equal to the identity.
\end{prop}
\begin{proof}
The functors $\pi^{-1}$ is exact because filtered inductive limits are exact in the category of sets. We have to check that $(\pi_*, \pi^{-1})$ is an adjoint pair of functors, \ie that given any couple of sheaves $\sF \in \ob(\wtilde{\cC})$ and $\sG \in \ob(\wtilde{\cC}_\bPro)$ there is a bijection
\[ \Hom_{\wtilde{\cC}_\bPro} (\pi^{-1} \sF, \sG) \cong \Hom_{\wtilde{\cC}} (\sF, \pi_* \sG). \]
On the one hand, any map $\phi: \pi^{-1} \sF \to \sG$ induces a morphism $\sF \to \pi_* \sG$ simply by restricting $\phi$ to trivial projective systems. On the other hand, given any map
$\phi: \sF \to \pi_* \sG$ and any $``\underset{i \in I}\limpro" X_i$ we can define the map
\[ \limind_{i \in I} \phi: \pi^{-1} \sF (``\limpro_{i \in I}" X_i) \to \pi^{-1} \pi_* \sG (``\limpro_{i \in I}" X_i) \]
because $\pi^{-1} \sF (``\underset{i \in I}\limpro" X_i) = \underset{i \in I}\limind \sF(X_i)$ and $\pi^{-1} \pi_* \sG (``\underset{i \in I}\limpro" X_i) = \underset{i \in I}\limind \sG(X_i)$ and because the direct limit is a functor. So, there is a canonical morphism
\[ \limind_{i \in I} \sG ( X_i ) \to \sG (``\limpro_{i \in I}" X_i) \]
by mean of which we get the required morphism of sheaves $\pi^{-1} \sF \to \sG$.
\end{proof}

\begin{prop}
For any $\sF \in \ob(\wtilde{\cC})$ the canonical morphism $\sF \to \pi_* \pi^{-1} \sF$ is an isomorphisms. In particular $\pi^{-1}$ is fully faithful.
\end{prop}
\begin{proof}
For any $C \in \ob(\cC)$ the morphism $\sF(C) \to \pi_* \pi^{-1} \sF(C)$ is an isomorphism by the definitions of $\pi^{-1}$ and $\pi_*$. Therefore, $\pi^{-1}$ is fully faithful by abstract non-sense.
\end{proof}

Let now $(X, \cO_X)$ be a ringed site. By mean of the canonical morphism of sites $\pi: X_\bPro \to X$ we can pullback the structural sheaf $\cO_X$ to a sheaf of rings
$\pi^{-1} \cO_X$ on $X_\bPro$ which, in this way, can be equipped canonically with a structure of ringed site. The morphism $\pi: X_\bPro \to X$ then becomes a 
morphism of ringed sites and also $(\pi_*, \pi^{-1}): \wtilde{\cC}_\bPro \to \wtilde{\cC}$ becomes a morphism of ringed topoi.

The association $\cC \mapsto \bPro(\cC)$ is functorial in the category of all (small) categories. In fact, if $F: \cC \to \cD$ is a functor then for any $\cC$-pro-object 
$``\underset{i \in I}\limpro" X_i$ we can associate the $\cD$-pro-object $``\underset{i \in I}\limpro" F(X_i)$ and moreover given a morphism $f: ``\underset{i \in I}\limpro" X_i \to ``\underset{j \in J}\limpro" Y_j$ then
we can associate to it the morphism $F f: ``\underset{i \in I}\limpro" F(X_i) \to ``\underset{j \in J}\limpro" F(Y_j)$, because $f$ is a morphism of systems, hence it is defined as a morphism of diagrams of $\cC$.

\index{cofinite diagram}
\begin{defn} \label{defn:cofinite_diagram}
\begin{itemize}
\item A category $I$ is said \emph{loopless} if it has no non-identity endomorphisms.
\item A category $I$ is said \emph{cofinite} if it is small, loopless, and for every object $i \in \ob(I)$, the set of arrows in $I$ with domain $i$ is finite.
\item A diagram $I \to \cC$ is said \emph{cofinite} if the category $I$ is cofinite.
\end{itemize} 
\end{defn}

\begin{prop} \label{prop:pro_site_functorial}
Let $\cC$ and $\cD$ be two sites. Suppose that both $\cC$ and $\cD$ admits all finite limits and that $f: \cC \to \cD$ is a morphism of sites. Then, $f$ induces a morphism of sites $f_\bPro: \cC_\bPro \to \cD_\bPro$.
\end{prop}
\begin{proof}
Since $\cC$ and $\cD$ have fiber products, then every Grothendieck topology on them can be defined by a pre-topology, see \cite{SGA4} 1.3.1, expos\'e II. So, We suppose having pre-topologies that define the topologies of $\cC$ and $\cD$. By definition a family of morphisms $\{ U_i \to U \}_{i \in I}$ in $\cC_\bPro$ is a covering if and only if for every sheaf on $\cC$ the sequence
\begin{equation} \label{eqn:exact_pro}
\pi_{\cC}^{-1} \sF(U) \to \prod_i \pi_{\cC}^{-1} \sF(U_i) \rightrightarrows \prod_{i,j} \pi_{\cC}^{-1} \sF(U_i \times_U U_j)
\end{equation} 
is exact, because of definition \ref{defn_pro_site}. Since the diagram $\{ U_i \to U \}_{i \in I}$ is cofinte, we can apply theorem 3.3 of \cite{ISAK} to deduce a level representation\footnote{See the second section of \cite{ISAK} for the precise definition of level representation. Here we only give the rough idea that it means that we can find an equivalent diagram to the given one where all pro-objects are index by the same category and all morphisms can be expressed as system morphisms between these uniformly indexed systems. We avoid to enter in a detailed explanation of what an equivalence of diagrams is.}
\[ \{ ``\limpro_{k \in K}" V_{i, k} \to ``\limpro_{k \in K}" W_k \}_{i \in I} \]
for a small filtered category $K$, which does not depend on $i \in I$, where $U = ``\underset{k \in K}\limpro" W_k$ and $U_i = ``\underset{k \in K}\limpro" V_{i,k}$. Using this level representation and the definition of the pullback $\pi_{\cC}^{-1}$, the equation (\ref{eqn:exact_pro}) becomes
\[ \limind_{k \in K} \sF(W_k) \to \prod_{i \in I} \limind_{k \in K} \sF(V_{i,k}) \rightrightarrows \prod_{i,j \in I} \limind_{k \in K} \sF(V_{i,k} \times_{W_k} V_{j,k}). \]
Let $\ol{f}: \cD \to \cC$ be the continuous functor associated to $f$. As remarked so far, $\ol{f}$ induces a functor $\ol{f}_{\bPro}: \cD_{\bPro} \to \cC_{\bPro}$,
and since $\ol{f}$ is continuous for any covering $\{ U_i \to U \}_{i \in I}$ in $\cD$ the family $\{ \ol{f}( U_i) \to \ol{f} (U) \}_{i \in I}$ is a covering of $\cC$. Hence, considering a covering $\{ U_i \to U \}$ of $\cD_\bPro$ and $\sF \in \ob(\wtilde{\cC})$, then the exactness of
\[ \pi_{\cC}^{-1} \sF(\ol{f}_\bPro (U)) \to \prod_{i \in I} \pi_{\cC}^{-1} \sF(\ol{f}_\bPro (U_i)) \rightrightarrows \prod_{i,j \in I} \pi_{\cC}^{-1} \sF(\ol{f}_\bPro(U_i) 
   \times_{\ol{f}_\bPro(U)} \ol{f}_\bPro(U_j)) \]
is equivalent to the exactness of
\[ \limind_{k \in K} \sF(\ol{f}(W_k)) \to \prod_{i \in I} \limind_{k \in K} \sF(\ol{f}(V_{i,j})) \rightrightarrows 
   \prod_{i,j \in I} \limind_{k \in K} \sF(\ol{f}(V_{i,k}) \times_{\ol{f}(W_k)} \ol{f}(V_{j,k})) \]
which is equivalent to
\[ \limind_{k \in K} \sF(f(W_k)) \to \prod_{i \in I} \limind_{k \in K} \sF(f(V_{i,j})) \rightrightarrows 
   \prod_{i,j \in I} \limind_{k \in K} \sF(f(V_{i,k} \times_{W_k} V_{j,k})) \]
which is exact by the definition of pro-site.

To conclude the proof we need to show that the pullback functor $f_\bPro^{-1}$ is exact. Without lost of generality we can consider a finite loopless diagram $\{U_i\}_{i \in I}$ in $\cD_\bPro$, because every finite diagram as a cofinal loopless diagram. We can apply theorem 3.1 of \cite{ISAK} to get a level representation for the diagram $\{U_i\}_{i \in I}$. Then, by the definition of $f_\bPro$ we can apply it to the level representation of $\{U_i\}_{i \in I}$ and obtain the exactness of $f_\bPro$ from the levelwise exactness given by the fact that $f$ is a morphism of sites.
\end{proof}

\begin{rmk}
	We expect that proposition \ref{prop:pro_site_functorial} holds without any restrictions on the categories $\cC$ and $\cD$, but the hypothesis that they have pullbacks simplify a lot the proof and it is not a restrictive hypothesis for applications.
\end{rmk}

We now describe some examples of the definitions given so far.

\begin{prop} \label{prop:pro_opne}
If $X$ is a $T_1$ topological spaces then there is an inclusion functor $\sP(X) \rhook \bPro(\bOuv(X))$.
\end{prop}
\begin{proof}
If $X$ is a $T_1$ topological space, then points are closed subsets. Hence, given any subset $S \subset X$ and any open neighborhood $U \subset X$ of $S$ and any finite sets of points $s_1, \ldots, s_n \in X - S$ the set $U - \{s_1, \ldots, s_n\}$ is open in $X$. This proves that the intersection of the system of neighborhoods of $S$ is equal to $S$, which is therefore a 
``pro-open'' set.
\end{proof}

Therefore, every subset of a $T_1$ topological space becomes a pro-open set in the pro-site associated to it.

\begin{exa}
\ben
\item Consider the Zariski site on $\Max(\C [T]) \cong \C$ and the subset $D = \{ z \in \C | |z| \le 1 \} \subset \Max(\C [T]) \C$. Since the Zariski topology on $\C$ is $T_1$, 
      then $D$ is a pro-open set and $\cO_{Zar}^{\bPro}(D)$ (the canonical pro-ringed sites structure on $\Max(\C [T])$ as described so far) coincides with the rational functions in $\C (T)$ with poles outside $D$. 
\item In analogy with the category of pro-analytic spaces we can construct the category of pro-algebraic variety and pro-schemes.
\item The pro-analytic site of a complex analytic space gives a systematic way to deal with the operation of taking germs of analytic functions, as dagger affinoid algebras do.
\item Another interesting example, not strictly linked with the scope of this thesis, is the application of the construction of the pro-site to the \'etale site of the category of
      schemes or to the \'etale site over a scheme. In fact, there is already in literature a construction called the ``pro-\'etale'' site, see \cite{STACK} for an introduction to
      the topic. In particular by the results of section $17$ of \cite{STACK}, we see that there is a canonical morphism of sites from the pro-site of the \'etale site, in our terminology, to the pro-\'etale site of \cite{STACK}. It is not clear whether these two definitions agree.
\item Consider the topological space $X = [0, 1]$ equipped with its usual topology. The subintervals $[0, \frac{1}{n}], \ldots, [\frac{n -1}{n}, n]$ form an admissible covering for the pro-site topology
      on $X_\bPro$. This is due to the fact that $\limind$ is an exact functor that commutes with finite limits. This explains in a general way why different constructions scattered in literature have to
      required some finiteness conditions on the coverings: the compact Stein site, the rigid weak $G$-topology, the dagger weak $G$-topology are the main examples. 
      The requirement of finiteness on covering is necessary because the functor $\limind$ is exact but does not commute in general with all limits. This fact naturally bound the coverings of the pro-site. Indeed, since $X$ is Hausdorff any subset of $X$ is pro-open (cf. \ref{prop:pro_opne}) and we can cover $X$ with the following covering
      \[ \bigcup_{x \in [0, 1]} \{ x \} = [0, 1]. \]
      All the subset $\{x\}$ are admissible opens for the pro-site associated to $X$ but is easy to see that the covering is not an admissible covering for the topology of $X_\bPro$.
\een
\end{exa}

Finally we remark that $(\bAn_k)_\bPro$ is equivalent to the category of pro-analytic spaces defined by Berkovich as a bare category. So, we are adding more structure on it by putting a natural topology induced by the ones on $\bAn_k$ and consider it as a site.
It seems that is worth to study these notions in more details in the future and their relations between other known notions like the pro-\'etale site.

We conclude this appendix with a counter-example. We have seen how the notion of pro-open generalize the notion of open space of a geometric space. This generalization is very broad and allows pathological examples. The next counter-example shows how the theory of dagger analytic spaces is precisely what is needed for discerning the good from the non-good pro-open subset.

\begin{exa}
	Consider the complex compact Stein space $K'$ described after remark \ref{rmk:gg_end}. This space has the very unpleasant property that its associated algebra of germs of analytic functions is non-Noetherian. Using Siu's theorem (cf. \ref{thm:siu}) it is not hard to see that $K'$ is not a $\C$-dagger affinoid space. Moreover, it can be proved (not as easily) that conversely every compact Stein space whose associated algebra of germs of analytic functions is Noetherian is a $\C$-dagger analytic space (we will work out the details of this issue in the future work \cite{BA2}).	
\end{exa}

So, the theory of $\C$-dagger analytic spaces gives a natural theoretical framework for explaining the classical interests of complex geometers for compact Stein spaces whose associated algebra of germs of analytic functions is Noetherian in comparison with the ones for which it is non-Noetherian: The former ones are $\C$-dagger analytic spaces whereas the latters are not. Finally, we remark that what we just described for $\C$ holds also over non-archimedean fields and it is explained (among other results) in the work of Liu \cite{Liu}. In that article it is shown, in the context of classical rigid spaces, that compact Stein spaces over $k$ can be described as finite unions of affinoid spaces over $k$ and their algebras of analytic functions are Noetherian. The theory of $\C$-dagger analytic spaces permits to make the results of Liu uniform over any base field, in the overconvergent setting, and to uniformize the notions of compact Stein space.

\chapter{Generalized rings of Durov} \label{app_durov}

This appendix is a recall of the work done by Durov in his Ph.D. thesis \cite{DUR}. We recall only the definitions and results that are used in this work and add some remarks in the final part of the appendix.

\section{Monads and generalized rings}

\index{monad}
\begin{defn} \label{monad}
A \emph{monad} over a category $\bC$ is a triple $\Sigma= (\Sigma, \mu, \epsilon)$ where $\Sigma: \bC \longrightarrow \bC$ is
an endofunctor, $\mu: \Sigma^2 \longrightarrow \Sigma$ (where $\Sigma^2$ denote the composition $\Sigma \circ \Sigma$), 
$\epsilon: \id_\bC \longrightarrow \Sigma$ are natural transformations, called \emph{multiplication} and \emph{identity}, such that $\mu$ 
and $\epsilon$ respect the axioms of associativity and of unit, \ie for any $X \in \ob(\bC)$ the following diagrams
\[
\begin{tikzpicture}
\matrix(m)[matrix of math nodes,
row sep=2.6em, column sep=2.8em,
text height=1.5ex, text depth=0.25ex]
{ \Sigma^3(X) & \Sigma^2(X) \\
  \Sigma^2(X) & \Sigma(X) \\};
\path[->,font=\scriptsize]
(m-1-1) edge node[auto] {$\mu_{\Sigma(X)}$} (m-2-1);
\path[->,font=\scriptsize]
(m-1-1) edge node[auto] {$\Sigma(\mu_X)$} (m-1-2);
\path[->,font=\scriptsize]
(m-1-2) edge node[auto] {$\mu_X$}  (m-2-2);
\path[->,font=\scriptsize]
(m-2-1) edge node[auto] {$\mu_X$} (m-2-2);
\end{tikzpicture} ,
\begin{tikzpicture}
\matrix(m)[matrix of math nodes,
row sep=2.6em, column sep=2.8em,
text height=1.5ex, text depth=0.25ex]
{ \Sigma(X)  & \Sigma^2(X)  & \Sigma(X) \\
             & \Sigma(X) \\};
\path[->,font=\scriptsize]
(m-1-1) edge node[auto] {$\epsilon_{\Sigma(X)}$} (m-1-2);
\path[->,font=\scriptsize]
(m-1-1) edge node[below] {$\id_{\Sigma(X)}$} (m-2-2);
\path[->,font=\scriptsize]
(m-1-2) edge node[auto] {$\mu_X$}  (m-2-2);
\path[->,font=\scriptsize]
(m-1-3) edge node[auto] {$\id_{\Sigma(X)}$}  (m-2-2);
\path[->,font=\scriptsize]
(m-1-3) edge node[above] {$\Sigma(\epsilon_X)$}  (m-1-2);
\end{tikzpicture}
\]
commutes.
 
A \textit{morphism of monads} over $\bC$, $\varphi: \Sigma_1 \longrightarrow \Sigma_2$, is a natural transformation between the endofunctor which
defines $\Sigma_1$ and $\Sigma_2$ such that
\begin{align*}\label{morphmonad}
&\varphi_Y \circ \Sigma_1(f)=\Sigma_2(f)\circ \varphi_X, \quad  \varphi_X\circ\epsilon_{\Sigma_1 (X)}=\epsilon_{\Sigma_2(X)} \ \ \text{ and }\nonumber \\
\mu_{\Sigma_2(X)}\circ&\Sigma_2(\varphi_X) \circ \varphi_{\Sigma_1(X)}=\varphi_X\circ\mu_{\Sigma_1(X)}=
\mu_{\Sigma_2(X)}\circ\varphi_{\Sigma_2(X)} \circ \Sigma_1(\varphi_X)
\end{align*}
for each couple $X,Y \in \ob(\bC)$ and any morphism $f:X \longrightarrow Y$. We can restate the last conditions requiring the following diagrams to commute
\[
\begin{tikzpicture}
\matrix(m)[matrix of math nodes,
row sep=2.6em, column sep=2.8em,
text height=1.5ex, text depth=0.25ex]
{ \Sigma_1(X) & \Sigma_2(X) \\
  \Sigma_1(Y) & \Sigma_2(Y) \\};
\path[->,font=\scriptsize]
(m-1-1) edge node[auto] {$\Sigma_1(f)$} (m-2-1);
\path[->,font=\scriptsize]
(m-1-1) edge node[auto] {$\varphi_X$} (m-1-2);
\path[->,font=\scriptsize]
(m-2-1) edge node[auto] {$\varphi_Y$}  (m-2-2);
\path[->,font=\scriptsize]
(m-1-2) edge node[auto] {$\Sigma_2(f)$} (m-2-2);
\end{tikzpicture} ,
\begin{tikzpicture}
\matrix(m)[matrix of math nodes,
row sep=2.6em, column sep=2.8em,
text height=1.5ex, text depth=0.25ex]
{   & X \\
\Sigma_1(X) &  & \Sigma_2(X) \\};
\path[->,font=\scriptsize]
(m-1-2) edge node[auto] {$\epsilon_{\Sigma_1(X)}$} (m-2-1);
\path[->,font=\scriptsize]
(m-2-1) edge node[auto] {$\varphi_X$} (m-2-3);
\path[->,font=\scriptsize]
(m-1-2) edge node[auto] {$\epsilon_{\Sigma_2(X)}$}  (m-2-3);
\end{tikzpicture}
\]
\[
\begin{tikzpicture}
\matrix(m)[matrix of math nodes,
row sep=2.6em, column sep=2.8em,
text height=1.5ex, text depth=0.25ex]
{                     & \Sigma_1(\Sigma_1(X)) & \\
\Sigma_1(\Sigma_2(X)) &                       & \Sigma_2(\Sigma_1(X))   \\
                      &      \Sigma_1(X)     &                         \\
\Sigma_2(\Sigma_2(X)) &                       & \Sigma_2(\Sigma_2(X))   \\
                      &      \Sigma_2(X)     &    \\};
\path[->,font=\scriptsize]
(m-1-2) edge node[auto] {$\Sigma_1(\varphi_X)$} (m-2-1);
\path[->,font=\scriptsize]
(m-1-2) edge node[auto] {$\varphi_{\Sigma_1(X)}$} (m-2-3);
\path[->,font=\scriptsize]
(m-1-2) edge node[auto] {$\mu_{\Sigma_1, X}$}  (m-3-2);
\path[->,font=\scriptsize]
(m-2-1) edge node[auto] {$\varphi_{\Sigma_2(X)}$} (m-4-1);
\path[->,font=\scriptsize]
(m-2-3) edge node[auto] {$\Sigma_2(\varphi_X)$} (m-4-3);
\path[->,font=\scriptsize]
(m-4-1) edge node[auto] {$\mu_{\Sigma_2, X}$} (m-5-2);
\path[->,font=\scriptsize]
(m-3-2) edge node[auto] {$\varphi_X$} (m-5-2);
\path[->,font=\scriptsize]
(m-4-3) edge node[auto] {$\mu_{\Sigma_2, X}$} (m-5-2);
\end{tikzpicture}.
\]

\end{defn}

The class of all monads over $\bC$ is therefore endowed with a structure of category and it is denoted by $\bMon(\bC)$.

\index{field without elements}
\begin{defn}
$\bMon(\bSet)$ has an initial object (the identity functor) which is denoted by $\F_\void$ and it is called \emph{the field without elements}.
\end{defn}

We recall the following definition of category theory.

\index{sub-functor}
\begin{defn}
Let $F: \bC \to \bSet$ be a functor. We say that $G: \bC \to \bSet$ is a \emph{sub-functor} of $F$, denoted by $G \subset F$, if 
\ben
\item for all $X \in \ob(\bC)$, $G(X) \subset F(X)$ and
\item for every morphism $f: X \longrightarrow Y$, $G(f)=F(f)_{|G(X)}$.
\een
\end{defn}

We can go on with our exposition of the theory of monads.

\index{sub-monad}
\begin{defn}
Let $\Sigma \in \ob(\bMon(\bSet))$.
A \emph{sub-monad} $\Sigma'$ of $\Sigma$ is a subfunctor $\Sigma' \subset \Sigma$ stable for $\mu$ and $\epsilon,$ \ie a monad $\Sigma'$ such that
\ben
\item $\Sigma'(X) \subset \Sigma(X)$ for all $X \in \ob(\bSet)$; 
\item the natural transformation $\Sigma'\longrightarrow \Sigma$ induced by the inclusion is a morphism of monads.
\een
\end{defn}

Given two submonads $\Sigma_1, \Sigma_2 \subset \Sigma$, one defines their intersection as the fiber product:
\[ \Sigma_1 \cap \Sigma_2 \doteq \Sigma_1 \times_\Sigma \Sigma_2. \]

\index{algebraic functor}
\begin{defn}
An endofunctor $\Sigma: \bC \to \bC$ is said \emph{algebraic} if it commutes with filtered colimits. 
A monad $\Sigma= (\Sigma, \mu, \epsilon)$ over $\bC$ is said \emph{algebraic} if the functor $\Sigma$ is algebraic.
\end{defn}

The composition of two algebraic functors is an algebraic functor and an algebraic functor $\Sigma: \bSet \to \bSet$ is uniquely determined 
by its values on finite sets, because any set is isomorphic to the filtered direct limit of its finite subsets.

\begin{notation}
We denote by $\ul{\N}$ the full subcategory of $\bSet$ formed by the \emph{standard finite sets} \ie the sets 
$\{1, 2, ..., n\}$ for all $n \in \N$. Moreover, we denote objects of $\ul{\N}$ by bold letter like $\bn, \bm, \bone, \btwo$.
\end{notation}

Thus, the restriction $\Sigma \mapsto \Sigma|_{\ul{\N}}$ defines an equivalence of categories between the category
of algebraic endofunctor over $\bSet$ and $\Hom_{\bCat} (\ul{\N}, \bSet)$.
Henceforth all our monads are defined over $\bSet$, if not differently stated.

\begin{notation}
Let $\Sigma \in \ob( \bMon(\bSet))$ be algebraic. We denote:
\ben
\item $|\Sigma| \doteq \Sigma (\bone)$ and we call it \emph{the underlying monoid} of $\Sigma$;
\item $\| \Sigma \| \doteq \underset{\bn \in \ul{\N}}\coprod \Sigma(\bn)$;
\item we call an element $t \in \Sigma(\bn)$ an \emph{$n$-ary operation}.
\een
\end{notation}

We recall that the data of an algebraic monad $\Sigma$ is equivalent to the following set of data (cf. \cite{DUR} section 4.1 and 4.3):

\ben
\item A collection of sets $\{\Sigma(\bn)\}_{n \in \N}$ and maps $\Sigma(\phi):\Sigma(\bm)\to\Sigma(\bn)$, defined for any map $\phi: \bm \to \bn$, 
      subject to the conditions $\Sigma(\id_\bn) = \id_{\Sigma(\bn)}$ and $\Sigma(\psi \circ \phi)=\Sigma(\psi) \circ \Sigma(\phi)$.
      
\item A collection of ``multiplication'' or ``evaluation'' maps $\mu_n^{(k)}: \Sigma(\bk) \times \Sigma(\bn)^k \to \Sigma(\bn)$, subject to the conditions 
      $\mu_n^{(k)} \circ (\id_{\Sigma(\bk)} \times \Sigma(\phi)^n) = \mu_n^{(k')} \circ (\Sigma(\phi) \times \id_{\Sigma(\bn)^{k'}})$ for any 
      $\phi:\bk \to \bk'$, and $\mu_n^{(k)} \circ (\id_{\Sigma(k)} \times \Sigma(\psi)^k) = \mu_m^{(k)}$ for any $\psi: \bm \to \bn$.
      We write $[t]_{\Sigma(\bn)}(x_1, x_2, \ldots, x_k) = t(x_1, x_2, \ldots, x_k)$ instead of $\mu_n^{(k)}(t; x_1, x_2, \ldots, x_k)$ for any $t \in \Sigma(k)$ and 
      $x_1$, $\ldots$, $x_k \in \Sigma(n)$. In this case the above requirements can be rewritten as 
      \[ \phi_*(t)(x_1,\ldots,x_{k'})= t(x_{\phi(1)},\ldots,x_{\phi(k)}) \]
      and 
      \[ \psi_* \l ( t (x_1, \ldots, x_k) \r ) = t (\psi_*x_1, \ldots, \psi_*x_k). \]
\item An element $\bee \in\Sigma(1)$, called the \emph{identity} of $\Sigma$, such that $\bee(x) = x$ for any $x \in \Sigma(\bn)$.
\een

Similarly, a pre-action $\alpha: \Sigma(X) \to X$ is given by a collection of maps $\a^{(n)}: \Sigma(n) \times X^n \to X$, 
      subject relations analogous to (2) and we will write $[t]_\alpha(x_1,\ldots,x_n) = [t]_X(x_1,\ldots,x_n) = t(x_1,\ldots,x_n)$
      instead of $\alpha^{(n)}(t;x_1,\ldots,x_n)$. Notice that $\mu_X: \Sigma^2(X) \to \Sigma(X)$ 
      defines a pre-action of $\Sigma$ on $\Sigma(X)$, so the notation $[t]_{\Sigma(X)}(x_1,\ldots,x_n)$ makes sense for any $t \in \Sigma(n)$ 
      and $x_1$, $\ldots$, $x_n\in\Sigma(X)$. If we take $X=\bn$, we recover again the maps $\mu_n^{(k)}$, so the notation $[t]_{\Sigma(n)}$ 
      is consistent.

Let $\Sigma$ be an algebraic monad and $U \subset ||\Sigma||$, one can consider the smallest submonad of $\Sigma$ containing $U$,
\ie the intersection of the family of submonads $\Sigma_\alpha$ of $\Sigma$ which contain $U$, and denote it by
\[ \langle U \rangle \doteq  \bigcap_{U\subset \|\Sigma_\alpha\|} \Sigma_\alpha. \]

\begin{prop}  \label{prop:existence_submonad}
Let $U \subset ||\Sigma||$ be a subset of a monad $\Sigma$. Then $\Sigma':= \langle U \rangle$ is the submonad of $\Sigma$ whose sets $\Sigma'(\bn)$ are obtained applying a finite number of times the following rules:
\ben
\item $\{k\}_\bn \in \Sigma'(\bn)$ for all $1\leq k \leq n.$
\item (Replacement property) Let $u \in U$ be a $k$-operation (\ie an element of $\Sigma(\bk)$). Then, for all $t_1, \ldots, t_k \in \Sigma'(\bn)$, one has
      \[ [u]_{\Sigma(\bn)}(t_1, \ldots, t_k) \in \Sigma'(\bn). \]
\een
\end{prop}
\begin{proof} 
Cf. \cite{DUR} 4.5.2.
\end{proof}

We can generalize the previous definition in the following way. Let $\rho: \Sigma_0 \longrightarrow \Sigma$ be a morphism of monads and 
$U \subset \|\Sigma\|$, one calls the \emph{sub-monad generated by $U$ over $\Sigma_0$} the monad 
\[ \Sigma_0\langle U \rangle= \langle U \cup ||\rho(\Sigma_0)||\rangle. \]

\index{monad of finite type}
\begin{defn}
Let $\rho: \Sigma_0 \longrightarrow \Sigma$ be a morphism of monads.
An algebraic monad $\Sigma$ is said of \emph{finite type} over $\Sigma_0$ if there exists a finte set $U$ of $\|\Sigma\|$ such that
$\Sigma = \Sigma_0 \langle U \rangle$. 

A monad is said \emph{absolutely of finite type} if it is of finite type over $\F_\void.$ 
\end{defn}

\index{commutative monad}
\begin{defn}
Let $\Sigma$ be an algebraic monad. Two operations $t \in \Sigma(\bn)$ and $s \in \Sigma(\bm)$ are said to \emph{commute} if for any
$\Sigma$-module $X$ and any $x_{i, j} \in X^{n m}$, with $1 \le i \le n, 1 \le j \le m$ we have
\[ s (t(x_{1, 1}, \cdots, x_{n, 1}), \cdots, t(x_{1, m}, \cdots, x_{n, m})) = t(s(x_{1, 1}, \cdots, x_{1, m}), \cdots,s(x_{m, 1}, \cdots, x_{m, n})). \]
An algebraic monad monad is said \emph{commutative} if any couple of elements of $\|\Sigma\|$ commutes.

A commutative algebraic monad over $\bSet$ is called a \emph{generalized ring} (in Durov terminology).
\end{defn}

The terminology "generalized ring" is justified by the fact that the category of rings embeds fully faithfully in the category of generalized rings, which also contains as a full sub-category the category of monoids and much more objects. We refer to \cite{DUR} for the explicit description of how to see rings and monoids as algebraic monads.

\section{Monads given by multiplicative submonoids of rings}

We recall the following construction from \cite{DUR}, 4.3.8 and 5.1.3.

\begin{defn}
Let $X$ be a set we define the \emph{absolute endomorphism ring} of $X$ as the monad $\End [X]$ given by the following algebraic monad:
\[ \End [X](\bn) \doteq \Hom_\bSet(X^n, X), \ \ \ \bee_{\End [X]} = \Id_X \in \Hom_\bSet(X, X) \]
with the multiplication maps
\[ \mu_n^{(k)}: \Hom_\bSet(X^k, X) \times \Hom_\bSet(X^n, X)^k  \to \Hom_\bSet(X^n, X) \]
given by the usual composition of maps, \ie 
\[ (g, f_1, \ldots, f_k) \mapsto g \circ (f_1 \times \ldots \times f_k). \]
\end{defn}

\begin{rmk}
Let $\Sigma$ be an algebraic monad. To give an action $\a: \Sigma(X) \to X$ of $\Sigma$ on $X$ is equivalent to give a morphism of monads
$\rho: \Sigma \to \End[X]$.
\end{rmk}

\begin{defn}
Let $X$ be a set and $Y \subset X$. We define the \emph{relative endomorphism ring of $Y$ in $X$} as the monad $\End_X [Y]$ given by the submonad
of $\End [X]$ of all the maps that stabilize $Y$, \ie
\[ \End_X [Y](\bn) \doteq \{ f \in \Hom_\bSet(X^n, X) | f(Y^n) \subset Y \}. \]
\end{defn}

It is readily verified that $ \End_X [Y]$ is a submonad of $ \End [X]$.
In the case when $Y$ is a submonoid of a ring, the last definition is equivalent to the following.

\begin{defn} \label{defn:power_bounded_monad}
Let $A$ be a ring and $M \subset A$ be a multiplicative submonoid of $A$ with $0 \in M$. For any $n \in \N$ we define
\[ \cR^n(M) \doteq \l \{ (\la_1, \cdots, \la_n) \in A^n | \sum_{i = 1}^n \la_i x_i \in M, \forall x_i \in M \r \} \]
and
\[ \cR(M) \doteq \bigcup_{n \in \N} \cR^n(M). \]
We associate to $M$ the functor $\Sigma_M \in \bMon(\bSet)$ by taking finite linear combinations of the form
\[ \Sigma_M(X) = \l \{ \sum_{x \in X} m_x x | x \in X, (m_x) \in \cR(M) \r \}. \]
We also define the map $\epsilon_X: X \to X$ to be the identity and $\mu_X: \Sigma_M(\Sigma_M(X)) \to \Sigma_M(X)$ as
\[ \mu_X \l ( \sum_{j = 1}^m \mu_j \l ( \sum_{i = 1}^n \la_{i, j} x_i \r ) \r ) = 
   \sum_{i = 1}^n \l ( \sum_{j = 1}^m \mu_j \la_{i, j} \r ) x_i. \]
\end{defn}

\begin{rmk}
The sets $\cR^n(M) \subset M^n$ are ideals of $M^n$ equipped with the direct product monoid structure. Indeed, if 
$(\la_1, \cdots, \la_n) \in \cR^n(M)$, $(a_1, ..., a_n) \in M^n$ and $(x_1, \cdots, x_n) \in M^n$ then
\[ \sum_{i = 1}^n \la_i a_i x_i = \sum_{i = 1}^n \la_i (a_i x_i) \in M, \]
hence $(\la_1 a_1, \cdots, \la_n a_n) \in \cR^n(M)$.
\end{rmk}

\begin{prop}
For a multiplicative monoid in a ring $M \subset A$ the triple defined above $(\Sigma_M, \mu, \epsilon)$ is an
algebraic monad over the category of sets.
\end{prop}
\begin{proof}
The only non-trivial fact to check is that $\mu_X$ satisfies the associativity law. We have to check that given $\a \in \Sigma_M^3(X)$ then the equality
\[ (\mu_X \circ \mu_{\Sigma_M(X)})(\a) = (\mu_X \circ \Sigma_M(\mu_X))(\a) \]
holds. The elements of $\Sigma_M^3(X)$ are of the form
\[ \sum_{l = 1}^p \sigma_l \l ( \sum_{j = 1}^m \mu_{j, l} \l ( \sum_{i = 1}^n \la_{i, j, l} x_i \r ) \r ) = \a. \]
So,
\[ \mu_{\Sigma_M(X)}(\a) = \sum_{j = 1}^m \l ( \sum_{l = 1}^p \sigma_l \mu_{j, l} \r ) \l ( \sum_{i = 1}^n \la_{i, j, l} x_i  \r ) \]
and
\[ (\mu_X \circ \mu_{\Sigma_M(X)})(\a) = 
   \sum_{i = 1}^n \l ( \sum_{j = 1}^m \l ( \sum_{l = 1}^p \sigma_l \mu_{j, l} \la_{i, j, l} \r ) \r ) x_i . \]
On the other hand, 
\[ \Sigma_M(\mu_X)(\a) = \sum_{l = 1}^p \sigma_l \l ( \sum_{i = 1}^n \l ( \sum_{j = 1}^m \mu_{j, l} \la_{i, j, l} \r ) x_i  \r ) \]
\[ (\mu_X \circ \Sigma_M(\mu_X))(\a) =   \sum_{i = 1}^n \l ( \sum_{j = 1}^m \l ( \sum_{l = 1}^p \sigma_l \mu_{j, l} \la_{i, j, l} \r ) \r ) x_i. \]
The algebraicity follows directly from the definition of $\Sigma_M$ because only finite sums are used.
\end{proof}

\begin{prop} \label{prop:monad_functorial}
Let $A, B$ be commutative rings, $M_A \subset A, M_B \subset B$ submonoids, and $f: A \to B$ a morphism of rings such that $f(M_A) \subset M_B$.
Then, there is a commutative diagram of monad morphisms
\[
\begin{tikzpicture}
\matrix(m)[matrix of math nodes,
row sep=2.6em, column sep=2.8em,
text height=1.5ex, text depth=0.25ex]
{ \Sigma_A     & \Sigma_B     \\
  \Sigma_{M_A} & \Sigma_{M_B} \\};
\path[->,font=\scriptsize]
(m-1-1) edge node[auto] {$f$} (m-1-2);
\path[->,font=\scriptsize]
(m-2-1) edge node[auto] {$f^\circ$} (m-2-2);
\path[right hook->,font=\scriptsize]
(m-2-1) edge node[auto] {} (m-1-1);
\path[right hook->,font=\scriptsize]
(m-2-2) edge node[auto] {} (m-1-2);
\end{tikzpicture}.
\]
\end{prop}
\begin{proof}
It is easy to check that the monads $\Sigma_{M_A}$ and $\Sigma_{M_B}$ associated to $M_A \subset A$ and $M_B \subset B$ are sub-monads of $\Sigma_A$ and $\Sigma_B$. Therefore, by proposition \ref{prop:existence_submonad} there exists a smallest sub-monad of $\Sigma_B$ which contains $f(M_A)$, which is necessarily contained in $\Sigma_{M_B}$. Therefore, the restriction of $f$ to $\Sigma_{M_A}$ (\ie the pullback with respect to the inclusion $\Sigma_{M_A} \to \Sigma_A$) must factor through $\Sigma_{M_B}$.
\end{proof}

\chapter{The bornology of coefficientwise convergence} \label{app_coeff_wise}

In this appendix we define a bornology on the ring $R \ldbrack X_1, ..., X_n \rdbrack$ for any bornological ring $R$. This bornology, which will be called the bornology of coefficientwise convergence, is used for proving one of our main result, theorem \ref{thm_W}. 

\section{The bornology of coefficientwise convergence}

Let $R$ be a bornological ring and let $R \ldbrack X_1, ..., X_n \rdbrack$ be the ring of formal power-series on $R$. The bornology that we will define on
$R \ldbrack X_1, ..., X_n \rdbrack$ is the bornological analogous of the topology of coefficientwise convergence defined in the first chapter \cite{GR}. 

For all $k \in \N^n$ we denote the projection map $\pi_k: R \ldbrack X_1, ..., X_n \rdbrack \to R$ defined
\[ \pi_k(\sum_{i \in \N^n} a_i X^i) = a_k. \]
This is an $R$-linear map of $R$-modules. We can consider the $R$-modulus homomorphism
\[ \pi: \prod_{i \in \N^n} \pi_v: R \ldbrack X_1, ..., X_n \rdbrack \to R^{\N^n} \]
and consider on $R^{\N^n}$ the product bornology. On $R \ldbrack X_1, ..., X_n \rdbrack$ we can canonically put 
the biggest bornology which makes the map $\pi$ bounded, and with this bornology $\pi$ becomes an isomorphism of bornological modules.
We call it the \emph{bornology of coefficientwise boundedness} on $R \ldbrack X_1, ..., X_n \rdbrack$.

\begin{prop}
Let $R$ be a bornological ring and let $R \ldbrack X_1, ..., X_n \rdbrack$ be equipped with the bornology of coefficientwise boundedness, then
$R \ldbrack X_1, ..., X_n \rdbrack$ is a bornological ring.
\end{prop}
\begin{proof}
On $R^{\N^n}$ the sum and the product by scalars are bounded operations and $\pi: R \ldbrack X_1, ..., X_n \rdbrack \to R^{\N^n}$ is an isomorphism of
$R$-modules, hence also the addition and the scalar product on 
$R \ldbrack X_1, ..., X_n \rdbrack$ are bounded. It remains to check that the multiplication map of $R \ldbrack X_1, ..., X_n \rdbrack$ is 
a bounded map. 

Let 
\[ \mu: R \ldbrack X_1, ..., X_n \rdbrack \times R \ldbrack X_1, ..., X_n \rdbrack \to R \ldbrack X_1, ..., X_n \rdbrack  \]
denote the multiplication map. By the bornological identification $R^{\N^n} \cong R \ldbrack X_1, ..., X_n \rdbrack$ of modules we have to check 
that the map
\[ \pi \circ \mu: R \ldbrack X_1, ..., X_n \rdbrack \times R \ldbrack X_1, ..., X_n \rdbrack \to R^{\N^n} \]
is bounded, which is true if and only if all maps
\[ \pi_k \circ \mu: R \ldbrack X_1, ..., X_n \rdbrack \times R \ldbrack X_1, ..., X_n \rdbrack \to R \]
are bounded. Given $f, g \in R \ldbrack X_1, ..., X_n \rdbrack$, with $f = \underset{j \in \N^n}\sum a_j X^j, g = \underset{j \in \N^n}\sum b_j X^j$ 
\[ \pi_k \circ \mu(f, g) = \sum_{i + j = k} \pi_i(f) \pi_j(g) = \sum_{i + j = k} a_i b_j \]
is a finite composition of the multiplication and the addition maps of $R$. Hence $\pi_k \circ \mu$ is bounded because is a finite composition of bounded maps.
\end{proof}

Specializing the result when $k$ is a valued field, we obtain.

\begin{prop}
$k \ldbrack X_1, ..., X_n \rdbrack$ equipped with the bornology of coefficientwise boundedness is a bornological algebra whose underlying bornological vector space is of convex type.
\end{prop}
\begin{proof}
$k \ldbrack X_1, ..., X_n \rdbrack$ is a bornological algebra as proved so far. The underlying bornological vector space of $k \ldbrack X_1, ..., X_n \rdbrack$ is of convex type because it is isomorphic to the direct product of a family of
bornological vector spaces of convex type.
\end{proof}

We recall some definitions given by Grauert and Remmert in \cite{GR}.
For any $\rho = (\rho_1, ..., \rho_n) \in \R_+^n$, we define the algebra
\[ B_\rho = \l \{ \sum_{i \in \N^n} a_i X^i | \sum_{i \in \N^n} |a_i| \rho^i < \infty \r \} \subset k \ldbrack X_1 ,..., X_n \rdbrack \]
for any complete valued field $k$, archimedean or not. In this case the summation symbol in the formula $\underset{i \in \N^n}\sum |a_i| \rho^i$ always means the usual sum of real numbers. This definition does not fit well with the discussion we developed so far (in particular in chapter \ref{sec:spectra}). Therefore, we change the definition of $B_\rho$ taking
\[ B_\rho \doteq k \lt \rho^{-1} X \gt, \]
where $k \lt \rho^{-1} X \gt$ is defined as in equation \ref{eqn:polydisk_alg} of the beginning of section \ref{sec:over_pw_series}. 

\begin{rmk}
Our change in the definition of $B_\rho$ does not affect the archimedean side of the theory. For the non-archimedean side of the theory, nothing substantial changes. Indeed, we would like to study direct limits of the algebras $B_\rho$ for $\rho \to r$, for some $r$ which can also be $0$. In turns out that both definitions give the same results to the limits. See section 6 of \cite{BABE} for a proof of this fact. 

Therefore, we use the definition of $B_\rho$ which is compatible with what we explained so far.
\end{rmk}

Now, following \cite{GR}, we call the algebra
\[ \limind_{\rho > 0} B_\rho = K_n \]
the $n$-dimensional \emph{stellen algebra} on $k$. $B_\rho$ are $k$-Banach algebras when endowed with the usual norm, hence $K_n$ becomes 
naturally a complete bornological m-algebra when equipped with the direct limit bornology.

\index{bornology of coefficientwise boundedness}
\begin{defn}
We define the \emph{bornology of coefficientwise boundedness} on $K_n$ as the bornology induced by the inclusion 
$K_n \subset k \ldbrack X_1, ..., X_n \rdbrack$, giving to $k \ldbrack X_1, ..., X_n \rdbrack$ the bornology of coefficientwise boundedness.
\end{defn}

It is a classical consequence of Weierstrass preparation theorem that $K_n$ is a local Noetherian algebra. See the first chapter of \cite{GR} for a proof of this fact.

\begin{prop}  \label{prop_coeff_bound_univ_prop}
Let $\fm \subset K_n$ denotes the maximal ideal of $K_n$. 
Then, the coefficientwise boundedness bornology is the weakest for which all the projections 
\[ \a_e: K_n \to \frac{K_n}{\fm^e} \]
for $e \ge 1$ are bounded.
\end{prop}
\begin{proof}
We have the isomorphism
\[ \frac{K_n}{\fm^e} \cong \l \{ \sum_{|i| < e} a_i X^i \in K_n \r \} \cong k^{\binom{n + e - 1}{e -1}} \]
of bornological vector spaces, where the finite dimensional $k$-vector space $\frac{K_n}{\fm^e}$ is equipped with the quotient bornology and $k^{\binom{n + e - 1}{e -1}}$ is equipped with the product bornology of finitely many copies of $k$. This follows from the commutative diagram
\[
\begin{tikzpicture}
\matrix(m)[matrix of math nodes,
row sep=2.6em, column sep=2.8em,
text height=1.5ex, text depth=0.25ex]
{  K_n       & k \ldbrack X_1, ..., X_n \rdbrack  \\
   \frac{K_n}{\fm^e}  & k^{\binom{n + e - 1}{e -1}} \\};
\path[->,font=\scriptsize]
(m-1-1) edge node[auto] {} (m-1-2);
\path[->,font=\scriptsize]
(m-1-1) edge node[auto] {$\a_e$} (m-2-1);
\path[->,font=\scriptsize]
(m-2-1) edge node[auto] {$\cong$} (m-2-2);
\path[->,font=\scriptsize]
(m-1-2) edge node[auto] {$\pi_e$} (m-2-2);
\end{tikzpicture}
\]
where both vertical arrows are quotient maps and $k \ldbrack X_1, ..., X_n \rdbrack$ is equipped by definition with the projective limit bornology induced by the system of maps $\pi_e$ for $e \in \N$.

Suppose that $K_n$ is equipped with a bornology such that all $\a_e$ are bounded, then there exists a unique bounded map
$K_n \to k \ldbrack X_1, ..., X_n \rdbrack$ which makes the diagram
\[
\begin{tikzpicture}
\matrix(m)[matrix of math nodes,
row sep=2.6em, column sep=2.8em,
text height=1.5ex, text depth=0.25ex]
{  K_n       & k \ldbrack X_1, ..., X_n \rdbrack  \\
   \frac{K_n}{\fm^e} &  \frac{k \ldbrack X_1, ..., X_n \rdbrack}{(X_1, ..., X_n)^e} \\};
\path[->,font=\scriptsize]
(m-1-1) edge node[auto] {} (m-1-2);
\path[->,font=\scriptsize]
(m-1-1) edge node[auto] {$\a_e$} (m-2-1);
\path[->,font=\scriptsize]
(m-2-1) edge node[auto] {$\cong$} (m-2-2);
\path[->,font=\scriptsize]
(m-1-2) edge node[auto] {} (m-2-2);
\end{tikzpicture}
\]
commutative. Set-theoretically this map is the injection $K_n \rhook k \ldbrack X_1, ..., X_n \rdbrack$. So, the bornology on $K_n$ must be stronger than the bornology induced on $K_n$ by the coefficientwise boundedness bornology of $k \ldbrack X_1, ..., X_n \rdbrack$.
\end{proof}

We need to recall the following consequence of Krull intersection theorem.

\begin{lemma} \label{lemma:krull}
Let $A$ be a Noetherian local ring with maximal ideal $\fm$ and $J,I \subset A$ two proper ideals of $A$ with $J \ne 0$. If
$I \subset J + \fm^e$ for all $e \ge 1$ then $I \subset J$.
\begin{proof}
By Krull intersection lemma we know that there exists a $k \ge 1$ such that $\fm^k \subset J$. Then, for all $e \ge k$ we have that
$J + \fm^e = J$, hence
\[ I \subset J + \fm^e, \forall e \ge 1 \then I \subset \bigcap_{e \ge 1} (J + \fm^e) = J. \]
\end{proof}
\end{lemma}

\begin{prop}
Every ideal of $K_n \subset k \ldbrack X_1 ,..., X_n \rdbrack$ is bornologically closed for the bornology of coefficientwise boundedness.
\end{prop}
\begin{proof}
By definition \ref{defn:born_closed} to say that a subspace of a bornological vector space $F \subset E$ is bornologically closed is the same to say that all sequences
$\{ f_n \}_{n \in \N} \subset F$, which converges bornologically to a limit in $E$, the limit $\underset{n \to \infty} \lim f_n \in F$.
Let $I \subset K_n$ be an ideal and let $\{ f_n \}_{n \in \N} \subset I \subset K_n$ be a sequence of elements which converges bornologically to an element $f = \underset{n \to \infty} \lim f_n \in K_n$ (see definition \ref{defn:born_closed}). Since $\a_e: K_n \to \frac{K_n}{\fm^e}$ are bounded maps then all the sequences $\{ \a_e(f_n) \}_{n \in \N}$ are convergent in 
$\frac{K_n}{\fm^e} \cong k^{\binom{n + e - 1}{e -1}}$ which is finite dimensional over $k$ and equipped
with the direct product bornology. Thus, necessarily $\a_e(I)$ is closed in $\frac{K_n}{\fm^e}$, hence $\a_e(f) \in \a_e(I)$. 
This is true for each $e$ hence
\[ K_n f \subset I + \fm^e \]
for any $e$ and applying lemma \ref{lemma:krull} we deduce that $f \in I$.
\end{proof}

\begin{cor}
All the ideals of $K_n$ are bornologically closed for the direct limit bornology induced by the direct limit $K_n = \underset{\rho >0} \limind B_\rho$.
\end{cor}
\begin{proof}
For the direct limit bornology all the maps
\[ K_n \to \frac{K_n}{\fm^e} \]
for $e \ge 1$ are bounded when $\frac{K_n}{\fm^e}$ is equipped with the bornology of $k^{\binom{n + e - 1}{e -1}}$. Therefore, this bornology must be finer the bornology of coefficientwise boundedness. By proposition \ref{prop_closed_subspaces} the preimages by bounded maps of bornologically closed subset are bornologically closed, and we get the corollary.
\end{proof}

\begin{cor}
For any ideal $I \subset K_n$ the ideal $I \cap B_\rho$ is closed in $B_\rho$.
\end{cor}
\begin{proof}
The canonical injection $\phi_\rho: B_\rho \rhook K_n$ is a bounded map, hence 
\[ \phi_\rho^{-1}(I) = I \cap B_\rho \]
is a bornologically closed set in $B_\rho$ by proposition \ref{prop_closed_subspaces}.
\end{proof}

We notice that algebraically $K_n \cong \cO_{k^n, 0}$, the ring of germs of analytic functions at $0 \in k^n$, when $k = \R, \C$, and it is 
isomorphic to its "rigid-analytic" counterpart if $k$ is non-archimedean.
Moreover, given any $x \in k^n$ we can consider a neighborhood basis of $x$ by polydisks centred in $x$ and obtain $K_n(x) \cong \cO_{k^n, x}$ as a 
direct limit of $k$-Banach algebras over these polydisks. We call $K_n(x)$ the $n$-dimensional $k$-stellen algebra centred in $x$.
It is clear that there is an isomorphism $K_n \cong K_n(x)$ of algebras and we can put on $K_n(x)$ the bornology of coefficientwise boundedness induced by this isomorphism. We have the following proposition.

\begin{prop} \label{prop:coeff_closed_ideals}
Any ideal of $K_n(x)$ is bornologically closed for the bornology of coefficientwise boundedness and this bornology is the smallest for which the maps
\[ \a_e: K_n(x) \to \frac{K_n(x)}{\fm_x^e} \]
are all bounded, where $\fm_x$ is the maximal ideal of $K_n(x)$.
\end{prop}
\begin{proof}
We have the following commutative diagram
\[
\begin{tikzpicture}
\matrix(m)[matrix of math nodes,
row sep=2.6em, column sep=2.8em,
text height=1.5ex, text depth=0.25ex]
{K_n    & \frac{K_n}{\fm^e} \\
 K_n(x) & \frac{K_n(x)}{\fm_x^e} \\};
\path[->,font=\scriptsize]
(m-1-1) edge node[auto] {} (m-1-2);
\path[->,font=\scriptsize]
(m-1-1) edge node[auto] {} (m-2-1);
\path[->,font=\scriptsize]
(m-2-1) edge node[auto] {} (m-2-2);
\path[->,font=\scriptsize]
(m-1-2) edge node[auto] {} (m-2-2);
\end{tikzpicture}
\]
where $\frac{K_n}{\fm^e} \to \frac{K_n(x)}{\fm_x^e}$ is an isomorphism of bornological vector spaces since is an algebraic isomorphism of finite
dimensional vector spaces equipped with the direct product bornology. Hence all the properties claimed for $K_n(x)$ follow from the respective properties for $K_n$.
\end{proof}

\pagestyle{empty}

\printindex

\end{document}